%% file: ThomasIrredSubgp.tex
\documentclass{memo-l}

\usepackage{amssymb}
\usepackage{float}
\usepackage{latexsym}
\usepackage{longtable}
\usepackage{enumitem}
\usepackage{array}
\usepackage{tikz}
\usepackage{pdflscape}
\usepackage{calc}
\usepackage{adjustbox}
\usepackage[bf,calcwidth=.9\linewidth,hang]{caption}

\newcommand\gtwo[1]{{#1}}
\newcommand\ffour[1]{{#1}}
\newcommand\esix[1]{{#1}}
\newcommand\eseven[1]{{#1}}
\newcommand\eeight[1]{{#1}}

\newtheorem{thm*}{Theorem}
\newtheorem{thm}{Theorem}[chapter]

\newtheorem{lem}[thm]{Lemma}
\newtheorem{cor}[thm]{Corollary}
\newtheorem{cor*}{Corollary}

\theoremstyle{definition}
\newtheorem{defn}[thm]{Definition}

\theoremstyle{remark}

\newenvironment{pf}[1]{\begin{trivlist}
\item[\hskip \labelsep {\bfseries #1}]}{\end{trivlist}}

\newcommand\modcounter{%
  \addtocounter{table}{-1}
  \renewcommand{\thetable}{\arabic{table}A}%
}

\newcommand\unmodcounter{%
  \renewcommand{\thetable}{\arabic{table}}%
}

\numberwithin{section}{chapter}
\numberwithin{equation}{chapter}

\makeindex

\begin{document}

\frontmatter

\title{The Irreducible Subgroups of Exceptional Algebraic Groups}


\author{Adam R. Thomas}
\address{School of Mathematics, University of Bristol, Bristol, BS8 1TW, UK, and Heilbronn Institute for Mathematical Research, Bristol, UK}
\curraddr{}
\email{adamthomas22@gmail.com}
\thanks{The author is indebted to Prof.\ M.\ Liebeck for his help in producing this paper. He would also like to thank Dr\ A.\ Litterick and Dr\ T.\ Burness for their comments on previous versions of this paper. Finally, the author would like to thank the anonymous referee for their careful reading of this paper and many insightful comments and corrections. }

\date{}

\subjclass[2010]{Primary 20G41, 20G07}

\keywords{algebraic groups, exceptional groups, complete reducibility, G-irreducible subgroups, subgroup structure}



\begin{abstract}
This paper is a contribution to the study of the subgroup structure of exceptional algebraic groups over algebraically closed fields of arbitrary characteristic. Following Serre, a closed subgroup of a semisimple algebraic group $G$ is called irreducible if it lies in no proper parabolic subgroup of $G$. In this paper we complete the classification of irreducible connected subgroups of exceptional algebraic groups, providing an explicit set of representatives for the conjugacy classes of such subgroups. Many consequences of this classification are also given. These include results concerning the representations of such subgroups on various $G$-modules: for example, the conjugacy classes of irreducible connected subgroups are determined by their composition factors on the adjoint module of $G$, with one exception.

A result of Liebeck and Testerman shows that each irreducible connected subgroup $X$ of $G$ has only finitely many overgroups and hence the overgroups of $X$ form a lattice. We provide tables that give representatives of each conjugacy class of connected overgroups within this lattice structure. We use this to prove results concerning the subgroup structure of $G$: for example, when the characteristic is $2$, there exists a maximal connected subgroup of $G$ containing a conjugate of every irreducible subgroup $A_1$ of $G$. 
\end{abstract}

\maketitle

\tableofcontents

\mainmatter
\include{Introduction}

\include{Notation}
\include{Prelims}
\include{Strat}
\include{G2}
\include{F4}
\include{E6}

\include{E7}
\include{E8}

\include{Corollaries}
\include{TablesThms}

\include{CompFactorsIrred}

\include{CompFactorsLevi}

\backmatter
\bibliographystyle{amsalpha}
\bibliography{biblio}
\printindex

\end{document}

%% file: Introduction.tex
\chapter{Introduction} \label{intro}

Let $G$ be a reductive connected algebraic group over an algebraically closed field $K$ of characteristic $p$. A closed subgroup $X$ of $G$ is said to be \emph{$G$-completely reducible} (or \emph{$G$-cr} for short) if, whenever it is contained in a parabolic subgroup $P$ of $G$, it is contained in a Levi subgroup of $P$. This definition, due to Serre in \cite{ser04}, generalises the standard notion of a completely reducible subgroup of $\text{GL}(V)$. Indeed, if $G = \text{GL}(V)$, a subgroup $X$ is $G$-completely reducible if and only if $X$ acts completely reducibly on $V$. 

The concept of $G$-complete reducibility offers a bridge between many different branches of mathematics. These include the theory of buildings, Kac-Moody groups, geometric invariant theory, Lie algebras and representation theory. The definition for spherical buildings was introduced in the original paper of Serre \cite{ser04}. This was extended to subgroups of Kac-Moody groups by Caprace in \cite{Cap09} and more generally for twin buildings in \cite{Daw13}. Complete reducibility has been studied from a geometric point of view in a series of works by Bate, Martin, R{\"o}hrle, et al. in which the authors apply ideas from geometric invariant theory (see for example \cite{BGM16}, \cite{bmr}, \cite{BMRT13}). There is also a natural generalisation of complete reducibility to subalgebras of Lie algebras of algebraic groups, as introduced by McNinch in \cite{McN07}, and studied further in \cite{BMRT11} and \cite{StewTho}. 

The notion of complete reducibility for subgroups of algebraic groups is very familiar in characteristic $0$. Indeed, in characteristic $0$ a subgroup $X$ of $G$ is $G$-cr if and only if $X$ is reductive \cite[Proposition 4.2]{ser04}. In positive characteristic, a $G$-cr subgroup is still reductive \cite[Proposition 4.1]{ser04} but the converse need not be true. However, non-$G$-cr reductive connected subgroups are inherently a low characteristic phenomenon: all reductive connected subgroups are $G$-cr when $p \geq a(G)$, where $a(G)$ is equal to $\textrm{rank}(G) + 1$ if $G$ is simple and equal to the supremum of $(1,a(G_1), \ldots, a(G_r))$ otherwise, where $G_i$ are all the simple quotients of $G$, as proved by Jantzen, McNinch and Liebeck--Seitz (see \cite[Theorem 4.4]{ser04}). There are also results concerning the complete reducibility of finite subgroups of $G$, see \cite[Theorem A]{gur99}, \cite[Theorem 1.9]{GHT16}, and \cite[Corollary~5]{Lit16}.

This paper is concerned with an important subset of $G$-cr subgroups, namely the $G$-irreducible subgroups. The definition, again introduced by Serre in \cite{ser04}, is as follows. A closed subgroup $X$ of $G$ is called \emph{$G$-irreducible} if it is not contained in any proper parabolic subgroup of $G$. We also say that the subgroup is \emph{irreducible} if $G$ is clear from the context. It is immediate from the definition that when $G = \text{GL}(V)$, a subgroup $X$ is $G$-irreducible if and only if $V$ restricted to $X$ is an irreducible $X$-module. In \cite{LT}, Liebeck and Testerman studied $G$-irreducible connected subgroups when $G$ is semisimple. They showed, amongst other things, that all irreducible connected subgroups are semisimple and have only a finite number of overgroups in $G$.

The $G$-irreducible connected subgroups play an important role in determining both the $G$-cr and non-$G$-cr connected subgroups of $G$. The $G$-cr subgroups of $G$ are simply the $L$-irreducible subgroups of $L$ for each Levi subgroup $L$ of $G$ (noting that $G$ is a Levi subgroup of itself). To determine the non-$G$-cr subgroups of $G$, one strategy is as follows. Let $P$ be a proper parabolic subgroup with unipotent radical $Q$ and Levi complement $L$. Then for each $L$-irreducible subgroup $X$ of $L$, determine the complements to $Q$  in $Q X$ that are not $Q$-conjugate to $X$ (if any exist). Any non-$G$-cr connected subgroup will be of this form for some $L$-irreducible connected subgroup $X$. This strategy has been used in \cite{g2dav}, \cite{dav} and more recently in \cite{littho}. 

We now restrict our attention to the case where $G$ is a simple algebraic group over an algebraically closed field $K$ of characteristic $p$ (setting $p= \infty$ for characteristic 0). If $G$ is of classical type then determining the $G$-irreducible subgroups reduces to representation-theoretic considerations by \cite[Lemma 2.2]{LT}. In particular, when $(G,p)$ is not ($D_n,2)$, the action of a subgroup $X$ on the natural module for $G$ determines whether or not $X$ is $G$-irreducible. In any case, when $G$ is of small rank (at most $8$ suffices for the purpose of this paper) it is possible to determine the conjugacy classes of $G$-irreducible connected subgroups. 

Now suppose that $G$ is of exceptional type. The simple $G$-irreducible connected subgroups have already been classified through a series of works by various authors. Firstly, work of Liebeck--Seitz in \cite[Theorem~1]{LS3} shows that all reductive connected subgroups $X$ are $G$-cr under the assumption that $p > N(X,G)$, where $N(X,G)$ is a prime depending on the type of $X$ and $G$ and always at most $7$. They use this result to classify the simple connected subgroups of rank at least 2 and Lawther--Testerman used this in \cite{Law} to classify the subgroups of type $A_1$, both when $p > N(X,G)$. In these cases, the $G$-irreducible subgroups are those $G$-cr subgroups with trivial connected centraliser and so one can find all of the $G$-irreducible simple connected subgroups under their assumptions on $p$. In \cite{dav}, Stewart classified the $F_4$-irreducible simple connected subgroups of $F_4$ of rank at least 2 without any assumption on $p$. Amende in \cite{bon} determined the $G$-irreducible subgroups of type $A_1$ when $G$ is not of type $E_8$. Finally, work of the author in \cite{tho1}, \cite{tho2} completed the classification of the simple $G$-irreducible connected subgroups of $G$.

This paper completes the classification of $G$-irreducible connected subgroups of $G$. Our main theorem is the following (here $\text{Aut}(G)$ denotes the group of algebraic automorphisms of $G$).  

\begin{thm*} \label{MAINTHM}
Let $G$ be a simple exceptional algebraic group and $X$ be a $G$-irreducible connected subgroup of $G$. Then $X$ is $\textrm{Aut}(G)$-conjugate to exactly one subgroup in Tables \ref{G2tab}--\ref{E8tab} and each subgroup in the tables is $G$-irreducible.  
\end{thm*}

We also determine the composition factors of each irreducible connected subgroup in the action on the adjoint and minimal modules for $G$; these can be found in Tables \ref{G2tabcomps}--\ref{E8tabcomps}. By ``minimal module'' we mean the smallest dimensional non-trivial module for $G$ (which coincides with the adjoint module when $G$ is of type $E_8$). The dimensions of such a module are $7$ ($6$ if $p=2$), $26$ ($25$ if $p=3$), $27$ and $56$ for $G = G_2$, $F_4$, $E_6$ and $E_7$, respectively. 

We explain how to read Tables \ref{G2tab}--\ref{E8tab} in Section \ref{thmtabs}, where they are presented, but let us make some important remarks about them now. The irreducible connected subgroups $X$ of $G$ listed in the tables are given an identification number $G(\#n)$ (ID number for short). In reference to Theorem \ref{MAINTHM}, we only count each subgroup with ID number $n$ once in Tables \ref{G2tab}--\ref{E8tab}, even if it appears multiple times; the ID number $n$ appears in italics each time the subgroup is repeated. These repeats are necessary to give the lattice structure of the connected overgroups of irreducible subgroups of $G$ and we discuss this further in Section \ref{thmtabs}. We also provide Tables \ref{G2tabaux}--\ref{E8tabaux} in Section \ref{thmtabs} to help recover this lattice structure. They give the conjugacy classes of immediate connected overgroups for irreducible subgroups of $G$. In particular, they make it easier to find all the repetitions of a subgroup $X$ in Tables \ref{G2tab}--\ref{E8tab}. Another important remark to make is that we list large collections of diagonal irreducible connected subgroups in separate tables. We do this to improve the readability of Tables \ref{G2tab}--\ref{E8tab} and in order to condense the presentation of the large number of conjugacy classes of such diagonal subgroups.  

A large part of this paper is devoted to proving Theorem \ref{MAINTHM}. We do this by proving Theorems \ref{G2THM}--\ref{E8THM} which classify the irreducible subgroups of $G_2$--$E_8$, respectively. As mentioned, the author already completed the classification of the simple $G$-irreducible connected subgroups in \cite{tho1}, \cite{tho2} and so the main focus is on the non-simple irreducible connected subgroups. We discuss the strategy used for the proofs in detail in Section \ref{strat}. The strategy involved is different to that used in [Loc. cit.] with the main distinction being our methods for finding the irreducible subgroups of each maximal connected subgroup of $G$. This difference allows us to study the lattice structure of the connected overgroups of each irreducible connected subgroup; the overgroups can be read off from Tables \ref{G2tab}--\ref{E8tab} as aforementioned. This lattice structure allows us to prove Corollaries \ref{A1overgroups} and \ref{A2overgroups} below, as well as Corollary \ref{corvarstein} in Section \ref{steinberg}. Moreover, we believe that the presentation of the lattice structure of all $G$-irreducible connected subgroups will be beneficial to future readers.    

In the remainder of this introduction we present many corollaries of Theorem \ref{MAINTHM}. To do this we require notation used throughout the paper to describe representations of algebraic groups, diagonal subgroups, identification of irreducible subgroups etc. This is explained in Section \ref{nota}. 

For the first two of these corollaries we need the following definition.  Let $X$ and $Y$ be semisimple subgroups of a semisimple algebraic group $G$ and let $V$ be a $G$-module. Then we say that $X$ and $Y$ \emph{have the same composition factors on $V$} if there exists an isomorphism from $X$ to $Y$ sending the set of composition factors of $V \downarrow X$ to the set of composition factors $V\downarrow Y$ (counted with multiplicity).

The first of our corollaries shows that if $G$ is a simple exceptional algebraic group then, with one exception, conjugacy between $G$-irreducible connected subgroups is determined by their composition factors on the adjoint module for $G$, which we denote by $L(G)$.

\begin{cor*} \label{comps}
Let $G$ be a simple exceptional algebraic group and $X$ and $Y$ be $G$-irreducible connected subgroups of $G$. If $X$ and $Y$ have the same composition factors on $L(G)$ then either:
\begin{enumerate}[leftmargin=*,label=\normalfont(\arabic*)]
\item $X$ is conjugate to $Y$ in Aut$(G)$, or
\item $G = E_8$, $X, Y$ are of type $A_2$, $p \neq 3$, $X \hookrightarrow A_2^2 < \bar{D}_4^2$ via $(10,10^{[r]})$ and $Y \hookrightarrow A_2^2 < \bar{D}_4^2$ via $(10,01^{[r]})$ (or vice versa) where $r \neq 0$ and $A_2^2$ is irreducibly embedded in $\bar{D}_4^2$. In the notation of Table \ref{E8tab} the subgroup $X = E_8(\#\eeight{53})$ and the subgroup $Y = E_8(\#\eeight{54})$. 
\end{enumerate}
\end{cor*}

We also deduce that for $G$ not of type $E_8$, the $\text{Aut}(G)$-conjugacy class of a $G$-irreducible connected subgroup of $G$ is determined by its composition factors on the minimal module for $G$.

\begin{cor*} \label{corconjmin}
Let $G$ be a simple exceptional algebraic group not of type $E_8$, and let $X$ and $Y$ be $G$-irreducible connected subgroups of $G$. If $X$ and $Y$ have the same composition factors on a minimal module for $G$ then $X$ is conjugate to $Y$ in $\text{Aut}(G)$. 
\end{cor*}

From the lists of composition factors provided in Tables \ref{G2tabcomps}--\ref{E8tabcomps}, one can determine the $G$-irreducible connected subgroups $X$ for which $V \downarrow X$ is multiplicity-free when $V$ is either the minimal or adjoint module for $G$. This is a specific case of a more general project of Liebeck, Seitz and Testerman; see \cite{LST15} for further details.   

The next corollary highlights interesting subgroups that are not $G$-irreducible but are $M$-irreducible for some reductive, maximal connected subgroup $M$. When we say a reductive, maximal connected subgroup we mean a reductive subgroup that is maximal among all closed connected subgroups; these have been classified by Liebeck--Seitz and are listed in Theorem \ref{maximalexcep}. Here ``interesting'' means that the $M$-irreducible subgroup is not $M_1$-reducible for some other reductive, maximal connected subgroup $M_1$ nor contained in a proper Levi subgroup of $G$. 

We explain some notation we only use in Table \ref{cortab}. A subgroup $X$ is said to be ``embedded via $\lambda_i$'' in a simple classical group $M$ if $V_{M}(\lambda_1) \downarrow X = V_{X}(\lambda_i)$. This determines $X$ up to $M$-conjugacy unless $M$ is of type $D_n$, in which case there may be two classes. Indeed, this happens for each of the subgroups of $D_8$ given in Table \ref{cortab}. However, we distinguish between the two classes in each case by their composition factors on $V_{D_8}(\lambda_7)$ which leads to the definition of $B_4(\ddagger)$ and $A_1 C_4 (\ddagger)$ in Section \ref{D8inE8}. Every subgroup $X$ of $D_8$ listed in Table \ref{cortab} is contained in $B_4(\ddagger)$ or $A_1 C_4 (\ddagger)$ and hence the conjugacy class of $X$ is uniquely determined.     

\begin{cor*} \label{nongcr}
Let $G$ be a simple exceptional algebraic group and $X$ be a connected subgroup of $G$. Suppose that whenever $X$ is contained in a reductive, maximal connected subgroup $M$ of $G$ it is $M$-irreducible and assume that such an overgroup $M$ exists. Assume further that $X$ is not contained in a proper Levi subgroup of $G$. Then either: 

\begin{enumerate}[leftmargin=*,label=\normalfont(\arabic*)]
\item $X$ is $G$-irreducible, or
\item $X$ is $\text{Aut}(G)$-conjugate to a subgroup in Table \ref{cortab}. Such $X$ are non-$G$-cr and satisfy the above hypothesis.  
\end{enumerate}
\end{cor*}

\begin{longtable}{p{0.07\textwidth - 2\tabcolsep}p{0.12\textwidth - 2\tabcolsep}p{0.06\textwidth - 2\tabcolsep}>{\raggedright\arraybackslash}p{0.75\textwidth-\tabcolsep}@{}}

\caption{Non-$G$-cr subgroups that are irreducible in every (and at least one) maximal, reductive overgroup. \label{cortab}} \\

\hline \noalign{\smallskip}
 
 $G$ & Max. $M$ & $p$ & $M$-irreducible subgroup $X$ \\
\hline \noalign{\smallskip} 

$G_2$ & $A_1 \tilde{A}_1$ & $2$ & $A_1 \hookrightarrow M$ via $(1,1)$ \\

\hline \noalign{\smallskip} 

$F_4$ & $A_2 \tilde{A}_2$ & $3$ & $A_2 \hookrightarrow M$ via $(10,01)$ \\
\endfirsthead
\hline \noalign{\smallskip}

$E_6$ & $A_2 G_2$ & $3$ & $A_2 \hookrightarrow A_2 A_2 = E_6(\#\esix{48})$ via $(10,10)$ \\

\hline \noalign{\smallskip}

$E_7$ & $A_1 G_2$ & $7$ & $A_1 \hookrightarrow A_1 A_1= E_7(\#\eseven{334})$ via $(1,1)$ \\

& $\bar{A}_2 A_5$ & $3$ & $A_2 \hookrightarrow \bar{A}_2 A_2 = E_7(\#\eseven{303})$ via $(10,10)$ \\

& $A_7$ & $2$ & $C_4$ embedded via $\lambda_1$ \\

& & & $D_4$ embedded via $\lambda_1$ \\

& & & $B_3$ embedded via $001$ \\

& & & $A_1 B_2$ embedded via $(1,01)$ \\

& & & $A_2$ embedded via $11$ \\

& & & $A_1^3$ embedded via $(1,1,1)$ \\

& & & $A_1 A_1 \hookrightarrow A_1^3$ via $(1_a,1_a^{[r]},1_b)$ $(r \neq 0)$ \\

& & & $A_1 \hookrightarrow A_1^3$ via $(1,1^{[r]},1^{[s]})$ $(0 < r < s)$ \\

& $G_2 C_3$ & $2$ & $G_2 \hookrightarrow G_2 G_2 = E_7(\#\eseven{320})$ via $(10,10)$ \\

\hline \noalign{\smallskip}

$E_8$ & $D_8$ & $2$ & $B_4(\ddagger)$ embedded via $\lambda_4$ \\

& & & $B_2^2 < B_4(\ddagger)$ embedded via $(01,01)$ \\

& & & $B_2 \hookrightarrow B_2^2$ via $(10,10^{[r]})$ $(r \neq 0)$, $(10,02)$ or $(10,02^{[r]})$ $(r \neq 0)$ \\

& &  & $A_1^2 B_2 <  B_4(\ddagger)$ embedded via $(1,1,01)$ \\

& & & $A_1 B_2 \hookrightarrow A_1^2 B_2$ via $(1,1^{[r]},10)$ ($r \neq 0$) \\

& & & $A_1^4 < B_4(\ddagger)$ embedded via $(1,1,1,1)$ \\

& & & $A_1 A_1^2 \hookrightarrow A_1^4$ via $(1_a,1_a^{[r]},1_b,1_c)$ $(r \neq 0)$ \\

& & & $A_1 A_1 \hookrightarrow A_1^4$ via $(1_a,1_a^{[r]},1_a^{[s]},1_b)$ $(0 < r < s)$ \\

& & & $A_1 A_1 \hookrightarrow A_1^4$ via $(1_a,1_a^{[r]},1_b,1_b^{[s]})$ $(0 < r \leq s)$ \\

& & & $A_1 \hookrightarrow A_1^4$ via $(1,1^{[r]},1^{[s]},1^{[t]})$ $(0 < r < s < t)$ \\

& & & $A_1 C_4 (\ddagger)$ embedded via $(1,\lambda_1)$ \\

& &  & $A_1 D_4 < A_1 C_4 (\ddagger)$ embedded via $(1,\lambda_1)$ \\

& &  & $A_1 B_3 < A_1 C_4 (\ddagger)$ embedded via $(1,001)$ \\

& & & $A_1 A_2 < A_1 C_4 (\ddagger)$ embedded via $(1,11)$ \\

& $A_8$ & $3$ & $A_2^2$ embedded via $(10,10)$ \\

& & & $A_2 \hookrightarrow A_2^2$ via $(10,10^{[r]})$ $(r \neq 0)$ or $(10,01^{[r]})$ $(r \neq 0)$  \\

& $G_2 F_4$ & $7$ & $G_2 \hookrightarrow G_2 G_2 = E_8(\#\eeight{1035})$ via $(10,10)$ \\

\hline

\end{longtable}

The corollary gives examples showing that one cannot generalise \cite[Theorem 3.26]{bmr} in certain ways for exceptional algebraic groups. This theorem states that in good characteristic every $H$-cr subgroup is $G$-cr for a regular reductive subgroup $H$ of $G$, where a subgroup is regular if it is normalised by a maximal torus of $G$. We say that a prime $p$ is \emph{good} for a simple exceptional algebraic group $G$ if $p \geq 5$ for $G$ of type $G_2$, $F_4$, $E_6$ or $E_7$, and $p \geq 7$ for $G$ of type $E_8$. We say a prime is \emph{bad} for $G$ if it is not good. 

Firstly, one cannot allow arbitrary bad characteristics. The subgroup $A_1 < G_2$ given in Table \ref{cortab} is $\bar{A}_1 \tilde{A}_1$-irreducible yet non-$G_2$-cr when $p=2$, and $\bar{A}_1 \tilde{A}_1$ is a regular subgroup of $G_2$. 

Secondly, if one considers reductive, maximal connected subgroups of $G$, many of these are regular. However, we do not have such a result for arbitrary reductive, maximal connected subgroups. For example, the subgroup $A_1 \hookrightarrow A_1 A_1 < A_1 G_2 < E_7$ via $(1,1)$ from Table \ref{cortab} is $A_1 G_2$-irreducible yet non-$E_7$-cr when $p=7$, which is even a good characteristic for $E_7$.

It is natural to ask whether $G$-irreducible subgroups of a certain type exist, especially in small characteristics. When $G$ is a simple exceptional algebraic group it is shown in \cite[Theorem 2]{LT} (corrected in \cite[Theorem 7.4]{bon}) that $G$-irreducible connected subgroups of type $A_1$ exist, except for $G = E_6$ when $p=2$. We extend this result to subgroups of type $A_1^n$. 

\begin{cor*} \label{A1subgroups}
Let $G$ be a simple exceptional algebraic group of rank $l$. Then $G$ contains a $G$-irreducible connected subgroup of type $A_1^n$ for all $n \leq l$, unless $G = E_6$. For $G = E_6$, there exists a $G$-irreducible subgroup of type $A_1^n$ if and only if $p \neq 2$ and $n \leq 3$. 
\end{cor*}

Given the existence of irreducible subgroups of type $A_1^n$, we study their overgroups. The next result shows the existence of a reductive, maximal connected subgroup that contains representatives of each conjugacy class of $G$-irreducible subgroups of type $A_1^n$ in small characteristics, with one exception. 

\vspace{0.1cm}

\begin{cor*} \label{A1overgroups}
Let $G$ be a simple exceptional algebraic group and $p=2$ or $3$. Then there exists a reductive, maximal connected subgroup $M$ containing representatives of every $\text{Aut}(G)$-conjugacy class of $G$-irreducible subgroups of type $A_1^n$, unless $G=F_4$ and $p=3$ (in which case two reductive, maximal connected subgroups are required). The following table provides examples of such overgroups $M$.   
\end{cor*}
\setlength\LTleft{0.31\textwidth}
\vspace{0.1cm}
\begin{longtable}{p{0.08\textwidth - 2\tabcolsep}>{\raggedright\arraybackslash}p{0.19\textwidth-2\tabcolsep}>{\raggedright\arraybackslash}p{0.11\textwidth-\tabcolsep}@{}}

\caption{Maximal connected overgroups for $G$-irreducible subgroups of type $A_1^n$. \label{A1overgroupstab}} \\

\hline

$G$ & $p=3$ & $p=2$ \\

\hline

$G_2$ & $\bar{A}_1 \tilde{A}_1$ & $\bar{A}_1 \tilde{A}_1$ \\

$F_4$ & $B_4$ and $A_1 C_3$ & $B_4$ \\

$E_6$ & $C_4$ & \textbf{---} \\

$E_7$ & $\bar{A}_1 D_6$ & $\bar{A}_1 D_6$ \\

$E_8$ & $D_8$ & $D_8$ \\
 \endfirsthead
\hline

\end{longtable}

\setlength\LTleft{0pt}
\vspace{0.1cm}
We also prove similar results for $G$-irreducible subgroups of type $A_2^n$. 
\vspace{0.1cm}
\begin{cor*} \label{A2subgroups}
Let $G$ be a simple exceptional algebraic group of rank $l$. Then for $G$ not of type $E_7$, there exists a $G$-irreducible connected subgroup of type $A_2^n$ if and only if $n \leq \frac{l}{2}$. For $G$ of type $E_7$, there exists a $G$-irreducible subgroup of type $A_2^n$ if and only if $p \neq 2$ and $n \leq 2$.  
\end{cor*}
\vspace{0.1cm}
\begin{cor*} \label{A2overgroups}
Let $G$ be a simple exceptional algebraic group. Then there exists a reductive, maximal connected subgroup $M$ containing representatives of every $\text{Aut}(G)$-conjugacy class of $G$-irreducible subgroups of type $A_2^n$, unless $(G,p)$ is one of the following:  $(G_2,3)$, $(E_6, p \neq 2)$, $(E_7, p \geq 5)$ or $(E_8, p \neq 3)$ (in all cases at most three reductive, maximal connected subgroups are required). The following table provides examples of such overgroups $M$.   
\end{cor*}
\setlength\LTleft{0.14\textwidth}

\begin{longtable}{p{0.08\textwidth - 2\tabcolsep}>{\raggedright\arraybackslash}p{0.2\textwidth-2\tabcolsep}>{\raggedright\arraybackslash}p{0.25\textwidth-2\tabcolsep}>{\raggedright\arraybackslash}p{0.19\textwidth-\tabcolsep}@{}}

\caption{Maximal connected overgroups for $G$-irreducible subgroups of type $A_2^n$. \label{A2overgroupstab}} \\

\hline

$G$ & $p \geq 5$ & $p=3$ & $p=2$ \\

\hline

$G_2$ & $\bar{A}_2$ &  $\bar{A}_2$ and $\tilde{A}_2$ & $\bar{A}_2$ \\

$F_4$ & $\bar{A}_2 \tilde{A}_2$ & $\bar{A}_2 \tilde{A}_2$ & $\bar{A}_2 \tilde{A}_2$ \\

$E_6$ & $\bar{A}_2^3$ and $A_2$ & $\bar{A}_2^3$, $A_2 G_2$ and $G_2$ & $\bar{A}_2^3$  \\

$E_7$ & $\bar{A}_2 A_5$ and $A_2$ & $\bar{A}_2 A_5$ & \textbf{---}\\

$E_8$ & $\bar{A}_2 E_6$ and $D_8$ & $\bar{A}_2 E_6$ & $\bar{A}_2 E_6$ and $D_8$ \\
 \endfirsthead
\hline

\end{longtable}

\setlength\LTleft{0pt}

%% file: Notation.tex
\chapter{Notation} \label{nota}

In this section we present the notation used throughout the paper. At many points in this paper we use the results of \cite{tho1} and \cite{tho2} and therefore have tried to be consistent with the notation in those papers, where possible. In particular, when we come to describe the identification number given to an irreducible connected subgroup $X$ we have chosen to keep the same identification number given to $X$ in \cite{tho2} when $X$ is of type $A_1$. 

Firstly we note that by a subgroup of an algebraic group we always mean a closed subgroup. Similarly, all representations of algebraic groups are assumed to be rational. 

Let $G$ be a simple algebraic group over an algebraically closed field $K$. Let $\Phi$ be the root system of $G$ and $\Phi^+$ be a fixed set of positive roots in $\Phi$. Write $\Pi = \{ \alpha_1, \ldots, \alpha_l \}$ for the simple roots of $G$ and $\lambda_1, \ldots, \lambda_l$ for the fundamental dominant weights of $G$, both with respect to the ordering of the Dynkin diagram as given in \cite[p.\ 250]{bourbaki}. We sometimes use $a_1 a_2 \ldots a_l$ to denote a dominant weight $a_1 \lambda_1 + a_2 \lambda_2 + \cdots + a_l \lambda_l$. We denote by $V_G(\lambda)$ (or just $\lambda$) the irreducible $G$-module of dominant high weight $\lambda$. Similarly, the Weyl module of high weight $\lambda$ is denoted by $W(\lambda) = W_G(\lambda)$ and the tilting module of high weight $\lambda$ is denoted by $T(\lambda)$. Another module we refer to frequently is the adjoint module for $G$; we recall that we denote this by $L(G)$. We let $$V_7 : = W_{G_2}(\lambda_1), \ V_{26}:=W_{F_4}(\lambda_4), \ V_{27} := V_{E_6}(\lambda_1), \ V_{56}:=V_{E_7}(\lambda_7).$$ For $G$-modules $V$ and $W$ we write $V + W$ for the module $V \oplus W$ and let $V^*$ denote the dual module of $V$. If $Y = Y_1 Y_2 \ldots Y_k$, a commuting product of simple algebraic groups, then $(V_1, \ldots, V_k)$ denotes the $Y$-module $V_1 \otimes \dots \otimes V_k$, where each $V_i$ is an irreducible $Y_i$-module. The notation $\bar{X}$ denotes a subgroup of $G$ that is generated by long root subgroups of $G$. If the root system of $G$ has short roots then $\tilde{X}$ denotes a subgroup generated by short root subgroups of $G$.

Suppose that char($K)=p < \infty$, recalling our convention that $p = \infty$ represents characteristic 0. Let $F: G \rightarrow G$ be the standard Frobenius endomorphism (acting on root groups $U_\alpha = \{ u_\alpha(c) \mid c \in K\}$ by $u_{\alpha}(c) \mapsto u_{\alpha}(c^p)$) and $V$ be a $G$-module afforded by a representation $\rho: G \rightarrow \text{GL}(V)$. If $r$ is a positive integer then $V^{[r]}$ denotes the module afforded by the representation $\rho^{[r]} : = \rho \circ F^r$. Let $M_1, \ldots, M_k$ be $G$-modules and $n_1, \ldots, n_k$ be positive integers. Then $M_1^{n_1} / \ldots / M_k^{n_k}$ denotes a $G$-module having the same composition factors as $M_1^{n_1} + \dots + M_k^{n_k}$. Furthermore, $V = M_1 | \ldots | M_k$ denotes a $G$-module with a socle series as follows: $M_k \cong \text{Soc}(V) = \text{Soc}^1(V)$ and for $i > 0$, the module $M_{k-i}$ is isomorphic to $\text{Soc}^{i+1}(V) = \text{Soc}(V/N_{i})$ where $N_{i}$ is the inverse image in $V$ of $\text{Soc}^{i}(V)$ under the quotient mapping $V \rightarrow V / N_{i-1}$ (so $N_0 = 0$ and $N_1 = M_k$). Sometimes, to make things clearer, we will use a tower of modules $$\cfrac{M_1}{\cfrac{M_2}{M_3}}$$ to denote $V = M_1 | M_2 | M_3$.

We need a notation for diagonal subgroups of $Y = H_1 H_2 \ldots H_k$, a commuting product with all of the subgroups $H_i$ simple and of the same type. Let $H$ be a simply connected algebraic group of type $H_1$ and $\hat{Y} = H \times H \times \cdots \times H$, the direct product of $k$ copies of $H$. Then we may regard $Y$ as $\hat{Y} / Z$, where $Z$ is a subgroup of the centre of $\hat{Y}$, and $H_i$ is then regarded as the image of the $i$th projection map. A diagonal subgroup of $\hat{Y}$ is a subgroup $\hat{X} \cong H$ of the following form: $\hat{X} = \{ (\phi_1(h), \ldots, \phi_k(h)) \mid h \in H \}$ where each $\phi_i$ is a surjective endomorphism of $H$. A diagonal subgroup $X$ of $Y$ is the image of a diagonal subgroup of $\hat{Y}$ under the natural map $\hat{Y} \rightarrow Y$. To describe such a subgroup it therefore suffices to give a surjective endomorphism, $\phi_i$, of $H$ for each $i$. By \cite[Section 1.15]{gls3}, $\phi_i = \alpha_i \theta_i F^{r_i}$ where $\alpha_i$ is an inner automorphism, $\theta_i$ is a graph automorphism and $F^{r_i}$ is a power of the standard Frobenius endomorphism. We only wish to distinguish these diagonal subgroups up to conjugacy. Therefore, we assume each $\alpha_i$ is trivial and give a (possibly trivial) graph automorphism $\theta_i$ of $H$ and a non-negative integer $r_i$, for each $1 \leq i \leq k$. 

Such a diagonal subgroup $X$ is denoted by \begin{align*} X \hookrightarrow Y \text{ via } (\lambda_1^{[\theta_1 r_1]}, \lambda_1^{[\theta_2 r_2]},  \ldots, \lambda_1^{[\theta_k r_k]}),\end{align*} where $\lambda_1$ is the first dominant weight of $X$. We often abbreviate this to \begin{align*} X \text{ via } (\lambda_1^{[\theta_1 r_1]}, \ldots, \lambda_1^{[\theta_k r_k]}) \end{align*} if the group $Y$ is clear. Unless $X$ is of type $D_n$ ($n \geq 4$), a graph automorphism is uniquely determined (up to conjugacy) by the image of $\lambda_1$ (including the exceptional graph automorphisms of $B_2$, $F_4$ when $p=2$ and $G_2$ when $p=3$, which takes $\lambda_1$ to $2\lambda_2$, $2\lambda_4$ and $3 \lambda_2$, respectively). In these cases, instead of writing $\lambda_1^{[\theta_i r_i]}$ we write $\mu^{[r_i]}$ where $\mu$ is the image of $\lambda_1$ under $\theta_i$. The only time we need a diagonal subgroup of a product of subgroups of type $D_n$ is when dealing with $Y = \bar{D}_4^2$. We give a notation for the standard graph automorphisms of $\bar{D}_4$: we let $\tau$ denote a standard triality automorphism induced by the permutation $(\alpha_1,\alpha_3,\alpha_4)$ and let $\iota$ denote a standard involutory automorphism induced by the permutation $(\alpha_3,\alpha_4)$. 

For clarity, note that field twists $r, s, t, \ldots$ are not assumed to be distinct. This is consistent with \cite{tho2} but not with \cite{tho1}. 

We extend this notation to describe certain semisimple subgroups of the form $X = X_1 X_2 \ldots X_n$ of $Y =  H_1 H_2 \ldots H_k$ ($n < k$) where each $X_i$ is of type $H$ and the projection of $X$ to each $H_i$ is surjective. Any such subgroup is a commuting product of diagonal subgroups of distinct subsets of the $H_i$. For this reason, we extend our use of the term ``diagonal subgroup'' to include such subgroups $X$. For example, consider diagonal subgroups isomorphic to $A_1^2$ contained in $A_1^4$. They are either a commuting product of one $A_1$ factor and a diagonal subgroup of $A_1^3$, or a commuting product of a diagonal subgroup of $A_1^2$ and another diagonal subgroup of the other $A_1^2$. Therefore, our notation needs to distinguish which of the subgroups $H_i$ each of the simple factors of $X$ project non-trivially to. We give the first factor of $X$ the label $a$, the second factor of $X$ the label $b$ and so on. Then for each $i$ such that $X_1$ has non-trivial projection to $H_i$ we give a subscript $a$ to $\lambda_1^{[\theta_i r_i]}$. For each $j$ such that $X_2$ has non-trivial projection to $H_j$ we give a subscript $b$ to $\lambda_1^{[\theta_j r_j]}$ and so on. For example, consider $X = A_1 A_1$ and $Y = A_1^4$ with the first $A_1$ factor of $X$ embedded diagonally in the first two factors of $Y$ (with field twists $0$ and $r$) and the second $A_1$ factor embedded diagonally in the last two factors of $Y$ (with field twists 0 and $s$). Then we write $X \hookrightarrow Y$ via $(1_a,1^{[r]}_a,1_b,1^{[s]}_b)$. 

We make another natural extension of this notation. Let $X$ and $Y$ be as above and suppose that $S = S_1 \ldots S_l$ is a semisimple group with no factors isomorphic to $H$. Then a subgroup $X S$ of $Y S$ is denoted by \begin{align*} X S \hookrightarrow Y S \text{ via } (\lambda_1^{[\theta_1 r_1]}, \lambda_1^{[\theta_2 r_2]},  \ldots, \lambda_1^{[\theta_k r_k]}, \nu^1_1, \nu^2_1, \ldots, \nu^l_1) \end{align*} where $\nu^i_1$ denotes the first fundamental dominant weight of $S_i$. Again, we will still refer to such subgroups as ``diagonal subgroups''. For example, consider a subgroup $X Z = A_1 B_2$ of $Y Z = A_1^2 B_2$ where $X = A_1$ is a simple diagonal subgroup of $Y = A_1^2$ with field twists 0 and $r$. Then we write $A_1 B_2 \hookrightarrow A_1^2 B_2$ via $(1,1^{[r]},10)$. 

Finally, we show how to combine all of these notations in the most general setting. Suppose that $A$ is a diagonal subgroup of $Y = H_1 \dots H_n$ and $B$ is a diagonal subgroup of $Z = J_1 \dots J_m$, where all the $H_i$ (resp.\ $J_i$) are simple and of the same type $H$ (resp.\ $J$) and $H$ is not isomorphic to $J$. Suppose also that $S$ is a semisimple group with no factors of type $H$ or $J$. Then we naturally concatenate the notations above to give a notation for the subgroup $ABS$ of $YZS$. For example, consider a subgroup $ABS = A_1^3 B_2^2 B_3$ of $YZS = A_1^4 B_2^3 B_3$ where $A$ is a diagonal subgroup of $Y$ via $(1_a,1_b,1_b^{[r]},1_c)$ and $B$ is a diagonal subgroup of $Z$ via $(10_a,10_b^{[s]},10_b^{[t]})$. Then we write $A_1^3 B_2^2 B_3 \hookrightarrow A_1^4 B_2^3 B_3$ via $(1_a,1_b,1_b^{[r]},1_c,10_a,10_b^{[s]},10_b^{[t]},100)$.   

When determining the different conjugacy classes of diagonal subgroups of a given group $X$ we need to understand certain automorphisms of $X$ induced from $G$. We define the following notation to describe the graph automorphisms of a subgroup $X$ of $G$ induced by $N_{G}(X)$. We let $\text{Aut}_G(X) = N_{G}(X) / C_G(X)$ and $\text{Inn}_G(X) = X / Z(X) \cong (X C_G(X)) / C_G(X)$. Therefore $\text{Aut}_G(X) / \text{Inn}_G(X) \cong N_{G}(X) / (X C_G(X))$, the group of graph automorphisms of $X$ induced by $N_G(X)$ and we let $\text{Out}_G(X) = \text{Aut}_G(X) / \text{Inn}_G(X)$. 

Now let $G$ be of exceptional type. In Tables \ref{G2tab}--\ref{E8tab} we give an identification number to each of the conjugacy classes of $G$-irreducible connected subgroups arising in Theorems \ref{G2THM}--\ref{E8THM}. The notation $G(\#a)$ (or simply $a$ if $G$ is clear from the context) means the $G$-irreducible subgroup corresponding to the ID number $a$. We set $G(\#0)$ to be $G$ itself. Sometimes $G(\#a)$ will refer to infinitely many conjugacy classes of $G$-irreducible subgroups. This situation only occurs for diagonal subgroups, where the conjugacy class will depend on field twists $r_1, \ldots, r_k$ and graph automorphisms $\theta_1, \ldots, \theta_k$. Sometimes we refer to a subset of the conjugacy classes that $G(\#a)$ represents; we only do this when all graph automorphisms are trivial. These are thus described by an ordered set of field twists $r_1, \ldots, r_k$ and are denoted by $G(\#a^{\{r_1, \ldots, r_k\}})$. Let us give a concrete example to make this clearer. Consider $G_2(\#\gtwo{1})$, the conjugacy classes of diagonal subgroups $A_1 \hookrightarrow \bar{A}_1 \tilde{A}_1$ via $(1^{[r]},1^{[s]})$ ($rs=0$; $r \neq s$) (see Table \ref{G2tab}). Then the notation $G_2(\#\gtwo{1}^{\{r,0\}})$ refers to the conjugacy classes with $s=0$ and the notation $G_2(\#\gtwo{1}^{\{1,0\}})$ refers to the single conjugacy class $A_1 \hookrightarrow \bar{A}_1 \tilde{A}_1$ via $(1^{[1]},1)$.

In the tables in Section \ref{thmtabs} we use a shorthand for $n^{\{s_1, \ldots, s_k\}}$ in certain situations. The notation $n^{\{\underline{0}\}}$ simply means that each $s_i$ is equal to $0$. The notation $n^{\{\delta_j\}}$ means $s_j = 1$ and $s_i = 0$ for all $i \neq j$.  

We have chosen these identification numbers to be consistent with those given to $G$-irreducible subgroups of type $A_1$ in \cite[Tables~4--8]{tho2}. After the first $m$ identification numbers have accounted for the subgroups of type $A_1$ we then give the next set of identification numbers to the simple subgroups of rank at least 2 from \cite[Tables~3--7]{tho1}. The next set is given to the remaining non-simple reductive, maximal connected subgroups and after this they 
are given in order. 

In Tables \ref{G2tab}--\ref{E8tab} we also need a notation to be able to describe conjugacy classes of $M$-irreducible connected subgroups $X$ of reductive, maximal connected subgroups $M$ of $G$. Suppose that $M = M_1 M_2 \ldots M_r$. If all of the factors are simple classical algebraic groups then we define $$V_M := V_{M_1}(\lambda_1) \otimes V_{M_2}(\lambda_1) \otimes \cdots \otimes V_{M_r}(\lambda_1)$$ and let $V_M \downarrow X$ be the usual restriction of the $M$-module $V_M$ to $X$. 

Now suppose that some $M_i$ is of exceptional type. We do not wish to list composition factors for the action of $X$ on a minimal module for $M_i$ as this will make things more difficult to read. Instead, we use the fact that the projection of $X$ to $M_i$ will be $M_i$-irreducible (see Lemma \ref{easy}) and therefore $X$ has a unique ID number. In this case we give the ID number rather than the composition factors of the restriction of any module. So, if $M = M_1$ then we define the notation $V_{M} \downarrow X$ to simply be $M(\#a)$ where $a$ is the ID number of the subgroup $X$ of $M$. The final possibility is $M = M_1 M_2$ with at least one $M_i$ of exceptional type. In this case, we combine the previous two notations by denoting $V_M \downarrow X = (V_{M_1} \downarrow X_1, V_{M_2} \downarrow X_2)$.  

As above, let $G$ be a simple exceptional algebraic group. Many of the $G$-irreducible connected subgroups have simple factors $\bar{X}$ generated by long root subgroups of $G$. In the following cases all subgroups of the given type are generated by long root subgroups of $G$ and we will therefore omit the bar.  
\begin{align*}
G = E_6, \hspace{5pt} & X  = A_{4}, A_{5}, D_{5}, E_{6}. \\
G = E_7, \hspace{5pt} & X  = A_{4}, A_{5}, A_{6}, A_{7}, D_{5}, D_{6}, E_{6}, E_{7}. \\
G = E_8, \hspace{5pt} & X  = A_{5}, A_{6}, A_{7}, A_{8}, D_{5}, D_{6}, D_{7}, D_{8} , E_{6}, E_{7}, E_{8}.
\end{align*} 
Finally, we note the standard notation we will use for certain finite groups. A symmetric group acting on a finite set of size $m$ will be denoted by $S_m$. A dihedral group of order $2n$ will be denoted by $\text{Dih}_{2n}$.  

%% file: Prelims.tex
\chapter{Preliminaries} \label{prelims}

To prove our main theorems we require a number of preliminary results, which we record in this section. The first result is the starting point for our strategy. Recall that when we say a reductive, maximal connected subgroup we mean a reductive subgroup that is maximal among all connected subgroups. 

\begin{thm}[{\cite[Corollary 2]{LS1}}] \label{maximalexcep}
Let $G$ be a simple exceptional algebraic group. Let $M$ be a reductive, maximal connected subgroup of $G$. Then $M$ is $\text{Aut}(G)$-conjugate to precisely one subgroup $X$ as follows, where each isomorphism type of $X$ denotes one $G$-conjugacy class of subgroups.

\begin{longtable}{>{\raggedright\arraybackslash}p{0.05\textwidth - 1\tabcolsep}>{\raggedright\arraybackslash}p{0.95\textwidth - 2\tabcolsep}}
\\ \hline $G$ & $X$ \\ \hline
$G_2$ & $\bar{A}_2$, $\tilde{A}_2$ $(p =3)$, $\bar{A}_1 \tilde{A}_1$, $A_1$ $(p \geq 7)$ \\
$F_4$ & $B_4$, $C_4$ $(p=2)$, $\bar{A}_1 C_3$ $(p \neq 2)$, $A_1 G_2$ $(p \neq 2)$, $\bar{A}_2 \tilde{A}_2$, $G_2$ $(p = 7)$, $A_1$ $(p \geq 13)$ \\
$E_6$ & $\bar{A}_1 A_5$, $\bar{A}_2^3$, $F_4$, $C_4$ $(p \neq 2)$, $A_2 G_2$, $G_2$ $(p \neq 7$), $A_2$ $(p \geq 5)$ \\
$E_7$ & $\bar{A}_1 D_6$, $\bar{A}_2 A_5$, $A_7$, $G_2 C_3$, $A_1 F_4$, $A_1 G_2$ $(p \neq 2)$, $A_2$ $(p \geq 5)$, $A_1 A_1$ $(p \geq 5)$, $A_1$ $(p \geq 17)$, $A_1$ $(p \geq 19)$ \\
$E_8$ & $D_8$, $\bar{A}_1 E_7$, $\bar{A}_2 E_6$, $A_8$, $\bar{A}_4^2$, $G_2 F_4$, $B_2$ $(p \geq 5)$, $A_1 A_2$ $(p \geq 5)$, $A_1$ $(p \geq 23)$, $A_1$ $(p \geq 29)$, $A_1$ $(p \geq 31)$ \\ \hline
\end{longtable}
\end{thm}

Let $G$ be a classical simple algebraic group, which we refer to as $\text{Cl}(V) = \text{SL}(V)$, $\text{Sp}(V)$ or $\text{SO}(V)$ for some finite-dimensional vector space $V$. We need to determine the $G$-conjugacy classes of reductive, maximal connected subgroups of $G$, when $G$ is of small rank. Firstly, we need part of a theorem of Liebeck and Seitz concerning the maximal subgroups of classical algebraic groups. Let $H$ be a subgroup of $G$. We introduce the following classes of subgroups (which is a subset of those from \cite{LS2}).

\vspace{0.1cm}
{\bf Class} $\mathcal{C}_1$: \emph{Subspace stabilisers}. Here $H \in \mathcal{C}_1$ if $H = \text{Stab}_G(W)$ where $W$ is either a non-degenerate subspace of $V$, or $(G,p) = (\text{SO}(V),2)$ and $W$ is a non-singular subspace of dimension 1.

\vspace{0.1cm}
{\bf Class} $\mathcal{C}_4$: \emph{Tensor product subgroups}. Suppose that $V=V_1 \otimes V_2$ with $\text{dim} V_i > 1$. Then $H \in \mathcal{C}_4$ if $H = \text{Cl}(V_1) \circ \text{Cl}(V_2)$ which acts naturally on $V$ as follows: $(g_1,g_2) (v_1 \otimes v_2) := (g_1 v_1) \otimes (g_2 v_2)$. The tensor product subgroups occurring are: \begin{align*} \text{SL} \otimes \text{SL} < \text{SL}, \hspace{1cm} & \hspace{1cm} \text{Sp} \otimes \text{SO} < \text{Sp} \hspace{0.1cm} (p \neq 2), \\ \text{Sp} \otimes \text{Sp} < \text{SO}, \hspace{1cm} & \hspace{1cm} \text{SO} \otimes \text{SO} < \text{SO} \hspace{0.1cm} (p \neq 2). \end{align*} 

The following theorem can be immediately deduced from {\cite[Theorem 1]{LS2}}. 

\begin{thm} \label{maximal}
Let $G$ be a classical simple algebraic group. Suppose that $M$ is a reductive, maximal connected subgroup of $G$. Then one of the following holds:
\begin{enumerate}[label=\normalfont(\roman*),leftmargin=*,widest=iii, align=left]
\item \label{C1lab} $M$ belongs to $\mathcal{C}_1$;
\item \label{C4lab} $M$ belongs to $\mathcal{C}_4$;
\item \label{simplelab} $M$ is a simple algebraic group and $V \downarrow M$ is irreducible and restricted.
\end{enumerate}
\end{thm} 
\vspace{0.1cm}
In the following lemma we now apply this theorem to certain classical simple groups of small rank; the ones we treat are those which will arise in the work of Sections 5--9. 
\vspace{0.1cm}
\begin{lem} \label{maxclassical}
Suppose that $G$ is a classical group of type $A_n$ for $n \in \{2,3,4,5,7,8\}$, $B_n$ for $n \in \{2,3,4,5,6,7\}$, $C_n$ for $n \in \{3,4\}$, or  $D_n$ for $n \in \{4,5,6,8\}$. Then the following table gives all $G$-conjugacy classes of reductive, maximal connected subgroups of $G$. 
\end{lem}

\begin{longtable}{>{\raggedright\arraybackslash}p{0.06\textwidth - 2\tabcolsep}>{\raggedright\arraybackslash}p{0.19\textwidth - 2\tabcolsep}>{\raggedright\arraybackslash}p{0.4\textwidth - 2\tabcolsep}>{\raggedright\arraybackslash}p{0.35\textwidth-\tabcolsep}@{}}

\caption{The maximal subgroups of certain low rank classical algebraic groups. \label{maxclassicaltab}} \\

\hline

$G$ & Max. sub. $M$ & $V_{G}(\lambda_1) \downarrow M$ & Comments\\ 

\hline

$A_2$ & $A_1$ $(p \neq 2)$ & $2$  \\

\hline

$B_2$ & $\bar{A}_1^2$ & $(1,1) + (0,0)$  \\

& $\tilde{A}_1^2$ $(p=2)$ & $(2,0) + (0,2)$  \\ 

& $A_1$ $(p \geq 5)$ & $4$ \\

\hline
\endfirsthead
$A_3$ & $B_2$ & $01$  \\

& $A_1^2$ $(p \neq 2)$ & $(1,1)$  \\

\hline

$B_3$ & $\bar{A}_3$ & \multicolumn{2}{l}{$010 + 000$ $(p \neq 2)$ or $010$ $(p =2)$} \\

& $\bar{A}_1^2 \tilde{A}_1$ $(p \neq 2)$ & $(1,1,0) + (0,0,2)$  \\

& $G_2$ & $10$  \\

& $\tilde{A}_1 B_2$ $(p=2)$ & $(2,00) + (0,10)$  \\

\hline

$C_3$ & $\bar{A}_1 C_2$ & $(1,00) + (0,10)$  \\

& $A_1 A_1$ $(p \neq 2)$ & $(2,1)$  \\

& $\tilde{A}_3$ $(p=2)$ & $010$  \\

& $G_2$ $(p=2)$ & $10$  \\

& $A_1$ $(p \geq 7)$ & $5$  \\

\hline

$A_4$ & $B_2$ & $10$  \\

\hline

$B_4$ & $\bar{D}_4$ & $\lambda_1 + 0$ $(p \neq 2)$ or $\lambda_1$ $(p=2)$ & \\

& $\tilde{A}_1 \bar{A}_3$ $(p \neq 2)$ & $(2,000) + (0,010)$ & \\

& $\bar{A}_1^2 B_2$ $(p \neq 2)$ & $(1,1,00) + (0,0,10)$  & \\

& $A_1^2$ $(p \neq 2)$ & $(2,2)$ &  \\

& $A_1$ $(p \geq 11)$ & $8$ & \\

& $\tilde{A}_1 B_3$ $(p=2)$ & $(2,000) + (0,100)$ & \\

& $B_2^2$ $(p=2)$ & $(10,00) + (00,10)$ & \\

\hline

$C_4$ & $C_2^2$ & $(10,00) + (00,10)$  \\

& $\bar{A}_1 C_3$ & $(1,000) + (0,100)$  \\

& $A_1^3$ $(p \neq 2)$ & $(1,1,1)$  \\

& $A_1$ $(p \geq 11)$ & $7$  \\

& $\tilde{D}_4$ $(p=2)$ & $\lambda_1$  \\

\hline

$D_4$ & $B_3$ (3 classes) & $100 + 000$ $(p \neq 2)$ or $000|100|000$ $(p =2)$ \newline $001$& classes are permuted by $\text{Out}(D_4) \cong S_3$  \\

& & $001$  \\

& $A_1 B_2$ $(p \neq 2)$ \newline (3 classes) & $(2,00) + (0,10)$ \newline $(1,01)$  & classes are permuted by $\text{Out}(D_4) \cong S_3$   \\

& & $(1,01)$   \\

& $\bar{A}_1^4$ & $(1,1,0,0) + (0,0,1,1)$  \\

& $A_2$ $(p \neq 3)$ \newline & $11$ & involutory graph aut. of $D_4$ induces graph aut. of $A_2$ \\

\hline

$A_5$ & $A_1 A_2$ & $(1,10)$ \\

& $C_3$ & $100$ \\

& $A_3$ $(p \neq 2)$ & $010$ \\

& $A_2$ $(p \neq 2)$ & $20$ \\

\hline

$B_5$ & $\bar{D}_5$ & $\lambda_1 + 0$ $(p \neq 2)$ or $\lambda_1$ $(p = 2)$ & \\

& $\tilde{A}_1 \bar{D}_4$ $(p \neq 2)$ & $(2,0) + (0,\lambda_1)$ & \\

& $\bar{A}_1^2 B_3$ $(p \neq 2)$ &$(1,1,000) + (0,0,100)$ & \\

& $B_2 \bar{A}_3$ $(p \neq 2)$ & $(10,000) + (00,010)$ & \\

& $A_1$ $(p \geq 11)$ & $10$ & \\

& $\tilde{A}_1 B_4$ $(p =2)$ & $(2,0) + (0,\lambda_1)$ & \\

& $B_2 B_3$ $(p =2)$ & $(10,000) + (00,100)$ & \\

\hline

$D_5$ & $B_4$ & $\lambda_1 + 0$ $(p \neq 2)$ or $0 | \lambda_1 | 0$ $(p=2)$ &  \\

& $A_1 B_3$ $(p \neq 2)$ & $(2,000) + (0,100)$ & \\

& $\bar{A}_1^2 \bar{A}_3$ & $(1,1,000) + (0,0,010)$ & \\

& $B_2^2$ $(p \neq 2)$ & $(10,00) + (00,10)$ & \\

& $B_2$ $(p \neq 2)$ \newline  (2 classes) & $02$ \newline $02$ & classes permuted by involutory graph aut. \\

\hline

$B_6$ & $\bar{D}_6$ & $\lambda_1 + 0$ $(p \neq 2)$ or $\lambda_1$ $(p = 2)$ & \\

& $\tilde{A}_1 \bar{D}_5$ $(p \neq 2)$ & $(2,0) + (0,\lambda_1)$ & \\

& $\bar{A}_1^2 B_4$ $(p \neq 2)$ &$(1,1,0) + (0,0,\lambda_1)$ & \\

& $B_2 \bar{D}_4$ $(p \neq 2)$ & $(10,0) + (00,\lambda_1)$ & \\

& $\bar{A}_3 B_3$ $(p \neq 2)$ & $(010,000) + (000,100)$ & \\

& $A_1$ $(p \geq 13)$ & $12$ & \\

& $B_2$ $(p=5)$ & $20$ & \\

& $C_3$ $(p=3)$ & $010$ & \\

& $\tilde{A}_1 B_5$ $(p =2)$ & $(2,0) + (0,\lambda_1)$ & \\

& $B_2 B_4$ $(p =2)$ & $(10,0) + (00,\lambda_1)$ & \\

&$B_3^2$ $(p=2)$ &  $(100,000) + (000,100)$ & \\

\hline

$D_6$ & $B_5$ & $\lambda_1 + 0$ $(p \neq 2)$ or $0 | \lambda_1 | 0$ $(p=2)$ &  \\

& $A_1 B_4$ $(p \neq 2)$ & $(2,0) + (0,\lambda_1)$ & \\

& $\bar{A}_1^2 \bar{D}_4$ & $(1,1,0) + (0,0,\lambda_1)$ & \\

& $B_2 B_3$ $(p \neq 2)$ & $(10,000) + (00,100)$ & \\

& $\bar{A}_3^2$ & $(010,000) + (000,010)$ & \\

& $A_1 C_3$ (2 classes) & $(1,100)$ \newline $(1,100)$ & classes permuted by involutory graph aut. \\

\hline

$A_7$ & $C_4$ & $\lambda_1$ & \\

& $D_4$ $(p \neq 2)$ & $\lambda_1$ & \\

& $A_1 A_3$ & $(1,100)$ & \\

\hline

$B_7$ & $\bar{D}_7$ & $\lambda_1 + 0$ $(p \neq 2)$ or $\lambda_1$ $(p = 2)$ & \\

& $\tilde{A}_1 \bar{D}_6$ $(p \neq 2)$ & $(2,0) + (0,\lambda_1)$ & \\

& $\bar{A}_1^2 B_5$ $(p \neq 2)$ &$(1,1,0) + (0,0,\lambda_1)$ & \\

& $B_2 \bar{D}_5$ $(p \neq 2)$ & $(10,0) + (00,\lambda_1)$ & \\

& $\bar{A}_3 B_4$ $(p \neq 2)$ & $(010,0) + (000,\lambda_1)$ & \\

& $B_3 \bar{D}_4$ $(p \neq 2)$ & $(010,0) + (000,\lambda_1)$ & \\

& $A_3$ $(p \neq 2)$ & $101$ & \\

& $A_1 B_2$ $(p \neq 2)$ & $(2,10)$ & \\

& $A_1$ $(p \geq 17)$ & $14$ & \\

& $\tilde{A}_1 B_6$ $(p =2)$ & $(2,0) + (0,\lambda_1)$ & \\

& $B_2 B_5$ $(p =2)$ & $(10,0) + (00,\lambda_1)$ & \\

& $B_3 B_4$ $(p=2)$ &  $(100,0) + (000,\lambda_1)$ & \\

\hline

$A_8$ & $B_4$ $(p \neq 2)$ & $\lambda_1$ & \\

& $A_2^2$ & $(10,10)$ & \\

\hline

$D_8$ & $B_7$ & $\lambda_1 + 0$ $(p \neq 2)$ or $0 | \lambda_1 | 0$ $(p=2)$ &  \\

& $A_1 B_6$ $(p \neq 2)$ & $(2,0) + (0,\lambda_1)$ & \\

& $\bar{A}_1^2 \bar{D}_6$ & $(1,1,0) + (0,0,\lambda_1)$ & \\

& $B_2 B_5$ $(p \neq 2)$ & $(10,0) + (00,\lambda_1)$ & \\

& $\bar{A}_3 \bar{D}_5$ & $(010,000) + (000,\lambda_1)$ & \\

& $B_3 B_4$ $(p \neq 2)$ & $(100,0) + (000,\lambda_1) $ & \\

& $\bar{D}_4^2$ & $(\lambda_1,0) + (0,\lambda_1)$ & \\

& $B_2^2$ ($p \neq 2$) \newline (2 classes) & $(01,01)$ \newline $(01,01)$ & classes permuted by involutory graph aut. \\

& $A_1 C_4$ \newline (2 classes) & $(1,\lambda_1)$ \newline $(1,\lambda_1)$ & classes permuted by involutory graph aut. \\

& $B_4$ (2 classes) & $\lambda_4$ \newline $\lambda_4$ & classes permuted by involutory graph aut. \\

\hline

\end{longtable}

\begin{proof}

We will give the details of how to apply Theorem \ref{maximal} when $G$ is of type $A_n$ ($2 \leq n \leq 5$) or of type $D_4$. The other types are similar. The strategy is to find all possible subgroups of $G$ in $\mathcal{C}_1$, $\mathcal{C}_4$ and \ref{maximal} \ref{simplelab}, and then check whether there are any containments amongst them.  

Suppose that $G$ is of type $A_n$. We apply Theorem \ref{maximal}, considering $A_n$ as $\text{SL}(V)$ where $V$ is of dimension $n+1$ and equipped with the $0$-form. Firstly, there are no subgroups in $\mathcal{C}_1$ since all subspaces are degenerate. Next, we consider $\mathcal{C}_4$. If $n+1$ is prime then there are no subgroups. For $G$ of type $A_3$ we obtain the subgroup $A_1^2$ acting as $(1,1)$ on $V$. This is only maximal when $p \neq 2$ because this subgroup is contained in $B_2$ when $p=2$ (the subgroup $B_2$ comes from \ref{maximal} \ref{simplelab} as explained below). For $G$ of type $A_5$ we obtain the subgroup $A_1 A_2$ acting as $(1,2)$. 

Now consider subgroups from \ref{maximal} \ref{simplelab}. We use \cite{lubeck} to find the simple groups $S$ with an $(n+1)$-dimensional irreducible restricted representation. Any such representation will embed $S$ into $G$ so it remains to determine the maximality. When $G$ is of type $A_2$ the only possibility for $S$ is a subgroup $A_1$ when $p \neq 2$ acting as the symmetric square of its natural representation. Since there are no other possible subgroups it follows that $S$ is maximal. Now let $G$ be of type $A_3$. The possibilities for $S$ are a subgroup $B_2$ acting as $V_{B_2}(01)$ and a subgroup $A_1$ $(p \geq 5)$ acting as $V_{A_1}(4)$. We claim the subgroup $A_1$ is contained in $B_2$. This follows since for any group of type $A_n$ ($n \geq 3$), the subgroup $A_1$ acting as $V_{A_1}(n)$ is contained in $C_{(n+1)/2}$ ($n$ odd) or $B_{n/2}$ ($n$ even) because it fixes a symplectic ($n$ odd) or orthogonal ($n$ even) form. Therefore only the subgroup $B_2$ is maximal. 

Let $G$ be of type $A_4$. Then the only possibility for $S$ is a subgroup $B_2$ acting as $V_{B_2}(10)$. Since there are no other possible subgroups this $B_2$ is maximal. Finally, let $G$ be of type $A_5$. Then the possibilities for $S$ are $C_3$ acting via $V_{C_3}(100)$, $A_3$ acting via $V_{A_3}(010)$, $A_2$ $(p \neq 2)$ acting via $V_{A_2}(20)$ and $A_1$ $(p \geq 7)$. The subgroup $A_1$ $(p \geq 7)$ is contained in $C_3$ by the argument in the previous paragraph. Also, the subgroup $A_3$ is contained in $C_3$ when $p=2$ ($\text{SO}_6 < \text{Sp}_6$). There are no further containments. 

Now let $G$ be of type $D_4$, which we consider as $\text{SO}_8$ with natural module $V$. The subgroups in $\mathcal{C}_1$ are $B_3$ ($\text{SO}_7$), $A_1 B_2$ $(p \neq 2)$ ($\text{SO}_3 \text{SO}_5$) and $A_1^4$ $(\text{SO}_4 \text{SO}_4 < \text{SO}_8)$. The tensor decomposition $V = V_4 \otimes V_2$ (with $V_i$ of dimension $i$) gives two $D_4$-conjugacy classes of maximal subgroups $A_1 B_2$ ($\text{Sp}_2 \otimes \text{Sp}_4$) when $p \neq 2$. Indeed, there is only one $\text{GO}_8$-conjugacy class but this splits in the index two subgroup $\text{SO}_8$ because any element $g \in \text{GO}_8 \setminus \text{SO}_8$ normalising $A_1 B_2$ centralises it (since the outer automorphism group of $A_1 B_2$ is trivial) and is hence contained in the centre of $\text{GO}_8$. However, the centre of $\text{GO}_8$ is contained in $\text{SO}_8$ so no such $g$ exists and $N_{\text{GO}_8}(A_1 B_2) < \text{SO}_8$. It follows that there is an involution in the automorphism group of $D_4$ which swaps the two $\text{SO}_8$-conjugacy classes. Furthermore, they are both mapped to the subgroup $A_1 B_2$ from $\mathcal{C}_1$ by appropriate triality automorphisms since neither is centralised by any triality automorphisms (\cite[Table~4.7.1]{gls3}). 

Now we consider simple groups with an $8$-dimensional irreducible restricted representation from \cite{lubeck}. The group $A_2$ has such a representation when $p \neq 3$, namely the adjoint module $V_{A_2}(11)$. Since it is the adjoint module for $A_2$, we know $A_2$ preserves the Killing form on it. The Killing form is a bilinear symmetric form and therefore $V_{A_2}(11)$ is an orthogonal module and $A_2 < D_4$ when $p>3$. When $p=2$ the proof of \cite[Proposition 2.3.3]{K} shows that $A_2$ preserves a quadratic form on $V_{A_2}(11)$ and thus $A_2$ is a subgroup of $D_4$. Moreover, there is an element of $\text{GO}_8 \setminus \text{SO}_8$ acting as a graph automorphism of $A_2$, again by the proof of \cite[Proposition 2.3.3]{K}. Therefore, $N_{\text{GO}_8}(A_2) \not < \text{SO}_8$ and there is just one $D_4$-conjugacy class of maximal subgroup $A_2$. We also note that we have shown that $A_2.2 < D_4.2$ but $A_2.2 \not < D_4$. The module $V_{A_1}(7)$ is a self-dual, restricted, 8-dimensional representation for $A_1$ when $p \geq 11$. However, this module is symplectic and therefore there is no maximal subgroup $A_1$ of $D_4$. The group $B_3$ has a self-dual, restricted, 8-dimensional, representation $V_{B_3}(001)$. This is indeed an orthogonal representation by \cite[Proposition 5.4.9]{KL} (the argument holds for algebraic groups also) and hence yields a subgroup $B_3$ of $D_4$. Since the outer automorphism group of $B_3$ is trivial, it follows as before that there are two $D_4$-conjugacy classes and they are swapped by an involution in the automorphism group of $D_4$; they are mapped to the subgroup $B_3$ from $\mathcal{C}_1$ under triality automorphisms. Finally, when $p=2$ the group $C_3$ has an 8-dimensional self-dual irreducible representation obtained by the special isogeny map from $C_3$ to $B_3$. Therefore this module is also orthogonal, but the image of this $\text{Sp}_6$ inside $D_4$ is a $B_3$, otherwise a triality automorphism conjugates a subgroup $C_3$ to one of type $B_3$, which is absurd since they are not even isomorphic as algebraic groups.    
\end{proof}

For the following lemmas let $G$ be a semisimple connected algebraic group. We state some elementary results about $G$-irreducible subgroups. 

\begin{lem} [{\cite[Lemma 2.1]{LT}}] \label{semirr}
If $X$ is a $G$-irreducible connected subgroup of $G$, then $X$ is semisimple and $C_G(X)$ is finite.   
\end{lem}

\begin{lem} [{\cite[Lemma 3.6]{tho1}}] \label{easy}
Suppose that a $G$-irreducible subgroup $X$ is contained in $K_1 K_2$, a commuting product of connected non-trivial subgroups $K_1$, $K_2$ of $G$. Then $X$ has non-trivial projection to both $K_1$ and $K_2$. Moreover, each projection is a $K_i$-irreducible subgroup. 
\end{lem}

\begin{lem}[{\cite[Lemma 2.2]{LT}}] \label{class}
Suppose that $G$ is a classical simple algebraic group, with natural module $V = V_G(\lambda_1)$. Let $X$ be a semisimple connected subgroup of $G$. If $X$ is $G$-irreducible then one of the following holds:

\begin{enumerate}[label=\normalfont(\roman*),leftmargin=*,widest=iii, align=left]

\item $G = A_n$ and $X$ is irreducible on $V$; 

\item $G = B_n, C_n$ or $D_n$ and $V \downarrow X = V_1 \perp \ldots \perp V_k$ with the $V_i$ all non-degenerate, irreducible and inequivalent as $X$-modules;

\item $G = D_n$, $p=2$, $X$ fixes a non-singular vector $v \in V$, and $X$ is a $G_v$-irreducible subgroup of $G_v = B_{n-1}$.
\end{enumerate}
\end{lem}

The next lemma and corollary are used in the proofs of Theorems \ref{G2THM} to \ref{E8THM} to show that an $M$-irreducible subgroup is $G$-irreducible, where $M$ is a reductive, maximal connected subgroup of a simple exceptional algebraic group $G$. 

\begin{lem}  [{\cite[Lemma 3.8]{tho1}}] \label{wrongcomps}
Let $X$ be a semisimple connected subgroup of $G$ and let $V$ be a $G$-module. Suppose that $X$ does not have the same composition factors as any semisimple $L'$-irreducible connected subgroup $H$ of the same type as $X$ for some proper Levi subgroup $L$ of $G$. If $X$ is of type $B_n$ and $p=2$ then assume further that there is no subgroup $H$ of type $C_n$ with $H \leq L'$ and $L$-irreducible, for some Levi subgroup $L$ of $G$, such that there is an isogeny $\phi: X \rightarrow H$ inducing a mapping which takes the composition factors of $V \downarrow X$ to those of $V \downarrow H$. Then $X$ is $G$-irreducible.  
\end{lem}

\begin{cor}  [{\cite[Corollary 3.9]{tho1}}] \label{notrivs}
Suppose that $X < G$ is semisimple and $L(G) \downarrow X$ has no trivial composition factors. Then $X$ is $G$-irreducible. 
\end{cor}

The following two well known results will be used throughout the proof of Theorem \ref{MAINTHM} when showing that a $G$-irreducible subgroup contained in a maximal subgroup $M_1$ is conjugate to a $G$-irreducible subgroup of a maximal subgroup $M_2$. 

\begin{lem} [{\cite[Lemma 11.13]{LS9}}] \label{centralisers}
Let $G$ be an adjoint simple algebraic group of exceptional type and let $s \in G$ be an element of prime order $r \neq p$. Then $C_{G}(s)^\circ$ is semisimple if and only if $C_{G}(s)$ is listed below. 
\end{lem}

\setlength{\LTleft}{0.18\textwidth}
\begin{longtable}{>{\raggedright\arraybackslash}p{0.06\textwidth - 2\tabcolsep}>{\raggedright\arraybackslash}p{0.36\textwidth - 2\tabcolsep}>{\raggedright\arraybackslash}p{0.24\textwidth - 1\tabcolsep}}
\\ \hline $G$ & $C_G(s)$ & $r$ \\ \hline
$G_2$ & $\bar{A}_1 \tilde{A}_1$, $\bar{A}_2$ & 2, 3 (resp.) \\
$F_4$ & $\bar{A}_1 C_3$, $B_4$, $\bar{A}_2 \tilde{A}_2$ & 2, 2, 3 (resp.) \\
$E_6$ & $\bar{A}_1 A_5$, $\bar{A}_2^3 . 3$ & 2, 3 (resp.)  \\
$E_7$ & $\bar{A}_1 D_6$, $A_7 . 2$, $\bar{A}_2 A_5$ & 2, 2, 3 (resp.) \\
$E_8$ & $\bar{A}_1 E_7$, $D_8$, $\bar{A}_2 E_6$, $A_8$, $\bar{A}_4^2$ & 2, 2, 3, 3, 5 (resp.) \\ \hline
\end{longtable} 

\begin{lem} [{\cite[Section 4]{LS6}}] \label{rootycentralisers}
Let $G$ be a simple algebraic group of exceptional type and $X$ be a subgroup of type $A_1$ or $A_2$. If $X$ is generated by long root subgroups of $G$ then $C_G(X)^\circ$ is given in the table below.
\end{lem}

\setlength{\LTleft}{0.35\textwidth}
\begin{longtable}{>{\raggedright\arraybackslash}p{0.06\textwidth - 2\tabcolsep}>{\raggedright\arraybackslash}p{0.12\textwidth - 2\tabcolsep}>{\raggedright\arraybackslash}p{0.12\textwidth - 1\tabcolsep}}
\\ \hline $G$ & $C_G(\bar{A}_1)^\circ$ & $C_G(\bar{A}_2)^\circ$ \\ \hline
$G_2$ & $\tilde{A}_1$ & $1$ \\
$F_4$ & $C_3$ & $\tilde{A}_2$ \\
$E_6$ & $A_5$ & $\bar{A}_2^2$  \\
$E_7$ & $D_6$ & $A_5$ \\
$E_8$ & $E_7$ & $E_6$ \\ \hline
\end{longtable} 

%% file: Strat.tex
\chapter{Strategy for the proofs of Theorems \ref{G2THM}--\ref{E8THM}} \label{strat}

To prove Theorem \ref{MAINTHM} we prove Theorems \ref{G2THM}--\ref{E8THM} in Sections \ref{secG2}--\ref{secE8}, respectively. In this section we describe the strategy used in proving those theorems.

Let $G$ be a simple exceptional algebraic group over an algebraically closed field of characteristic $p$. Suppose that $X$ is a $G$-irreducible connected subgroup of $G$. Then $X$ is contained in a maximal connected subgroup $M$ of $G$. Since $X$ is $G$-irreducible, $M$ is reductive. Furthermore, $X$ is $M$-irreducible as any parabolic subgroup of $M$ is contained in a parabolic subgroup of $G$ by the Borel-Tits Theorem \cite[Th\'{e}or\`{e}me 2.5]{BT}. Therefore, $X$ is $M$-irreducible in some reductive, maximal connected subgroup $M$ of $G$ and the following strategy will find all such $X$. Not only do we wish to classify all $G$-irreducible connected subgroups of $G$, we also want to describe the lattice structure of connected overgroups. This last point explains why the strategy we describe below is different to that used in \cite{tho1} and \cite{tho2}. 

Take the first reductive, maximal connected subgroup $M$ from Theorem \ref{maximalexcep} (the ordering is chosen to make the proof and resulting tables easier to follow; of course, one could use any ordering). We iterate the following schematic process, explaining each stage below. Find all reductive, maximal connected subgroups of $M$; these subgroups will be denoted by $M_1$ in each case. All such subgroups will clearly be $M$-irreducible. Now find all of the reductive, maximal connected subgroups of each $M_1$, in turn. At this point it is not necessarily the case that all such subgroups will be $M$-irreducible and we only wish to keep the irreducible ones. Continuing in this way will lead to a set containing every $M$-irreducible subgroup; we want to consider the $G$-conjugacy classes of these subgroups and so we specify when a class is repeated in our list. Moreover, for a given representative of a $G$-conjugacy class of irreducible subgroups we know all of its connected overgroups and thus understand the lattice structure of the irreducible connected subgroups of $M$. 

To find the reductive, maximal connected subgroups of a subgroup of $M$ we repeatedly use Theorem \ref{maximalexcep} and Lemma \ref{maxclassical}. To check whether a subgroup $X$ of $M$ is $M$-irreducible we first note that by Lemma \ref{easy}, the projection of $X$ to each simple factor of $M$ must be irreducible. We then use Lemma \ref{class} for classical factors of $M$; if $M$ has a factor of exceptional type then this will have rank less than that of $G$ and so we may use Theorem \ref{MAINTHM} inductively, since we will prove this in ascending order of rank. We also need to consider $G$-conjugacy rather than just $M$-conjugacy of all of the $M$-irreducible connected subgroups. There are a large number of cases where $G$ does not fuse any $M$-classes together. Given two $M$-conjugacy classes of irreducible connected subgroups, with representatives $X_1$ and $X_2$, say, it is easy to check that $X_1$ is not $G$-conjugate to $X_2$ by comparing the composition factors of their actions on the adjoint module. If the composition factors are the same then we will give an argument showing either that $X_1$ is $G$-conjugate to $X_2$, or that we are in the case of the exception given in Corollary \ref{comps} with $G = E_8$, $p=3$, and $X_1$ and $X_2$ are of type $A_2$. Many of the arguments are based on information from \cite{car} on $\text{Out}_G(H)$ for certain maximal rank subgroups $H$ of $G$. We give each $G$-conjugacy class of $M$-irreducible connected subgroups a unique identification number $n$, and the class (or representative subgroup) is denoted by $G(\#n)$ (see Section \ref{nota}). 

We now consider the second reductive, maximal connected subgroup of $G$ from Theorem \ref{maximalexcep}. In the proofs and tables, this will also be denoted by $M$, but to be clear in our explanation in this paragraph we denote it by $N$ so as to be able to differentiate between the first and second reductive, maximal connected subgroups. We use the same method as above to find the lattice structure of $G$-conjugacy classes of $N$-irreducible connected subgroups but with an important extra step. For each $N$-irreducible connected subgroup $X$ we need to determine whether $X$ is conjugate to an $M$-irreducible connected subgroup of $M$. To show that this does not happen is straightforward. Indeed, we can compare the composition factors of $X$ on $L(G)$ with those of each $M$-irreducible subgroup of the same type as $X$. This routine but tedious check will not be explicitly mentioned. When $X$ has the same composition factors as an $M$-irreducible subgroup $Y$ then we need to prove that $X$ is conjugate to $Y$, since we also wish to prove Corollary \ref{comps} (when $X$ is simple this has already been proved in \cite[Corollary 1]{tho1}, \cite[Corollary 3]{tho2}). Proving that $X$ and $Y$ are conjugate requires results from Section \ref{prelims} such as Lemmas \ref{centralisers}, \ref{rootycentralisers} and some ad-hoc methods. The main task is showing that $X$ is contained in some conjugate of $M$. The fact that $X$ and $Y$ are conjugate then follows from the fact that composition factors on $L(G)$ determine conjugacy of $M$-irreducible subgroups. For example, suppose $G$ is of type $E_7$. If a subgroup $X$ of $N$ has an $\bar{A}_1$ factor generated by long root subgroups of $G$ then $X$ will be contained in some conjugate of $M = \bar{A}_1 D_6$ by Lemma \ref{rootycentralisers}. For a specific example, consider $N = G_2 C_3$ and $X = \bar{A}_1 A_1 C_3 = E_7(\#\eseven{42})$. 

We repeat this process until we have considered all of the reductive, maximal connected subgroups from Theorem \ref{maximalexcep}. The information obtained from this method is displayed in Tables \ref{G2tab}--\ref{E8tab}, except for the $M$-irreducible subgroups which are $G$-reducible. We explained the notation used in these tables in Section~\ref{nota}, and we explain how to read the tables at the start of Section~\ref{thmtabs}. In particular, when the identification number $n$ for a subgroup $X$ is written in italics it means that $X$ is listed elsewhere in the table and so should be discounted if one wants exactly one conjugate of each $G$-irreducible connected subgroup, as in Theorem~\ref{MAINTHM}. 

We now describe the final part of the strategy: determining whether an $M$-irreducible connected subgroup is $G$-irreducible. Let $X$ be an $M$-irreducible connected subgroup of $G$. If $X$ is simple then we use \cite{tho1} and \cite{tho2} to check if $X$ is irreducible. If $X$ contains a simple $G$-irreducible subgroup then $X$ is of course itself $G$-irreducible. In both of these cases we will not explicitly mention this in the proofs. Now suppose that $X$ is non-simple and contains no simple $G$-irreducible subgroup. Then to prove that $X$ is $G$-irreducible we use  Lemma \ref{wrongcomps} and Corollary \ref{notrivs}, for which we require the composition factors for the action of $X$ on the minimal or adjoint module. These can be found by restricting the composition factors of $M$ to $X$ and the composition factors for the $G$-irreducible ones can be found in Section \ref{tab:compositionfactors}. To apply Lemma \ref{wrongcomps} we also need the composition factors for the action of the Levi subgroups of $G$ on the minimal and adjoint modules. These can be found in Section \ref{levicomps}. In most cases, a non-simple $M$-irreducible connected subgroup is $G$-irreducible. Corollary \ref{nongcr} lists the connected subgroups $X$ (simple or not) for which $X$ is irreducible in every reductive, maximal connected overgroup, yet $G$-reducible. To prove that $X$ is $G$-reducible requires different methods and we explain these as and when we use them. 

We modify the above approach when considering diagonal subgroups. Suppose that $X$ is an $M$-irreducible connected subgroup of $M$ of the form $A^n B$ for some $n \geq 2$ with $A$ and $B$ of different types. Then $X$ has maximal diagonal subgroups of the form $A^{n-1} B$. We do not want to then consider all of their maximal subgroups including those of the form $A^{n-2} B$. Instead, we list all diagonal subgroups of $X$ immediately and do not list their subgroups in the table. Doing this significantly reduces the size of the tables, and it does not mean that we miss any $M$-irreducible subgroups. It does however mean that some additional combinatorial work is required to recover the lattice of overgroups of certain diagonal subgroups. In light of the large number of diagonal subgroups, we often produce a complete set of non-conjugate classes of diagonal subgroups of an irreducible subgroup $X$ in a supplementary table. This allows for easier reading of the tables referenced in Theorem \ref{MAINTHM}.  

There is another situation where we deviate slightly from the above approach, but only to alter the order in which we study the subgroups in question. Let us explain this with an example. Suppose that $M = \text{SO}_8$ of type $D_4$. When $p \neq 2$ there is a maximal connected subgroup $\text{SO}_5 \text{SO}_3$. However, when $p=2$ this subgroup still exists and is still $\text{SO}_8$-irreducible but is now contained in $\text{SO}_7$. In view of this we use the same identification number $n$ for both groups but use $n$a for the subgroup when $p \neq 2$ and $n$b for the subgroup when $p=2$. Moreover, we wish to study the subgroups of $n$a and $n$b together. In the tables, if we first arrive at the subgroup labelled by $n$b we postpone listing its subgroups until we come to the subgroup labelled by $n$a. 

This generalises to any situation where a subgroup $X$ occurs for $p \neq m$ somewhere in the lattice and elsewhere when $p = m$. In particular, we use this only if the composition factors on both the minimal and adjoint module for $n$b are the same as those for $n$a when $p = m$.

%% file: G2.tex
\chapter{Irreducible subgroups of $G_2$} \label{secG2}

In this section we deduce Theorem \ref{MAINTHM} when $G$ is of type $G_2$. It follows immediately from Theorem \ref{maximalexcep} that the only non-simple $G_2$-irreducible connected subgroup is the maximal subgroup $A_1 \tilde{A}_1$. The simple $G_2$-irreducible connected subgroups of rank at least $2$ are given in \cite[Lemma 3.3]{tho1}. The $G_2$-irreducible subgroups of type $A_1$ are given in \cite[Theorem 2]{tho2}, where it is also proved that $A_1 < A_2$ embedded via the representation $V_{A_1}(2)$ is conjugate to $A_1 \hookrightarrow A_1 \tilde{A}_1$ via $(1,1)$ when $p \neq 2$. Combining these results we deduce the following theorem, concluding the case where $G$ is of type $G_2$. 

\begin{thm} \label{G2THM}
Let $X$ be a $G_2$-irreducible connected subgroup of $G_2$. Then $X$ is conjugate to exactly one subgroup of Table \ref{G2tab} and each subgroup in Table \ref{G2tab} is $G_2$-irreducible. Moreover, Table \ref{G2tab} gives the lattice structure of the irreducible connected subgroups of $G_2$.   
\end{thm}

%% file: F4.tex

\chapter{Irreducible subgroups of $F_4$} \label{secF4}

In this section we use the strategy described in Section \ref{strat} to prove Theorem \ref{MAINTHM} when $G$ is of type $F_4$. In doing this we prove the following theorem. 

\begin{thm} \label{F4THM}
Let $X$ be an $F_4$-irreducible connected subgroup of $F_4$. Then $X$ is conjugate to exactly one subgroup of Table \ref{F4tab} and each subgroup in Table \ref{F4tab} is $F_4$-irreducible. Moreover, Table \ref{F4tab} gives the lattice structure of the irreducible connected subgroups of $F_4$.   
\end{thm}

As mentioned in Section \ref{strat}, throughout the proof we will use \cite[Theorem 3.4]{tho1} and \cite[Theorem 3]{tho2}, which classify the simple $F_4$-irreducible connected subgroups, without reference. In addition, we will implicitly use the fact that any subgroup that contains an $F_4$-irreducible subgroup is itself $F_4$-irreducible. 

By Theorem~\ref{maximalexcep}, the reductive, maximal connected subgroups of $F_4$ are $B_4$, $C_4$ $(p=2)$, $\bar{A}_1 C_3$ $(p \neq 2)$, $A_1 G_2$ $(p \neq 2)$, $\bar{A}_2 \tilde{A}_2$, $G_2$ $(p=7)$ and $A_1$ $(p \geq 13)$. The maximal connected subgroup $A_1$ contains no proper irreducible connected subgroups and so requires no further consideration. In the following sections we consider each of the remaining reductive, maximal connected subgroups $M$ in turn.  

\section{$M = B_4$ $(F_4(\#\ffour{12}))$}      

Here we will only treat the case $p \neq 2$; the $p=2$ case will follow from Section \ref{F4C4} below, where we consider the subgroups of $C_4$, by applying an exceptional graph automorphism of $F_4$. The reductive, maximal connected subgroups of $B_4$ are given by Lemma \ref{maxclassical} and the possibilities are $\bar{D}_4$, $\tilde{A}_1 \bar{A}_3$, $\bar{A}_1^2 B_2$, $A_1^2$ and $A_1$ $(p \geq 11)$. There are no proper non-simple $M$-irreducible connected subgroups contained in either $A_1^2$ or $A_1$ $(p \geq 11)$. We consider the remaining maximal connected subgroups $M_1$ in the following sections. 

\subsection{$M_1 = \bar{D}_4$ $(F_4(\#\ffour{13}))$}

From Lemma \ref{maxclassical}, the reductive, maximal connected subgroups of $\bar{D}_4$ are $B_3$, $A_1 B_2$ ($3$ conjugacy classes of each), $A_2$ $(p \neq 3)$ and $A_1^4$. We have $\text{Out}_{F_4}(\bar{D}_4) \cong S_3$ by \cite[Table 8]{car} and therefore there is only one $F_4$-conjugacy class of maximal subgroups of type $B_3$ and $A_1 B_2$. The subgroup $B_3$ is $F_4$-reducible, since it is a Levi subgroup of $F_4$ (recalling $p \neq 2$). The subgroup $A_1 B_2$ is a diagonal subgroup of $\bar{A}_1^2 B_2 = F_4(\#\ffour{28})$ since there is only one $\text{SO}_9$-conjugacy class of subgroups acting as $\text{SO}_3 \text{SO}_5$ on the natural 9-dimensional module. Therefore $A_1 B_2$ is conjugate to $F_4(\#\ffour{40}^{\{0\}})$ and no further consideration is given to it here. The maximal subgroup $A_2$ is conjugate to $F_4(\#\ffour{19})$. It remains to consider the $M$-irreducible connected subgroups of $\bar{A}_1^4$, all of which are diagonal. From \cite[Table 8]{car}, we see $\text{Out}_{F_4}(\bar{A}_1^4) \cong S_4$, acting naturally on the four $\bar{A}_1$ factors. The conjugacy classes of diagonal subgroups are hence as listed in Table \ref{F4tab} with the following exception: the subgroup $\bar{A}_1^2 \tilde{A}_1 \hookrightarrow \bar{A}_1^4$ via $(1_a,1_b,1_c,1_c)$ is not $M$-irreducible (and hence none of its subgroups is) as it is contained in the $F_4$-reducible subgroup $B_3$. This explains the $(r \neq 0)$ condition on $F_4(\#\ffour{33})$ and subsequent conditions on $F_4(\#\ffour{34})$, $F_4(\#\ffour{35})$ and $F_4(\#\ffour{1})$.    

\subsection{$M_1 = \tilde{A}_1 \bar{A}_3$ $(F_4(\#\ffour{27}))$}

The reductive, maximal connected subgroups of $\tilde{A}_1 \bar{A}_3$ are of the form $\tilde{A}_1 Y$ where $Y$ is a reductive, maximal connected subgroup of $\bar{A}_3$. From Lemma \ref{maxclassical} we see the possibilities for $Y$ are $B_2$ and $\tilde{A}_1^2$ (the subgroups of type $A_1$ are generated by short root subgroups of $B_4$ and hence of $F_4$). The subgroup $\tilde{A}_1 B_2$ is contained in $\bar{A}_1^2 B_2$ and conjugate to $F_4(\#\ffour{40}^{\{0\}})$. The diagonal subgroups of $\tilde{A}_1^3$ follow from the fact that $\text{Out}_{F_4}(\tilde{A}_1^3) \cong S_3$, acting naturally on the three $\tilde{A}_1$ factors. Note that the subgroup $A_1 \tilde{A}_1 \hookrightarrow \tilde{A}_1^3$ via $(1_a,1_a,1_b)$ is $M$-reducible by Lemma \ref{class} since it acts as $(2,0)^2 + (0,2)$ on $V_{B_4}(\lambda_1)$.     

\subsection{$M_1 = \bar{A}_1^2 B_2$ $(F_4(\#\ffour{28}))$} \label{a1a1b2inf4}

The reductive, maximal connected subgroups of $\bar{A}_1^2 B_2$ are either of the form $\bar{A}_1^2 Y$ where $Y$ is a reductive, maximal connected subgroup of $B_2$, or a diagonal subgroup $A_1 B_2 \hookrightarrow \bar{A}_1^2 B_2$ via $(1_a,1_a^{[r]},10)$. By Lemma \ref{maxclassical}, the possibilities for $Y$ are $\bar{A}_1^2$ and $A_1$ $(p \geq 5)$ (recalling our assumption that $p \neq 2$). The subgroup $\bar{A}_1^4$ is conjugate to $F_4(\#\ffour{32})$ since there is only one class of such subgroups in $B_4$. The diagonal subgroups of $X = \bar{A}_1^2 A_1$ $(p \geq 5)$ are as in Table \ref{F4tab}, noting that the normaliser of $X$ contains an involution swapping the two $\bar{A}_1$ factors. As explained in Section \ref{strat}, we do not consider the reductive, maximal connected subgroups of $A_1 B_2 \hookrightarrow \bar{A}_1^2 B_2$ separately. 

This completes the $M = B_4$ case. 

\section{$M = C_4$ ($p=2$) $(F_4(\#\ffour{14}))$} \label{F4C4}   
By Lemma \ref{maxclassical}, the reductive, maximal connected subgroups of $C_4$ when $p=2$ are $B_2^2$, $\bar{A}_1 C_3$ and $\tilde{D}_4$ (this subgroup $\tilde{D}_4$ is generated by short root subgroups of $F_4$ because it is the image of $\bar{D}_4$ under the exceptional graph automorphism of $F_4$). The subgroup $\bar{A}_1 C_3$ is a maximal subgroup of $F_4$ when $p \neq 2$. We therefore consider the subgroups of $\bar{A}_1 C_3$ when $p=2$ in the $\bar{A}_1 C_3$ section below (as indicated by the ID number $24$b in Table \ref{F4tab}). The remaining maximal connected subgroups $M_1$ are considered in the following sections. 

\subsection{$M_1 = B_2^2$ $(F_4(\#\ffour{30}))$} \label{b2b2inf4}

The only reductive, maximal connected subgroups of $B_2$ when $p=2$ are $\bar{A}_1^2$ and $\tilde{A}_1^2$. Almost everything follows from this, by working out the conjugacy classes of diagonal $C_4$-irreducible subgroups of $\bar{A}_1^2 \tilde{A}_1^2$, $\tilde{A}_1^4$ and $\bar{A}_1^4$. We note that there is only one $F_4$-conjugacy class of $\bar{A}_1^4$ and that it is contained in $\bar{D}_4 < B_4$ even when $p=2$. Applying the exceptional graph automorphism of $F_4$ shows that $\text{Out}_{F_4}(\bar{A}_1^4) \cong S_4$ implies that $\text{Out}_{F_4}(\tilde{A}_1^4) \cong S_4$. The diagonal $C_4$-irreducible subgroups of $\tilde{A}_1^4$ now follow. The normalisers of both $\bar{A}_1^2$ and $\tilde{A}_1^2$ in $B_2$ contain an involution swapping the two $A_1$ factors and hence $\text{Out}_{F_4}(\bar{A}_1^2 \tilde{A}_1^2) \cong S_2 \times S_2$. Again, the diagonal $C_4$-irreducible subgroups now follow. 

\subsection{$M_1 = \tilde{D}_4$ $(F_4(\#\ffour{15}))$}

By Lemma \ref{maxclassical}, the reductive, maximal connected subgroups of $\tilde{D}_4$ when $p=2$ are $B_3$ ($3$ conjugacy classes), $\tilde{A}_1^4$ and $A_2$. Since $\text{Out}_{F_4}(\tilde{D}_4) \cong S_3$, it follows that there is only one $F_4$-conjugacy class of maximal subgroups of type $B_3$. Therefore Lemma \ref{class} shows that $B_3$ is $C_4$-reducible since one of the $\tilde{D}_4$ classes acts as $000|100|000$ on $V_{\tilde{D}_4}(\lambda_1)$. The subgroup $A_2$ is conjugate to $F_4(\#\ffour{18}^{\{1,0\}})$. Finally, consider $\tilde{A}_1^4$. This subgroup is contained in $B_2^2$ because there is only one $C_4$-conjugacy class of subgroups of type $A_1^4$ acting as $(1,1,0,0) + (0,0,1,1)$ on $V_{C_4}(\lambda_1)$. Therefore $\tilde{A}_1^4$ is conjugate to $F_4(\#\ffour{51})$. 

This completes the $C_4$ $(p=2)$ case. 

\section{$M = \bar{A}_1 C_3$ $(F_4(\#\ffour{24}))$}

The subgroup $\bar{A}_1 C_3$ is maximal when $p \neq 2$ and contained in $C_4$ when $p=2$. We consider both cases in this section. By Lemma \ref{maxclassical}, the reductive, maximal connected subgroups of $C_3$ are $\bar{A}_1 B_2$, $A_1 A_1$ $(p \neq 2)$, $\tilde{A}_3$ $(p=2)$, $G_2$ $(p=2)$ and $A_1$ $(p \geq 7)$. This yields the reductive, maximal connected subgroups of $\bar{A}_1 C_3$, as listed in Table \ref{F4tab}. 

By \cite[p.333, Table 2]{LS6}, we have $C_{F_4}(B_2) = \bar{A}_1^2$ $(p \neq 2)$ and $B_2$ $(p = 2)$. Thus, the subgroup $\bar{A}_1^2 B_2$ has already been considered in Sections \ref{a1a1b2inf4} and \ref{b2b2inf4}, respectively and in particular, is conjugate to $F_4(\#\ffour{28})$. The only thing noteworthy about the diagonal subgroups of $\bar{A}_1 A_1 A_1$ is that $X = \bar{A}_1 A_1 \hookrightarrow \bar{A}_1 A_1 A_1$ via $(1_a,1_b,1_b)$ is contained in $\bar{A}_1^2 B_2$ when $p \geq 5$ since ${2 \otimes 1} = 3 + 1$ and it follows that $X$ is conjugate to $F_4(\#\ffour{42}^{\{0\}})$. When $p=3$, we have ${2 \otimes 1} = T(3) = 1 | 3 | 1$ and hence $X$ is $M$-reducible by Lemma \ref{class}.  

The only reductive, maximal connected subgroup of $\tilde{A}_3$ when $p=2$ is $B_2$, which is $C_3$-reducible by Lemma \ref{class}. Therefore $\bar{A}_1 \tilde{A}_3$ contains no proper $F_4$-irreducible subgroups. However, since it has maximal rank it is clearly $F_4$-irreducible. 

Finally, we consider the subgroups contained in $\bar{A}_1 G_2$ $(p=2)$. The $G_2$-irreducible subgroups are given by Theorem \ref{G2THM}. The maximal connected subgroup $A_2 < G_2$ is $C_3$-reducible by Lemma \ref{class} and so $\bar{A}_1 \tilde{A}_2$ is $M$-reducible. The subgroup $\bar{A}_1 \tilde{A}_1 A_1 < \bar{A}_1 G_2$ is also contained in $\bar{A}_1^2 \tilde{A}_1^2 < \bar{A}_1^2 B_2 < C_4$ and conjugate to $F_4(\#\ffour{44}^{\{0,1\}})$. Indeed, the subgroup $\tilde{A}_1 A_1 < G_2 < C_3$ is conjugate to $\tilde{A}_1 A_1 \hookrightarrow \tilde{A}_1^2 \bar{A}_1 < C_3$ via $(1_a,1_b,1_b^{[1]})$ since both act as $(1,1) + (0,2)$ on $V_{C_3}(\lambda_1)$.  

\section{$M = A_1 G_2$ $(p \neq 2)$ $(F_4(\#\ffour{25}))$}

The lattice structure of $M$-irreducible connected subgroups of $M$ follows from Theorem \ref{G2THM}. The $G_2$ factor of $M$ is contained in $\bar{D}_4$ as seen from the construction in \cite[3.9]{se2} and therefore subgroups generated by long root subgroups of $G_2$ are generated by long root subgroups of $F_4$. It follows that $A_1 \bar{A}_1 A_1$ is conjugate to $F_4(\#\ffour{57})$, since $C_{F_4}(\bar{A}_1)^\circ = C_3$ by Lemma \ref{rootycentralisers}.

\section{$M = \bar{A}_2 \tilde{A}_2$ $(F_4(\#\ffour{26}))$}

By Lemma \ref{maxclassical}, the only reductive, maximal connected subgroup of $A_2$ is $A_1$ $(p \neq 2)$. The reductive, maximal connected subgroups of $M$ are thus as in Table \ref{F4tab}. We note that $\bar{A}_2 A_1$ $(p \neq 2)$ is conjugate to $F_4(\#\ffour{64}) < A_1 G_2$, since the $\bar{A}_2 < G_2$ is generated by long root subgroups of $F_4$ and $C_{F_4}(\bar{A}_2) = \tilde{A}_2$ by Lemma \ref{rootycentralisers}. It subsequently follows that $A_1 A_1 < \bar{A}_2 \tilde{A}_2$ is conjugate to $F_4(\#\ffour{62}^{\{0,0\}})$ since $A_1 A_1$ is a subgroup of $\bar{A}_2 A_1$.    

\section{$M = G_2$ $(p=7)$ $(F_4(\#\ffour{16}))$}

The lattice structure of $M$-irreducible subgroups is given by Theorem \ref{G2THM}. We show that all reductive, maximal connected subgroups of $M$ are conjugate to subgroups already considered. For the maximal connected subgroups $A_2$ and $A_1$ this follows from the proofs of \cite[Theorem 3.4]{tho1} and \cite[Theorem 3]{tho2}. Now consider $X = A_1 A_1$; we claim that $X$ is  conjugate to $Y = F_4(\#\ffour{61}^{\{0,0\}}) < \bar{A}_1 C_3$. To do this it suffices to show that $A_1 A_1$ is contained in $\bar{A}_1 C_3$ as comparing composition factors on $V_{26}$ shows that $X$ is then conjugate to $Y$. We know by Lemma \ref{centralisers} that $X$ is the centraliser in $G_2$ of a semisimple element of order 2, call it $t$. Again by Lemma \ref{centralisers}, we know that the centraliser in $F_4$ of $t$ is $B_4$ or $\bar{A}_1 C_3$ and by \cite[Proposition 1.2]{LS8} the trace on $V_{26}$ is $-6$ or $2$, respectively. We calculate the trace of $t$ on $V_{26}$ from the restriction $V_{26} \downarrow X = (2,2) + (1,1) + (1,3) + (0,4)$, noting that the element $t$ can be seen as minus the identity in both $A_1$ factors and hence has trace $2$ on $V_{26}$. Therefore $X$ is contained in $\bar{A}_1 C_3$, proving the claim.     

This completes the proof of Theorem \ref{F4THM}.

%% file: E6.tex
\chapter{Irreducible subgroups of $G = E_6$} \label{secE6}

In this section we use the strategy described in Section \ref{strat} to prove the following theorem. 

\begin{thm} \label{E6THM}
Let $X$ be an $E_6$-irreducible connected subgroup of $E_6$. Then $X$ is $\text{Aut}(E_6)$-conjugate to exactly one subgroup of Table \ref{E6tab} and each subgroup in Table \ref{E6tab} is $E_6$-irreducible. Moreover, Table \ref{E6tab} gives the lattice structure of the irreducible connected subgroups of $E_6$.   
\end{thm}

Recall that throughout the proof we will use \cite[Theorem 1]{tho1} and \cite[Theorem 4]{tho2}, which classify the simple $E_6$-irreducible connected subgroups, without reference. 

By Theorem~\ref{maximalexcep}, the reductive, maximal connected subgroups of $E_6$ are $\bar{A}_1 A_5$, $\bar{A}_2^3$, $A_2 G_2$, $F_4$, $C_4$ $(p \neq 2)$, $G_2$ $(p \neq 7)$ and $A_2$ $(p \geq 5)$. The only irreducible connected subgroup contained in $A_2$ is $A_1$ $(p \neq 2)$. Since this is a simple subgroup, the maximal connected subgroup $A_2$ $(p \geq 5)$ requires no further consideration. In the following sections we consider each of the remaining reductive, maximal connected subgroups $M$ in turn. 

\section{$M = \bar{A}_1 A_5$ $(E_6(\#\esix{24}))$}
 
The reductive, maximal connected subgroups of $A_5$ are $C_3$, $A_1 A_2$, $A_3$ $(p \neq 2)$ and $A_2$ $(p \neq 2)$ by Lemma \ref{maxclassical}. Therefore the reductive, maximal connected subgroups of $M$ are as in Table \ref{E6tab}, and we consider them in the following sections. 

\subsection{$M_1 = \bar{A}_1 C_3$  $(E_6(\#\esix{27}))$}

Using Lemma \ref{maxclassical}, we find that the reductive, maximal connected subgroups of $\bar{A}_1 C_3$ are $\bar{A}_1^2 C_2$, $\bar{A}_1 A_1 A_1$ $(p \neq 2)$, $\bar{A}_1 A_1$ $(p \geq 7)$, $\bar{A}_1 A_3$ $(p=2)$ and $\bar{A}_1 G_2$ $(p=2)$. Since the subgroup $\bar{A}_1 C_2 < C_3$ acts reducibly on $V_{C_3}(100)$, it follows that $\bar{A}_1^2 C_2$ is $M$-reducible. 

The diagonal subgroups of $\bar{A}_1 A_1 A_1$ $(p \neq 2)$ and $\bar{A}_1 A_1$ $(p \geq 7)$ are easily seen to be as in Table \ref{E6tab}, noting that $\bar{A}_1 A_1 \hookrightarrow \bar{A}_1 A_1 A_1$ via $(1_a,1_b,1_b)$ is $M$-reducible by Lemma \ref{class}. Neither $X = \bar{A}_1 A_3$ $(p=2)$ nor $Y = \bar{A}_1 G_2$ $(p=2)$ has a proper $M$-irreducible connected subgroup, by Lemma \ref{class}. 

It remains to prove that $X$ and $Y$ are $E_6$-irreducible. Suppose that $X$ is $E_6$-reducible. Then by Lemma \ref{wrongcomps}, there exists a subgroup $Z$ of type $A_1 A_3$ contained $L$-irreducibly in a Levi subgroup $L$, such that $X$ and $Z$ have the same composition factors on $V_{27}$. It follows that $L'$ has type $D_5$ or $A_1 A_3$. From Table \ref{E6tabcomps}, the dimensions of the $X$-composition factors on $V_{27}$ are $14,12,1$, whereas Table \ref{levie6} shows that $D_5$ and $A_1 A_3$ have composition factors of dimensions $16,10,1$ and $8,6,4^2,2^2,1$, respectively.  Thus $X$ and $Z$ do not have the same composition factors on $V_{27}$, a contradiction. Hence $X$ is $E_6$-irreducible. 

The same argument applies to $Y$ since the composition factors of $V_{27} \downarrow Y$ also have dimensions $14, 12, 1$ and only a Levi subgroup of type $D_5$ contains an irreducible subgroup of type $A_1 G_2$ when $p=2$. 

\subsection{$M_1 = \bar{A}_1 A_1 A_2$  $(E_6(\#\esix{28}))$} \label{E6A1A1A2}

The only reductive, maximal connected subgroup of $A_2$ is $A_1$ $(p \neq 2)$. Therefore the lattice structure of $M_1$-irreducible subgroups is as given in Table \ref{E6tab} along with the subgroup $Y = A_1 A_2 \hookrightarrow \bar{A}_1 A_1 A_2$ via $(1,1,10)$ when $p=2$. There is only one $A_5$-conjugacy class of subgroups $A_1 A_1$ acting as $(2,1)$ on $V_{A_5}(\lambda_1)$ and so $\bar{A}_1 A_1 A_1$ is contained in $\bar{A}_1 C_3$ and thus conjugate to $E_6(\#\esix{31})$. 

It remains to prove that when $p=2$ the subgroup $X = A_1 A_2 \hookrightarrow \bar{A}_1 A_1 A_2$ via $(1^{[r]},1^{[s]},10)$ $(rs=0; r \neq s)$ is $E_6$-irreducible and that $Y$, defined as above, is $E_6$-reducible. To prove that $X$ is $E_6$-irreducible we use Lemma \ref{wrongcomps}. Suppose that $Z$ is an $L$-irreducible subgroup $A_1 A_2$ of a Levi subgroup $L$, such that $X$ and $Z$ have the same composition factors on $V_{27}$. Since $L'$ contains an irreducible subgroup $A_1 A_2$ it follows from Lemma \ref{class} that $L'$ is of type $A_5$, $\bar{A}_1 \bar{A}_2^2$, $\bar{A}_1^2 A_2$ or $\bar{A}_1 A_2$. The $X$-composition factors of $V_{27}$ are $(1^{[r]} \otimes 1^{[s]},10) /$ $\!\! (2^{[s]},10) /$ $\!\! (0,02) /$ $\!\! (0,10)^2$. Now it follows that $Z$ is not contained in $A_5$ since $V_{27} \downarrow A_5 = \lambda_1^2 / \lambda_4$. Similarly, since $Z$ has a composition factor of dimension 12, we find that $Z$ is not a subgroup of $\bar{A}_1 \bar{A}_2^2$, $\bar{A}_1^2 A_2$ or $\bar{A}_1 A_2$. This is a contradiction and thus $X$ is $E_6$-irreducible.

Now consider $Y$. It is shown in Section \ref{e6a2g2} below that $M_1$ is conjugate to $A_2 \bar{A}_1 A_1 < A_2 G_2$. When $p=2$, the subgroup $A_1 \hookrightarrow \bar{A}_1 A_1 < G_2$ via $(1,1)$ is $G_2$-reducible by Theorem \ref{G2THM}. Therefore $Y$ is $A_2 G_2$-reducible and thus $E_6$-reducible.  

\subsection{$M_1 = \bar{A}_1 A_3$ $(p \neq 2)$  $(E_6(\#\esix{29}\mathrm{a}))$ and $ \bar{A}_1 A_2$ $(p \neq 2)$  $(E_6(\#\esix{30}))$}

There are no proper $M$-irreducible connected subgroups of either $\bar{A}_1 A_3$ $(p \neq 2)$ or $\bar{A}_1 A_2$ $(p \neq 2)$ since any proper connected subgroup of $A_3$ or $A_2$ is $A_5$-reducible by Lemma \ref{class}. It remains to prove that they are both $E_6$-irreducible. The dimensions of the composition factors on $V_{27}$ are $15,12$ for both subgroups, as seen from Table \ref{E6tab}. An easy application of Lemma \ref{wrongcomps}, as in the previous sections, shows that both are $E_6$-irreducible.

\section{$M = \bar{A}_2^3$  $(E_6(\#\esix{25}))$}

The only reductive, maximal connected subgroup of $A_2$ is $A_1$ $(p \neq 2)$. We have $\text{Out}_{\text{Aut}(E_6)}(\bar{A}_2^3) \cong S_3 \times S_2$ where the central involution acts as a graph automorphism on each $A_2$ factor; a two-cycle in the $S_3$ direct factor swaps two of the $A_2$ factors whilst also inducing a graph automorphism on them; a three-cycle in the $S_3$ direct factor acts naturally as a three-cycle on the three $A_2$ factors. The lattice structure of the $M$-irreducible subgroups then follows. It remains for us to show that $X = A_1 A_2 \hookrightarrow A_1 \bar{A}_2^2$ via $(1,10,01)$ is conjugate to $Y = E_6(\#\esix{37}^{\{0,0\}})$ when $p \neq 2$. Let $Z = E_6(\#\esix{28}) = \bar{A}_1 A_1 A_2$ so $Y \hookrightarrow Z$ via $(1,1,10)$. It suffices to prove that $Y$ is contained in $\bar{A}_2^3$. It is shown in Section \ref{e6a2g2} below that $Z$ is conjugate to $A_2 \bar{A}_1 A_1 < A_2 G_2$. Therefore $Y$ is conjugate to $A_2 A_1 \hookrightarrow A_2 \bar{A}_1 A_1 < A_2 G_2$ via $(10,1,1)$. By Theorem \ref{G2THM} this implies that $Y$ is also conjugate to $A_2 A_1 < A_2 \bar{A}_2 < A_2 G_2$. Using Lemma \ref{rootycentralisers}, we see that $A_2 \bar{A}_2$ is contained in $\bar{A}_2^3$, and hence $Y$ is contained in $\bar{A}_2^3$, as required.     

\section{$M = A_2 G_2$  $(E_6(\#\esix{26}))$} \label{e6a2g2}

By Lemma \ref{maxclassical}, the only proper irreducible connected subgroup of $A_2$ is $A_1$ $(p \neq 2)$; the $G_2$-irreducible connected subgroups are given by Theorem \ref{G2THM}. Using this we find that the lattice structure of $M$-irreducible connected subgroups is as given in Table \ref{E6tab}, noting that the $G_2$ factor is contained in a Levi subgroup $D_4$ and hence the long root subgroups of $G_2$ are long root subgroups of $E_6$.

Since $C_{E_6}(\bar{A}_2) = \bar{A}_2^2$ by Lemma \ref{rootycentralisers}, the subgroup $A_2 \bar{A}_2$ is contained in $\bar{A}_2^3$ and comparing composition factors shows that it is conjugate to $E_6(\#\esix{40}^{\{0\}})$. Similarly, the subgroup $A_2 \bar{A}_1 A_1$ is conjugate to $E_6(\#\esix{28}) < \bar{A}_1 A_5$. All other conjugacies follow from these two facts. 

\section{$M = F_4$  $(E_6(\#\esix{7}))$}

Theorem \ref{F4THM} gives the lattice structure of the $F_4$-irreducible subgroups. The maximal subgroup $F_4$ of $E_6$ is the centraliser of a standard graph automorphism of $E_6$ by \cite[Theorem 15.1]{se2}. Therefore the maximal subgroup $B_4$ of $F_4$ is contained in a $D_5$ Levi subgroup. It also follows that the maximal connected subgroup $\bar{A}_2 A_2$ is conjugate to $E_6(\#\esix{39}^{\{0\}})$.   

The subgroup $C_4$ is a maximal subgroup of $E_6$ when $p \neq 2$; we therefore consider the subgroups of the maximal subgroup $C_4$ $(p=2)$ in the next section. Similarly, the subgroup $G_2$  is a maximal subgroup of $E_6$ when $p \neq 7$; we therefore consider the subgroups of the maximal subgroup $G_2$ $(p=7)$ in Section \ref{e6g2}.   

The maximal subgroup $\bar{A}_1 C_3$ $(p \neq 2)$ is contained in $\bar{A}_1 A_5$ since $C_{E_6}(\bar{A}_1)^\circ = A_5$ by Lemma \ref{rootycentralisers} and hence $\bar{A}_1 C_3$ is conjugate to $E_6(\#\esix{27})$. Finally, we consider the maximal subgroup $X = A_1 G_2$ $(p \neq 2)$. By \cite[p.333, Table 3]{LS6}, we have $C_{E_6}(G_2)^\circ = A_2$ and thus $X < A_2 G_2$. It follows that $X$ is conjugate to $E_6(\#\esix{47})$. 

\section{$M = C_4$  $(E_6(\#\esix{8}))$} \label{e6c4}

The subgroup $C_4$ is maximal when $p \neq 2$ and contained in $F_4$ when $p=2$. We consider both cases in this section. Lemma \ref{maxclassical} yields the reductive, maximal connected subgroups of $C_4$. They are $B_2^2$, $\bar{A}_1 C_3$, $A_1^3$ $(p \neq 2)$, $A_1$ $(p \geq 11$) and $D_4$ $(p =2)$. The subgroup $B_2^2$ is $E_6$-reducible by Lemma \ref{semirr} since $C_{E_6}(B_2)^\circ = B_2 T_1$ by \cite[p.333, Table 2]{LS6}. The subgroup $\bar{A}_1 C_3$ is conjugate to $E_6(\#\esix{27})$, by the argument given in the previous section. 

Next, we prove that the maximal connected subgroup $X = A_1^3$ $(p \neq 2)$ is conjugate to $Y = E_6(\#\esix{44}) < \bar{A}_2^3$. It suffices to show that $Y$ is contained in $C_4$. Consider the standard graph automorphism of $E_6$, call it $\tau$ and let $w_0$ be the longest word of the Weyl group. Then $w_0 = - \tau$ and so $t := \tau w_0$ acts by inversion on a maximal torus of $E_6$. Therefore $t$ induces a graph automorphism on each $A_2$ factor of $\bar{A}_2^3$. Thus $Y < C_{E_6}(t)^\circ$ since the irreducible subgroup $A_1$ in a subgroup $A_2$ is centralised by the graph automorphism of $A_2$. Finally, we check that $\mathrm{dim}(C_{L(E_6)}(t)) = 36$ and so $\mathrm{dim}(C_{E_6}(t)) = 36$ (by \cite[9.1]{borel}, since $t$ is semisimple). Therefore, $C_{E_6}(t)^\circ = C_4$ by \cite[Table 4.3.1]{gls3}.       

Finally, we consider the subgroups of the maximal subgroup $D_4$ $(p=2)$. In this case $M$ is contained in $F_4$ since $p=2$ and so we need only consider the $F_4$-irreducible subgroups contained in $D_4$, which are $A_2$ and $\bar{A}_1^4$. Moreover, Theorem \ref{F4THM} shows that $A_2$ is a subgroup of $\bar{A}_2 A_2 < F_4$ and hence conjugate to $E_6(\#\esix{15}^{\{1\}})$. It also shows that $\bar{A}_1^4$ is a subgroup of $B_2^2$ and thus $E_6$-reducible.   

\section{$M = G_2$ $(E_6(\#\esix{10}))$} \label{e6g2}

The subgroup $G_2$ is maximal when $p \neq 7$ and contained in $F_4$ when $p=7$. In this section we consider both cases. Theorem \ref{G2THM} gives the lattice structure of $M$-irreducible connected subgroups. We need only consider the maximal connected subgroup $X = A_1 A_1$. We claim that if $p \neq 2$ then $X$ is conjugate to $E_6(\#\esix{34}^{\{0,0\}})$ but when $p=2$ it is $E_6$-reducible. 

First suppose that $p \neq 2$. It suffices to prove that $X$ is contained in $\bar{A}_1 A_5$. By Lemma \ref{centralisers}, $X$ is the centraliser in $G_2$ of a semisimple element of order 2, call this $t$. Moreover, since $X$ contains a simple $E_6$-irreducible subgroup of type $A_1$ it is $E_6$-irreducible. Therefore, the centraliser in $E_6$ of $t$ is $\bar{A}_1 A_5$ by Lemma \ref{centralisers}.  

Now let $p=2$. To prove that $X$ is $E_6$-reducible we consider the action of $X$ on $L(E_6)$. By \cite[Table~10.1]{LS1}, we have $L(E_6) \downarrow G_2 = 11 + 01$. In Table \ref{G2tabcomps}, the composition factors of $V_{G_2}(01) \downarrow X$ are given and moreover, $V_{G_2}(01) \downarrow X = ((0,0) | ((2,0) + (0,2)) | (0,0)) + (1,3)$. Therefore, $X$ fixes a non-trivial vector of $L(E_6)$. By \cite[Lemma~1.3]{se2}, it follows that $X$ is contained in either a parabolic subgroup, or $\bar{A}_1 A_5$ or $\bar{A}_2^3$. The $\bar{A}_1 A_5$ and $\bar{A}_2^3$ sections show that neither $\bar{A}_1 A_5$ nor $\bar{A}_2^3$ contains an $E_6$-irreducible subgroup $A_1 A_1$ when $p=2$. Therefore $X$ is $E_6$-reducible, as claimed. 

This completes the proof of Theorem \ref{E6THM}.

%% file: E7.tex
\chapter{Irreducible subgroups of $G = E_7$} \label{secE7}

In this section we prove Theorem \ref{MAINTHM} when $G$ is of type $E_7$ by proving the following theorem. 

\begin{thm} \label{E7THM}
Let $X$ be an $E_7$-irreducible connected subgroup of $E_7$. Then $X$ is conjugate to exactly one subgroup of Table \ref{E7tab} and each subgroup in Table \ref{E7tab} is $E_7$-irreducible. Moreover, Table \ref{E7tab} gives the lattice structure of the irreducible connected subgroups of $E_7$.   
\end{thm}

As in the previous sections we consider each of the reductive, maximal connected subgroups of $E_7$ in turn. The simple $E_7$-irreducible connected subgroups are classified in \cite[Theorem~2]{tho1} and \cite[Theorem~5]{tho2} and we will use these without reference throughout. By Theorem \ref{maximalexcep}, the reductive, maximal connected subgroups of $E_7$ are $\bar{A}_1 D_6$, $\bar{A}_2 A_5$, $A_7$, $G_2 C_3$, $A_1 F_4$, $A_1 G_2$ $(p \neq 2)$, $A_1 A_1$ $(p \geq 5)$, $A_2$ $(p \geq 5)$ and $A_1$ (2~classes, $p \geq 17, 19$). 

There are no proper non-simple irreducible subgroups of $A_1 A_1$ $(p \geq 5)$, $A_2$ $(p \geq 5)$, $A_1$ $(p \geq 17)$ or $A_1$ $(p \geq 19)$ so these require no examination. We treat the remaining cases in the following sections. 

\section{$M = \bar{A}_1 D_6$  $(E_7(\#\eseven{30}))$} 

Using Lemma \ref{maxclassical} we find that the reductive, maximal connected subgroups of $M$ are $\bar{A}_1^3 \bar{D}_4$, $\bar{A}_1 A_1 B_4$ $(p \neq 2)$, $\bar{A}_1 B_2 B_3$ $(p \neq 2)$, $\bar{A}_1 \bar{A}_3^2$, $\bar{A}_1 B_5$ and $\bar{A}_1 A_1 C_3$ (2 classes). We consider each of these in the following sections. 

\subsection{$M_1 = \bar{A}_1^3 \bar{D}_4$ $(E_7(\#\eseven{36}))$} 

We first note that $\text{Out}_{E_7}(M_1) \cong S_3$ where the $S_3$ acts simultaneously as the outer automorphism group of both $\bar{A}_1^3$ and $\bar{D}_4$. Therefore, by Lemma \ref{maxclassical}, the $E_7$-conjugacy classes of reductive, maximal connected subgroups of $M_1$ are $\bar{A}_1^7$, $\bar{A}_1^3 B_3$, $\bar{A}_1^3 A_1 B_2$ $(p \neq 2)$, $\bar{A}_1^3 A_2$ $(p \neq 3)$ and maximal diagonal subgroups. As explained in Section \ref{strat}, we do not explicitly consider the subgroup lattice of the maximal diagonal subgroups and we write down all classes of diagonal subgroups of $M_1$ in the same place as the reductive, maximal connected subgroups in Table \ref{E7tab}. 

All $M_1$-irreducible subgroups of $\bar{A}_1^7$ are diagonal and are given in Table \ref{A17E7diags}. They follow in the same way that the diagonal subgroups of type $A_1$ are found in the proof of \cite[Theorem~5]{tho2}. Indeed, we have chosen the same isomorphism of $\text{PSL}(2,7) \cong \text{Out}_{E_7}(\bar{A}_1^7)$ with a subgroup of $S_7$ acting on the seven factors, namely that the generators are mapped to $(1,2,3) (5,6,7)$ and $(2,4) (3,5)$.

By Lemmas \ref{maxclassical} and \ref{class}, the $M$-irreducible, reductive, maximal connected subgroups of $\bar{A}_1^3 B_3$ are $\bar{A}_1^5 A_1$ $(p \neq 2)$, $\bar{A}_1^3 G_2$, $\bar{A}_1^3 A_1 B_2$ $(p=2)$ and maximal diagonal subgroups. The diagonal subgroups follow from noting that $\text{Out}_{E_7}(\bar{A}_1^3 B_3) \cong S_2$ where the involution swaps the second and third $\bar{A}_1$ factors. The subgroup $\bar{A}_1^3 A_1 B_2$ is a maximal connected subgroup of $\bar{A}_1^3 \bar{D}_4$ when $p \neq 2$ and we will consider it below (in the notation of Table \ref{E7tab}, the subgroup $\bar{A}_1^3 A_1 B_2$ is $E_7(\#\eseven{45}\mathrm{b})$). The subgroup $\bar{A}_1^5 A_1$ is contained in $\bar{A}_1^7$ since the sixth $A_1$ factor is a subgroup of $\bar{A}_1^2$ diagonally embedded via $(1,1)$. 

Now consider $X = \bar{A}_1^3 G_2$. The diagonal $M$-irreducible subgroups of $X$ follow from Lemma \ref{class} and the fact that $\text{Out}_{E_7}(X) \cong S_3$ acting naturally on the three $\bar{A}_1$ factors. Theorem \ref{G2THM} yields the remaining $X$-irreducible connected subgroups. By Lemma \ref{class}, the subgroup $\bar{A}_1^3 \bar{A}_2$ is $M$-reducible. We also see that the subgroup $\bar{A}_1 A_1 < G_2 < \bar{D}_4$ is a diagonal subgroup of $\bar{A}_1^4 < \bar{D}_4$ via $(1_a,1_b,1_b,1_b)$ and hence $\bar{A}_1^4 A_1$ is a subgroup of $\bar{A}_1^7$. The subgroup $\bar{A}_1^3 A_2$ $(p =3)$ is a maximal connected subgroup of $\bar{A}_1^3 \bar{D}_4$ when $p \neq 3$ and we will consider the subgroups for all $p$ below. Finally, the diagonal subgroups of $\bar{A}_1^3 A_1$ $(p \geq 7)$ follow in the same way as the diagonal subgroups of $X$.   

Next, we consider the subgroup $X = \bar{A}_1^3 A_1 B_2$ which is a maximal connected subgroup of $M_1$ when $p \neq 2$ and contained in $\bar{A}_1^3 B_3$ when $p=2$. The diagonal subgroups follow from the fact that $\text{Out}_{E_7}(X) \cong S_2$ where the involution swaps the second and third $\bar{A}_1$ factors. Using Lemma \ref{maxclassical}, we see the reductive, maximal connected subgroups of $X$ are $\bar{A}_1^5 A_1$, $\bar{A}_1^3 A_1^3$ $(p=2)$ and $\bar{A}_1^3 A_1 A_1$ $(p \geq 5)$. The first subgroup is contained in $\bar{A}_1^7$ and conjugate to $E_7(\#\eseven{49}^{\{0\}})$. The diagonal subgroups of $Y_1 =\bar{A}_1^3 A_1^3$ are found in Table \ref{A13A13E7diags} and are deduced from the fact that $\text{Out}_{E_7}(Y_1) \cong S_2 \times S_3$, where the $S_3$ direct factor acts naturally on $A_1^3$ and the central involution swaps the second and third $\bar{A}_1$ factors. Similarly, the diagonal subgroups of $Y_2 = \bar{A}_1^3 A_1 A_1$ can be found in Table \ref{A13A1A1E7diags}. Here they follow from the fact that $\text{Out}_{E_7}(Y_2) \cong S_2$, where the involution swaps the second and third $\bar{A}_1$ factors. However, the normaliser of $\bar{A}_1^3 A_1 \hookrightarrow Y_2$ via $(1_a,1_b,1_c,1_d,1_d)$ contains an $S_3$ acting naturally on the three $\bar{A}_1$ subgroups. Indeed, since a subgroup $A_1$ acting as $4 + 2$ on $V_{\bar{D}_4}(\lambda_1)$ is contained in both a $\bar{D}_4$-irreducible $A_2$ and $A_1 B_2$ it is centralised by both a triality and involutory automorphism  of $\bar{D}_4$. 

Finally, consider $\bar{A}_1^3 A_2$ which is a maximal connected subgroup of $M_1$ when $p \neq 3$ and contained in $\bar{A}_1^3 G_2$ when $p=3$. Other than diagonal subgroups, the only reductive, maximal connected subgroup is $\bar{A}_1^3 A_1$ $(p \neq 2)$. When $p \geq 5$ this is a subgroup of $\bar{A}_1^3 A_1 A_1 = E_7(\#\eseven{87})$ and when $p=3$ the maximal subgroup $A_1$ of $A_2 < G_2$ is also contained in $\bar{A}_1 A_1 < G_2$ and hence $\bar{A}_1^3 A_1$ is a subgroup of $\bar{A}_1^7$. 

It remains to prove that when $p=2$ the diagonal subgroups $E_7(\#\eseven{184})$ are $E_7$-irreducible. Let $X = A_1 A_2 \hookrightarrow \bar{A}_1^3 A_2 = E_7(\#\eseven{46})$ via $(1,1^{[r]},1^{[s]},10)$ $(0 < r < s)$ and suppose that $X$ is $E_7$-reducible. From Table \ref{E7tabcomps}, we have $V_{56} \downarrow \bar{A}_1^3 A_2 = (1,1,1,00) /$ $\!\! (1,0,0,11) /$ $\!\! (0,1,0,11) /$ $\!\! (0,0,1,11)$. Therefore, $X$ has three $16$-dimensional composition factors on $V_{56}$. By Lemma \ref{wrongcomps}, there exists a subgroup $Z$ of type $A_1 A_2$ contained $L$-irreducibly in a Levi subgroup $L$, such that $X$ and $Z$ have the same composition factors on $V_{56}$. Since $Z$ has three $16$-dimensional composition factors on $V_{56}$, it follows from Table \ref{levie7} that $L'$ is of type $A_1 D_5$. Using Lemma \ref{class}, we see that $L' = A_1 D_5$ does not contain an $L'$-irreducible subgroup $A_1 A_2$, a contradiction. Hence $X$ is $E_7$-irreducible.       
 
\subsection{$M_1 = \bar{A}_1 A_1 B_4$ $(E_7(\#\eseven{37}))$} 

In this section we consider $\bar{A}_1 A_1 B_4$ both when it is a maximal connected subgroup of $M$ $(p \neq 2)$ and when it is contained in $\bar{A}_1 B_5$ $(p =2)$. By Lemma \ref{maxclassical}, the reductive, maximal connected subgroups of $\bar{A}_1 A_1 B_4$ are those listed in Table \ref{E7tab}. The subgroups $\bar{A}_1 A_1 \bar{D}_4$ and $\bar{A}_1^3 A_1 B_2$ $(p \neq 2)$ are both contained in $\bar{A}_1^3 \bar{D}_4$, as seen by considering their action on $V_{M}$. 

We now consider $\bar{A}_1 A_1^2 \bar{A}_3$, which is maximal connected when $p \neq 2$ and contained in $\bar{A}_1 A_1^2 B_3$ when $p=2$. The reductive, maximal connected subgroups of $\bar{A}_3$ are $B_2$ and $A_1^2$ $(p \neq 2)$. However, when $p=2$ the subgroup $\bar{A}_1 A_1^2 B_2$ is $M$-reducible by Lemma \ref{class}. When $p \neq 2$ the subgroup $\bar{A}_1 A_1^2 B_2$ is contained in $\bar{A}_1^3 \bar{D}_4$ and thus conjugate to $E_7(\#\eseven{91}^{\{0\}})$. The diagonal subgroups of $Y = \bar{A}_1 A_1^4$ follow from $\text{Out}_{E_7}(Y) \cong S_4$, acting naturally on the last four factors. The diagonal $M$-irreducible subgroups of $\bar{A}_1 A_1^2 \bar{A}_3$, as listed in Table \ref{E7tabcomps}, follow from noting that there is an involution in the Weyl group of $D_6$ swapping the second and third $A_1$ factors, and the fact that the field twists corresponding to the second and third factors must be distinct for the subgroup to be $M$-irreducible.  

Now consider $X = \bar{A}_1 A_1 A_1^2$ $(p \neq 2)$. Since $N_{B_4}(A_1^2)$ contains an involution swapping the two $A_1$ factors we have $\text{Out}_{E_7}(X) \cong S_2$. Also, the subgroup $A_1 \hookrightarrow A_1^2 < B_4$ via $(1,1)$ is contained in $\bar{D}_4$ when $p \geq 5$ and is $B_4$-reducible when $p=3$. Indeed, this follows from Lemma \ref{class}, since $2 \otimes 2 = 4 + 2 + 0$ when $p \geq 5$ and $2 \otimes 2 = (0|4|0) + 2$ when $p=3$. The non-simple diagonal subgroups of $X$ follow from this and are found in Table \ref{A1A1A12E7diags}. 

Now consider $X = \bar{A}_1 A_1 B_2^2$ $(p=2)$. The reductive, maximal connected subgroups of $B_2$ when $p=2$ are $\bar{A}_1^2$ and $A_1^2$. The maximal connected subgroup $\bar{A}_1^3 A_1 B_2$ is contained in $\bar{A}_1^3 \bar{D}_4$ and the diagonal subgroups of $X$ follow easily, noting that there is an involution in $N_{B_4}(B_2^2)$ swapping the $B_2$ factors and that $\bar{A}_1 A_1 B_2\hookrightarrow X$ via $(1,1,10,10)$ is $\bar{A}_1 D_6$-reducible by Lemma \ref{class}. That leaves us to consider the subgroups of $Y = \bar{A}_1 A_1^3 B_2$.    

The diagonal subgroups of $Y$ follow as usual, noting that $\text{Out}_{M}(Y) \cong S_3$, acting naturally on the three conjugate $A_1$ factors. The non-diagonal reductive, maximal connected subgroups of $Y$ are $\bar{A}_1^3 A_1^3$ and $\bar{A}_1 A_1^5$. The first subgroup is contained in $\bar{A}_1^3 \bar{D}_4$ and has already been considered. Let $Z = \bar{A}_1 A_1^5$. Then $\text{Out}_{M}(Z) \cong S_5$ acting naturally on the last five $A_1$ factors, and the diagonal $M$-irreducible subgroups follow using Lemma \ref{class}. 

Finally, it remains to consider $X = \bar{A}_1 A_1^2 B_3$ $(p=2)$. The reductive, maximal connected subgroups of $B_3$ when $p=2$ are $A_1 B_2$, $\bar{A}_3$ and $G_2$. The maximal subgroup $\bar{A}_1 A_1^3 B_2$ is contained in $\bar{A}_1 A_1 B_2^2$, as may be seen by considering its action on $V_{M}$ and thus has already been considered. The diagonal irreducible subgroups of $X$ are as given in Table \ref{E7tab}, noting that there is an involution in the Weyl group of $D_6$ swapping the second and third $A_1$ factors. The maximal subgroup $\bar{A}_1 A_1^2 \bar{A}_3$ was considered above and we are left to consider $Y = \bar{A}_1 A_1^2 G_2$. 

The reductive, maximal connected subgroups of $G_2$ when $p=2$ are $A_2$ and $\bar{A}_1 A_1$, by Theorem \ref{G2THM}. The subgroup $A_2$ is $B_3$-reducible by Lemma \ref{class} and hence $\bar{A}_1 A_1^2 A_2$ is $M$-reducible. The subgroup $\bar{A}_1 A_1$ is contained in $A_1 \bar{A}_1^2 < A_1 B_2 < B_3$ and it follows that $\bar{A}_1^2 A_1^2 A_1$ is contained in $\bar{A}_1^3 \bar{D}_4$ and conjugate to $E_7(\#\eseven{134}^{\{0,0\}})$. Finally, the diagonal irreducible subgroups of $Y$ follow in the same way as those of $X$. 

\subsection{$M_1 = \bar{A}_1 B_2 B_3$ $(E_7(\#\eseven{38}))$} 

The subgroup $\bar{A}_1 B_2 B_3$ is maximal when $p \neq 2$ and contained in $\bar{A}_1 B_5$ when $p=2$. We consider both cases in this section. By Lemma \ref{maxclassical}, the reductive, maximal connected subgroups of $\bar{A}_1 B_2 B_3$ are those listed in Table \ref{E7tab}. The subgroups $\bar{A}_1^3 B_3$ and $\bar{A}_1^3 A_1 B_2$ $(p \neq 2)$ are both contained in $\bar{A}_1^3 \bar{D}_4$, as seen by considering their action on $V_{M}$. Similarly, the subgroups $\bar{A}_1 A_1^2 B_3$ $(p=2)$ and $\bar{A}_1 A_1 B_2^2$ $(p=2)$ are both contained in $\bar{A}_1 A_1 B_4$. We now consider the three remaining reductive, maximal connected subgroups. 

First, let $X = \bar{A}_1 A_1 B_3$ $(p \geq 5$). The reductive, maximal connected subgroups again follow from Lemma \ref{maxclassical}. The subgroups that require further comment are $Y_1 = \bar{A}_1 A_1 \bar{A}_3$ and $Y_2 = \bar{A}_1 A_1 G_2$. First note that the diagonal subgroups of $Y_1$ and $Y_2$ follow in the same way as the diagonal subgroups of $X$. 

The reductive, maximal connected subgroups of $\bar{A}_3$ when $p \geq 5$ are $B_2$ and $A_1^2$ by Lemma \ref{maxclassical}. The subgroup $\bar{A}_1 A_1 B_2$ of $Y_1$ is $M$-reducible by Lemma \ref{class}. The subgroup $\bar{A}_1 A_1 A_1^2$ of $Y_1$ is contained in $\bar{A}_1^3 \bar{D}_4$ and conjugate to $E_7(\#\eseven{102}^{\{0\}})$.

Now consider $Y_2$. The $G_2$-irreducible subgroups are given by Theorem \ref{G2THM}. The subgroup $\bar{A}_2 < G_2$ is $B_3$-reducible by Lemma \ref{class}. The subgroup $\bar{A}_1 A_1$ is contained in $\bar{A}_1^2 A_1 < B_3$ and therefore $\bar{A}_1^2 A_1 A_1$ is contained in $\bar{A}_1^3 \bar{D}_4$ and conjugate to $E_7(\#\eseven{103}^{\{0,0\}})$. The diagonal subgroups of $\bar{A}_1 A_1 A_1$ $(p \geq 7)$ are clear and finish the study of the subgroups of $X$. 

Next, we consider $\bar{A}_1 B_2 \bar{A}_3$. The reductive, maximal connected subgroups $\bar{A}_1^3 \bar{A}_3$ and $\bar{A}_1 B_2^2$ are both $M$-reducible. The other reductive, maximal connected subgroups of $\bar{A}_1 B_2 \bar{A}_3$ are listed in Table \ref{E7tab} and have been considered previously, as seen from their action on $V_{M}$. 

Finally, we consider $X = \bar{A}_1 B_2 G_2$. All of the reductive, maximal connected subgroups of the form $\bar{A}_1 Y G_2$, where $Y$ is a reductive, maximal connected subgroup of $B_2$, have been considered above and the same is true for $\bar{A}_1 B_2 \bar{A}_1 A_1$. The subgroup $\bar{A}_1 B_2 \bar{A}_2$ is $M$-reducible by Lemma \ref{class}. That leaves us to consider $\bar{A}_1 B_2 A_2$ $(p=3)$ and $\bar{A}_1 B_2 A_1$ $(p \geq 7)$. By Theorem \ref{G2THM}, the maximal subgroup $A_1 < A_2 < G_2$ when $p=3$ is contained in $\bar{A}_1 A_1$ and hence we have already considered all of the reductive, maximal connected subgroups of $\bar{A}_1 B_2 A_2$ $(p=3)$. The only reductive, maximal connected subgroups of $\bar{A}_1 B_2 A_1$ $(p \geq 7)$ we have not already considered are diagonal subgroups, which are listed in Table \ref{E7tab}.  

\subsection{$M_1 = \bar{A}_1 \bar{A}_3^2$ $(E_7(\#\eseven{39}))$} 

By \cite[Table~10]{car}, we have $\text{Out}_{E_7}(M_1) \cong S_2 \times S_2$, where one generator swaps the two $\bar{A}_3$ factors and another generator induces a graph automorphism on both of them. The diagonal irreducible subgroups of $M_1$ are therefore as in Table \ref{E7tab}. The reductive, maximal connected subgroups of $\bar{A}_3$ are $B_2$ and $A_1^2$ $(p \neq 2)$ and so the non-diagonal reductive, maximal connected subgroups of $\bar{A}_1 \bar{A}_3^2$ are $\bar{A}_1 B_2 \bar{A}_3$ and $\bar{A}_1 A_1^2 \bar{A}_3$. The first subgroup is contained in $\bar{A}_1 B_2 B_3$ and the second is contained in $\bar{A}_1 A_1 B_4$ and hence both have already been considered.    

\subsection{$M_1 = \bar{A}_1 B_5$ $(E_7(\#\eseven{40}))$} 

As before, we use Lemma \ref{maxclassical} to find that the reductive, maximal connected subgroups of $\bar{A}_1 B_5$ are as listed in Table \ref{E7tab}, as well as the $M$-reducible subgroup $\bar{A}_1 D_5$. The subgroups $\bar{A}_1 A_1 \bar{D}_4$ and $\bar{A}_1^3 B_3$ are contained in $\bar{A}_1^3 \bar{D}_4$, and the subgroup $\bar{A}_1 B_2 \bar{A}_3$ is contained in $\bar{A}_1 B_2 B_3$. The subgroups $\bar{A}_1 A_1 B_4$ and $\bar{A}_1 B_2 B_3$ are both maximal subgroups of $M$ when $p \neq 2$; we considered their subgroups in previous sections. 

The only subgroup left to consider is $\bar{A}_1 A_1 = E_7(\#\eseven{267})$ $(p \geq 11)$. All proper irreducible connected subgroups are diagonal and as listed in Table \ref{E7tab}. 

\subsection{$M_1 = \bar{A}_1 A_1 C_3$ $(E_7(\#\eseven{41}))$} 

We note that there are two classes of $\bar{A}_1 A_1 C_3$ and they are distinguished by their action on $V_{56}$, for example, as given in Table \ref{E7tabcomps}.

By Lemma \ref{maxclassical}, the reductive, maximal connected subgroups of $C_3$ are $\bar{A}_1 B_2$, $A_1 A_1$ $(p \neq 2)$, $A_1$ $(p \geq 7)$, $A_3$ $(p=2)$ and $G_2$ $(p=2)$. The subgroup $\bar{A}_1^2 A_1 B_2$ is contained in $\bar{A}_1^3 \bar{D}_4$, as seen from its action on $V_{M}$; there are two classes of such subgroups, namely $E_7(\#\eseven{90}^{\{0,0\}})$ and $E_7(\#\eseven{92}^{\{0,0\}})$. By considering composition factors on $V_{56}$, we see that $\bar{A}_1^2 A_1 B_2$ is conjugate to $E_7(\#\eseven{90}^{\{0,0\}})$.   

The irreducible subgroups of $X =\bar{A}_1 A_1 A_1 A_1$ $(p \neq 2)$ are all diagonal and these can be found in Table \ref{A1A1A1A1E7diags}. We note that if $Y_1 \hookrightarrow X$ via $(1_a,1_b,1_b,1_c)$ then $Y_1$ acts as $(1,3,1) + (1,1,1)$ on $V_{M}$ when $p \geq 5$ and as $(1,T(3),1)$ on $V_{M}$ when $p =3$. Therefore, $Y_1$ is contained in $\bar{A}_1^3 A_1 A_1$ and conjugate to $E_7(\#\eseven{119}^{\{0,0,0,0\}})$ when $p \geq 5$, but $Y_1$ is $M$-reducible by Lemma \ref{class} when $p=3$. Similarly, if $Y_2 \hookrightarrow X$ via $(1_a,1_b,1_c,1_c)$, then $Y_2$ is conjugate to $E_7(\#\eseven{117}^{\{0,0,0,0\}})$ when $p \geq 5$, but $M$-reducible when $p=3$. Finally, if $Y_3 \hookrightarrow X$ via $(1_a,1_b,1_c,1_b)$, then $Y_3$ acts as $(1,2,2) + (1,2,0)$ and is conjugate to $E_7(\#\eseven{207}^{\{0,0\}})$.    

The diagonal subgroups of $X = \bar{A}_1 A_1 A_1$ $(p \geq 7)$ are easily seen to be as in Table \ref{E7tab}. Since $1 \otimes 5 = 6 + 4$, it follows that $Y \hookrightarrow X$ via $(1_a,1_b,1_b)$ is conjugate to $E_7(\#\eseven{262}^{\{0,0\}})$. 

By Lemma \ref{class}, the only $\bar{A}_1 A_1 C_3$-irreducible subgroups contained in $\bar{A}_1 A_1 A_3$ $(p=2)$ are diagonal subgroups. Again by Lemma \ref{class}, the only $\bar{A}_1 A_1 C_3$-irreducible subgroups contained in $\bar{A}_1 A_1 G_2$ $(p=2)$ are $\bar{A}_1 A_1 A_1  A_1$ and diagonal subgroups. The subgroup $A_1 A_1 < G_2$ acts on $V_{C_3}(\lambda_1)$ as $(1,1) + (0,2)$ and is hence a subgroup of $\bar{A}_1 B_2 < C_3$. Therefore, we have already considered $\bar{A}_1 A_1 A_1 A_1 < \bar{A}_1 A_1 G_2$.  

Finally, we need to prove that when $p=2$ the diagonal subgroups $E_7(\#\eseven{289})$ are $E_7$-irreducible. Let $X = A_1 A_3 \hookrightarrow \bar{A}_1 A_1 A_3 = E_7(\#\eseven{270})$ via $(1^{[r]},1^{[s]},100)$ $(rs=0)$ and suppose that $X$ is $E_7$-reducible. From Table \ref{E7tabcomps}, we have $V_{56} \downarrow \bar{A}_1  A_1 A_3 = (1,1,010) /$ $\!\! (0,3,000) /$ $\!\! (0,1,101)$. Therefore, $X$ has a $28$-dimensional composition factor on $V_{56}$. By Lemma \ref{wrongcomps}, there exists a subgroup $Z$ of type $A_1 A_3$ contained $L$-irreducibly in a Levi subgroup $L$, such that $X$ and $Z$ have the same composition factors on $V_{56}$. Since $Z$ has a $28$-dimensional composition factor on $V_{56}$, it follows from Table \ref{levie7} that $L'$ is of type $D_6$. Using Lemma \ref{maxclassical}, we find that there are two conjugacy classes of $D_6$-irreducible subgroups of type $A_1 A_3$ when $p=2$ both acting as $(1,010)$ on $V_{D_6}(\lambda_1)$. The composition factors of $Z$ on $V_{56}$ are thus $(1,010)^2 /$ $\!\! (1,101) /$ $\!\! (3,000)$ or $(1,010)^2 /$ $\!\! (2,010) /$ $\!\! (0,200) /$ $\!\! (0,010)^2 /$ $\!\! (0,002)$. In both cases $Z$ does not have the same composition factors as $X$ on $V_{56}$, which is a contradiction. Hence $X$ is $E_7$-irreducible.       

\subsection{$M_1 = \bar{A}_1 A_1 C_3$ $(E_7(\#\eseven{42}))$} \label{sec:E7A1A1C3}

The reductive, maximal connected subgroups of $E_7(\#\eseven{42})$ follow in the same way as those of $E_7(\#\eseven{41})$ in the previous section. We note that the subgroup $\bar{A}_1^2 A_1 B_2$ is conjugate to $E_7(\#\eseven{92}^{\{0,0\}})$, which can be seen by considering its composition factors on $V_{56}$. There is only one $M$-conjugacy class of $\bar{A}_1 A_1 A_1 A_1$ acting as $(1,1,2,1)$ on $(1,\lambda_1)$, since the graph automorphism of $D_6$ swaps the second and fourth $A_1$ factors.

We claim that $X = A_1 C_3 \hookrightarrow \bar{A}_1 A_1 C_3$ via $(1,1,100)$ is $E_7$-reducible when $p=2$. It then follows that both $A_1 A_3 \hookrightarrow \bar{A}_1 A_1 A_3$ via $(1,1,100)$ and $A_1 G_2 \hookrightarrow \bar{A}_1 A_1 G_2$ via $(1,1,10)$ are $E_7$-reducible. To prove the claim we first show that $\bar{A}_1 A_1 C_3$ is contained in the maximal subgroup $G_2 C_3$. Indeed, the $C_3$ factor of $G_2 C_3$ is contained in $A_5 < D_6$ and is generated by long root subgroups of $E_7$. By \cite[p.333, Table 3]{LS6}, the centraliser of such a subgroup $C_3$ is either $\bar{A}_1 A_1$ or $G_2$. By considering the composition factors on $V_{56}$, we see that the $C_3$ factor of $E_7(\#\eseven{42})$ has connected centraliser $G_2$ and therefore $\bar{A}_1 A_1 C_3 < G_2 C_3$. Now, by Theorem \ref{G2THM}, the subgroup $A_1 \hookrightarrow \bar{A}_1 A_1 < G_2$ via $(1,1)$ is $G_2$-reducible when $p=2$ and hence $X$ is $G_2 C_3$-reducible, proving the claim.   

Finally, we prove the diagonal subgroups $X = E_7(\#\eseven{298})$ are $E_7$-irreducible when $p=2$. The dimensions of the $X$-composition factors of $V_{56}$ are $24, 12, 6^2, 4^2$. Using a similar argument to that given for $E_7(\#\eseven{289})$ in the previous section we find that no subgroup of a Levi factor has the same composition factors as $X$ on $V_{56}$. Therefore $X$ is $E_7$-irreducible by Lemma \ref{wrongcomps}.   

This completes the case $M = \bar{A}_1 D_6$. 

\section{$M = \bar{A}_2 A_5$  $(E_7(\#\eseven{31}))$}

By Lemma \ref{maxclassical}, the reductive, maximal connected subgroups of $A_5$ are $A_2 A_1$, $C_3$, $A_3$ $(p \neq 2)$ and $A_2$ $(p \neq 2)$. Similarly, by Lemma \ref{maxclassical}, the only reductive, maximal connected subgroup of $\bar{A}_2$ is $A_1$ $(p \neq 2)$. It follows that the reductive, maximal connected subgroups of $\bar{A}_2 A_5$ are as in Table \ref{E7tab}. We consider these in the following sections. Note that all reductive, maximal connected subgroups of $A_1 A_5$ $(p \neq 2)$ are contained in $\bar{A}_2 X$ for some reductive, maximal connected subgroups $X$ of $A_5$. We therefore give them no further consideration. 

\subsection{$M_1 = \bar{A}_2 A_2 A_1$  $(E_7(\#\eseven{300}))$} \label{sec:E7A2A2A1}

The $M$-irreducible subgroups contained in $M_1$ are straightforward to find and given in Table \ref{E7tab}, along with the subgroup $X = A_2 A_1 \hookrightarrow M_1$ via $(10,01,1)$ when $p=3$. We need to show that $X$ is conjugate to $E_7(\#\eseven{184}^{\{0,0,0\}}) < \bar{A}_1 D_6$ when $p \neq 3$ and $E_7$-reducible when $p=3$. We also show that $Y = A_1 A_1 A_1 < \bar{A}_2 A_2 A_1$ is conjugate to $E_7(\#\eseven{275}^{\{0,0\}}) < \bar{A}_1 D_6$.

We first note that $\bar{A}_2 A_2 A_1$ is a subgroup of $A_1 F_4$. Indeed, consider the maximal subgroup $\bar{A}_2 A_2$ of $F_4$. The $\bar{A}_2$ factor generated by long root subgroups of $F_4$ is generated by long root subgroups of $E_7$ and so $A_1 \bar{A}_2 A_2 < \bar{A}_2 C_{E_7}(\bar{A}_2)^\circ = \bar{A}_2 A_5$, by Lemma \ref{rootycentralisers}. Theorem \ref{F4THM} shows that $Z = A_2 \hookrightarrow \bar{A}_2 A_2 < F_4$ via $(10,01)$ is conjugate to $A_2 < \bar{D}_4$ embedded via $V_{A_2}(11)$ when $p \neq 3$ but $Z$ is $F_4$-reducible when $p=3$. It follows that $X < \bar{D}_4 C_{E_7}(\bar{D}_4)^\circ = \bar{A}_1^3 \bar{D}_4 < \bar{A}_1 D_6$ when $p \neq 3$ and $X$ is $A_1 F_4$-reducible when $p=3$. Comparing composition factors when $p \neq 3$ shows that $X$ is conjugate to $E_7(\#\eseven{184}^{\{0,0,0\}})$, as required. Also by Theorem \ref{F4THM}, the irreducible subgroup $A_1 A_1 < \bar{A}_2 A_2 < F_4$ is contained in $\bar{A}_1 C_3$. Therefore, $Y < \bar{A}_1 C_3 A_1 < F_4 A_1$. Comparing composition factors shows that $Y$ is conjugate to $E_7(\#\eseven{275}^{\{0,0\}})$.

When $p=2$ the diagonal subgroups $E_7(\#\eseven{307})$ and $E_7(\#\eseven{308})$ contain no proper $E_7$-irreducible subgroups and so we need to prove that they are $E_7$-irreducible. From Table \ref{E7tabcomps}, the composition factors of $E_7(\#\eseven{300}) = \bar{A}_2 A_2 A_1$ acting on $V_{56}$ are $(10,10,1) /$ $\!\! (01,01,1) /$ $\!\! (00,11,1) /$ $\!\! (00,00,3)$. Let $X$ be the diagonal subgroup embedded via $(10,10,1)$. We prove that $X$ is $E_7$-irreducible; the other cases are similar and easier. Suppose that $X$ is $E_7$-reducible. By Lemma \ref{wrongcomps}, there exists a subgroup $Z$ of type $A_1 A_2$ contained $L$-irreducibly in a Levi subgroup $L$, such that $X$ and $Z$ have the same composition factors on $V_{56}$. The $Z$-composition factors of $V_{56}$ are thus $(3,00) /$ $\!\! (1,20) /$ $\!\! (1,11) /$ $\!\! (1,10)^2 /$ $\!\! (1,02) /$ $\!\! (1,01)^2 $; in particular, the dimensions of the composition factors are $16, 6^6, 4$. Using Table \ref{levie7} we see that $L'$ has type $D_6$, $A_1 A_5$ or $A_5$. Using Lemma \ref{class} we find all of the $L'$-irreducible subgroups of type $A_1 A_2$. If $L'$ has type $D_6$ then $Z$ acts as $( (0,00) | (2,00) | (0,00) ) + (0,11)$ on $V_{D_6}(\lambda_1)$. Therefore $Z$ has a trivial composition factor on $V_{56}$, a contradiction. If $L'$ has type $A_1 A_5$ then $Z$ is contained in the subgroup $A_1 A_1 A_2$ where the $A_1 A_2 < A_5$ acts on $V_{A_5}(\lambda_1)$ as $(1,10)$. Therefore the $Z$-composition factors of $V_{56}$ are $(1^{[r]} \otimes 1^{[s]},10) /$ $\!\! (1^{[r]} \otimes 1^{[s]},01) /$ $\!\! (3^{[s]},00) /$ $\!\! (1^{[s]},11) /$ $\!\! (1^{[s]},10) /$ $\!\! (1^{[s]},01)$ for some $r,s$ with $rs=0$. These are not the same as the $X$-composition factors of $V_{56}$ for any $r,s$, which is a contradiction. Finally, suppose that $L'$ is of type $A_5$, of which there are two such conjugacy classes. Then $Z$ acts as $(1,10)$ on $V_{A_5}(\lambda_1)$. Using the previous case and Table \ref{levie7} we see that the $Z$-composition factors of $V_{56}$ are not the same as the $X$-composition factors. This final contradiction shows that $X$ is $E_7$-irreducible.  

\subsection{$M_1 = \bar{A}_2 C_3$  $(E_7(\#\eseven{301}))$}

The only irreducible subgroup of $\bar{A}_2$ is $A_1$ when $p \neq 2$. We claim that $A_1 C_3$ is contained in $\bar{A}_1 D_6$. To prove this, we first note that by \cite[p.333, Table~3]{LS6}, the centraliser in $E_7$ of the $C_3$ factor is $G_2$ and therefore $\bar{A}_2 C_3$ is contained in $G_2 C_3$. By Theorem \ref{G2THM}, the irreducible subgroup $A_1$ of $\bar{A}_2 < G_2$ is contained in $\bar{A}_1 A_1 < G_2$. It follows that $A_1 C_3$ is therefore contained in $\bar{A}_1 A_1 C_3 < \bar{A}_1 D_6$ and by considering composition factors we conclude that $A_1 C_3$ is conjugate to $E_7(\#\eseven{294}^{\{0,0\}})$. 

From Lemma \ref{maxclassical}, the reductive, maximal connected subgroups of $C_3$ are $A_1 B_2$, $A_1 A_1$ $(p \neq 2)$, $A_1$ $(p \geq 7)$, $A_3$ $(p =2)$ and $G_2$ $(p=2)$. The action of $A_1 B_2$ on $V_{A_5}(\lambda_1)$ is reducible and hence $\bar{A}_2 A_1 B_2$ is $M$-reducible by Lemma \ref{class}. We also note that the maximal subgroup $A_1 A_1$ $(p \neq 2)$ acts as $(2,1)$ on $V_{A_5}(\lambda_1)$ and therefore $\bar{A}_2 A_1 A_1$ is conjugate to $E_7(\#\eseven{306})$.  

The subgroup $\bar{A}_2 A_3$ is a maximal connected subgroup of $M$ when $p \neq 2$ and so we consider the subgroups of $\bar{A}_2 A_3$ $(p=2)$ in the next section.    

Finally, consider $X = \bar{A}_2 G_2$ $(p=2)$. The reductive, maximal connected subgroups of $G_2$ when $p=2$ are $M$-reducible by Lemma \ref{class} and there are no reductive, maximal connected subgroups of $\bar{A}_2$ when $p=2$. Thus there are no proper $M$-irreducible subgroups of $X$. It remains to prove that $X$ is $E_7$-irreducible. From Table \ref{E7tabcomps}, we have $V_{56} \downarrow X = (10,10) /$ $\!\! (01,10) /$ $\!\! (00,20) /$ $\!\! (00,10)^2 /$ $\!\! (00,00)^2$. As in previous cases, it is straightforward to show that there are no Levi subgroups with an irreducible subgroup $A_2 G_2$ having the same composition factors as $X$ on $V_{56}$. Therefore Lemma \ref{wrongcomps} implies that $X$ is $E_7$-irreducible.  

\subsection{$M_1 = \bar{A}_2 A_3$  $(E_7(\#\eseven{302}))$} \label{e7a2a3}

The subgroup $\bar{A}_2 A_3$ is a maximal connected subgroup of $M$ when $p \neq 2$ and contained in $\bar{A}_2 C_3$ when $p=2$. In this section we consider both cases.

Both reductive, maximal connected subgroups of $A_3$ act reducibly on $V_{A_3}(010)$ and hence on $V_{A_5}(\lambda_1)$. Therefore, the only proper irreducible subgroup to consider is $X = A_1 A_3$ $(p \neq 2)$, where the $A_1$ factor is irreducibly embedded in $A_2$. We need to prove that $X$ is $E_7$-irreducible when $p \neq  2$ and that $M_1$ is $E_7$-irreducible when $p=2$. 

First we consider $M_1$ when $p=2$. Suppose that $M_1$ is $E_7$-reducible. By Lemma \ref{wrongcomps}, there exists a subgroup $Z$ of type $A_2 A_3$ contained $L$-irreducibly in a Levi subgroup $L$, such that $X$ and $Z$ have the same composition factors on $V_{56}$. It follows that $L'$ has no simple factors of type $A_1$ and has rank at least $5$. Therefore $L'$ is of type $A_2 A_3$. Using Table \ref{E7tabcomps} and Table \ref{levie7}, we see that the composition factors of the Levi subgroup $A_2 A_3$ are not the same as the composition factors of $M_1$ on $V_{56}$. This contradiction proves that $M_1$ is $E_7$-irreducible.

Now we consider $X$. When $p \geq 5$, the action of $X$ on $L(E_7)$ has no trivial composition factors and so $X$ is $E_7$-irreducible by Corollary \ref{notrivs}. When $p=3$, Table \ref{E7tabcomps} shows that $V_{56} \downarrow X = (2,010)^2 /$ $\!\! (0,200) /$ $\!\! (0,002)$.  As in previous cases, it is straightforward to show that there are no Levi subgroups with a subgroup $A_1 A_3$ having the same composition factors as $X$ on $V_{56}$. Therefore  Lemma \ref{wrongcomps} implies that $X$ is $E_7$-irreducible.  

\subsection{$M_1 = \bar{A}_2 A_2$ $(p \neq 2)$  $(E_7(\#\eseven{303}))$}

The only non-simple proper irreducible connected subgroup to consider is $X = A_1 A_2 = E_7(\#\eseven{314})$, since the irreducibly embedded subgroup $A_1$ of the second $A_2$ factor acts reducibly on $V_{A_5}(\lambda_1)$. Since $X$ contains no proper $E_7$-irreducible subgroups we need to show that it is $E_7$-irreducible. 

When $p \geq 5$ we see from Table \ref{E7tabcomps} that there are no trivial $X$-composition factors of $L(E_7)$ and hence $X$ is $E_7$-irreducible by Corollary \ref{notrivs}. Now let $p=3$ and suppose that $X$ is $E_7$-reducible.  From Table \ref{E7tabcomps} we have $V_{56} \downarrow X = (2,20) /$ $\!\! (2,02) /$ $\!\! (0,30) /$ $\!\! (0,11)^2 /$ $\!\! (0,03)$.  By Lemma \ref{wrongcomps}, there exists a subgroup $Z$ of type $A_1 A_2$ contained $L$-irreducibly in a Levi subgroup $L$, such that $X$ and $Z$ have the same composition factors on $V_{56}$. The dimensions of the $Z$-composition factors of $V_{56}$ are $18^2, 7^2, 3^2$. Using Table \ref{levie7}, we deduce that there is no Levi subgroup containing such a subgroup $Z$. Therefore $X$ is $E_7$-irreducible.  

This completes the case $M = \bar{A}_2 A_5$. 

\section{$M = A_7$  $(E_7(\#\eseven{22}))$} \label{e7a7}

From Lemma \ref{maxclassical}, we find that the reductive, maximal connected subgroups of $A_7$ are $C_4$, $D_4$ $(p \neq 2)$ and $A_1 A_3$. The subgroup $C_4$ is $E_7$-reducible and hence we do not need to consider any of its subgroups.  

Next we consider $D_4$ $(p \neq 2)$. By \cite[Lemma 2.15]{clss}, we have $\text{Out}_{E_7}(D_4) \cong S_3$. It follows that the only reductive, maximal connected subgroup of $D_4$ which is $A_7$-irreducible is $A_2$ when $p \geq 5$, acting as $11$ on $V_{A_7}(\lambda_1)$. Such a subgroup $A_2$ has already been considered in \cite{tho1} and shown to be contained in $\bar{A}_2 A_5$. By comparing composition factors, we find that it is conjugate to $E_7(\#\eseven{24})$. The only proper irreducible connected subgroup of $A_2$ is $A_1$. Such a subgroup $A_1$ acts as $4 + 2$ on $V_{A_7}(\lambda_1)$ and is therefore $M$-reducible.

We now show that $X = A_1 A_3$ is contained in $\bar{A}_2 A_5$. Let $Y$ be the $A_3$ factor of $X$. Then $Y$ acts as $100^2$ on $V_{A_7}(\lambda_1)$. Thus $Y$ is contained in $\bar{A}_3^2$, a Levi subgroup of $A_7$. By \cite[p.333, Table~2]{LS6}, we have $C_{E_7}(\bar{A}_3)^\circ = \bar{A}_1 \bar{A}_3$ and it follows that $\bar{A}_3^2 < C_{E_7}(\bar{A}_1)^\circ = D_6$. By considering the composition factors of $Y$ on $V_{56}$ it follows that $Y$ acts as $010^2$ on $V_{D_6}(\lambda_1)$ and is hence contained in an $A_5$-Levi subgroup of $D_6$, acting as $010$ on $V_{A_5}(\lambda_1)$. There are two $E_7$-conjugacy classes of $A_5$-Levi subgroups and again by considering the composition factors of $Y$, it follows that $Y$ is contained in an $A_5$-Levi whose connected centraliser is $\bar{A}_2$. Therefore $X < \bar{A}_2 A_5$. 

When $p \neq 2$, we find that $X$ is conjugate to $E_7(\#\eseven{313})$, as may be seen by considering their composition factors on $V_{56}$. When $p=2$, there are no $\bar{A}_2 A_5$-irreducible subgroups of type $A_1 A_3$ contained in $\bar{A}_2 A_5$. Therefore $X$ is $\bar{A}_2 A_5$-reducible and hence $G$-reducible.   

\section{$M = G_2 C_3$  $(E_7(\#\eseven{32}))$}\label{e7g2c3}

The reductive, maximal connected subgroups of $G_2$ are given by Theorem \ref{G2THM} and the reductive, maximal subgroups of $C_3$ are given by Lemma \ref{maxclassical}. The lattice of $M$-irreducible subgroups now easily follows with further use of Lemma \ref{maxclassical} and Lemma \ref{class}. It remains for us to show the claimed conjugacies between subgroups of $M$ and subgroups of previously considered reductive, maximal connected subgroups of $E_7$. 

The $G_2$ factor of $G_2 C_3$ is contained in a Levi subgroup $D_4$ and hence the subgroups of $G_2$ generated by long root subgroups of $G_2$ are generated by long root subgroups of $E_7$. It follows that $\bar{A}_2 C_3$ is conjugate to $E_7(\#\eseven{301}) < \bar{A}_2 A_5$ and $\bar{A}_1 A_1 C_3$ is conjugate to $E_7(\#\eseven{42}) < \bar{A}_1 D_6$. By Theorem \ref{G2THM}, the irreducible subgroup $A_1$ of $A_2 < G_2$ $(p=3)$ is contained in $\bar{A}_1 A_1$ and so when $p=3$, the subgroup $A_1 C_3$ is conjugate to $E_7(\#\eseven{294}^{\{1,0\}})$. The subgroup $G_2 \bar{A}_1 B_2$ is contained in $\bar{A}_1 D_6$ and conjugate to $E_7(\#\eseven{253})$. 

 We need to prove that some diagonal subgroups of $X = A_1 A_1 A_1 = E_7(\#\eseven{323})$  $(p \geq 7)$ have also been seen previously. Firstly, the subgroup $A_1 A_1$ embedded via $(1_a,1_b,1_b)$ is contained in $A_1 \bar{A}_1 B_2$ since $2 \otimes 1 = 3 + 1$. Secondly, we claim the subgroup $Y = A_1 A_1$ embedded via $(1_a,1_a,1_b)$ is contained in $\bar{A}_1 D_6$ also. To prove this, first consider the maximal connected subgroup $A_1 A_1 G_2$ of $A_1 F_4$. Then $C_{E_7}(G_2)^\circ = C_3$ and hence $A_1 A_1 G_2 < G_2 C_3$. By considering composition factors on $V_{56}$ we see that the $A_1$ factor of $A_1 F_4$ is the third $A_1$ factor of $X$. Appealing to Theorem \ref{F4THM} shows that $A_1 \hookrightarrow A_1 A_1 < A_1 G_2 < F_4$ via $(1,1)$ is conjugate to $A_1 \hookrightarrow \bar{A}_1 A_1 < \bar{A}_1 C_3 < F_4$ via $(1,1)$ since both have identification number $F_4(\#\ffour{8}^{\{0,0\}})$. Therefore $Y$ is contained in $A_1 \bar{A}_1 C_3$, which is a subgroup of $\bar{A}_1 D_6$ by Lemma \ref{rootycentralisers}.  

\section{$M = A_1 F_4$  $(E_7(\#\eseven{33}))$}   

Theorem \ref{F4THM} gives the lattice structure of the $M$-irreducible subgroups of $M$. By \cite[p.333, Table~3]{LS6}, $C_{E_7}(B_4)^\circ = \bar{A}_1 A_1$ and hence $A_1 B_4$ is conjugate to $E_7(\#\eseven{190}^{\{0,0\}}) < \bar{A}_1 D_6$. Similarly, $A_1 \bar{A}_1 C_3$ is a subgroup of $\bar{A}_1 D_6$ and $A_1 \bar{A}_2 A_2$ is a subgroup of $\bar{A}_2 A_5$. In Section \ref{e7g2c3}, we also proved that $A_1 A_1 G_2$ is contained in $G_2 C_3$. We are therefore left to consider $A_1 C_4$ $(p=2)$, $A_1 G_2$ $(p=7)$ and $A_1 A_1$ $(p \geq 13)$. The diagonal subgroups of $A_1 A_1$ $(p \geq 13)$ follow immediately. 

From Theorem \ref{F4THM}, we see that the $M$-irreducible maximal connected subgroups of $A_1 C_4$ are $A_1 D_4$, $A_1 B_2^2$ and $A_1 \bar{A}_1 C_3$. The second subgroup is contained in $A_1 B_4$ and the third subgroup is contained in $\bar{A}_1 D_6$ as well. Hence both have already been considered. The $M$-irreducible maximal connected subgroups of $A_1 D_4$ are $A_1 A_1^4$ and $A_1 A_2$. The former is a subgroup of $A_1 B_4$ and the latter a subgroup of $A_1 \bar{A}_2 A_2$. Hence both have already been considered. 

Similarly, using Theorem \ref{F4THM} we see that all $F_4$-irreducible connected subgroups of $G_2$ $(p=7)$ are contained in either $\bar{A}_2 A_2$ or $\bar{A}_1 C_3$ and so we have already considered all $M$-irreducible subgroups of $A_1 G_2$ ${(p=7)}$. 
 
\section{$M = A_1 G_2$  $(p \neq 2)$ $(E_7(\#\eseven{34}))$}   

The lattice structure of $M$-irreducible subgroups follows from Theorem \ref{G2THM}. We claim that $A_1 A_2$ (where the factor $A_2$ is generated by long root subgroups of $G_2$) is contained in $\bar{A}_2 A_5$ and conjugate to $E_7(\#\eseven{309}^{\{0,0\}})$. To prove this, let $X$ be the $A_2$ factor. The $G_2$ factor of $M$ is contained in $A_6$, acting on $V_{A_6}(\lambda_1)$ as $V_{G_2}(10)$, by the construction in \cite[3.12]{se2}. It follows that $X$ acts as $10 + 01 + 00$ on $V_{A_6}(\lambda_1)$ and therefore $X$ is contained in a Levi subgroup $\bar{A}_2^2$ of $A_6$ embedded via $(10,01)$. Furthermore, $X$ is contained in $A_5$ acting as $10^2$, since $\bar{A}_2^2$ acts as $(10,00) + (00,01)$ on $V_{A_5}(\lambda_1)$. We now have $C_{E_7}(X)^\circ \geq \bar{A}_2 A_1$ and since $X \bar{A}_2 A_1$ is $E_7$-irreducible (subgroup $E_7(\#\eseven{300})$) this must be an equality. Therefore $A_1 X < \bar{A}_2 C_{E_7}(\bar{A}_2)^\circ = \bar{A}_2 A_5$ and comparing composition factors finishes the proof of the claim. 

Now consider $X = A_1 A_1 A_1$. The projection of $X$ to $G_2$ is the centraliser in $G_2$ of a semisimple involution. Therefore $X$ centralises a semisimple involution $t$ of $E_7$ and thus $X$ is contained in $\bar{A}_1 D_6$ or $A_7$ by Lemma \ref{centralisers}. Since $X$ contains a simple $E_7$-irreducible subgroup, $X$ is $E_7$-irreducible. As $A_7$ does not contain an $E_7$-irreducible subgroup of type $A_1 A_1 A_1$ it follows that $X$ is contained in $\bar{A}_1 D_6$. By considering the $X$-composition factors of $V_{56}$, we see that $X$ is conjugate to $E_7(\#\eseven{206}^{\{0,0\}})$. 

The only reductive, maximal connected subgroup of $A_1 A_2$ $(p=3)$ is $A_1 A_1$ and by Theorem \ref{G2THM} this is a subgroup of $A_1 A_1 A_1$ and has thus already been considered. 

This completes the proof of Theorem \ref{E7THM}.  

%% file: E8.tex
\chapter{Irreducible subgroups of $G = E_8$} \label{secE8}

In this section we prove Theorem \ref{MAINTHM} when $G$ is of type $E_8$; we classify the $E_8$-irreducible connected subgroups thus proving the following theorem. 

\begin{thm} \label{E8THM} 
Let $X$ be an $E_8$-irreducible connected subgroup of $E_8$. Then $X$ is conjugate to exactly one subgroup of Table \ref{E8tab} and each subgroup in Table \ref{E8tab} is $E_8$-irreducible. Moreover, Table \ref{E8tab} gives the lattice structure of the irreducible connected subgroups of $E_8$.   
\end{thm}

As in the previous sections we consider each of the reductive, maximal connected subgroups of $E_8$ in turn. The simple connected $E_8$-irreducible subgroups are classified in \cite[Theorem 3]{tho1} and \cite[Theorem 6]{tho2} and we will use these without reference throughout. By Theorem \ref{maximalexcep}, the reductive, maximal connected subgroups of $E_8$ are $D_8$, $\bar{A}_1 E_7$, $\bar{A}_2 E_6$, $A_8$, $\bar{A}_4^2$, $G_2 F_4$, $B_2$ $(p \geq 5)$, $A_1 A_2$ $(p \geq 5)$ and $A_1$ (3 classes, $p \geq 23, 29, 31$). 

We remind the reader that throughout the proof we will only make reference to the $E_8$-irreducibility of an $M$-irreducible subgroup when it does not properly contain an $E_8$-irreducible subgroup. Since the simple $E_8$-irreducible connected subgroups are known we only need to consider a small number of cases.

\section{$M = D_8$ $(E_8(\#\eeight{43}))$} \label{D8inE8}

The reductive, maximal connected subgroups of $M$ are given in Lemma \ref{maxclassical}. They are $\bar{A}_1^2 D_6$, $\bar{D}_4^2$, $\bar{A}_3 D_5$, $B_7$, $A_1 B_6$ $(p \neq 2)$, $B_2 B_5$ $(p \neq 2)$, $B_3 B_4$ $(p \neq 2)$, $B_2^2$ $(p \neq 2)$ (2 classes), $A_1 C_4$ (2 classes) and $B_4$ (2 classes). In the cases where there are two conjugacy classes these are distinguished by their action on $L(E_8)$, which is given in Table \ref{E8tabcomps}. 

Firstly we prove that two of the maximal connected subgroups are $E_8$-reducible when $p=2$. In \cite[Lemma~7.4]{tho1}, it is shown that one of the classes of $B_4$ subgroups is $E_8$-reducible, denoted $B_4(\ddagger)$. In the current notation this is subgroup $E_8(\#\eeight{46})$ for which $p \neq 2$. We also claim one of the classes of $A_1 C_4$ is $E_8$-reducible when $p=2$. Let $A_1 C_4 (\ddagger)$ denote the subgroup $A_1 C_4$ with composition factors $(4,0) /$ $\!\! (2,\lambda_2) /$ $\!\! (2,0)^2 /$  $\!\! (0,\lambda_2)^2 /$ $\!\! (0,\lambda_4) /$ $\!\! (0,0)^2$ on $V_{D_8}(\lambda_7)$ when $p=2$ (this is subgroup $E_8(\#\eeight{116})$ when $p \neq 2$). By \cite[p.333, Table~3]{LS6}, the connected centraliser of the $C_4$ factor of $A_1 C_4(\ddagger)$ is $\bar{A}_1 U_5$, where $U_5$ is a $5$-dimensional connected unipotent subgroup. Therefore $X < U_5 \bar{A}_1 C_4$, which by the Borel-Tits Theorem \cite[Th\'{e}or\`{e}me~2.5]{BT} is contained in a parabolic subgroup of $E_8$. Therefore $X$ is $E_8$-reducible. 

There are a very large number of $E_8$-irreducible subgroups contained in $D_8$. We therefore suppress many of the routine parts of constructing the lattice of irreducible connected subgroups since we have been more explicit in the proofs of the previous theorems, especially the case $\bar{A}_1 D_6 < E_7$. For example, the use of Lemma \ref{maxclassical} to find the lattice of $D_8$-irreducible connected subgroups is mainly omitted, as are most of the considerations of when two subgroups contained in different reductive, maximal connected subgroups of $D_8$ are conjugate. Our use of Lemma \ref{class} to remove any $D_8$-reducible classes inside reductive, maximal connected subgroups of $D_8$ will also be implicit.

When finding the $E_8$-conjugacy classes of diagonal subgroups contained in $X < D_8$ we give only the required information concerning $\text{Out}_{E_8}(X)$. On many occasions $\text{Out}_{E_8}(X) \cong \text{Out}_{D_8}(X)$; in such cases we will not be explicit. 

\subsection{$M_1 = \bar{A}_1^2 D_6$ $(E_8(\#\eeight{107}))$}  

By \cite[Table~11]{car}, we have $\text{Out}_{E_8}(M_1) \cong S_2$, where the involution acts simultaneously as a graph automorphism of $\bar{A}_1^2$ and $D_6$. Up to $E_8$-conjugacy we may therefore choose to take just one representative of each non-conjugate pair of $D_6$-irreducible subgroups which are conjugate in $D_6.2$. This leads to the reductive, maximal connected subgroups of $\bar{A}_1^2 D_6$ given in Table \ref{E8tab}. We will consider each of them in turn.  

\subsubsection{$M_2 = \bar{A}_1^4 \bar{D}_4$ $(E_8(\#\eeight{117}))$}

By \cite[Table~11]{car}, we have $\text{Out}_{E_8}(M_2) \cong S_4$ acting naturally on the four $A_1$ factors and inducing the full outer automorphism group of $\bar{D}_4$. Therefore, we have only one $E_8$-conjugacy class of each of the subgroups $\bar{A}_1^4 B_3$ and $\bar{A}_1^4 A_1 B_2$ contained in $M_2$. Moreover, the stabiliser of $B_3$ and $A_1 B_2$ under the action of $S_4$ is $\text{Dih}_8$. Therefore $\text{Out}_{E_8}(\bar{A}_1^4 B_3) \cong \text{Out}_{E_8}(\bar{A}_1^4 A_1 B_2) \cong \text{Dih}_8$. The irreducible subgroups $G_2$ and $A_2$ are normalised by all outer automorphisms of $\bar{D}_4$ and so $\text{Out}_{E_8}(\bar{A}_1^4 G_2) \cong \text{Out}_{E_8}(\bar{A}_1^4 A_2) \cong S_4$. Finding most of the diagonal subgroups contained in $M_2$ now follows; we give further details below where needed. 

From \cite[Table~10.3]{LS1} we see $\text{Out}_{E_8}(\bar{A}_1^8) \cong \text{AGL}_3(2)$. The action on the eight factors is described in the proof of \cite[Theorem~6]{tho2} and the non-simple diagonal subgroups follow using the same method as for the simple ones.    

We note two intricacies in finding the diagonal subgroups of $E_8(\#\eeight{189}) = \bar{A}_1^4 A_1 A_1$ $(p \geq 5)$. There are two classes of pairs $\bar{A}_1^2 < \bar{A}_1^4$, which may be represented by the first and second $A_1$ factors and the first and third $A_1$ factors, for example. When taking a diagonal subgroup, say $X$, in the first and second $A_1$ factors via $(1,1)$ it is then conjugate to the fifth $A_1$ factor, say $Y$. It then follows that there is an involution in $D_8$ that swaps $X$ and $Y$. This leads to the condition ``if $r=s$ then $r < t$'' in the definition of $E_8(\#\eeight{213})$ in Table \ref{barA14A1A1p5diags}. The second thing to note is that a diagonal subgroup of the fifth and sixth $A_1$ factors via $(1,1)$ is conjugate to a subgroup $A_1$ contained in $A_2 < \bar{D}_4$. Indeed, the irreducible subgroup $A_1$ of $A_2$ acts as $4 + 2$ on $V_{\bar{D}_4}(\lambda_1)$. Since this subgroup $A_2$ is centralised by a triality automorphism of $\bar{D}_4$ it follows that $\text{Out}_{E_8}(E_8(\#\eeight{211}^{\{0,0\}})) \cong S_4$ acting naturally on the first four $A_1$ factors.    

Finally, we consider the subgroup $E_8(\#\eeight{190}) = \bar{A}_1^4 A_1^3 < \bar{A}_1^4 A_1 B_2$ when $p =2$. Since $\text{Out}_{\bar{D}_4}(A_1^3) \cong S_3$ it follows that $\text{Out}_{E_8}(E_8(\#\eeight{190})) \cong S_4 \times S_3$, where the $S_4$ direct factor acts naturally on the first four factors and the $S_3$ direct factor acts naturally on the final three factors. 

\subsubsection{$M_2 = \bar{A}_1^2 B_5$ $(E_8(\#\eeight{118}))$}

Since $B_5$ is centralised by a graph automorphism of $D_6$ we have $\text{Out}_{E_8}(M_2) \cong S_2$. This is enough to determine the classes of diagonal subgroups of all subgroups of $M_2$.  

\subsubsection{$M_2 = \bar{A}_1^2 \bar{A}_3^2$ $(E_8(\#\eeight{119}))$}

Here $\text{Out}_{E_8}(M_2) \cong S_2 \times S_2$ by \cite[Table~11]{car}, with an involution swapping the two $\bar{A}_3$ factors and another involution simultaneously swapping the two $\bar{A}_1$ factors whilst fixing one $\bar{A}_3$ factor and acting as a graph automorphism on the other. The diagonal subgroups of $M_2$ now follow.    

Because $\text{Out}_{\bar{A}_3}(A_1^2) \cong S_2$ it follows that $\text{Out}_{E_8}(\bar{A}_1^2 A_1^2 \bar{A}_3) \cong S_2 \times S_2$, containing an involution swapping the first and second $\bar{A}_1$ factors whilst acting as a graph automorphism of the $\bar{A}_3$ factor and another involution swapping the third and fourth $A_1$ factors. For the diagonal subgroup $X = A_1^3 \bar{A}_3 \hookrightarrow \bar{A}_1^2 A_1^2 \bar{A}_3$ via $(1_a,1_a,1_b,1_c,100)$ we have $\text{Out}_{E_8}(X) \cong S_3$ acting naturally on the three $A_1$ factors. Finally, note that $\text{Out}_{E_8}(\bar{A}_1^2 A_1^4) \cong S_2 \times S_4$ and taking a diagonal subgroup $X = A_1^5$ via $(1_a,1_a,1_b,1_c,1_d,1_e)$ yields $\text{Out}_{E_8}(X) \cong S_5$. 

It remains to prove that the diagonal subgroups $E_8(\#\eeight{373})$ are irreducible when $p=2$ as they do not properly contain any $E_8$-irreducible subgroups. Let $X = E_8(\#\eeight{373}) = A_1 \bar{A}_3 \hookrightarrow \bar{A}_1^2 A_1^2 \bar{A}_3 = E_8(\#\eeight{355}) =  Y$ via $(1^{[r]},1^{[s]},1^{[t]},1^{[u]},100)$ where $rt=0$, $r < s$ and $t < u$. From Table \ref{E8tabcomps}, we have \begin{align*} L(E_8) \downarrow Y = \ & (2,0,0,0,000) / (1,1,2,0,000) / (1,1,0,2,000) / (1,1,0,0,010) / \\ \ & (1,1,0,0,000)^2 / (1,0,1,1,100) / (1,0,1,1,001) / (0,2,0,0,000) / \\ \ & (0,1,1,1,100) / (0,1,1,1,001) / (0,0,2,2,000) / (0,0,2,0,010) / \\ \ & (0,0,2,0,000)^2  / (0,0,0,2,010) / (0,0,0,2,000)^2 / (0,0,0,0,101) / \\ \ & (0,0,0,0,010)^2 / (0,0,0,0,000)^6. \end{align*} Therefore $X$, has at most seven trivial composition factors and at least two composition factors of dimension 32. Looking for a contradiction we suppose that $X$ is $E_8$-reducible. By Lemma \ref{wrongcomps}, there exists a subgroup $Z$ of type $A_1 A_3$ contained $L$-irreducibly in a Levi subgroup $L$, such that $X$ and $Z$ have the same composition factors on $L(E_8)$. By considering the number of trivial composition factors on $L(E_8)$ of each Levi subgroup that contains an irreducible subgroup $A_1 A_3$ (using Table \ref{levie8} and Lemma \ref{class}), it follows that $L'$ has type $E_7$, $A_1 E_6$, $D_7$, $D_6$, $A_1 D_5$, $A_7$, $A_1 A_5$ or $A_3^2$.

Suppose that $L'$ has type $E_7$. Then Theorem \ref{E7THM} shows that $Z$ is conjugate to $E_7(\#m)$ for $m = \eseven{194}$, $\eseven{265}$, $\eseven{266}$, $\eseven{289}$, or $\eseven{298}$. Using the restriction $L(E_8) \downarrow E_7 = V_{E_7}(\lambda_1) /$ $\!\! V_{E_7}(\lambda_7)^2 /$ $\!\! V_{E_7}(0)^4$ from Table \ref{levie8} and the composition factors given in Table \ref{E7tabcomps}, we calculate that $Z$ does not have the same composition factors as $X$ on $L(E_8)$, a contradiction. 

Now suppose that $L'$ has type $A_1 E_6$. Then by Theorem \ref{E6THM} the projection of $Z$ to $E_6$ is conjugate to $E_6(\#\esix{29})$. As before, we use the composition factors of the action of $A_1 E_6$ on $L(E_8)$ from Table \ref{levie8} and the composition factors of $E_6(\#\esix{29})$ acting on $V_{27}$ and $L(E_6)$ from Table \ref{E6tabcomps} to see that $Z$ does not have the same composition factors as $X$ on $L(E_8)$. 

In the remaining cases all of the simple factors of $L'$ are classical. We therefore use Lemma \ref{class} to find the irreducible subgroups of type $A_1 A_3$. Firstly, suppose that $L'$ has type $D_7$. Then by Lemma \ref{class}, $Z$ acts on $V_{D_7}(\lambda_1)$ as $((0,000) |((2,000) + (2^{[r]},000) + (2^{[s]},000))| (0,000)) + (0,010)$ with $0 < r < s$ since $Z$ is $D_7$-irreducible. Using the composition factors of $L(E_8) \downarrow D_7$ given in Table \ref{levie8}, we calculate that $Z$ does not have the same composition factors as $X$ on $L(E_8)$, a contradiction. Next suppose that $L'$ has type $D_6$ or $A_1 D_5$. Then the action of $Z$ on $\lambda_1$ (respectively $(0,\lambda_1)$) has two trivial composition factors. It then follows that $Z$ has more than seven trivial composition factors on $L(E_8)$, a contradiction.  

 Now suppose that $L'$ has type $A_7$. Then using Lemma \ref{class} we find that there is a unique conjugacy class of $A_7$-irreducible subgroups of type $A_1 A_3$; they act as $(1,100)$ on $V_{A_7}(\lambda_1)$. Using the restriction $L(E_8) \downarrow A_7$ from Table \ref{levie8} we find that $Z$ is not contained in $A_7$. Finally, suppose that $L'$ has type $A_1 A_5$ or $A_3^2$. Using the composition factors of $L(E_8) \downarrow L'$ given in Table \ref{levie8}, it follows that $Z$ has at most one composition factor of dimension at least 32. This is a contradiction since $X$ has at least two such composition factors.  

\subsubsection{$M_2 = \bar{A}_1^2 A_1 B_4$ $(E_8(\#\eeight{120}))$ or $\bar{A}_1^2 B_2 B_3$ $(E_8(\#\eeight{121}))$}

Since the subgroups $A_1 B_4$ and $B_2 B_3$ of $D_6$ are centralised by a graph automorphism of $D_6$ we have $\text{Out}_{E_8}(M_2) \cong S_2$ with the involution swapping the two $\bar{A}_1$ factors. We note that as in previous discussions the subgroup $X = A_1 \hookrightarrow \bar{A}_1^2 < D_8$ via $(1,1)$ is conjugate to $Y = A_1 < D_8$ acting via $2 + 0^{13}$ on $V_{D_8}(\lambda_1)$. This is enough to determine the classes of diagonal subgroups of all subgroups of $M_2$.

\subsubsection{$M_2 = \bar{A}_1^2 A_1 C_3$ $(E_8(\#\eeight{122}))$}

Firstly we note that $\text{Out}_{E_8}(M_2)$ is trivial, as can be seen from the $M_2$-composition factors of $L(E_8)$. The diagonal subgroups of $M_2$ therefore follow. We also note $X = A_1 \bar{A}_1 C_3 \hookrightarrow M_2$ via $(1_a,1_b,1_a,100)$ is contained in the $E_8$-reducible subgroup $A_1 C_4$ when $p=2$ and hence $X$ is $E_8$-reducible.  

Next, we consider $X = \bar{A}_1^2 A_1^2 A_1 = E_8(\#\eeight{558})$ $(p \neq 2)$. A graph automorphism of $D_6$ swaps the third and fourth $A_1$ factors and hence $\text{Out}_{E_8}(\bar{A}_1^2 A_1^2 A_1) \cong S_2$ where the involution simultaneously swaps the two $\bar{A}_1$ factors as well as swapping the third and fourth $A_1$ factors. The subgroup $A_1^2 A_1 < D_6$ acts as $(1,1,2)$ on $V_{D_6}(\lambda_1)$. Therefore the subgroups $\bar{A}_1^2 A_1 A_1 \hookrightarrow X$ via $(1_a,1_b,1_c,1_c,1_d)$, $(1_a,1_b,1_c,1_d,1_c)$ or $(1_a,1_b,1_c,1_d,1_d)$ are contained in previously considered subgroups of $\bar{A}_1^2 D_6$.  

It remains to prove that when $p=2$ the diagonal subgroups $X = E_8(\#\eeight{614}) = A_1 A_3 \hookrightarrow \bar{A}_1^2 A_1 A_3 = E_8(\#\eeight{560}) = Y$ via $(1^{[r]},1^{[s]},1^{[t]},100)$ $(rst=0; r \neq t)$ are $E_8$-irreducible. From Table \ref{E8tabcomps}, we have \begin{align*} L(E_8) \downarrow Y = \ & (2,0,0,000) /  (1,1,1,010) /  (1,0,3,000) /  (1,0,1,101) /  (0,2,0,000) / \\ \ & (0,1,2,010) /  (0,1,0,200) /   (0,1,0,010)^2 / (0,1,0,002) /  (0,0,2,101) / \\ \ & (0,0,2,000) /  (0,0,0,101)^2 /  (0,0,0,020) /  (0,0,0,000)^4. \end{align*} It follows that $X$ has exactly $4$ trivial composition factors as well as composition factors of dimensions $56$ and $48$. A routine use of Lemma \ref{wrongcomps} as before shows that $X$ is $E_8$-irreducible.

\subsection{$M_1 = \bar{D}_4^2$ $(E_8(\#\eeight{108}))$}  

It easily follows from \cite[Table~10.3]{LS1} that $\text{Out}_{E_8}(M_1) \cong S_2 \times S_3$, where the central involution swaps the two $\bar{D}_4$ factors, and the $S_3$ direct factor acts simultaneously as the group of graph automorphisms of each $\bar{D}_4$ factor. In particular, there is just one $E_8$-conjugacy class of each of the subgroups $B_3 \bar{D}_4$ and $A_1 B_2 \bar{D}_4$ contained in $M_1$. It then follows that there are two $E_8$-conjugacy classes of each of the subgroups $B_3^2$, $B_3 A_1 B_2$ and $A_1^2 B_2^2$. We make some remarks about them now. 

Firstly, one class of subgroups $B_3^2$ acts with composition factors $(W(100),000) /$ $\!\! (000,W(100)) /$ $\!\! (000,000)^2$ on $V_{D_8}(\lambda_1)$ and is hence $D_8$-reducible by Lemma \ref{class}. The other class is $D_8$-irreducible, namely $E_8(\#\eeight{623})$. One class of subgroups $A_1 B_2 B_3$ acts with composition factors $(W(2),00,000) /$ $\!\! (0,W(10),000) /$ $\!\! (0,00,W(100)) /$ $\!\! (0,00,000)$ on $V_{D_8}(\lambda_1)$ and is therefore $D_8$-irreducible if and only if $p \neq 2$. Moreover, when $p \neq 2$ this subgroup is contained in $\bar{A}_1^2 D_6$ and conjugate to $E_8(\#\eeight{533}^{\{0\}})$. Similarly, one class of subgroups $A_1^2 B_2^2$ acts  with composition factors $(W(2),0,00,00) /$ $\!\! (0,W(2),00,00) /$ $\!\! (0,0,W(10),00) /$ $\!\! (0,0,00,W(10))$ on $V_{D_8}(\lambda_1)$ and is therefore $D_8$-irreducible if and only if $p \neq 2$. This class is denoted $E_8(\#\eeight{677})$ and in this case $\text{Out}_{E_8}(A_1^2 B_2^2) = \text{Out}_{D_8}(A_1^2 B_2^2) \cong S_2 \times S_2$ where one involution swaps the two $A_1$ factors and another involution swaps the two $B_2$ factors. The second class of subgroups $A_1^2 B_2^2$ is denoted $E_8(\#\eeight{676})$ and in that case $\text{Out}_{E_8}(A_1^2 B_2^2) \cong S_2$ where the involution simultaneously swaps the two $A_1$ factors and the two $B_2$ factors. 

Most of the classes of diagonal subgroups contained in $M_1$ follow immediately. We will point out those that are not entirely obvious. First, consider $X = B_3 A_1^3 = E_8(\#\eeight{636})$. In this case $\text{Out}_{\bar{D}_4}(A_1^3) \cong S_3$ and so $\text{Out}_{E_8}(X) \cong S_3$. Next, we consider the subgroup $Y = A_1^4 B_2 = E_8(\#\eeight{638}) < B_3 A_1^3$. We note that $\text{Out}_{E_8}(Y) \cong S_3$, with the $S_3$ acting naturally on the last three $A_1$ factors. This yields the diagonal subgroups of $Y$. Finally, we consider the subgroup $Z = A_1^6 < B_3 A_1^3$. Here we have $\text{Out}_{E_8}(Z) \cong S_3 \times S_3 \times S_2$, where the first $S_3$ direct factor acts naturally on the first three $A_1$ factors, the second $S_3$ direct factor acts naturally on the final three $A_1$ factors and the central involution swaps the first three $A_1$ factors with the last three $A_1$ factors. 

We also note that $\text{Out}_{E_8}(A_2^2) \cong S_2 \times S_2$, where one involution simultaneously acts as a graph automorphism of both $A_2$ factors and another involution swaps the two $A_2$ factors. 

It remains to prove that the diagonal subgroups $X = E_8(\#\eeight{673}) = A_1 A_2 \hookrightarrow A_1^3 A_2 = E_8(\#\eeight{670}) = Y$ via $(1,1^{[r]},1^{[s]},10)$ $(0 < r < s)$ are $E_8$-irreducible when $p=2$. To do this we use a standard application of Lemma \ref{wrongcomps}. The $X$-composition factors of $L(E_8)$ can be found from those of $L(E_8) \downarrow Y$, which are given in Table \ref{E8tabcomps}. We note that there are six trivial $X$-composition factors of $L(E_8)$ as well as two $64$-dimensional $X$-composition factors. We omit the details of checking that no Levi subgroup $L$ contains an $L'$-irreducible subgroup of type $A_1 A_2$ having the same composition factors as $X$ on $L(E_8)$. 

\subsection{$M_1 = \bar{A}_3 D_5$ $(E_8(\#\eeight{109}))$ or $B_7$ $(E_8(\#\eeight{47}))$ or $A_1 B_6$ $(E_8(\#\eeight{110}))$ or $B_2 B_5$ $(E_8(\#\eeight{111}))$ or $B_3 B_4$ $(E_8(\#\eeight{112}))$}  

There is nothing for us to explicitly note here. In particular, all outer automorphisms acting on subgroups contained in $M_1$ are induced by elements of $D_8$.

\subsection{$M_1 = B_2^2$ (maximal if $p \neq 2$, contained in $B_4$ if $p=2$) $(E_8(\#\eeight{113}))$}  \label{b2b2ine8}

The $D_8$-irreducible subgroups contained in $M_1$ are straightforward to determine using Lemma \ref{maxclassical}. When $p=2$ the subgroup $A_1^4$ acting as $(1,1,1,1)$ on $V_{D_8}(\lambda_1)$ is also contained in $A_1 C_4 = E_8(\#\eeight{115})$. Indeed, there are just two $D_8$-classes of such $A_1^4$ subgroups and these can be distinguished by their composition factors on $L(E_8)$, as given in Table \ref{E8tabcomps}. We consider $A_1^4$ as a subgroup of $E_8(\#\eeight{115})$ in Section \ref{a1c4ine8}.   

\subsection{$M_1 = B_2^2$ $(p \neq 2)$ $(E_8(\#\eeight{114}))$} 

We note that this subgroup $B_2^2$ is contained in the $E_8$-reducible maximal subgroup $B_4$ of $D_8$ when $p=2$ and is therefore $E_8$-reducible itself. The subgroups of $M_1$ follow in a similarly straightforward manner to those of $E_8(\#\eeight{113})$ in Section \ref{b2b2ine8}.   

\subsection{$M_1 = A_1 C_4$ $(E_8(\#\eeight{115}))$}  \label{a1c4ine8}

By Lemma \ref{maxclassical}, the reductive, maximal connected subgroups of $M_1$ are $A_1 \bar{A}_1 C_3$, $A_1 B_2^2$, $A_1^4$ $(p \neq 2)$, $A_1 A_1$ $(p \geq 11)$ and $A_1 D_4$ $(p=2)$. The subgroup $A_1 \bar{A}_1 C_3$ is contained in $\bar{A}_1^2 A_1 C_3 = E_8(\#\eeight{122})$ and by considering composition factors of $L(E_8)$ we find that it is conjugate to $E_8(\#\eeight{564}^{\{0,0\}})$. The maximal connected subgroup $A_1 B_2^2$ acts as $(1,01,00) + (1,00,01)$ on $V_{D_8}(\lambda_1)$ and is thus a subgroup of either $E_8(\#\eeight{676})$ or $E_8(\#\eeight{677})$, both of which are of type $A_1^2 B_2^2$. It again follows from considering the composition factors of $L(E_8)$ that $A_1 B_2^2$ is conjugate to $E_8(\#\eeight{708}^{\{0\}})$. 

The maximal subgroup $A_1^4$ $(p \neq 2)$ is still $E_8$-irreducible when $p=2$, but it is now contained in $A_1^2 B_2 < A_1 B_3 < A_1 D_4$. The diagonal subgroups of $A_1^4$ follow easily since $\text{Out}_{D_8}(A_1^4) \cong S_4$. The diagonal subgroups of the maximal subgroup $A_1 A_1$ $(p \geq 11)$ are clear. 

Next, we consider the maximal subgroup $M_2 = A_1 D_4$. We note that since $\text{Out}_{C_4}(D_4) \cong S_2$ we have $\text{Out}_{D_8}(M_2) \cong S_2$, where the involution acts as an involutory graph automorphism of the $D_4$ factor. There are thus two classes of maximal subgroups of type $A_1 B_3$. One such class acts as $(1,000) | (1,100) | (1,000)$ on $V_{D_8}(\lambda_1)$ and is thus $D_8$-reducible by Lemma \ref{class} and hence $E_8$-reducible. The subgroup $A_1 A_1^4$ of $M_2$ is contained in $A_1 B_2^2 < M_1$ and hence conjugate to $E_8(\#\eeight{649}^{\{0,0\}})$. 

Finally, we prove that $X = A_1 A_2 = E_8(\#\eeight{872})$ is $E_8$-irreducible when $p=2$. From Table \ref{E8tabcomps}, we have \begin{align*} L(E_8) \downarrow X = \ & (3,11) / (2,30) / (2,11) / (2,03) / (2,00)^2 / (1,30)^2 / (1,22) / \\ \ & (1,03)^2 / (1,00)^4 / (0,30)^2 / (0,22) / (0,11)^2 / (0,03)^2 / (0,00)^4.\end{align*} Suppose that $X$ is $E_8$-reducible. By Lemma \ref{wrongcomps}, there exists a subgroup $Z$ of type $A_1 A_2$ contained $L$-irreducibly in a Levi subgroup $L$, such that $X$ and $Z$ have the same composition factors on $L(E_8)$. By considering the number of trivial composition factors on $L(E_8)$ of each Levi subgroup that contains an irreducible subgroup $A_1 A_2$ (using Table \ref{levie8} and Lemma \ref{class}) it follows that $L'$ has type $E_7$, $A_1 E_6$, $D_7$, $A_2 D_5$, $A_2 D_4$ or $A_1^2 A_2^2$. Suppose that $L'$ has type $E_7$. Then Theorem \ref{E7THM} shows that $Z$ is conjugate to $E_7(\#m)$ for $m = \eseven{184}, \eseven{307}$, or $\eseven{308}$. Using the restriction $L(E_8) \downarrow E_7 = V_{E_7}(\lambda_1) /$ $\!\! V_{E_7}(\lambda_7)^2 /$ $\!\! V_{E_7}(0)^4$ and the composition factors given in Table \ref{E7tabcomps}, we calculate that $Z$ does not have the same composition factors as $X$ on $L(E_8)$, a contradiction. Now suppose that $L'$ has type $A_1 E_6$. Then by Theorem \ref{E6THM} the projection of $Z$ to $E_6$ is conjugate to $E_6(\#m)$ for $m = \esix{13}, \esix{14}, \dots, \esix{23}$ or $\esix{37}$. As before, we calculate that $Z$ does not have the same composition factors as $X$ on $L(E_8)$. Suppose that $L'$ has type $D_7$. Then by Lemma \ref{class}, $Z$ acts on $V_{D_7}(\lambda_1)$ as $((0,00) |((2,00) + (2^{[r]},00))| (0,00)) + (0,11)$ with $r \neq 0$, since $Z$ is $D_7$-irreducible. Now, because $V_{D_7}(\lambda_1)$ occurs as composition factor of $L(E_8)$ with multiplicity two, it follows that $Z$ has at least six trivial composition factors on $L(E_8)$, which is a contradiction. Now suppose that $L'$ has type $A_2 D_5$ or $A_2 D_4$. Then the $A_2$ factor of $Z$ is conjugate to $\bar{A}_2$ and so in all cases the composition factors of $Z$ are not the same as those of $X$. Finally, suppose that $L'$ has type $A_1^2 A_2^2$.  Then $L'$ has only four composition factors of dimension at least $16$ on $L(E_8)$ and all four of them have dimension $18$. Since $X$ has four composition factors of dimension 18 and two composition factors of dimension 16, it follows that $X$ and $Z$ do not have the same composition factors on $L(E_8)$. This final contradiction shows that $X$ is $E_8$-irreducible. 

\subsection{$M_1 = A_1 C_4$ $(p \neq 2)$ $(E_8(\#\eeight{116}))$}  

Notice that $X = A_1 B_2^2 < M_1$ acts on $V_{D_8}(\lambda_1)$ as $(1,01,00) + (1,00,01)$ and is thus a subgroup of $E_8(\#\eeight{677}) = A_1^2 B_2^2 < \bar{D}_4^2$. By considering the $X$-composition factors of $L(E_8)$ we find that $X$ is conjugate to a $D_8$-reducible subgroup acting as $(2,00,00)^2 + (0,10,00) + (0,00,10)$. The rest of the lattice structure follows as for $E_8(\#\eeight{115})$ in Section \ref{a1c4ine8}.  

\subsection{$M_1 = B_4$ $(E_8(\#\eeight{45}))$}  

The reductive, maximal connected subgroups of $M_1$ are given by Lemma \ref{maxclassical}. It is then a routine calculation to find the action of each subgroup on $V_{D_8}(\lambda_1) \downarrow B_4 = V_{B_4}(\lambda_4)$, which is also required when finding the composition factors on $V_{26}$ of the subgroups of the maximal subgroup $B_4$ of $F_4$. 

Firstly, the maximal connected subgroup $D_4$ of $M_1$ acts as $\lambda_3 + \lambda_4$ on $V_{D_8}(\lambda_1)$. There are two $D_8$-conjugacy classes of $D_4$ subgroups acting in such a way. By restricting the action of $M_1$ on $L(E_8)$ to $D_4$ we find the composition factors of $D_4$ on $L(E_8)$. This shows that $D_4$ is conjugate to $E_8(\#\eeight{49})$. 

The maximal connected subgroup $A_1 A_3$ acts as $(1,100) + (1,001)$ on $V_{D_8}(\lambda_1)$. Therefore it is $D_8$-reducible by Lemma \ref{class}. The maximal connected subgroup $B_2^2$ $(p=2)$ acts as $(01,01)$ and is maximal when $p \neq 2$. Therefore it is given the identification number $E_8(\#\eeight{113}\mathrm{b})$ and has been considered previously. The action of the maximal subgroup $A_1^2 B_2$ $(p \neq 2)$ on $V_{D_8}(\lambda_1)$ is $(1,1,01)$ and therefore $A_1^2 B_2$ is conjugate to $E_8(\#\eeight{709}^{\{0\}})$. 

Next, we consider the maximal connected subgroup $A_1^2$. When $p \geq 5$ it acts as $(3,1) + (1,3)$ on $V_{D_8}(\lambda_1)$ and is conjugate to $E_8(\#\eeight{707}^{\{0,0,0,0\}})$. When $p=3$ it acts with composition factors $(3,1) / (1,1)^2 / (1,3)$ and is therefore $D_8$-reducible by Lemma \ref{class}. Finally, the maximal subgroup $A_1$ $(p \geq 11)$ acts as $10 + 4$ and is therefore conjugate to $E_8(\#\eeight{12}^{\{0,0\}})$.  

\subsection{$M_1 = B_4$ $(p \neq 2)$ $(E_8(\#\eeight{46}))$}  

As for $E_8(\#\eeight{45})$, we find the action of the maximal subgroups of $M_1$ on $V_{D_8}(\lambda_1)$. In this case the subgroups $D_4$, $A_1 A_3$, $A_1^2 B_2$ and $A_1^2$ $(p=3)$ are all conjugate to $D_8$-reducible subgroups. When $p \geq 5$, the maximal subgroup $A_1^2$ is conjugate to $E_8(\#\eeight{693}^{\{0,0,0,0\}})$. As in the previous case, the maximal subgroup $A_1$ $(p \geq 11)$ is conjugate to $E_8(\#\eeight{12}^{\{0,0\}})$. 

This completes the analysis of the case $M = D_8$.  

\section{$M = \bar{A}_1 E_7$ $(E_8(\#\eeight{102}))$}

In this section we find the lattice of $E_8$-irreducible subgroups contained in $M$. Theorem \ref{E7THM} allows us to immediately write down the lattice of $M$-irreducible subgroups. The subgroup $\bar{A}_1^2 D_6 < \bar{A}_1 E_7$ is conjugate to $E_8(\#\eeight{107}) < D_8$ and so any subgroup of $\bar{A}_1^2 D_6$ has already been considered. 

The remainder of this section deals with the question of whether an $M$-irreducible subgroup is $E_8$-irreducible. 

 There is one $M$-irreducible subgroup that we prove is $E_8$-reducible and hence all of its subgroups are thus $E_8$-reducible. Consider the subgroup $X = A_1 F_4 \hookrightarrow \bar{A}_1 A_1 F_4 = E_8(\#\eeight{881})$ via $(1,1,\lambda_1)$ when $p=2$. By \cite[p.333, Table~3]{LS6} the connected centraliser of the $F_4$ factor is $G_2$ and hence $\bar{A}_1 A_1 F_4 < G_2 F_4$. Theorem \ref{G2THM} implies that the projection of $X$ to $G_2$ is $G_2$-reducible since $p=2$. Therefore $X$ is $G_2 F_4$-reducible and thus $E_8$-reducible.  

Now we need to prove that each subgroup of $M$ in Table \ref{E8tab} is $E_8$-irreducible. Most subgroups either contain a simple $E_8$-irreducible subgroup or have already been considered as subgroups of $D_8$. The remaining cases are $E_8(\#n)$ for $n = \eeight{879}, \eeight{889}, \eeight{897}, \eeight{898}$ when $p=2$ and $n = \eeight{911}$, $\eeight{913}, \eeight{914}, \eeight{915}, \eeight{916}$ when $p \neq 2$. 

We start by considering the cases where $p \neq 2$. If $p \geq 5$ then the action of each subgroup on $L(E_8)$ has no trivial composition factors and hence all of the subgroups are $E_8$-irreducible by Corollary \ref{notrivs}. Since $E_8(\#\eeight{913})$ is only $A_1 E_7$-irreducible for $p \geq 5$, it requires no further consideration. 

Now let $p=3$ and $X = E_8(\#\eeight{911}) = A_1 A_3 \hookrightarrow \bar{A}_1 A_1 A_3 = E_8(\#\eeight{910}) = Y$ via $(1^{[r]},1^{[s]},100)$. Suppose that $X$ is $E_8$-reducible. By Lemma \ref{wrongcomps}, there exists a subgroup $Z$ of type $A_1 A_3$ contained $L$-irreducibly in a Levi subgroup $L$, such that $X$ and $Z$ have the same composition factors on $L(E_8)$. From Table \ref{E8tabcomps} we see that \begin{align*} L(E_8) \downarrow Y = \ & (2,0,000) / (1,2,010)^2 / (1,0,200) / (1,0,002) / (0,4,000) /  \\ \ & (0,2,101)^2 / (0,2,000) / (0,0,101) / (0,0,020) / (0,0,000)^2.\end{align*} Therefore, $X$ and $Z$ have only two trivial composition factors on $L(E_8)$ and so $L'$ has at most two trivial composition factors on $L(E_8)$. Using the restrictions in Table \ref{levie8}, we find that $L'$ has one of the following types: $A_1 E_6$, $D_7$, $A_2 D_5$, $A_7$, $A_3 A_4$, $A_1 A_6$, $A_1 A_2 A_4$, $A_1^2 A_4$, or $A_3^2$. Lemma \ref{class} immediately rules out $L'$ having type  $A_3 A_4$, $A_1 A_6$, $A_1 A_2 A_4$ or $A_1^2 A_4$ since such a Levi subgroup contains no irreducible subgroup $A_1 A_3$. We will now take the other five types and rule them out in turn. Suppose that $L'$ has type $A_1 E_6$. Then one $L'$-composition factor of $L(E_8)$ is isomorphic to $(1,0)$. Hence $Z$ has a $2$-dimensional composition factor on $L(E_8)$, a contradiction. Now suppose that $L'$ has type $A_7$. Then from Lemma \ref{class} we deduce that $Z$ acts on $V_{A_7}(\lambda_1)$ as $(1,100)$. It is straightforward to calculate the $Z$-composition factors of $L(E_8)$ from the $L'$-composition factors given in Table \ref{levie8}. These are not the same as those of $X$, a contradiction. Now suppose that $L'$ has type $A_2 D_5$ or $D_7$. Then the $A_3$ factor of $Z$ is conjugate to a Levi subgroup $\bar{A}_3$, acting as $010 + 000^4$ on $V_{D_5}(\lambda_1)$ and $010 + 000^8$ on $V_{D_7}(\lambda_1)$. In particular, there is no $Z$-composition factor of $L(E_8)$ isomorphic to $(0,020)$. Therefore the composition factors of $X$ and $Z$ are not the same on $L(E_8)$, a contradiction. Finally, suppose that $L'$ has type $A_3^2$. Then by Table \ref{levie8}, the largest dimension of any composition factor of $L(E_8) \downarrow L'$ is $24$. This is a contradiction since $X$ has a $45$-dimensional composition factor on $L(E_8)$, namely $(2^{[s]},101)$.  

A similar argument applies to show that $X = E_8(\#n)$ is $E_8$-irreducible for $n = \eeight{914}$, $\eeight{915}$ and $\eeight{916}$. In particular, there are only two trivial $X$-composition factors of $L(E_8)$.   
 
We now consider the cases where $p=2$. Firstly, since $\bar{A}_1 A_7 = E_8(\#\eeight{879})$ has rank $8$ it is clearly $G$-irreducible. Now let $X = \bar{A}_1 \bar{A}_2 A_3 = E_8(\#\eeight{889})$. Suppose that $X$ is $E_8$-reducible. By Lemma \ref{wrongcomps}, there exists a subgroup $Z$ of type $A_1 A_2 A_3$ contained $L$-irreducibly in a Levi subgroup $L$, such that $X$ and $Z$ have the same composition factors on $L(E_8)$. From Table \ref{E8tabcomps} we find that $X$ has just two trivial composition factors on $L(E_8)$. It thus follows that $L'$ is of type $A_2 D_5$. Using Lemma \ref{class}, we deduce that the $A_1 A_3$ factor of $Z$ is contained in the $D_5$ factor of $A_2 D_5$ and acts as $(1^{[r]} \otimes 1^{[s]},000) + (0,010)$ on $V_{D_5}(\lambda_1)$. From this we calculate the $Z$-composition factors of $L(E_8)$ and see that they are not the same as those of $X$, a contradiction. Therefore $X$ is $E_8$-irreducible. 

Finally, we need to consider the cases where $X$ is $E_8(\#\eeight{897})$ or $E_8(\#\eeight{898})$ and so $X = A_1 A_2 \hookrightarrow \bar{A}_1 \bar{A}_2 A_1 A_2 = E_8(\#\eeight{887}) = Y$ via $(1^{[r]},10^{[t]},1^{[s]},10^{[u]})$ $(rs=tu=0)$ or $(1^{[r]},10^{[t]},1^{[s]},01^{[u]})$ ($rs=tu=0$; $t \neq u$), respectively. Using the composition factors of $L(E_8) \downarrow Y$ given in Table \ref{E8tabcomps}, we find that $X$ has at most six trivial composition factors on $L(E_8)$. The case where $X$ has exactly six occurs when $X$ is embedded via $(1,10^{[1]},1,10)$. We prove that $X$ is $E_8$-irreducible in this case; the others are similar and easier. Suppose that $X$ is $E_8$-reducible. By Lemma \ref{wrongcomps}, there exists a subgroup $Z$ of type $A_1 A_2$ contained $L$-irreducibly in a Levi subgroup $L$, such that $X$ and $Z$ have the same composition factors on $L(E_8)$. Restricting the $Y$-composition factors of $L(E_8)$ to $X$ yields \begin{align*} L(E_8) \downarrow X = \ & (4,00) / (2,30)^2 / (2,11)^2 / (2,03)^2 / (2,00)^4 / (0,30)^4 / (0,22)^3 / \\ \ & (0,11)^4 / (0,03)^4 / (0,00)^6. \end{align*} By considering the number of trivial composition factors on $L(E_8)$ of each Levi subgroup containing an irreducible subgroup $A_1 A_2$ (using Table \ref{levie8} and Lemma \ref{class}), it follows that $L'$ has type $E_7$, $A_1 E_6$, $D_7$, $A_2 D_5$, $A_2 D_4$, $A_1 A_5$, $A_1 A_2 A_3$ or $A_1^2 A_2^2$. Suppose that $L'$ has type $E_7$ or $A_1 E_6$. Then Theorem \ref{E7THM} implies that $Z$ is conjugate to $E_7(\#m)$ for $m = \eseven{184}$, $\eseven{307}$, or $\eseven{308}$ and Theorem \ref{E6THM} implies that the projection of $Z$ to $E_6$ is conjugate to $E_6(\#m)$ for $m = \esix{13}, \esix{14}, \dots, \esix{23}$ or $\esix{37}$. Using the restrictions in Tables \ref{E6tabcomps}, \ref{E7tabcomps} and \ref{levie8} we calculate the $Z$-composition factors of $L(E_8)$; in all cases these are not the same as the $X$-composition factors, a contradiction. Now suppose that $L'$ has type $D_7$. Lemma \ref{class} implies that $Z$ acts on $V_{D_7}(\lambda_1)$ as $((0,00)|((2,00) + (2^{[r]},00))| (0,00)) + (0,11)$ with $r \neq 0$, since $Z$ is $D_7$-irreducible. From Table \ref{levie8}, the $D_7$-composition factor $V_{D_7}(\lambda_1)$ of $L(E_8)$ has multiplicity two. Hence $(2,00)$ and $(2^{[r]},00)$ both occur as $Z$-composition factors of $L(E_8)$ with multiplicity at least two. Comparing with the $X$-composition factors we see that $r=0$, a contradiction. Suppose that $L'$ has type $A_2 D_5$, $A_2 D_4$ or $A_1 A_2 A_3$. Then the $A_2$ factor of $Z$ is conjugate to $\bar{A}_2$ and so in all cases the composition factors of $Z$ are not the same as those of $X$. For the next case let $L'$ have type $A_1 A_5$. Then $L'$ has only five composition factors of dimension at least $16$ on $L(E_8)$: four of them have dimension $20$ and the other has dimension $23$. Since $X$ has four composition factors of dimension $18$ and two composition factors of dimension $16$, it follows that $X$ and $Z$ do not have the same composition factors on $L(E_8)$. Finally, we suppose that $L'$ has type $A_1^2 A_2^2$. Then $L'$ has only four composition factors of dimension at least $16$ on $L(E_8)$ and all four of them have dimension $18$. As in the previous case, it follows that $X$ and $Z$ do not have the same composition factors on $L(E_8)$. This final contradiction proves that $X$ is $E_8$-irreducible. 

This completes the analysis of the $E_8$-irreducible subgroups contained in $M = \bar{A}_1 E_7$.

\section{$M = \bar{A}_2 E_6$ $(E_8(\#\eeight{103}))$}

Theorem \ref{E6THM} and Lemma \ref{maxclassical} allow us to write down the lattice of $M$-irreducible subgroups. Firstly, we note that the maximal subgroup $\bar{A}_2 \bar{A}_1 A_5$ is contained in $\bar{A}_1 E_7$ and conjugate to $E_8(\#\eeight{878})$. We will now consider subgroups of $\bar{A}_2 Y$ for the remaining reductive, maximal connected subgroups $Y$ of $E_6$. In doing this we will consider all of the subgroups of $A_1 E_6 = E_8(\#\eeight{981})$ $(p \neq 2)$ in the following sections and thus give it no further consideration.  

\subsection{$M_1 = \bar{A}_2^4$ $(E_8(\#\eeight{975}))$}

As in \cite[Section 7.3]{tho1} we fix a conjugacy class representative of $M_1$ by giving the $M_1$-composition factors of $L(E_8)$ in Table \ref{E8tabcomps}. The diagonal subgroups follow from the action of $\text{Out}_{E_8}(M_1) \cong \textrm{GL}(2,3)$ (by \cite[Table~10.3]{LS1}) on the four $A_2$ factors.

We prove that $A_1^4 < M_1$ $(p \neq 2)$ is contained in $D_8$ and conjugate to $E_8(\#\eeight{873})$. Theorem \ref{E6THM} implies that $A_1^3 < \bar{A}_2^3 < E_6$ is contained in the maximal subgroup $C_4$ of $E_6$. When we consider $\bar{A}_2 C_4$ in Section \ref{a2c4ine8}, we prove that $A_1 C_4$ is conjugate to $E_8(\#\eeight{116})$. Hence $A_1^4$ is contained in $E_8(\#\eeight{116})$ and conjugate to $E_8(\#\eeight{873})$.

\subsection{$M_1 = \bar{A}_2 F_4$ $(E_8(\#\eeight{976}))$}

We first note $M_1$ is a subgroup of $G_2 F_4$. Thus, by Theorem \ref{G2THM}, the subgroup $A_1 F_4$ $(p \neq 2)$ is conjugate to $A_1 F_4 \hookrightarrow \bar{A}_1 A_1 F_4 = E_8(\#\eeight{881}) < \bar{A}_1 E_7$ via $(1,1,\lambda_1)$. Therefore we need only consider subgroups of the form $\bar{A}_2 X$ where $X$ is an $E_6$-irreducible connected subgroup of $F_4$. 

By Theorem \ref{E6THM}, we have that $\bar{A}_2 \bar{A}_1 C_3$ is a subgroup of $\bar{A}_2 \bar{A}_1 A_5$ and $\bar{A}_2^2 A_2$ is a subgroup of $\bar{A}_2^4$. We will consider $\bar{A}_2 C_4$ $(p=2)$ in Section \ref{a2c4ine8}, and similarly we consider $\bar{A}_2 G_2$ $(p=7)$ in Section \ref{a2g2ine8}.

\subsection{$M_1 = \bar{A}_2 C_4$ $(E_8(\#\eeight{977}))$} \label{a2c4ine8}

 The subgroup $M_1$ is maximal if $p \neq 2$ and contained in $\bar{A}_2 F_4$ when $p=2$. Firstly, we prove that $X = A_1 C_4$ $(p \neq 2)$ is contained in $D_8$ and conjugate to $Y = A_1 C_4 = E_8(\#\eeight{116})$. Indeed, since $p \neq 2$, \cite[Table 8.1]{LS3} shows that there are exactly two classes of subgroups of type $C_4$ in $E_8$; they are contained in the two classes of subgroups of type $A_7$. In particular, the $C_4$ factor of $Y$ is contained in $E_6$ and the connected centraliser is $\bar{A}_2$. It follows that $Y$ is contained in $\bar{A}_2 C_4$ and thus conjugate to $X$.  

We now need to consider subgroups of the form $X = \bar{A}_2 Y$ where $Y$ is an $E_6$-irreducible connected subgroup of $C_4$. If $Y$ is contained in $\bar{A}_1 C_3$ or $A_1^3$ then $X$ is contained in $\bar{A}_2 \bar{A}_1 A_5$ or $\bar{A}_2^4$, respectively. If $Y = A_1$ $(p \geq 11)$ then $X$ requires no further consideration. The final case is $Y = D_4$ $(p=2)$. In this case the only $E_6$-irreducible subgroup of $Y$ is $A_2$. This subgroup $A_2$ is contained in $\bar{A}_2^3$, again by Theorem \ref{E6THM}.

\subsection{$M_1 = \bar{A}_2 A_2 G_2$ $(E_8(\#\eeight{978}))$} \label{sec:E8A2A2G2}

Theorem \ref{E6THM} shows that $\bar{A}_2 A_1 G_2$ is contained in $\bar{A}_2 F_4$. Similarly, $\bar{A}_2 A_2 \bar{A}_1 A_1$ and $\bar{A}_2 A_2 \bar{A}_2$ are contained in $\bar{A}_2 \bar{A}_1 A_5$ and $\bar{A}_2^4$, respectively.  

The diagonal subgroups of $M_1$ follow by noting that $\bar{A}_2 A_2$ is a maximal subgroup of $F_4$ and thus there is an involution in $E_8$ acting as a graph automorphism on both $A_2$ factors. Moreover, by Theorem \ref{F4THM} the subgroup $A_2 G_2 \hookrightarrow M_1$ via $(10,10,10)$ is $G_2 F_4$-reducible when $p=3$. Note that looking at Theorem \ref{F4THM} would seem to suggest the subgroup $A_2 G_2 \hookrightarrow M_1$ via $(10,01,10)$ is the $G_2 F_4$-reducible subgroup. However, the given $M_1$-composition factors of $L(E_8)$ in Table \ref{E8tabcomps} show that a graph automorphism of the second $A_2$ factor has been introduced.   

The non-diagonal subgroups of $\bar{A}_2 A_2 A_2$ $(p=3)$ have all been covered previously and the simple diagonal subgroups are given by \cite[Lemma 7.11]{tho2}. The non-simple diagonal subgroups follow in the same way. Similarly, the non-diagonal subgroups of $\bar{A}_2 A_2 A_1$ $(p \geq 7)$ have already been considered and the diagonal subgroups follow in the same way as the diagonal subgroups of $M_1$.

\subsection{$M_1 = \bar{A}_2 G_2$ $(E_8(\#\eeight{979}))$} \label{a2g2ine8}

The subgroup $M_1$ is maximal if $p \neq 7$ and contained in $\bar{A}_2 F_4$ if $p=7$. By Theorem \ref{E6THM}, except for $A_2$ $(p=3)$, all reductive, maximal connected subgroups of the $G_2$ factor are contained in a previously considered maximal connected subgroup of $E_6$. It remains to consider subgroups of $\bar{A}_2 A_2$ $(p=3)$. We postpone this until Section \ref{a2a2inE8}.

\subsection{$M_1 = \bar{A}_2 A_2$ $(p \neq 2)$ $(E_8(\#\eeight{980}))$} \label{a2a2inE8}

The subgroup $M_1$ is maximal if $p \geq 5$ and contained in $\bar{A}_2 G_2$ if $p=3$. The non-diagonal reductive, maximal connected subgroups of $M_1$ are $\bar{A}_2 A_1$ and $A_1 A_2$. Theorem \ref{E6THM} implies that $\bar{A}_2 A_1$ is contained in $\bar{A}_2 \bar{A}_1 A_5$ and hence the subgroup $A_1 A_1$ of $A_1 A_2$ is also contained in $\bar{A}_2 \bar{A}_1 A_5$. 

The subgroup $X = A_1 A_2 = E_8(\#\eeight{1027})$ when $p=5$ is the only subgroup that does not properly contain an $E_8$-irreducible subgroup. The composition factors of $L(E_8) \downarrow X$ are given in Table \ref{E8tabcomps}. Since there are no trivial composition factors, Corollary \ref{notrivs} shows that $X$ is $E_8$-irreducible.  

This completes the case $M = \bar{A}_2 E_6$. 

\section{$M = A_8$ $(E_8(\#\eeight{62}))$}

By Lemma \ref{maxclassical}, the only reductive, maximal connected subgroups of $A_8$ are $B_4$ $(p \neq 2)$ and $A_2^2$. The subgroup $B_4$ is contained in $D_8$ and so its non-simple subgroups require no consideration. The subgroup $A_2^2$ was studied in \cite[Section 7.2]{tho2}. In particular, \cite[Lemma 7.10]{tho2} shows that $A_2^2$ is $E_8$-irreducible and conjugate to $E_8(\#\eeight{668})$ when $p \neq 3$ but $E_8$-reducible when $p=3$.

\section{$M = \bar{A}_4^2$ $(E_8(\#\eeight{104}))$}

By Lemma \ref{maxclassical}, the only reductive, maximal connected subgroup of $\bar{A}_4$ is $B_2$ $(p \neq 2)$. The only $\bar{A}_4$-irreducible maximal connected subgroup of $B_2$ is $A_1$ $(p \geq 5)$ embedded via the representation of high weight $4$. The final thing to note is that the subgroup $B_2^2$ of $\bar{A}_4^2$ is contained in $D_8$. Indeed, this is shown in \cite[p.63]{LS3}.

\section{$M = G_2 F_4$ $(E_8(\#\eeight{105}))$}

The irreducible connected subgroups of $G_2$ and $F_4$ are given by Theorems \ref{G2THM} and \ref{F4THM}, respectively. This allows us to write down all of the $M$-irreducible connected subgroups when noting the following details.

Firstly, the $G_2$ factor is contained in a Levi subgroup $D_4$ by \cite[3.16]{se2} and thus $\bar{A}_1 A_1 F_4$ is contained in $\bar{A}_1 E_7$ and $\bar{A}_2 F_4$ is contained in $\bar{A}_2 E_6$. Similarly, the subgroups $G_2 B_4$, $G_2 \bar{A}_1 C_3$ and $G_2 \bar{A}_2 A_2$ have all been considered before in the $D_8$, $\bar{A}_1 E_7$ and $\bar{A}_2 E_6$ cases, respectively.  

We now consider the remaining maximal connected subgroups of $M$ in turn. Once we have considered $A_2 F_4$ $(p=3)$ and $A_1 F_4$ $(p \geq 7)$ we will only consider subgroups of the form $G_2 X$ where $X$ is an $F_4$-irreducible subgroup, since all others will have been covered. In particular, we will not make any further mention of $G_2 A_1$ $(p \geq 13)$.

\subsection{$M_1 = A_2 F_4$ $(p=3)$ $(E_8(\#\eeight{1030}))$}

By Theorem \ref{G2THM}, the irreducible subgroup $A_1$ of the $A_2$ factor is contained in $\bar{A}_1 A_1$. We have hence already considered the maximal connected subgroup $A_1 F_4$ of $M_1$ and need only consider subgroups of the form $A_2 X$ where $X$ is an $F_4$-irreducible subgroup. Similarly, $A_2 B_4$, $A_2 \bar{A}_1 C_3$ and $A_2 \bar{A}_2 A_2$ have already been considered and so it remains to examine $A_2 A_1 G_2$. 

Theorem \ref{F4THM} implies that the only reductive, maximal connected subgroup of $A_1 G_2$ not previously considered is $A_1 A_2$ (where the subgroup $A_2$ of $G_2$ is generated by short root subgroups of $G_2$). This yields $A_2^2 A_1$, since the $G_2$ factor of $M$ and the $G_2$ factor of $A_1 G_2$ are conjugate and there exists an involution in the normaliser of $G_2 A_1 G_2$ swapping the two $G_2$ factors by the construction in \cite[p.39]{se2}. It remains to consider the diagonal subgroups of $A_2^2 A_1$. Each subgroup $G_2$ contains an element inducing a graph automorphism of $A_2$ and combined with the involution which swaps the $A_2$ subgroups we obtain the classes of diagonal subgroups $E_8(\#\eeight{1038})$.

\subsection{$M_1 = A_1 F_4$ $(p \geq 7)$ $(E_8(\#\eeight{1031}))$}

We need to consider the subgroups $A_1 A_1 G_2$, $A_1 A_1$ $(p \geq 13)$ and $A_1 G_2$ $(p=7)$ since the other reductive, maximal connected subgroups of $M_1$ have already been considered. Firstly, let $M_2 = A_1 A_1 G_2 = E_8(\#\eeight{1039})$. The diagonal subgroups of $M_2$ are as listed in Table \ref{E8tab}. We use Theorem \ref{F4THM} to prove that $X = A_1 G_2 \hookrightarrow M_2$ via $(1,1,10)$ is conjugate to $E_8(\#\eeight{952}^{\{0,0\}})$. We may consider $X$ as a subgroup of $G_2 A_1 A_1 < G_2 A_1 G_2 < G_2 F_4$ by swapping the two $G_2$ factors. Theorem \ref{F4THM} shows that $A_1 \hookrightarrow A_1 A_1 = F_4(\#\ffour{66}) < A_1 G_2 < F_4$ via $(1,1)$ is conjugate to $A_1 \hookrightarrow \bar{A}_1 A_1 = F_4(\#\ffour{58}) < \bar{A}_1 C_3 < F_4$ where the subgroup $A_1$ of $C_3$ is maximal. Therefore $X$ is a subgroup of $\bar{A}_1 G_2 C_3 < \bar{A}_1 E_7$ and comparing composition factors shows that $X$ is indeed conjugate to $E_8(\#\eeight{952}^{\{0,0\}})$.

Now, using Theorem \ref{F4THM}, we need only consider the reductive, maximal connected subgroup $A_1 A_1$ of $A_1 G_2 < F_4$, where the subgroup $A_1$ of $G_2$ is maximal. This yields $A_1^2 A_1 = E_8(\#\eeight{1042})$ since the $G_2$ factors are conjugate, as shown in the previous section. The classes of diagonal subgroups of $A_1^2 A_1$ then follow. We claim that $X_1 \hookrightarrow A_1^2 A_1$ via $(1_a,1_a,1_b)$ and $X_2 \hookrightarrow A_1^2 A_1$ via $(1_a,1_b,1_b)$ are contained in $\bar{A}_1 E_7$. To prove the claim for $X_1$, we consider the involution $t \in N_{E_8}(A_1^2 A_1)$ swapping the first two $A_1$ factors. The subgroup $X_1$ is centralised by $t$ and a routine check shows that the full centraliser of $t$ in $E_8$ is $\bar{A}_1 E_7$ rather than $D_8$. The claim for $X_2$ follows from Theorem \ref{F4THM}. Indeed, as above we have $A_1 \hookrightarrow A_1 A_1 < A_1 G_2 < F_4$ via $(1,1)$ is contained in $\bar{A}_1 C_3$ and thus $X_2$ is contained in $\bar{A}_1 G_2 C_3 < \bar{A}_1 E_7$, as claimed.  

The irreducible subgroups of $A_1 A_1$ $(p \geq 13)$ are all simple and thus require no further consideration. Finally, we use Theorem \ref{F4THM} to see that all subgroups of $A_1 G_2$ $(p=7)$ have already been considered.

\subsection{$M_1 = G_2 C_4$ $(p=2)$ $(E_8(\#\eeight{1032}))$}

We need only consider the maximal connected subgroup $G_2 D_4$ since $G_2 \bar{A}_1 C_3$ and $G_2 B_2^2$ are contained in $\bar{A}_1 E_7$ and $D_8$, respectively. Theorem \ref{F4THM} implies that the subgroups of the $D_4$ factor are contained in $B_4$ or $\bar{A}_2 A_2$ and so have already been considered.

\subsection{$M_1 = G_2^2 A_1 $ $(p \neq 2)$ $(E_8(\#\eeight{1033}))$}

As noted above, the $G_2$ factors of $M_1$ are $E_8$-conjugate and furthermore, there is an involution $t$ in $\text{Out}_{E_8}(G_2^2 A_1)$ that swaps the two factors. We need only consider the diagonal subgroups of $M_1$ since any other reductive, maximal connected subgroup of $M_1$ has been considered above. The diagonal subgroups of $M_1$ are as in Table \ref{E8tab}; we note that only when $p=3$ does there exist a special isogeny of $G_2$ yielding the subgroups $E_8(\#\eeight{1048})$, with the notation as in Section \ref{nota}. Finally, we note that $X = G_2 A_1 \hookrightarrow M_1$ via $(10,10,1)$ is centralised by $t$ and thus contained in $\bar{A}_1 E_7$, as above. By considering the $X$-composition factors of $L(E_8)$ we see that $X$ is conjugate to $E_8(\#\eeight{967}^{\{0,0\}})$.

\subsection{$M_1 = G_2 G_2 $ $(p = 7)$ $(E_8(\#\eeight{1035}))$}

From Theorem \ref{F4THM}, we see that the reductive, maximal connected subgroups of the second $G_2$ factor are contained in other reductive, maximal connected subgroups of $F_4$ and have thus already been considered.

This completes the case $M = G_2 F_4$.

\section{$M = B_2$ $(p \geq 5)$ $(E_8(\#\eeight{101}))$}

By Lemma \ref{maxclassical}, the only reductive, maximal connected subgroups of $B_2$ are $A_1^2$ and $A_1$. The subgroup $A_1^2$ is the centraliser in $B_2$ of a semisimple involution $t$ by \cite[Table 4.3.1]{gls3}. By Lemma \ref{centralisers}, the centraliser of $t$ in $E_8$ is either $D_8$ or $\bar{A}_1 E_7$. We have therefore considered the subgroup $A_1^2$ before and by considering composition factors we find that it is conjugate to $E_8(\#\eeight{707}^{\{0,0,0,0\}})$.

\section{$M = A_1 A_2$ $(p \geq 5)$ $(E_8(\#\eeight{106}))$}

The only reductive, maximal connected subgroup of $A_1 A_2$ is $A_1 A_1$, where the second $A_1$ factor is irreducibly embedded in $A_2$. By \cite[p. 31]{se2}, there exists an involution $t$ in $N_{E_8}(M)$ such that $t$ centralises the $A_1$ factor of $M$ and acts as a graph automorphism of the $A_2$ factor. Thus $A_1 A_1$ is contained in $C_{E_8}(t)^\circ$. As before, we calculate the connected centraliser of $t$ is $\bar{A}_1 E_7$ and hence $A_1 A_1$ is conjugate to $E_8(\#\eeight{973}^{\{0,0\}})$.

This completes the proof of Theorem \ref{E8THM} and thus the proof of Theorem \ref{MAINTHM}.

%% file: Corollaries.tex
\chapter{Corollaries} \label{cors}

In this section we give the proofs of Corollaries \ref{comps}--\ref{A2overgroups}, as well as giving further corollaries not mentioned in the introduction. Let $G$ be a simple exceptional algebraic group over an algebraically closed field of characteristic $p$. 

Corollaries \ref{A1subgroups} and \ref{A2subgroups} are immediate from Theorem \ref{MAINTHM}. Similarly, Corollaries \ref{A1overgroups} and \ref{A2overgroups} follow from the lattice structure given in Tables \ref{G2tab}--\ref{E8tab}. Further, the proofs of Corollaries \ref{comps} and \ref{corconjmin} follow from careful inspection of Tables \ref{G2tabcomps}--\ref{E8tabcomps}, recalling that one can deduce the composition factors of all $G$-irreducible connected subgroups from the composition factors given for the non-diagonally embedded subgroups. We are left to prove Corollary \ref{nongcr}. 

\begin{pf}{Proof of Corollary \ref{nongcr}} 
The strategy for the proof is as follows. For each simple exceptional algebraic group $G$ we find all $M$-irreducible connected subgroups that are not $G$-irreducible. Given such a subgroup $X$ we then check whether it satisfies the hypothesis of Corollary \ref{nongcr}. That is to say, we check whether $X$ is contained reducibly in another reductive, maximal connected subgroup, or if $X$ is contained in a Levi subgroup of $G$. To do this we use the composition factors of $X$ on the minimal or adjoint module for $G$, using restriction from $M$. Of course, since $X$ is $G$-reducible there exists some subgroup $Z$ of the same type as $X$ contained in a Levi factor $L'$ having the same composition factors as $X$. Therefore, we will require the exact module structure of $X$ acting on either the minimal or adjoint module for $G$ to prove that $X$ is not contained in $L'$. 

This has already been done when $X$ is simple in \cite[Corollary 5]{tho2} and \cite[Corollary 2]{tho1}. Moreover, suppose that $Y$ is a subgroup of $M$ containing such a simple subgroup $X$ and that $Y$ is $G$-reducible. Then $Y$ also satisfies the hypothesis of the corollary. 

By studying the proofs of Theorems \ref{G2THM}--\ref{F4THM} we find that there are no non-simple irreducible connected subgroups which are $M$-irreducible yet $G$-reducible when $G$ is of type $G_2$ or $F_4$. So we need only consider the cases where $G$ has type $E_6$, $E_7$ and $E_8$. 

Suppose that $G$ is of type $E_6$ and consider the $M$-irreducible non-simple subgroups that are $E_6$-reducible. These are all found in the proof of Theorem \ref{E6THM}; let $X$ be such a subgroup. Firstly, consider $M = \bar{A}_1 A_5$. Then $X$ is a subgroup of $A_1 A_2 \hookrightarrow \bar{A}_1 A_1 A_2 = E_6(\#\esix{28})$ via $(1,1,10)$. It is shown in Section \ref{E6A1A1A2} that $A_1 A_2$ is contained in $A_2 G_2$ and moreover, is $A_2 G_2$-reducible. Thus $X$ does not satisfy the hypothesis of the corollary. Now suppose that $M = F_4$ or $C_4$ $(p \neq 2)$. Then $X$ is a subgroup of $B_4$ or $B_2^2$, respectively, both of which are contained in a Levi subgroup of type $D_5$. Hence $X$ does not satisfy the hypothesis of Corollary \ref{nongcr}.   

The remaining case when $G$ is of type $E_6$ to consider is $X = A_1 A_1 < G_2 = M$ when $p=2$. We prove that $X$ does not satisfy the hypothesis of the corollary by showing that $X$ is contained reducibly in $\bar{A}_1 A_5$. Let $A B$ be the subgroup $A_1 A_1 < A_5$ acting as $(W(2),1)$ on $V_{A_5}(\lambda_1)$ and let $Y = A_1 A_1 \hookrightarrow \bar{A}_1 A B < \bar{A}_1 A_5$ via $(1_a,1_a,1_b)$. We claim that $X$ is conjugate to $Y$. 

Firstly, note that both $X$ and $Y$ are contained in parabolic subgroups of $E_6$. Moreover, $Y$ is contained in a parabolic subgroup of $\bar{A}_1 A_5$ with Levi factor $L'$ of type $\bar{A}_1 \bar{A}_1 A_3$ with the projection of $Y$ to $L'$ being $L'$-irreducible. The only Levi subgroup of $E_6$ containing $\bar{A}_1^2 A_3$ is $D_5$. Therefore $Y$ is contained in $P=QL$, a $D_5$-parabolic subgroup of $E_6$, with irreducible projection to $L' = D_5$. Specifically, the projection of $Y$ to $L'$ is conjugate to $Z = A_1 A_1$ acting on $V_{D_5}(\lambda_1)$ as $(1,1) + ((0,0)|((2,0) + (0,4))|(0,0))$. Moreover, using the action of $\bar{A}_1 A_5$ on $V_{27}$ we find that $Y$ and $Z$ are non-$\text{GL}_{27}(K)$-conjugate and thus $Y$ is non-$E_6$-cr. 

To prove that $X$ is conjugate to $Y$ we show that there are just two $\text{Aut}(E_6)$-conjugacy classes of subgroups of type $A_1 A_1$ contained in $Q Z$, namely the $E_6$-cr subgroup $Z$ and the non-$E_6$-cr subgroup $Y$. First note that $Q$ is abelian and by \cite{ABS}, $Q$ is an $L'$-module with high weight $\lambda_4$ (or $\lambda_5$ depending on the choice of $D_5$-parabolic subgroup). Moreover, the action of $Z$ on either spin module for $D_5$ is $(1,3) + (2, T(2))$. We need to calculate the Hochschild cohomology group $H^1(Z,Q)$. It is well known that $H^1(A_1,M) = 0$ for the modules $M = 1, T(2), 3$ and $H^1(A_1,2) \cong K$ so applying K\"{u}nneth's formula \cite[10.85]{rot} yields $H^1(Z,Q) \cong K$. Considering the non-trivial action of $Z(L)$ on $Q$ shows that there is just one conjugacy class of subgroups of type $A_1 A_1$ contained in $QZ$ which is not $Q$-conjugate to $Z$. This proves the claim. 

Now we need to prove that a conjugate of $X$ is contained in $P$ with its projection to $L'$ being conjugate to $Z$. From this it follows that $X$ does not satisfy the hypothesis of the corollary: $X$ is either contained in a Levi subgroup or contained reducibly in $\bar{A}_1 A_5$.  In fact, it follows that $X$ is conjugate to $Y$ by considering the action of $X$ on $V_{27}$ (using the action of the maximal subgroup $G_2$ given in \cite[Table~10.2]{LS1}). It is shown in the proof of Theorem \ref{E6THM} that $X$ is contained in a parabolic subgroup of $E_6$. The composition factors of the action of $X$ on $V_{27}$ are $(2,2) /$ $\!\! (2,0) /$ $\!\! (1,3) /$ $\!\! (1,1) /$ $\!\! (0,4) /$ $\!\! (0,2)^2 /$ $\!\! (0,0)^3$. By Lemma \ref{wrongcomps} there exists an irreducible subgroup of a Levi factor with the same composition factors as $X$. It is a routine calculation, using the composition factors in Table \ref{levie6}, to find that $D_5$ is the only Levi subgroup with an irreducible subgroup having the same composition factors as $X$, and that this irreducible subgroup is conjugate to $Z$. Therefore $X$ is contained in $QZ$, as required.     

We note that if one considers $E_6$-conjugacy rather than $\text{Aut}(E_6)$-conjugacy then there are two classes of $D_5$-parabolic subgroups, with representatives $P_1 = Q_1 L_1$ and $P_2 = Q_2 L_2$, say. Therefore, we have two classes of irreducible subgroups $A_1 A_1$, say $Z_1$ in $P_1$ and $Z_2$ in $P_2$. There are also two classes of maximal subgroup $G_2$ and so the class of $X$ splits into two classes, with representatives $X_1$ and $X_2$, say. Furthermore, there is no longer an element acting as a graph automorphism on $\bar{A}_1 A_5$ and thus the class of $Y$ splits into two classes, with representatives $Y_1$ and $Y_2$, say. One then finds that (up to reordering) $X_i$ is conjugate to $Y_i$, both being non-$G$-cr and contained in $Q_i Z_i$.  

Next, let $G$ be of type $E_7$ and consider the $M$-irreducible non-simple subgroups that are $E_7$-reducible. These are all found in the proof of Theorem \ref{E7THM} and we let $X$ be such a subgroup. First suppose that $M = \bar{A}_1 D_6$. Then $X$ is a subgroup of $A_1 C_3 \hookrightarrow \bar{A}_1 A_1 C_3 = E_7(\#\eseven{42})$ via $(1,1)$ with $p=2$. In Section \ref{sec:E7A1A1C3}, we proved that the subgroup $A_1 C_3$ is contained reducibly in $G_2 C_3$. Therefore $X$ does not satisfy the hypothesis of the corollary. 

Now let $M = \bar{A}_2 A_5$. Then $X$ is a subgroup of $A_2 A_1 \hookrightarrow \bar{A}_2 A_2 A_1$ via $(10,01,1)$ with $p=3$. It is shown in Section \ref{sec:E7A2A2A1} that $A_2 A_1$ is contained in the maximal subgroup $A_1 F_4$ and is $A_1 F_4$-reducible. Therefore $X$ does not satisfy the hypothesis of Corollary \ref{nongcr}.  

The last case for $G$ of type $E_7$ we need to consider is $M=A_7$. Here $X$ is any $A_7$-irreducible subgroup contained in $C_4$ when $p=2$. In particular, either $X$ contains the subgroup $Y = A_1$ acting as $1 \otimes 1^{[r]} \otimes 1^{[s]}$ on $V_{A_7}(\lambda_1)$ or $X$ is of type $A_2$ acting as $11$ on $V_{A_7}(\lambda_1)$. From \cite[Corollary~5]{tho2} we know that $Y$ satisfies the hypothesis of the corollary. It remains to consider $X = A_2$. The composition factors of $X$ on $V_{56}$ follow from $V_{56} \downarrow A_7 = V_{A_7}(\lambda_2)/$ $\!\!V_{A_7}(\lambda_6)$ and are thus $30^2 / 11^2 / 03^2 / 00^4$. By using them and considering the composition factors of Levi subgroups from Table \ref{levie7} and reductive, maximal connected subgroups from Table \ref{E7tabcomps}, we conclude that $X$ satisfies the hypothesis of the corollary or is contained in either $\bar{A}_2 A_5$ or a Levi subgroup $E_6$. Using the finite subgroup $S = A_2(4) < X$ and calculation in Magma \cite{magma} we find that $V_{56} \downarrow X = (00| (30 + 03) | 00)^2 + 11^2$. Since $X$ has two direct summands of dimension $20$ and no trivial direct summands, it follows that $X$ is not contained in $\bar{A}_2 A_5$ or $E_6$. Hence $X$ satisfies the hypothesis of the corollary.

Finally, suppose that $G$ is of type $E_8$. There are many subgroups of $D_8$ to consider when $p=2$. Specifically, we need to consider all $D_8$-irreducible subgroups of $A_1 C_4 (\ddagger)$ and $B_4 (\ddagger)$, as defined in Section \ref{D8inE8}. Using Lemma \ref{class} we find that all such subgroups are those given in Table \ref{cortab} as well as all $D_8$-irreducible subgroups of $Y = A_1 \bar{A}_1 C_3 < A_1 C_4 (\ddagger)$. Since $Y$ is a diagonal subgroup of $\bar{A}_1^2 A_1 C_3 = E_8(\#\eeight{122})$, it follows that $Y$ is a subgroup of $G_2 F_4$. Moreover, the projection of $Y$ to $G_2$ is conjugate to $A_1 \hookrightarrow \bar{A}_1 A_1 < G_2$ via $(1,1)$, as may be seen by noting that the $\bar{A}_1 C_3$ factor of $Y$ is contained in $F_4$. Thus Theorem \ref{G2THM} implies that $Y$ is $G_2 F_4$-reducible. Hence $Y$ does not satisfy the hypothesis of the corollary. 

We need to prove that the remaining $D_8$-irreducible subgroups of $A_1 C_4 (\ddagger)$ and $B_4 (\ddagger)$ satisfy the hypothesis of the corollary. In all but the case $X = A_1 A_2 < A_1 C_4$ such subgroups contain the $D_8$-irreducible subgroup of type $A_1$ acting as $1 \otimes 1^{[r]} \otimes 1^{[s]} \otimes 1^{[t]}$ $(0 < r < s < t)$. By \cite[Corollary~5]{tho2}, this subgroup $A_1$ satisfies the hypothesis of Corollary \ref{nongcr} and thus so does every subgroup containing it.  

Now we let $X = A_1 A_2 < A_1 C_4 (\ddagger)$ acting on $V_{D_8}(\lambda_1)$ as $(1,11)$. We consider the action of the finite subgroup $S = A_1(8) A_2(8) < X$ on $L(E_8)$. From \cite[Lemma~11.2]{LS9}, we have  $L(G) \downarrow D_8 = L(D_8) + V_{D_8}(\lambda_7)$. We use the inbuilt functionality of Magma \cite{magma} to find $L(D_8) \downarrow S$ and $V_{D_8}(\lambda_7) \downarrow S$. It follows that both modules have two direct summands, of dimensions $88, 32$ and $96, 32$, respectively. 

We conclude that $L(E_8) \downarrow X$ has a direct summand of dimension at least $96$, another of dimension at least $88$ and no trivial direct summands. It easily follows that $X$ is not contained in any Levi subgroup or another reductive, maximal connected subgroup of $G$. Since $X$ is $G$-reducible it follows that $X$ is non-$G$-cr and satisfies the hypothesis of Corollary \ref{nongcr}.

We now consider the $G$-reducible subgroups of the other reductive, maximal connected subgroups of $G$. Let $M = \bar{A}_1 E_7$ so $X$ is any $M$-irreducible subgroup of $A_1 F_4 \hookrightarrow \bar{A}_1 A_1 F_4 = E_8(\#\eeight{881})$ via $(1,1,\lambda_1)$ with $p=2$. Then $X$ does not satisfy the hypothesis of Corollary \ref{nongcr} because $A_1 F_4$ is $G_2 F_4$-reducible. Similarly, no subgroups of $M = \bar{A}_2 E_6$ satisfy the hypothesis of the corollary. Indeed, $X$ is any $M$-irreducible subgroup of $Y = A_2 G_2 \hookrightarrow \bar{A}_2 A_2 G_2 = E_8(\#\eeight{978})$ via $(10,10,10)$ with $p=3$, and $Y$ is $G_2 F_4$-reducible, as shown in Section \ref{sec:E8A2A2G2}. 

The last case to consider is $M = A_8$ and $X = A_2^2$ acting as $(10,10)$ on $V_{A_8}(\lambda_1)$. Since the diagonal $M$-irreducible subgroups of $X$ satisfy the hypothesis of the corollary it follows that $X$ itself does. \qed 
\end{pf}

\section{Variations of Steinberg's Tensor Product Theorem} \label{steinberg}

We need some background for the next corollary. Let $X$ be a simple, simply connected algebraic group over an algebraically closed field $K$ of characteristic $p < \infty$. We recall Steinberg's tensor product theorem \cite{stein}. It states that if $\phi \colon X \rightarrow \text{SL}(V)$ is an irreducible rational representation, then we can write $V = V_1^{[r_1]} \otimes \cdots \otimes V_k^{[r_k]}$, where the $V_i$ are restricted $X$-modules and the $r_i$ are distinct. The main result of \cite{LS4} generalises this conclusion to the situation where $\phi$ is a rational homomorphism from $X$ to an arbitrary simple algebraic group $G$. To describe this generalisation we need the following definition. Throughout this section we let $G$ be a simple exceptional algebraic group. 

\begin{defn} \textup{\cite[p. 263]{LS4}}
A simple, simply connected subgroup $X$ of $G$ is \emph{restricted} if all composition factors of $L(G) \downarrow X$ are restricted if $X$ is not of type $A_1$, and are of high weight at most $2p-2$ if $X$ is of type $A_1$.    
\end{defn}

\begin{thm} \textup{\cite[Corollary 1]{LS4}} \label{stein}
Assume $p$ is good for $G$. If $X$ is a connected simple $G$-cr subgroup of $G$, then there is a uniquely determined commuting product $E_1 \dots E_k$ with $X \leq E_1 \dots E_k \leq G$, such that each $E_i$ is a simple restricted subgroup of the same type as $X$, and each of the projections $X \rightarrow E_i/Z(E_i)$ is non-trivial and involves a different field twist. 
\end{thm}

Using our classification of $G$-irreducible connected subgroups, we investigate to what extent Theorem \ref{stein} is true in bad characteristics for simple $G$-irreducible connected subgroups. To save repeating ourselves, we say a subgroup $X$ satisfies the conclusion of Theorem \ref{stein} if there is a uniquely determined commuting product $E_1 \dots E_k$ with $X \leq E_1 \dots E_k \leq G$, such that each $E_i$ is a simple restricted subgroup of the same type as $X$, and each of the projections $X \rightarrow E_i/Z(E_i)$ is non-trivial and involves a different field twist. We have already considered the simple irreducible subgroups of rank at least $2$ in \cite[Section~9]{tho1} and so we need only prove the following results for subgroups of type $A_1$. The method is similar for all of them and so we give the proof for $G$ of type $E_8$ only.

\begin{cor} \label{corvarstein}
Let $G$ be a simple algebraic group of exceptional type and $X$ be a simple $G$-irreducible connected subgroup of $G$. Then either $X$ satisfies the conclusion of Theorem \ref{stein} or $p \leq 5$ and $X$ is conjugate to one of subgroups in the following table. We list irreducible subgroups by their identification number $n$ and any conditions given with $n$ refer to the field twists associated with $G(\#n)$.

\end{cor}
\setlength{\LTleft}{0\textwidth}

\begin{longtable}{>{\raggedright\arraybackslash}p{0.07\textwidth - 2\tabcolsep}>{\raggedright\arraybackslash}p{0.07\textwidth - 2\tabcolsep}>{\raggedright\arraybackslash}p{0.86\textwidth - \tabcolsep}@{}}

\caption{\label{varsteintab}} \\

\hline

$G$ & $p$ & $n$ \\

\hline

$G_2$ & $3$ & $\gtwo{5}$ \\

& $2$ & $\gtwo{1}$ \\
\endfirsthead
$F_4$ & $3$ & $\ffour{9}^{\{r,r,s\}}$ \\

& $2$ & $\ffour{2}^{\{0,r,r,s\}}$; $\ffour{2}^{\{r,s,0,r\}}$; $\ffour{14}$; $\ffour{15}$; $\ffour{17}$--$\ffour{20}$; $\ffour{22}$  \\

$E_6$ & $3$ & $\esix{2}^{\{r,r,s\}}$; $\esix{11}$; $\esix{12}$ \\

& $2$ & $\esix{8}$; $\esix{9}$; $\esix{13}$; $\esix{14}$; $\esix{15}$; $\esix{17}$; $\esix{19}$ \\

$E_7$ & $3$ & $\eseven{3}$ if $r=t$; $\eseven{10}^{\{r,r,r,s\}}$; $\eseven{10}^{\{r,r,s,r\}}$; $\eseven{26}$ \\

& $2$ & $\eseven{12}$ if the conditions of lines 4--7, 11, 14, 15 or 16 of Table \ref{condE712} are satisfied; $\eseven{13}$ if either $s$ or $t \in \{u,v,w\}$; $\eseven{28}$ \\

$E_8$ & $5$ & $\eeight{7}$; $\eeight{101}$ \\

& $3$ & $\eeight{6}$ if $v \in \{r,s\}$; $\eeight{24}$ if $t \in \{r,s\} \cap \{u,v\}$; $\eeight{90}$--$\eeight{99}$ \\

& $2$ & $\eeight{8}$; $\eeight{26}$ if the conditions of lines 3, 5, 8, 12, 15, 16 or 17 of Table \ref{E8cond26} are satisfied; $\eeight{28}$ if one of $t, u$ or $v \in \{0,r,s\}$; $\eeight{29}$ if one of $t,u,v,w$ or $x \in \{r,s\}$; $\eeight{53}$--$\eeight{57}$; $\eeight{66}$--$\eeight{81}$ \\

\hline

\end{longtable}

\begin{proof}
Suppose that $X$ is an $E_8$-irreducible subgroup of type $A_1$ not satisfying the conclusion of Theorem \ref{stein}. Then $p \leq 5$ and $X$ is given by Theorem \ref{E8THM}. Moreover, using the lattice structure given in Table \ref{E8tab} we find that $X$ is contained in at most one commuting product of restricted groups, of the same type as $X$, containing $X$ as a diagonal subgroup with distinct field twists. It remains to check whether each subgroup in the commuting product is restricted, using Table \ref{E8tabcomps}. We illustrate this with two examples; the remaining cases are all similar. 

Firstly, suppose that $X = E_8(\#\eeight{7})$ so $X \hookrightarrow A_1 A_1 = Y_1 Y_2 = E_8(\#\eeight{771})$ via $(1^{[r]},1^{[s]})$ $(rs = 0; r \neq s)$. From Table \ref{E8tabcomps} we have 
\begin{align*}
L(E_8) \downarrow Y_1 Y_2 = \ & (W(5),3) / (4,W(6)) / (4,2) / (3,W(7)) / (3,W(5)) / (3,1) / (2,W(8)) / \\ 
\ &  (2,4)^2 / (2,0) / (1,W(9)) / (1,W(5)) / (1,3) / (0,W(6)) / (0,2).
\end{align*}
Thus when $p=5$ we have that $Y_1$ is $p$-restricted but $Y_2$ is not since it has a composition factor of high weight $9$ which is greater than $2p -2 = 8$. Therefore $X$ does not satisfy the conclusion of Theorem \ref{stein}. 

Secondly, suppose that $X = E_8(\#\eeight{28})$ so $p=2$ and $X \hookrightarrow A_1^6 = E_8(\#\eeight{641})$ via $(1,1^{[r]},1^{[s]},1^{[t]},1^{[u]},1^{[v]})$ $(0 < r < s;  t < u < v$; if $t=0$ then $r \leq u$; if $t=0$ and $r=u$ then $s \leq v$). From Table \ref{E8tabcomps} we have 
\begin{align*}
L(E_8) \downarrow A_1^6 = \ & (2,2,0,0,0,0) / (2,0,2,0,0,0) / (2,0,0,1,1,1) / (2,0,0,0,0,0)^2 / \\ \ & (1,1,1,2,0,0) / (1,1,1,1,1,1) / (1,1,1,0,2,0) / (1,1,1,0,0,2) / \\ \ & (1,1,1,0,0,0)^2 / (0,2,2,0,0,0) / (0,2,0,1,1,1) / (0,2,0,0,0,0)^2 / \\ \ & (0,0,2,1,1,1) / (0,0,2,0,0,0)^2 / (0,0,0,2,2,0) / (0,0,0,2,0,2) / \\ \ & (0,0,0,2,0,0)^2 / (0,0,0,1,1,1)^2 /  (0,0,0,0,2,2) / (0,0,0,0,2,0)^2 / \\ \ & (0,0,0,0,0,2)^2 / (0,0,0,0,0,0)^8.
\end{align*} 
If $X$ has distinct field twists then it satisfies the conclusion of Theorem \ref{stein} since $2$ is the highest weight of any composition factor for an $A_1$ factor of $A_1^6$. If two of the field twists of the embedding of $X$ are equal then at least one of $t, u$ or $v \in \{0,r,s\}$. Since the composition factors $(1,1,1,2,0,0)$, $(1,1,1,0,2,0)$ and $(1,1,1,0,0,2)$ occur in $L(E_8) \downarrow A_1^6$ it follows that at least one of the subgroups of type $A_1$ in the commuting product containing $X$ will have a composition factor of high weight $3$ and therefore not be restricted.   
\end{proof}

%% file: TablesThms.tex
\chapter{Tables for Theorem \ref{MAINTHM}} \label{thmtabs}

In this section we present the tables referenced in Theorems \ref{G2THM}--\ref{E8THM}. The notation used is described in Section \ref{nota}. The tables contain all of the $G$-irreducible connected subgroups of each simple exceptional algebraic group $G$. They also contain enough information to find the lattice of connected overgroups of each $G$-irreducible subgroup. The structure of their presentation reflects the strategy described in Section \ref{strat} and the proofs of Theorems \ref{G2THM}--\ref{E8THM} in Sections \ref{secG2}--\ref{secE8}. Moreover, we give an auxiliary table for each exceptional algebraic group $G$, giving the immediate connected overgroups of every $G$-irreducible connected subgroup, except for certain diagonal subgroups. This is intended to help the reader construct the lattice of connected overgroups of each $G$-irreducible connected subgroup. We give further details regarding this lattice structure below.  

We start by explaining how to read Tables \ref{G2tab}--\ref{E8tab}. Let $G$ be a simple exceptional algebraic group. Each table is broken into sections, divided by pairs of horizontal lines, and within each section there are parts divided by single horizontal lines. There is one section for each reductive, maximal connected subgroup of $G$. There is then one part for each $G$-irreducible connected subgroup $X$ of $M$, such that $X$ contains a proper $G$-irreducible connected subgroup. This part gives the $G$-irreducible maximal connected subgroups of $X$, as well as all $G$-irreducible diagonal subgroups. There is a heading for each section and part that gives the type of the subgroup being considered, as well as the identification number and any restrictions on the characteristic $p$. There is one more piece of information in the heading which is the ``$M_i =$'' for $i = 0,1, \dots$ (where $M_0$ is simply written as $M$). This is intended to make it easier for the reader to follow the tables, and is explained below. 

We need to describe what the entries in each column of a generic row represent. If the first column is empty then this row is a heading of a section or part, as described above. There are two other types of rows, both of which correspond to $G$-irreducible subgroups $X$ contained in $Y$, where $Y$ is the subgroup of the current section or part, contained in the maximal connected subgroup $M$ (in this description we are including the case $Y = M$). In both cases, the entry in the first column gives the ID number for the conjugacy class or classes of such irreducible subgroups. The second and third column differ, depending on whether $X$ is a diagonal subgroup of $Y$ or not. If $X$ is a diagonal subgroup of $Y$ then the second and third columns are merged and the information given is of the form ``$X$ via \dots'', denoting the usual notation for an embedding of a diagonal subgroup of $Y$, as well as any restrictions on the characteristic $p$. If $X$ is not a diagonal subgroup of $Y$ then the second column gives the isomorphism type of $X$ and any restrictions on the characteristic $p$. The third column contains the description of $V_M \downarrow X$. Note that for a diagonal subgroup $X$ of $Y$, it is straightforward to work out $V_M \downarrow X$ from $V_M \downarrow Y$. We note that we do not repeat restrictions on $p$ as we list the subgroups of a $G$-irreducible connected subgroup $X$. So if, for example, the maximal connected subgroup $M$ exists only for $p \neq 2$ then we write $M$ $(p \neq 2)$ in the heading and it is assumed that any subgroup $X$ of $M$ inherits this restriction on $p$ without explicit labelling. However, when we consider the subgroups of $X$ we do explicitly repeat any restriction inherited from $M$ in the heading, for clarity.  

In Tables \ref{E7tab} and \ref{E8tab} there is one more possibility for a generic row. To make the tables easier to read we have moved large collections of diagonal subgroups to supplementary tables in Section \ref{sec:irreddiagsubs}. Thus the row ``See Table $x$'' means that all diagonal subgroups of the relevant subgroup $Y$ appear in Table $x$ in Section \ref{sec:irreddiagsubs}.  

We now explain how to find the lattice of connected overgroups of each $G$-irreducible connected subgroup given in Tables \ref{G2tab}--\ref{E8tab}. Each table starts with the $G$-irreducible subgroups that are maximal amongst reductive connected subgroups, with their identification number listed in the ID column. A pair of horizontal lines then indicates the end of that list and the beginning of the first section. We then write ``In $M = H_1 H_2 \dots $ ($G(\#n_1)$)'' where $n_1$ is the identification number for the first reductive, maximal connected subgroup of type $H_1 H_2 \dots $, and will include any restrictions on the characteristic $p$. We then list the $G$-irreducible maximal connected subgroups of $M$ as well as all diagonal connected subgroups of $M$, if there are any, not just the maximal ones. Recall that we will not explicitly consider the proper subgroups of any diagonal connected subgroups, as discussed in Section \ref{strat}. A horizontal line then indicates the end of this list. The next row will be a heading ``In $M_1 = X_1 X_2 \dots$ ($G(\#n_2)$)'', where $X_1 X_2 \dots$ is the first $G$-irreducible maximal connected subgroup of $M$. The ``$M_1=$'' tells the reader that we are now listing the subgroups of a maximal subgroup of a maximal connected subgroup. We then repeat the process, listing the $G$-irreducible maximal connected subgroups of $M_1$ and all diagonal subgroups of $M_1$. The next heading could be ``In $M_2 = Y_1 Y_2 \dots $ ($G(\#n_2)$)'', where $Y_1 Y_2 \dots$ is a maximal connected subgroup of $M_1$ or it could be ``In $M_1 = Z_1 Z_2 \dots $ ($G(\#n_3)$)'', where $Z_1 Z_2 \dots$ is the second maximal connected subgroup of $M$. This will depend on whether $M_1$ has any proper $G$-irreducible connected subgroups that need considering or not. Once all $G$-irreducible connected subgroups of $M$ have been listed in this way, a pair of horizontal lines indicates the end of the subgroups contained in the first reductive, maximal connected subgroup of $G$. The next heading will be ``In $M = K_1 K_2 \dots$ ($G(\#n_4)$)'' and we repeat the process again for the second reductive, maximal connected subgroup $K_1 K_2 \dots$ of $G$. We iterate this process until we have considered all of the $G$-irreducible subgroups contained in the final reductive, maximal connected subgroup. 

There is an important deviation from the process described above. Suppose that $X$ is a representative for a conjugacy class of $G$-irreducible connected subgroups with more than one conjugacy class of proper immediate overgroups. Then after one occurrence of $X$ in the table we list all repetitions of $X$ with the ID number in italics and do not reconsider its subgroups at those points in the table. Moreover, if $X$ occurs more than once in the same reductive, maximal connected subgroup $M$ then we only list $V_M \downarrow X$ in the third column once. We note that the first occurrence of a subgroup in the table can sometimes be in italics because it is clearer to consider its subgroups later in the table. For example, the subgroup $X = \bar{A}_1^2 B_2 \bar{A}_3 = E_8(\#\eeight{351})$ is listed in italics the first time we reach it in Table \ref{E8tab}, as a subgroup of $Y = \bar{A}_1^2 B_5 = E_8(\#\eeight{118})$. This is because the subgroup $X$ is defined for all $p$, whereas it only occurs as a maximal connected subgroup of $Y$ when $p \neq 2$. The second time we reach $X$ is as a subgroup of $Z = \bar{A}_1^2 \bar{A}_3^2 = E_8(\#\eeight{119})$. This time the ID number is not in italics, $V_M \downarrow X$ is given in the third column and we consider the subgroups of $X$ for all $p$ in the next part. 

There is another important example where we do not immediately consider the subgroups contained in $X$, listed as a subgroup of $Y_1$, say, even if the ID number for $X$ is not in italics nor is $X$ a diagonal subgroup of $Y_1$. In this case the ID number will be $n$b and the subgroup $X$ is defined for all $p \geq k$ but only a maximal connected subgroup of $Y_1$ for some prime $l \geq k$. At some point later in the table $X$ will be defined for all $p \geq k$ except $l$ and given ID number $n$a. At this point we will consider the subgroups of $X$ for all $p \geq k$ together. There are instances where the subgroup $X$ occurs again in the table. If $X$ is listed for $p = l$, we write $\mathit{n}$b in the ID column. Similarly, if $X$ is listed for $p \geq k$ except $l$, we write $\mathit{n}$a. It may be that $X$ is listed for all $p \geq k$ in a later reductive, maximal connected subgroup and in this case we simply write $\mathit{n}$ in the ID column. This could cause some confusion. However, upon finding any subgroup $X = G(\#n)$ we can look up $X$ in the auxiliary table for $G$. There are two possibilities. Either $n$ occurs in the first column of a row, in which case the third column of that row gives any characteristic restriction on $X$ in full generality; or $n$ does not occur in the first column of the table, in which case $X$ occurs only in the place one has found it, leaving no confusion about the characteristic restrictions.  

There are many $G$-irreducible connected subgroups occurring multiple times, especially when $G = E_7$ or $E_8$. For this reason we have provided an auxiliary table for each $G$; these are Tables \ref{G2tabaux}--\ref{E8tabaux}. Each row of the table gives a conjugacy class of $G$-irreducible connected subgroups, with the first column giving the ID number, the second column the isomorphism type of a representative $X$ and the third column gives any restrictions on the characteristic $p$. The fourth column gives all conjugacy classes of immediate connected overgroups of $X$. Note that we use the notation $X [\#n]$ to denote the subgroup $X = G(\#n)$, as a shorthand only in these five tables.  

If $X$ is a diagonal subgroup of $Y$ and all immediate connected overgroups of $X$ are also subgroups of $Y$ (and hence diagonal subgroups of $Y$ or just $Y$ itself) then $X$ does not appear in the auxiliary table. In this case there will be only one appearance of $X$ in Tables \ref{G2tab}--\ref{E8tab} (or the supplementary tables in Section \ref{sec:irreddiagsubs}) and all of its immediate overgroups can be straightforwardly computed from the other diagonal subgroups that will appear just above $X$ in the table. For this reason, the phrase ``immediate overgroups in $Y[\#n]$'' sometimes appears in the fourth column of a row corresponding to the irreducible subgroup $X$. This means that the subgroup $X$ is a diagonal subgroup of $Y$ and the immediate connected overgroups of $X$ in $Y$ are included in the list but not explicitly calculated.  

\begin{figure}
\centering
\begin{tikzpicture}

\node at (0,0) {};

\node at (\textwidth,0) {};

\node at (0,-4) {};

\node at (\textwidth,-4) {};

\draw  (0.5\textwidth,0.5) node[anchor=north,align=center]{$G_2 = G_2(\#\gtwo{0})$}; 

\draw (0,-2) node[anchor=west,align=center]{$\bar{A}_2 = G_2(\#\gtwo{4})$ \\ ${}$};

\draw (0.25\textwidth,-2) node[anchor=west,align=center]{$\bar{A}_1 \tilde{A}_1 = G_2(\#\gtwo{6})$ \\ ${}$};

\draw (0.55\textwidth,-2) node[anchor=west,align=center]{$\tilde{A}_2 = G_2(\#\gtwo{5})$ \\ $(p=3)$};

\draw (0.8\textwidth,-2) node[anchor=west,align=center]{$A_1 = G_2(\#\gtwo{3})$ \\ $(p \geq 7)$};

\draw (0,-4.5) node[anchor=west,align=center]{$A_1 = G_2(\#\gtwo{2})$ \\ $p\neq2$};

\draw (0.6\textwidth,-4.5) node[align=center,anchor=west]{$A_1 = G_2(\#\gtwo{1}^{\{1,0\}})$ \\ $(p = 3)$};

\draw (0.2\textwidth,-4.5) node[align=center,anchor=west]{$A_1 = G_2(\#\gtwo{1})$ \\ (if $p=3$ then $(r,s) \neq (1,0)$)};

\draw (0.1\textwidth,-1.5)--(0.5\textwidth,0);

\draw (0.35\textwidth,-1.5)--(0.5\textwidth,0);

\draw (0.65\textwidth,-1.5)--(0.5\textwidth,0);

\draw (0.9\textwidth,-1.5)--(0.5\textwidth,0);

\draw (0.1\textwidth,-2.2)--(0.1\textwidth,-3.9);

\draw (0.35\textwidth,-2.2)--(0.1\textwidth,-3.9);

\draw (0.35\textwidth,-2.2)--(0.35\textwidth,-3.9);

\draw (0.35\textwidth,-2.2)--(0.67\textwidth,-3.9);

\draw (0.67\textwidth,-2.5)--(0.67\textwidth,-3.9);









\end{tikzpicture}
\vspace{5pt}
\caption{The lattice of $G_2$-irreducible connected subgroups. \label{picg2subs}} 

\end{figure}
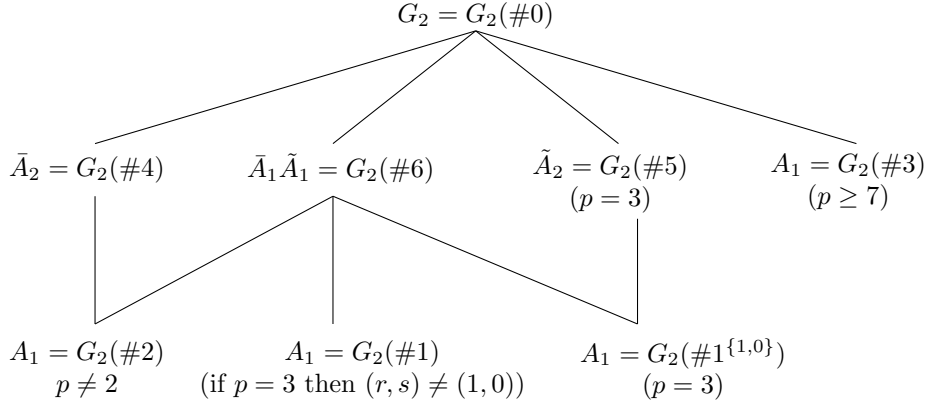

We give some examples of how to use the information in the two sets of tables to recover the lattice of overgroups for a given $G$-irreducible connected subgroup $X$. They are chosen to highlight as many different scenarios as possible. Firstly, let $G = G_2$. In this case we can easily recover the lattice of overgroups for all $G$-irreducible connected subgroups and even present this in a small diagram, see Figure \ref{picg2subs}. Table \ref{G2tab} contains precisely the same information as this diagram and in this case Table \ref{G2tabaux} is unnecessary, but we include it for completeness.

Now suppose $G = E_8$. For the first examples, we highlight the fact that a condition on the characteristic $p$ of a subgroup $X = G(\#n)$ contained in $Y_1$ is based on that specific embedding of $X$ into $Y_1$ and thus the subgroup $X$ may be considered for other excluded characteristics in a different overgroup $Y_2$.  For example, the subgroup $X = A_1 B_6 = E_8(\#\eeight{110})$ is defined for all $p$. However, we first arrive at it in Table \ref{E8tab} as a maximal connected subgroup of $D_8 = E_8(\#\eeight{43})$, with ID number $110$a and a restriction ``$(p \neq 2)$''. In this case, we know that the ``a'' implies that $A_1 B_6$ must be a maximal connected subgroup of at least one other irreducible connected subgroup of $G$ when $p=2$. To find all of the extra immediate overgroups we go to the row with $110$ in the first column in Table \ref{E8tabaux}. We find that $A_1 B_6$ is also a maximal connected subgroup of $B_7 = E_8(\#\eeight{44})$ when $p=2$. Indeed, when we come to the $M_1 = B_7$ part of the table we see that $X$ appears with ID number $110$b. 

For another example, suppose $X = \bar{A}_1^6 A_1 = E_8(\#\eeight{132}^{\{0\}})$. Then $X$ is a maximal connected subgroup of $\bar{A}_1^4 B_3 = E_8(\#\eeight{125})$ and $\bar{A}_1^4 A_1 B_2 = E_8(\#\eeight{126})$, and in both cases the ID number is listed in italics. The subgroup is not listed again in Table \ref{E8tab}. Therefore $X$ must be a diagonal subgroup listed in one of the supplementary tables in Section \ref{sec:irredE8diags}. Using Table \ref{E8tabaux} we immediately find that $X$ is a maximal connected subgroup of $\bar{A}_1^8 = E_8(\#\eeight{124})$.  

Finally, let us suppose that we are interested in subgroups of the reductive, maximal connected subgroup $A_8 = E_8(\#\eeight{62})$. We see from the $M = A_8$ section of Table \ref{E8tab} that there are two irreducible maximal connected subgroups, namely $X_1 = B_4 = E_8(\#\eeight{46})$ when $p \neq 2$, and $X_2 = A_2^2 = E_8(\#\eeight{669})$ when $p \neq 3$. Both ID numbers are in italics and therefore both $X_1$ and $X_2$ are contained in other irreducible connected subgroups of $G$. Looking at Table \ref{E8tabaux} we find that $X_1$ is a maximal connected subgroup of $D_8 = E_8(\#\eeight{43})$ also when $p \neq 2$, and indeed $X_1$ is only to be considered when $p \neq 2$, as per the information in the table. Similarly, we find that $X_2$ is a maximal connected subgroup of $A_2 \bar{D}_4 = E_8(\#\eeight{621})$ when $p \neq 3$ and again $X_2$ is only considered when $p \neq 3$. We can then continue looking up subgroups in Table \ref{E8tabaux} to find that when $p \neq 3$ the subgroup $A_2 \bar{D}_4$ is a maximal connected subgroup of $\bar{D}_4^2 = E_8(\#\eeight{108})$, which in turn is a maximal connected subgroup of $D_8 = E_8(\#\eeight{43})$.       

\setlength{\LTleft}{0.215\textwidth}

\begin{longtable}{>{\raggedright\arraybackslash}p{0.07\textwidth - 2\tabcolsep}>{\raggedright\arraybackslash}p{0.35\textwidth - 2\tabcolsep}>{\raggedright\arraybackslash}p{0.15\textwidth - \tabcolsep}@{}}

\caption{The lattice structure of irreducible connected subgroups of $G_2$. \label{G2tab}} \\

\hline

ID & Irreducible subgroup $X$ & $V_M \downarrow X$ \\

\hline

\gtwo{4} & $\bar{A}_2$ & \\

\gtwo{5} & $\tilde{A}_2$ $(p=3)$ & \\
\endfirsthead
\gtwo{6} & $\bar{A}_1 \tilde{A}_1$ & \\

\gtwo{3} & $A_1$ $(p \geq 7)$ &   \\

\hline \hline 

& In $M = \bar{A}_2$ ($G_2(\#\gtwo{4})$) & \\

\textit{\gtwo{2}} & $A_1$ $(p \neq 2)$ & $2$ \\

\hline \hline

& \multicolumn{2}{l}{In $M = \tilde{A}_2$ ($p=3$) ($G_2(\#\gtwo{5})$)} \\

$\mathit{\gtwo{1}^{\{\delta_1\}}}$ & $A_1$ & $2$ \\

\hline \hline

& In $M = \bar{A}_1 \tilde{A}_1$ ($G_2(\#\gtwo{6})$) & \\

\gtwo{1} & \multicolumn{2}{l}{$A_1$ via $(1^{[r]},1^{[s]})$ $(rs=0; r \neq s)$}  \\

$\gtwo{2}$ & $A_1$ via $(1,1)$ $(p \neq 2)$ & \\

\hline \hline

\end{longtable}

\setlength{\LTleft}{0.17\textwidth}
\modcounter

\begin{longtable}{p{0.07\textwidth - 2\tabcolsep}>{\raggedright\arraybackslash}p{0.08\textwidth - 2\tabcolsep}>{\raggedright\arraybackslash}p{0.08\textwidth - 2\tabcolsep}>{\raggedright\arraybackslash}p{0.43\textwidth-\tabcolsep}@{}}

\caption{The connected overgroups for $G$-irreducible subgroups of $G_2$. \label{G2tabaux}} \\

\hline

ID & $X$ & $p$ & Immediate connected overgroups \\

\hline

$\gtwo{1}^{\{\delta_1\}}$ & $A_1$ &  & $\tilde{A}_2[\#\gtwo{5}]$ $(p=3)$, $\bar{A}_1 \tilde{A}_1 [\#\gtwo{6}]$ \\ 

\gtwo{2} & $A_1$ & $ \neq 2$ & $\bar{A}_2 [\#\gtwo{4}]$, $\bar{A}_1 \tilde{A}_1 [\#\gtwo{6}]$ \\
\endfirsthead
\gtwo{3} & $A_1$ & $\geq 7$ & $G_2 [\#\gtwo{0}]$ \\

\gtwo{4} & $\bar{A}_2$ & & $G_2 [\#\gtwo{0}]$ \\

\gtwo{5} & $\tilde{A}_2$ & $3$ & $G_2 [\#\gtwo{0}]$ \\

\gtwo{6} & $\bar{A}_1 \tilde{A}_1$ & & $G_2 [\#\gtwo{0}]$ \\

\hline

\end{longtable}

\unmodcounter
\setlength{\LTleft}{0pt}

\begin{longtable}{>{\raggedright\arraybackslash}p{0.08\textwidth - 2\tabcolsep}>{\raggedright\arraybackslash}p{0.27\textwidth - 2\tabcolsep}>{\raggedright\arraybackslash}p{0.65\textwidth-\tabcolsep}@{}}

\caption{The lattice structure of irreducible connected subgroups of $F_4$. \label{F4tab}} \\

\hline

ID & Irred. subgroup $X$ & $V_M \downarrow X$ \\

\hline

\ffour{12} & $B_4$ &  \\

\ffour{14} & $C_4$ $(p=2)$ &  \\
\endfirsthead
\ffour{24}a & $\bar{A}_1 C_3$ $(p \neq 2)$ &  \\

\ffour{25} & $A_1 G_2$ $(p \neq 2)$ &  \\

\ffour{26} & $\bar{A}_2 \tilde{A}_2$ &  \\

\ffour{16} & $G_2$ $(p=7)$ &  \\

\ffour{10} & $A_1$ $(p \geq 13)$ & \\

\hline \hline 

& \multicolumn{2}{l}{In $M = B_4$ $(F_4(\#\ffour{12}))$} \\

\ffour{13} & $\bar{D}_4$ & $\lambda_1 / 0$ \\

\ffour{27}a & $\tilde{A}_1 \bar{A}_3$ $(p \neq 2)$ & $(2,000) / (0,010) $ \\
\hline \hline 
\end{longtable}

\setlength{\LTleft}{0.05\textwidth}
\modcounter

\begin{longtable}{p{0.1\textwidth - 2\tabcolsep}>{\raggedright\arraybackslash}p{0.12\textwidth - 2\tabcolsep}>{\raggedright\arraybackslash}p{0.1\textwidth - 2\tabcolsep}>{\raggedright\arraybackslash}p{0.58\textwidth-\tabcolsep}@{}}

\caption{The connected overgroups for irreducible subgroups of $F_4$. \label{F4tabaux}} \\

\hline

ID & $X$ & $p$ & Immediate connected overgroups   \\

\hline

$\ffour{6}^{\{0\}}$ & $A_1$ & $\geq 5$ & $A_1^2 [\#\ffour{29}]$, $A_1 A_1 [\#\ffour{41}^{\{0\}}]$, $A_1 \bar{A}_1 [\#\ffour{42}^{\{\underline{0}\}}]$ \\

$\ffour{7}$ & $A_1$ & $\geq 11$ & $B_4[\#\ffour{12}]$ \\
\endfirsthead
$\ffour{8}^{\{\underline{0}\}}$ & $A_1$ & $\geq 7$ & $\bar{A}_1 A_1 [\#\ffour{58}]$, $A_1 A_1[\#\ffour{66}]$ \\

$\ffour{8}^{\{\delta_1\}}$ & $A_1$ & $\geq 7$ & $G_2[\#\ffour{16}]$ $(p=7)$, $\bar{A}_1 A_1 [\#\ffour{58}]$ \\

$\ffour{10}$ & $A_1$ & $\geq 13$ & $F_4[\#\ffour{0}]$ \\

$\ffour{12}$ & $B_4$ & & $F_4[\#\ffour{0}]$ \\

$\ffour{13}$ & $\bar{D}_4$ & & $B_4[\#\ffour{12}]$ \\

$\ffour{14}$ & $C_4$ & $2$ & $F_4[\#\ffour{0}]$ \\

$\ffour{15}$ & $\tilde{D}_4$ & $2$ & $C_4[\#\ffour{14}]$ \\

$\ffour{16}$ & $G_2$ & $7$ & $F_4[\#\ffour{0}]$ \\
\hline

\end{longtable}

\unmodcounter

\setlength{\LTleft}{0pt}

\begin{longtable}{>{\raggedright\arraybackslash}p{0.1\textwidth - 2\tabcolsep}>{\raggedright\arraybackslash}p{0.5\textwidth - 2\tabcolsep}>{\raggedright\arraybackslash}p{0.4\textwidth-\tabcolsep}@{}}

\caption{The lattice structure of irreducible connected subgroups of $E_6$. \label{E6tab}} \\

\hline

ID & Irred. subgroup $X$ & $V_{M} \downarrow X$ \\

\hline

\esix{24} & $\bar{A}_1 A_5$ &  \\

\esix{25} & $\bar{A}_2^3$ & \\
\endfirsthead
\esix{26} & $A_2 G_2$ & \\

\esix{7} & $F_4$ & \\

\esix{8}a & $C_4$ $(p \neq 2)$ & \\

\esix{10}a & $G_2$ $(p \neq 7)$ & \\

\esix{11}a & $A_2$ $(p \geq 5)$ & \\

\hline \hline 

& In $M = \bar{A}_1 A_5$ ($E_6(\#\esix{24})$) & \\

\esix{27} & $\bar{A}_1 C_3$ & $(1,100)$ \\

\esix{28} & $\bar{A}_1 A_1 A_2$ & $(1,1,10)$ \\
\hline \hline 

\end{longtable}

\modcounter

\setlength{\LTleft}{0pt}

\begin{longtable}{p{0.08\textwidth - 2\tabcolsep}>{\raggedright\arraybackslash}p{0.12\textwidth - 2\tabcolsep}>{\raggedright\arraybackslash}p{0.09\textwidth - 2\tabcolsep}>{\raggedright\arraybackslash}p{0.71\textwidth-\tabcolsep}@{}}

\caption{The connected overgroups for irreducible subgroups of $E_6$. \label{E6tabaux}} \\

\hline

ID & $X$ & $p$ &  Immediate connected overgroups  \\

\hline

$\esix{1}^{\{\underline{0}\}}$ & $A_1$ & $\geq 7$ & $A_2 [\#\esix{11}]$, $\bar{A}_1 A_1[\#\esix{32}]$, $A_1 A_1 [\#\esix{51}]$ \\

$\esix{1}^{\{\delta_1\}}$ & $A_1$ & $\geq 7$ & $G_2[\#\esix{10}]$ $(p=7)$, $\bar{A}_1 A_1[\#\esix{32}]$ \\
\endfirsthead
$\esix{2}^{\{\delta_3\}}$ & $A_1$ & $\neq 2$ & $A_2 [\#\esix{11}]$ $(p=3)$, $A_1 A_1 [\#\esix{34}^{\{\underline{0}\}}]$, $A_1 A_1 [\#\esix{35}^{\{\delta_2\}}]$, $\bar{A}_1 A_1 [\#\esix{36}^{\{\delta_2\}}]$   \\

$\esix{5}$ & $A_1$ & $\geq 13$ & $F_4[\#\esix{7}]$ \\

$\esix{6}$ & $A_1$ & $\geq 11$ & $C_4 [\#\esix{8}]$, $G_2 [\#\esix{10}]$ \\

$\esix{7}$ & $F_4$ & & $E_6[\#\esix{0}]$ \\

$\esix{8}$ & $C_4$ & & $E_6 [\#\esix{0}]$ $(p \neq 2)$, $F_4 [\#\esix{7}]$ $(p=2)$ \\

$\esix{9}$ & $D_4$ & $2$ & $C_4 [\#\esix{8}]$ \\

$\esix{10}$ & $G_2$ & & $E_6 [\#\esix{0}]$ $(p \neq 7)$, $F_4 [\#\esix{7}]$ $(p=7)$ \\

$\esix{11}$ & $A_2$ & $\neq 2$ & $E_6 [\#\esix{0}]$ $(p \geq 5)$, $G_2 [\#\esix{10}]$ $(p=3)$ \\
\hline

\end{longtable}

\unmodcounter
\begin{longtable}{>{\raggedright\arraybackslash}p{0.1\textwidth - 2\tabcolsep}>{\raggedright\arraybackslash}p{0.27\textwidth - 2\tabcolsep}>{\raggedright\arraybackslash}p{0.63\textwidth-\tabcolsep}@{}}

\caption{The lattice structure of irreducible connected subgroups of $E_7$. \label{E7tab}} \\

\hline

ID & Irred. subgroup $X$ & $V_{M} \downarrow X$ \\

\hline

\eseven{30} & $\bar{A}_1 D_6$  & \\

\eseven{31} & $\bar{A}_2 A_5$ & \\
\endfirsthead
\eseven{22} & $A_7$ & \\

\eseven{32} & $G_2 C_3$ & \\

\eseven{33} & $A_1 F_4$ & \\

\eseven{34} & $A_1 G_2$ $(p \neq 2)$ & \\

\eseven{35} & $A_1 A_1$ ($p \geq 5$) & \\

\eseven{29} & $A_2$ ($p \geq 5$) & \\

\eseven{20} & $A_1$ ($p \geq 17$) & \\

\eseven{21} & $A_1$ ($p \geq 19$) & \\
\hline \hline 

\end{longtable}

\modcounter
\begin{longtable}{p{0.09\textwidth - 2\tabcolsep}>{\raggedright\arraybackslash}p{0.14\textwidth - 2\tabcolsep}>{\raggedright\arraybackslash}p{0.09\textwidth - 2\tabcolsep}>{\raggedright\arraybackslash}p{0.68\textwidth-\tabcolsep}@{}}

\caption{The connected overgroups for irreducible subgroups of $E_7$. \label{E7tabaux}} \\

\hline

ID & $X$ & $p$ & Immediate connected overgroups \\

\hline

$\eseven{5}^{\{\underline{0}\}}$ & $A_1$ & $\geq 11$ & immediate overgroups in $\bar{A}_1 A_1 A_1 [\#\eseven{187}]$, $A_1 A_1 [\#\eseven{324}]$, $A_1 A_1 [\#\eseven{334}]$ \\

$\eseven{20}$ & $A_1$ & $\geq 17$ & $E_7[\#\eseven{0}]$ \\
\endfirsthead
$\eseven{21}$ & $A_1$ & $\geq 19$ & $E_7[\#\eseven{0}]$ \\

$\eseven{22}$ & $A_7$ & & $E_7[\#\eseven{0}]$ \\

$\eseven{23}$ & $D_4$ & $\neq 2$ & $A_7[\#\eseven{22}]$ \\

$\eseven{24}$ & $A_2$ & $\geq 5$ & $D_4 [\#\eseven{23}]$, $\bar{A}_2 A_2 [\#\eseven{303}]$ \\

$\eseven{29}$ & $A_2$ & $\geq 5$ & $E_7[\#\eseven{0}]$ \\

$\eseven{30}$ & $\bar{A}_1 D_6$ & & $E_7[\#\eseven{0}]$ \\

$\eseven{31}$ & $\bar{A}_2 A_5$ & & $E_7[\#\eseven{0}]$ \\

$\eseven{32}$ & $G_2 C_3$ & & $E_7[\#\eseven{0}]$ \\
\hline

\end{longtable}

\unmodcounter

\begin{longtable}{>{\raggedright\arraybackslash}p{0.1\textwidth - 2\tabcolsep}>{\raggedright\arraybackslash}p{0.25\textwidth - 2\tabcolsep}>{\raggedright\arraybackslash}p{0.65\textwidth-\tabcolsep}@{}}

\caption{The lattice structure of irreducible connected subgroups of $E_8$. \label{E8tab}} \\

\hline

ID & Irred. subgroup $X$ & $V_M \downarrow X$ \\

\hline

\eeight{43} & $D_8$ &  \\

\eeight{102} & $\bar{A}_1 E_7$ & \\
\endfirsthead
\eeight{103} & $\bar{A}_2 E_6$ & \\

\eeight{62} & $A_8$ & \\

\eeight{104} & $\bar{A}_4^2$ & \\

\eeight{105} & $G_2 F_4$ & \\

\eeight{101} & $B_2$ $(p \geq 5)$ & \\

\eeight{106} & $A_1 A_2$ $(p \geq 5)$ & \\

\eeight{40} & $A_1$ $(p \geq 23)$  & \\

\eeight{41} & $A_1$ $(p \geq 29)$ & \\
\hline \hline 

\end{longtable}

\modcounter

\begin{longtable}{p{0.11\textwidth - 2\tabcolsep}>{\raggedright\arraybackslash}p{0.15\textwidth - 2\tabcolsep}>{\raggedright\arraybackslash}p{0.09\textwidth - 2\tabcolsep}>{\raggedright\arraybackslash}p{0.65\textwidth-\tabcolsep}@{}}

\caption{The connected overgroups for irreducible subgroups of $E_8$. \label{E8tabaux}} \\

\hline

ID & $X$ & $p$ & Immediate connected overgroups \\

\hline

$\eeight{10}$ & $A_1$ & $\geq 17$ & $B_7 [\#\eeight{44}]$ \\

$\eeight{11}^{\{\underline{0}\}}$ & $A_1$ & $\geq 13$ & $A_1 A_1 [\#\eeight{772}]$, $A_1 A_1 [\#\eeight{1040}]$ \\
\endfirsthead
$\eeight{12}^{\{\underline{0}\}}$ & $A_1$ & $\geq 11$ & $B_4 [\#\eeight{45}]$, $B_4 [\#\eeight{46}]$, $B_2 [\#\eeight{101}]$, $A_1 A_1 [\#\eeight{838}]$ \\

$\eeight{14}^{\{\underline{0}\}}$ & $A_1$ & $\geq 11$ & $A_1 A_1 [\#\eeight{854}]$, $A_1 A_1 [\#\eeight{866}]$, $A_1 A_1 [\#\eeight{874}]$ \\

$\eeight{17}^{\{\underline{0}\}}$ & $A_1$ & $\geq 7$ & immediate overgroups in $\bar{A}_1^2 A_1 A_1 [\#\eeight{538}]$, $A_1 A_1 [\#\eeight{771}]$, $A_1^2 [\#\eeight{862}]$, $A_1^2 [\#\eeight{864}]$ \\ 

$\eeight{17}^{\{s,s,t,r\}}$ & $A_1$ & $\geq 7$ & immediate overgroups in $\bar{A}_1^2 A_1 A_1 [\#\eeight{538}]$, $A_1 A_1 A_1 [\#\eeight{630}]$ \\

$\eeight{40}$ & $A_1$ & $\geq 23$ & $E_8 [\#\eeight{0}]$ \\

$\eeight{41}$ & $A_1$ & $\geq 29$ & $E_8 [\#\eeight{0}]$ \\

$\eeight{42}$ & $A_1$ & $\geq 31$ & $E_8 [\#\eeight{0}]$ \\

$\eeight{43}$ & $D_8$ & & $E_8 [\#\eeight{0}]$ \\
\hline

\end{longtable}

\unmodcounter

\section{Irreducible diagonal subgroups} \label{sec:irreddiagsubs}

In this section we give the tables of diagonal subgroups referred to in Tables \ref{E7tab} and \ref{E8tab}. The first column gives the ID number, as in the previous tables, and the second column gives the embedding of the diagonal subgroups. To describe the embeddings we use a slightly modified notation, to shorten the tables. Specifically, we introduce a shorthand for diagonal subgroups of $A_1^n Z$, where $Z$ has no simple factor of type $A_1$. For example, instead of writing $A_1^3 B_2 \hookrightarrow A_1^4 B_2$ via $(1_a^{[r]},1_a^{[s]},1_b,1_c,10)$ we just write $(a^{[r]},a^{[s]},b,c,10)$; from any such vector it is easy to recover the isomorphism type of the diagonal subgroup. Similarly, in Tables \ref{barA12A1B23diags}, \ref{barA24diags} and \ref{barA2A2A2diags} the usual notation for diagonal subgroups is used but we again omit the isomorphism type of each diagonal subgroup as they too can be easily recovered from the listed embedding. 

There are further tables which give the extra restrictions on the field twists in certain diagonal embeddings. These restrictions ensure there is no repetition of conjugacy classes and further, that each conjugacy class is $G$-irreducible. The restrictions are given in rows of the tables: the first column lists all permitted equalities amongst certain subsets of the field twists; the second column lists any further requirements. So an ordered set $\{0, r, \ldots \}$ is permitted if it satisfies the conditions in the first and second column of a row of the table. We note that a set of field twists satisfies the conditions of at most one row. We emphasise that an ordered set may be excluded either because it yields a $G$-reducible subgroup, or because it yields a repeated diagonal subgroup.

We give an example to illustrate this. Let $X = A_1 A_1 = E_8(\#\eeight{165})$, so that $X$ is a diagonal subgroup of $E_8(\#\eeight{124}) = \bar{A}_1^8 < D_8$ via $(1_a^{[r]}, 1_a^{[s]}, 1_a^{[t]}, 1_a^{[u]}, 1_b^{[v]},1_b^{[w]},1_b^{[x]},1_b^{[y]})$ with $rstu = vwxy=0$. Table \ref{condE8223} gives the extra conditions that an ordered set $r, s, \ldots, y$ needs to satisfy. The conditions in the first column only restrict the equalities allowed among elements from $r,s,t,u$ and separately the equalities allowed among elements from $v,w,x,y$. So in every row of the table $r$ is permitted to be equal to $v$, for example. 

In particular, the ordered set $0,1,2,3,0,3,1,2$ satisfies the conditions of the first row, as does $0,1,2,3,0,1,2,3$. We also note that the ordered set $0,0,1,2,0,1,2,3$ satisfies the conditions of the second row, whereas $2,1,0,0,0,1,2,3$ and $0,1,2,3,0,0,1,2$ do not satisfy the conditions of any of the rows. This is because these three ordered sets of field twists yield the same conjugacy class of $G$-irreducible subgroups and thus only one of them can be permitted.    

\subsection{Irreducible diagonal subgroups contained in $E_7$} \label{sec:irredE7diags}
\leavevmode

\begin{longtable}{>{\raggedright\arraybackslash}p{0.05\textwidth - 2\tabcolsep}>{\raggedright\arraybackslash}p{0.95\textwidth-\tabcolsep}@{}}

\caption{Irreducible diagonal subgroups of $\bar{A}_1^7 = E_7(\#\eseven{43})$.  \label{A17E7diags}} \\
\hline
ID & Embedding \\

\hline

\eseven{49} & $(a,a^{[r]},b,c,d,e,f)$  \\

\eseven{50} & $(a,a^{[r]},a^{[s]},b,c,d,e)$ $(r \leq s)$ \\
\endfirsthead
\eseven{51} & $(a,a^{[r]},b,a^{[s]},c,d,e)$ $(r \leq s)$  \\

\eseven{52} & $(a,a^{[r]},b^{[s]},b^{[t]},c,d,e)$ $(st=0)$  \\

\eseven{53} & $(a^{[r]},a^{[s]},b,c,c^{[t]},d,e)$ $(rs=0)$  \\

\eseven{54} & $(a,a^{[r]},b,c,d,e,c^{[s]})$ $(r \leq s; s \neq 0)$ \\

\eseven{55} & $(a^{[r]},a^{[s]},a^{[t]},a^{[u]},b,c,d)$ $(ru=0; r \leq s \leq t)$  \\

\eseven{56} & $(a,a^{[r]},b,a^{[s]},c,d,a^{[t]})$ $(r \leq s \leq t$; if two of $0,r,s,t$ are equal then the other two are not equal)  \\

\eseven{57} & $(a^{[r]},a^{[s]},a^{[t]},b,b^{[u]},c,d)$ $(rs=0; s \leq t;$ if $u=0$ then $s < t$)  \\

\eseven{58} & $(a^{[r]},a^{[s]},b,a^{[t]},b^{[u]},c,d)$ $(rst=0;$ if $u=0$ then $s < t$)  \\
\hline
\end{longtable}

\setlength{\LTleft}{0.09\textwidth}
\begin{longtable}{>{\raggedright\arraybackslash}p{0.37\textwidth-2\tabcolsep}>{\raggedright\arraybackslash}p{0.45\textwidth-\tabcolsep}@{}}
\caption{Conditions on field twists for $E_7(\#\eseven{71})$. \label{condE771}} \\

\hline \noalign{\smallskip}

Equalities among $0, r, \ldots, v$ & Further requirements on $0, r, \ldots, v$ \\

\hline 

none & $r < \text{min}\{s,t,u\}$ \\

$r=s$ & $t<u$ \\
\endfirsthead
$r=t$ & $s<u$ \\

$r=u$ & $s<t$ \\

$r=v$ & none \\

$r=s=t$ & none \\

$r=s=v$ & $t < u$ \\

$r=t=v$ & $s < u$ \\

$r=u=v$ & $s < t$ \\

$r=s=t=v$ & none \\

$r=0$ & $u < v$ \\

$r=0$; $s=t$ & $u < v$ \\

$r=0$; $s=u$ & none \\

$r=0$; $t=u$ & none \\

$r=0$; $s=t=u$ & none \\

$v=0$ & $r < \text{min}\{s,u\}$; $s < t$ \\

$v=0$; $r=s$ & $t < u$ \\

$r=s=0$ & $t < u < v$ \\

$r=t=0$ & $s < u < v$ \\

$r=u=0$ & $s < t$ \\

$r=u=0$; $s=v$ & none \\

$r=s=t=0$ & $u < v$ \\
\hline
\end{longtable}

\setlength{\LTleft}{0.11\textwidth}
\begin{longtable}{>{\raggedright\arraybackslash}p{0.33\textwidth-2\tabcolsep}>{\raggedright\arraybackslash}p{0.45\textwidth-\tabcolsep}@{}}
\caption{Conditions on field twists for $E_7(\#\eseven{73})$. \label{condE773}} \\

\hline \noalign{\smallskip}

Equalities among $0, v, w$ & Further requirements on $0,r, \ldots, w$ \\

\hline 

none & $v < w$  \\

$v=0$ & $r < s$ \\
\endfirsthead
$v=w=0$ & $\eseven{0} < r < s$ \\
\hline
\end{longtable}

\setlength{\LTleft}{0.1\textwidth}

\begin{longtable}{>{\raggedright\arraybackslash}p{0.35\textwidth-2\tabcolsep}>{\raggedright\arraybackslash}p{0.45\textwidth-\tabcolsep}@{}}
\caption{Conditions on field twists for $E_7(\#\eseven{74})$. \label{condE774}} \\

\hline \noalign{\smallskip}

Equalities among $0, r, s, t$ & Further requirements on $0, r, \ldots, w$ \\
\hline 

none & $\eseven{0} < r < s < t$ \\

$r=0$ & $s < t$; $v < w$ \\
\endfirsthead
$r=s=0$ & $v < w$ \\
\hline
\end{longtable}

\setlength{\LTleft}{0.09\textwidth}
\begin{longtable}{>{\raggedright\arraybackslash}p{0.37\textwidth-2\tabcolsep}>{\raggedright\arraybackslash}p{0.45\textwidth-\tabcolsep}@{}}
\caption{Conditions on field twists for $E_7(\#\eseven{12})$ and $E_8(\#\eeight{162})$. \label{condE712}} \\

\hline \noalign{\smallskip}

Equalities among $0, r, \ldots, w$ & Further requirements on $0, r, \ldots, w$ \\

\hline 

none & $r < \text{min}\{s,t\}$; $t<\text{min}\{u,v,w\}$ \\

$r=s$ & $t<\text{min}\{u,v\}$; $v<w$ \\
\endfirsthead
$r=t$ & $s < v$ \\

$r=s=t$ & none \\

$r=t=w$ & $s < u,v$ \\

$r=s=t=v$ & $u < w$ \\

$r=t$; $u=v=w$ & none \\

$r = 0$ & $t < \text{min}\{u,v\}$; $v < w$ \\

$r=0$; $s=t$ & $u < v$ \\

$r=0$; $t=u$ & $v < w$ \\

$r=0$; $s=t=u$ & $v < w$ \\

$r=0$; $s=t$; $u=w$ & none \\

$r=s=0$ & $t < u < v < w$ \\

$r=t=0$ & $s < u < v$ \\

$r=t=0$; $s=w$ & $u<v$ \\

$r=s=t=0$ & $u < v < w$ \\
\hline
\end{longtable}

\setlength{\LTleft}{0pt}
\begin{longtable}{>{\raggedright\arraybackslash}p{0.05\textwidth - 2\tabcolsep}>{\raggedright\arraybackslash}p{0.95\textwidth-\tabcolsep}@{}}

\caption{Diagonal irreducible subgroups of $\bar{A}_1^3 A_1 B_2 = E_7(\#\eseven{45})$.  \label{A13A1B2E7diags}} \\
\hline
ID & Embedding \\
\hline

\eseven{89} & $(a^{[r]},a^{[s]},b,c,10)$ $(rs=0)$  \\

\eseven{90} & $(a^{[r]},b,c,a^{[s]},10)$ $(rs=0)$ \\
\endfirsthead
\eseven{91} &  $(a,b,b^{[r]},c,10)$ (if $p=2$ then $r \neq 0$) \\

\eseven{92} &  $(a,b^{[r]},c,b^{[s]},10)$ $(rs=0)$  \\

\eseven{93} & $(a^{[r]},a^{[s]},a^{[t]},b,10)$  $(rs=0; s \leq t$; if $p=2$ then $s < t$)  \\

\eseven{94} & $(a,b^{[r]},b^{[s]},b^{[t]},10)$ $(rt=0; r \leq s$; if $p=2$ then $r < s$; if $r=s$ then $r < t$)  \\

\eseven{95} & $(a^{[r]},b,a^{[s]},a^{[t]},10)$ $(rst=0)$  \\

\eseven{96} & $(a^{[r]},a^{[s]},b^{[t]},b^{[u]},10)$  $(rs=tu=0)$  \\

\eseven{97} & $(a^{[r]},b,b^{[t]},a^{[s]},10)$ $(rs=0$; if $p=2$ then $t \neq 0$)   \\

\eseven{98} & $(a^{[r]},a^{[s]},a^{[t]},a^{[u]},10)$ $(rsu=0; s \leq t$; if $p=2$ then $s < t$; if $s=t$ then $s < u$)   \\
\hline
\end{longtable}

\begin{longtable}{>{\raggedright\arraybackslash}p{0.06\textwidth - 2\tabcolsep}>{\raggedright\arraybackslash}p{0.94\textwidth-\tabcolsep}@{}}

\caption{Diagonal irreducible subgroups of $\bar{A}_1^3 A_1 A_1 = E_7(\#\eseven{87})$ $(p \geq 5)$.  \label{A13A1A1E7diags}} \\
\hline
ID & Embedding \\
\hline

\eseven{99} & $(a^{[r]},a^{[s]},b,c,d)$ $(rs = 0)$ \\  

\eseven{100} & $(a^{[r]},b,c,a^{[s]},d)$ $(rs = 0)$ \\
\endfirsthead
\eseven{101} & $(a^{[r]},b,c,d,a^{[s]})$ $(rs = 0)$ \\

\eseven{102} & $(a,b,b^{[r]},c,d)$  \\

\eseven{103} & $(a,b^{[r]},c,b^{[s]},d)$ $(rs = 0)$  \\

\eseven{104} & $(a,b^{[r]},c,d,b^{[s]})$ $(rs = 0)$ \\

\eseven{105} & $(a,b,c,d^{[r]},d^{[s]})$ $(rs = 0)$  \\

\eseven{106} & $(a^{[r]},a^{[s]},a^{[t]},b,c)$  $(rst=0; s \leq t)$  \\  

\eseven{107} & $(a^{[r]},a^{[s]},b,a^{[t]},c)$ $(rst=0)$   \\  

\eseven{108} & $(a^{[r]},a^{[s]},b,c,a^{[t]})$ $(rst=0)$   \\
\hline
\end{longtable}

\begin{longtable}{>{\raggedright\arraybackslash}p{0.06\textwidth - 2\tabcolsep}>{\raggedright\arraybackslash}p{0.94\textwidth-\tabcolsep}@{}}

\caption{Diagonal irreducible subgroups of $\bar{A}_1^3 A_1^3 = E_7(\#\eseven{88})$ $(p=2)$.  \label{A13A13E7diags}} \\
\hline
ID & Embedding \\

\hline

\eseven{133} &  $(a^{[r]},a^{[s]},b,c,d,e)$ $(rs=0)$ \\

\eseven{134} &   $(a^{[r]},b,c,a^{[s]},d,e)$ $(rs=0)$ \\

\eseven{135} &   $(a,b,b^{[r]},c,d,e)$ $(r \neq 0)$ \\

\eseven{136} &   $(a,b^{[r]},c,b^{[s]},d,e)$ $(rs=0)$ \\ 
\endfirsthead
\eseven{137} &   $(a,b,c,d,d^{[r]},e)$ $(r \neq 0)$ \\

\eseven{138} &   $(a^{[r]},a^{[s]},a^{[t]},b,c,d)$ $(rs=0; s<t)$ \\

\eseven{139} &   $(a^{[r]},a^{[s]},b,a^{[t]},c,d)$ $(rst=0)$ \\

\eseven{140} & $(a^{[r]},b,c,a^{[s]},a^{[t]},d)$ $(rs=0; s<t)$ \\ 

\eseven{141} &   $(a,b^{[r]},b^{[s]},b^{[t]},c,d)$ $(rt=0; r < s)$ \\

\eseven{142} &  $(a,b,c^{[r]},c^{[s]},c^{[t]},d)$ $(rs = 0; s < t)$ \\

\eseven{143} &   $(a,b,c,d,d^{[r]},d^{[s]})$ $(0 < r < s)$ \\

\eseven{144} &   $(a^{[r]},a^{[s]},b^{[t]},b^{[u]},c,d)$ $(rs = tu =0)$ \\
\hline
\end{longtable}

\setlength{\LTleft}{0.17\textwidth}
\begin{longtable}{>{\raggedright\arraybackslash}p{0.06\textwidth - 2\tabcolsep}>{\raggedright\arraybackslash}p{0.6\textwidth-\tabcolsep}@{}}

\caption{Diagonal irreducible subgroups of $\bar{A}_1 A_1^4 = E_7(\#\eseven{191})$ $(p \neq 2)$.  \label{A1A14E7diags}} \\
\hline
ID & Embedding \\

\hline

\eseven{195} &  $(a^{[r]},a^{[s]},b,c,d)$ $(rs=0)$ \\

\eseven{196} &   $(a,b,b^{[r]},c,d)$ $(r \neq 0)$ \\

\eseven{197} &  $(a^{[r]},a^{[s]},a^{[t]},b,c)$ $(rs=0; s<t)$ \\ 

\eseven{198} &  $(a,b,b^{[r]},b^{[s]},c)$ $(0< r < s)$ \\ 

\eseven{199} & $(a^{[r]},a^{[s]},b,b^{[t]},c)$ $(rs=0 \neq t)$ \\

\eseven{200} & $(a,b,b^{[r]},c,c^{[s]})$ $(0 < r \leq s)$ \\

\eseven{201} & $(a^{[r]},a^{[s]},a^{[t]},a^{[u]},b)$ $(rs=0; s<t<u)$ \\

\endfirsthead 

\eseven{202} & $(a,b,b^{[r]},b^{[s]},b^{[t]})$ $(0<r<s<t)$ \\

\eseven{203} &  $(a^{[r]},a^{[s]},a^{[t]},b,b^{[u]})$ $(rs=0 \neq u; s<t)$ \\

\eseven{204} &  $(a^{[r]},a^{[s]},b,b^{[t]},b^{[u]})$ $(rs=0; 0<t<u)$ \\

\eseven{11} & $(a^{[r]},a^{[s]},a^{[t]},a^{[u]},a^{[v]})$ $(rs=0; s < t < u < v)$ \\
\hline
\end{longtable}

\setlength{\LTleft}{0pt}
\begin{longtable}{>{\raggedright\arraybackslash}p{0.09\textwidth - 2\tabcolsep}>{\raggedright\arraybackslash}p{0.39\textwidth - 2\tabcolsep}>{\raggedright\arraybackslash}p{0.06\textwidth - 2\tabcolsep}>{\raggedright\arraybackslash}p{0.46\textwidth-\tabcolsep}@{}}

\caption{Diagonal irreducible subgroups of $\bar{A}_1 A_1 A_1^2 = E_7(\#\eseven{186})$ $(p \neq 2)$.  \label{A1A1A12E7diags}} \\
\hline
ID & Embedding & ID & Embedding \\
\hline

\eseven{205} & $(a^{[r]},a^{[s]},b,c)$ $(rs=0)$ & \eseven{210} & $(a^{[r]},b,a^{[s]},a^{[t]})$ $(rst=0; s<t)$ \\

\eseven{206} & $(a^{[r]},b,a^{[s]},c)$ $(rs=0)$ & \eseven{211} & $(a,b^{[r]},b^{[s]},b^{[t]})$ $(rs=0; s <t)$ \\
\endfirsthead
\eseven{207} & $(a,b^{[r]},b^{[s]},c)$ $(rs=0)$ & \eseven{212} & $(a^{[r]},a^{[s]},b,b^{[t]})$ $(rs=0 \neq t)$ \\

\eseven{208} & $(a,b,c,c^{[r]})$ $(r \neq 0)$ & \eseven{213} & $(a^{[r]},b^{[t]},a^{[s]},b^{[u]})$ $(rs=tu=0)$ \\

$\mathit{\eseven{115}^{\{\underline{0}\}}}$ &  $(a,b,c,c)$ $(p \geq 5)$ & \eseven{10} & $(a^{[r]},a^{[s]},a^{[t]},a^{[u]})$ $(rst=0; t < u)$ \\

\eseven{209} & $(a^{[r]},a^{[s]},a^{[t]},b)$ $(rst=0)$ \\
\hline
\end{longtable}

\setlength{\LTleft}{0.12\textwidth}
\begin{longtable}{>{\raggedright\arraybackslash}p{0.06\textwidth - 2\tabcolsep}>{\raggedright\arraybackslash}p{0.7\textwidth-\tabcolsep}@{}}

\caption{Diagonal irreducible subgroups of $\bar{A}_1 A_1^5 = E_7(\#\eseven{223})$ $(p=2)$.  \label{A1A15E7diags}} \\
\hline
ID & Embedding \\

\hline

\eseven{230} &  $(a^{[r]},a^{[s]},b,c,d,e)$ $(rs=0)$ \\

\eseven{238} &   $(a^{[r]},a^{[s]},a^{[t]},b,b^{[u]},c)$ $(rs= 0 \neq u; s < t)$ \\

\eseven{231} &   $(a,b,b^{[r]},c,d,e)$ $(r \neq 0)$ \\

\eseven{239} &  $(a^{[t]},b,b^{[r]},b^{[s]},a^{[u]},c)$ $(0 < r < s; tu=0)$ \\

\endfirsthead
\eseven{232} &   $(a^{[r]},a^{[s]},a^{[t]},b,c,d)$ $(rs=0; s < t)$ \\ 

\eseven{240} &   $(a,b,b^{[r]},b^{[s]},c,c^{[t]})$ $(0 < r < s; t \neq 0)$ \\

\eseven{233} &   $(a,b,b^{[r]},b^{[s]},c,d)$ $(0 < r < s)$ \\ 

\eseven{241} &  $(a^{[r]},a^{[s]},b,b^{[t]},c,c^{[u]})$ $(rs=0; 0 < t \leq u)$ \\

\eseven{234} &  $(a^{[r]},a^{[s]},b,b^{[t]},c,d)$ $(rs=0 \neq t)$ \\ 

\eseven{242} &   $(a^{[r]},a^{[s]},a^{[t]},a^{[u]},a^{[v]},b)$ $(rs=0; s < t < u < v)$ \\

\eseven{235} &   $(a,b,b^{[r]},c,c^{[s]},d)$ $(0 < r \leq s)$ \\ 

\eseven{243} &   $(a,b,b^{[r]},b^{[s]},b^{[t]},b^{[u]})$ $(0 < r < s < t < u)$ \\

\eseven{236} & $(a^{[r]},a^{[s]},a^{[t]},a^{[u]},b,c)$ $(rs=0; s < t < u)$ \\

\eseven{237} & $(a,b,b^{[r]},b^{[s]},b^{[t]},c)$ $(0 < r < s < t)$ \\ 

\eseven{14} & $(a^{[r]},a^{[s]},a^{[t]},a^{[u]},a^{[v]},a^{[w]})$ $(rs=0$; $s < t < u < v < w)$ \\

\hline
\end{longtable}

\setlength{\LTleft}{0.17\textwidth}
\begin{longtable}{>{\raggedright\arraybackslash}p{0.1\textwidth - 2\tabcolsep}>{\raggedright\arraybackslash}p{0.56\textwidth-\tabcolsep}@{}}

\caption{Diagonal irreducible subgroups of $\bar{A}_1 A_1 A_1 A_1 = E_7(\#\eseven{268})$ ${(p \neq 2)}$.  \label{A1A1A1A1E7diags}} \\
\hline
ID & Embedding \\

\hline

\eseven{273} &  $(a^{[r]},a^{[s]},b,c)$  $(rs=0)$ \\

\eseven{274} &   $(a^{[r]},b,a^{[s]},c)$  $(rs=0)$ \\

\eseven{275} &  $(a^{[r]},b,c,a^{[s]})$  $(rs=0)$ \\ 
\endfirsthead
\eseven{276} &  $(a,b^{[r]},b^{[s]},c)$  $(rs=0; r \neq s)$ \\ 

$\mathit{\eseven{119}^{\{\underline{0}\}}}$ &  $(a,b,b,c)$ $(p \geq 5)$  \\

\eseven{277} &   $(a,b^{[r]},c,b^{[s]})$  $(rs=0; r \neq s)$ \\

$\mathit{\eseven{207}^{\{\underline{0}\}}}$ &  $(a,b,c,b)$  \\

\eseven{278} &   $(a,b,c^{[r]},c^{[s]})$  $(rs=0; r \neq s)$ \\ 

$\mathit{\eseven{117}^{\{\underline{0}\}}}$ & $(a,b,c,c)$ $(p \geq 5)$ \\

\eseven{279} &   $(a^{[r]},a^{[s]},a^{[t]},b)$  $(rst=0; s \neq t)$ \\

\eseven{280} &  $(a^{[r]},a^{[s]},b,a^{[t]})$  $(rst=0; s \neq t)$ \\

\eseven{281} &  $(a^{[r]},b,a^{[s]},a^{[t]})$  $(rst=0; s \neq t)$ \\

\eseven{282} &  $(a,b^{[r]},b^{[s]},b^{[t]})$  $(rst=0$; $r,s,t$ distinct) \\

\eseven{283} &  $(a^{[r]},a^{[s]},b^{[t]},b^{[u]})$  $(rs=tu=0; t \neq u)$ \\

\eseven{284} & $(a^{[r]},b^{[t]},a^{[s]},b^{[u]})$  $(rs=tu=0; t \neq u)$ \\

\eseven{285} &  $(a^{[r]},b^{[t]},b^{[u]},a^{[s]})$  $(rs=tu=0; t \neq u)$ \\

\eseven{3} & $(a^{[r]},a^{[s]},a^{[t]},a^{[u]})$ $(rstu=0; s,t,u$ distinct) \\

\hline
\end{longtable}

\subsection{Irreducible diagonal subgroups contained in $E_8$}  \label{sec:irredE8diags}
\leavevmode

\setlength{\LTleft}{0pt}
\begin{longtable}{>{\raggedright\arraybackslash}p{0.06\textwidth - 2\tabcolsep}>{\raggedright\arraybackslash}p{0.94\textwidth-\tabcolsep}@{}}

\caption{Irreducible diagonal subgroups of $\bar{A}_1^8 = E_8(\#\eeight{124})$.  \label{barA18diags}} \\
\hline
ID & Embedding \\

\hline

\eeight{132} & $(a,a^{[r]},b,c,d,e,f,g)$ \\ 

\eeight{133} & $(a,a^{[r]},a^{[s]},b,c,d,e,f)$ $(0 \leq r \leq s)$ \\
\endfirsthead

\eeight{134} & $(a,a^{[r]},b,b^{[s]},c,d,e,f)$ $(r \leq s$; $s \neq 0$) \\

\eeight{135} & $(a,a^{[r]},b,c,b^{[s]},d,e,f)$ $(r \leq s)$ \\

\eeight{136} & $(a,a^{[r]},a^{[s]},a^{[t]},b,c,d,e)$ $(0 \leq r \leq s \leq t;$ if two of $0,r,s,t$ are equal then the other two are not equal) \\ 

\eeight{137} & $(a,a^{[r]},a^{[s]},b,a^{[t]},c,d,e)$ $(0 \leq r \leq s \leq t)$ \\  

\eeight{138} & $(a,a^{[r]},a^{[s]},b^{[t]},b^{[u]},c,d,e)$ $(0 \leq r \leq s; tu=0)$ \\

\eeight{139} & $(a^{[r]},a^{[s]},a^{[t]},b,c,c^{[u]},d,e)$ $(rt=0; r \leq s$; if $u=0$ then $r < s$) \\

\eeight{140} & $(a,a^{[r]},b,b^{[s]},c,c^{[t]},d,e)$ $(r \leq s \leq t$; if $r=0$ then $s \neq 0$) \\

\eeight{141} & $(a,a^{[r]},b,c^{[t]},b^{[s]},c^{[u]},d,e)$ ($tu=0$; $s \leq t + u$; if $s=0$ then $t < u$) \\
\hline
\end{longtable}

\setlength{\LTleft}{0.09\textwidth}
\begin{longtable}{>{\raggedright\arraybackslash}p{0.37\textwidth-2\tabcolsep}>{\raggedright\arraybackslash}p{0.45\textwidth-\tabcolsep}@{}}
\caption{Conditions on field twists for $E_8(\#\eeight{154})$. \label{condE8t1a5b2}} \\

\hline \noalign{\smallskip}

Equalities among $0, r, \ldots, v$ & Further requirements on $0, r, \ldots, v$ \\

\hline 

none & $s < \text{min}\{t,u\}$; $u < v$ \\

$r=s$ & $u < v$ \\
\endfirsthead
$s=t$ & $u < v$ \\

$s=u$ & $t < v$ \\

$r=s=t$ & $u < v$ \\

$r=s=u$ & $t < v$ \\

$s=t=u$ & none \\

$r=s=t=u$ & none \\

$r=0$ & $s < \text{min}\{t,u\}$; $u < v$ \\

$r=0$; $s=u$ & $t < v$ \\

$s=0$ & $r < t$; $u < v$ \\

$s=0$; $r=u$ & none \\

$s=0$; $u=v$ & $r < t$ \\

$s=0$; $r=u=v$ & none \\

$r=s=0$ & $u < v$ \\

$r=s=0$; $t=u$ & none \\

$s=t=0$ & $u < v$ \\

$s=t=0$; $r=u$ & none \\

$s=u=0$ & $r < t < v$ \\

$s=t=u=0$ & $r< v$ \\
\hline
\end{longtable}

\setlength{\LTleft}{0.1\textwidth}
\begin{longtable}{>{\raggedright\arraybackslash}p{0.35\textwidth-2\tabcolsep}>{\raggedright\arraybackslash}p{0.45\textwidth-\tabcolsep}@{}}
\caption{Conditions on field twists for $E_8(\#\eeight{156})$. \label{condE8t1a5b4}} \\

\hline \noalign{\smallskip}

Equalities among $r,s,t,u$ & Further requirements on $0,r, \ldots, w$ \\

\hline 

none & $r=0$;  $s < t < u$ \\

$r=s$ & $t < u$; $v < w$ \\

$r=s=t$ & $v < w < x$ \\
\hline
\endfirsthead
\end{longtable}

\setlength{\LTleft}{0.12\textwidth}
\begin{longtable}{>{\raggedright\arraybackslash}p{0.33\textwidth-2\tabcolsep}>{\raggedright\arraybackslash}p{0.43\textwidth-\tabcolsep}@{}}
\caption{Conditions on field twists for $E_8(\#\eeight{157})$. \label{condE8t1a5b5}} \\

\hline \noalign{\smallskip}

Equalities among $v,w,x$ & Further requirements on $r, \ldots, x$ \\

\hline 

none & $\eeight{0} < v < w$  \\

$v=w$ & $t < u$ \\

$v=w=x$ & $s < t < u$ \\
\hline
\end{longtable}

\setlength{\LTleft}{0.09\textwidth}
\begin{longtable}{>{\raggedright\arraybackslash}p{0.37\textwidth-2\tabcolsep}>{\raggedright\arraybackslash}p{0.45\textwidth-\tabcolsep}@{}}
\caption{Conditions on field twists for $E_8(\#\eeight{163})$.\label{condE8221}} \\

\hline \noalign{\smallskip}

Equalities among $0, r, \ldots, x$ & Further requirements on $0, r, \ldots, x$ \\

\hline 

none & $s < \text{min}\{t,u\}$; $u < v$ \\

$r=s$ & $u < v$ \\

$s=t$ & $u < v$; $w < x$ \\

$s=u$ & $t < v$ \\

$r=s=t$ & $u < v$ \\

$r=s=u$ & $t < v$ \\

$s=t=u$ & none \\

$r=s=t=u$ & none \\

$r=0$ & $s < \text{min}\{t,u\}$; $u < v$; $w < x$ \\

$r=0$; $s=u$ & $t < v$; $w < x$ \\

$s=0$ & $r < t$; $u < v$ \\

$s=0$; $r=u$ & none \\

$s=0$; $u=v$ & $r < t$; $w < x$ \\

$s=0$; $r=u=v$ & $w < x$ \\

$r=s=0$ & $u < v$; $w < x$ \\

$r=s=0$; $t=u$ & $w < x$ \\

$s=t=0$ & $u < v$; $w < x$ \\

$s=t=0$; $r=u$ & $w < x$ \\

$s=u=0$ & $r < t < v$ \\

$s=t=u=0$ & $r< v$; $w < x$ \\
\hline
\end{longtable}

\setlength{\LTleft}{0.1\textwidth}
\begin{longtable}{>{\raggedright\arraybackslash}p{0.37\textwidth-2\tabcolsep}>{\raggedright\arraybackslash}p{0.43\textwidth-\tabcolsep}@{}}
\caption{Conditions on field twists for $E_8(\#\eeight{164})$. \label{condE8222}} \\

\hline \noalign{\smallskip}

Equalities among $r,s,t,u,v$ & Further requirements on $r, \ldots, y$ \\

\hline 

none & $r < s < t < u$   \\
\endfirsthead
$r=s$ & $t < u$; $x < y$ \\

$r=v$ & $s < t < u$ \\

$r=s=t$ & $w < x < y$ \\

$r=s=v$ & $t < u$; $x < y$ \\

$r=s=t=v$ & $w < x < y$ \\
\hline
\end{longtable}

\setlength{\LTleft}{0pt}
\begin{longtable}{>{\raggedright\arraybackslash}p{0.3\textwidth-2\tabcolsep}>{\raggedright\arraybackslash}p{0.7\textwidth-\tabcolsep}@{}}
\caption{Conditions on field twists for $E_8(\#\eeight{165})$. \label{condE8223}} \\

\hline \noalign{\smallskip}

Equalities among $r, s, t, u$ and equalities among $v, w, x, y$ & Further requirements on $r, \ldots, y$ \\

\hline 

none & $r=v=0$;  $s < t < u$; let $i < j < k$ represent the order of $w,x,y$: $s \leq i$; if $s=i$ then $t \leq j$; if $s=i$ and $t=j$ then $u \leq k$ \\

$r=s$ & $t < u$; $v=0$; $x<y$ \\

\endfirsthead

$r=s=t$ & $v=0$; $w < x < y$ \\

$r=s$; $v=x$ &  $t < u$; $w < y$; $r \leq v$; if $r=v$ then either $t < w$ or $t=w$ and $u \leq y$  \\
\hline
\end{longtable}

\setlength{\LTleft}{0pt}
\begin{longtable}{>{\raggedright\arraybackslash}p{0.35\textwidth-2\tabcolsep}>{\raggedright\arraybackslash}p{0.65\textwidth-\tabcolsep}@{}}
\caption{Conditions on field twists for $E_8(\#\eeight{166})$. \label{condE8224}} \\

\hline \noalign{\smallskip}

Equalities among $r,s,t,u$ and equalities among $v, w, x, y$ & Further requirements on $r, \ldots, y$ \\

\hline 

none & $r=0$;  $s < t < u$; let $0=i < j < k < l$ represent the order of $v,w,x,y$: $s \leq j$; if $s=j$ then $t \leq k$; if $s=j$ and $t=k$ then $u \leq l$  \\

$r=s$ & $t<u$; $x < y$ \\

\endfirsthead

$r=s=t$ & $w<x<y$ \\

$r=s$; $v=w$ & $t < u$; $x < y$; $r \leq v$; if $r=v$ then $t \leq x$; if $r=v$ and $t=x$ then $u \leq y$ \\

$r=s$; $v=x$ & $r \leq v$; if $r=v$ then $t \leq w$; if $r=v$ and $t=w$ then $u \leq y$ \\

$r=s$; $v=w=x$ & $t < u$ \\

$r=s=t$; $v=w$ & $x < y$ \\
\hline
\end{longtable}

\setlength{\LTleft}{0.06\textwidth}
\begin{longtable}{>{\raggedright\arraybackslash}p{0.38\textwidth-2\tabcolsep}>{\raggedright\arraybackslash}p{0.5\textwidth-\tabcolsep}@{}}
\caption{Conditions on field twists for $E_8(\#\eeight{26})$. \label{E8cond26}} \\

\hline \noalign{\smallskip}

Equalities among $0, r, \ldots, w$ & Further requirements on $0, r, \ldots, w$ \\

\hline \noalign{\smallskip}

none & $r < s$; $s < \text{min}\{t, u\}$; $u < \text{min}\{v, w, x\}$ \\

$r=s$ & $u<v, w, x$; $w<x$ \\
\endfirsthead
$r=s=t$ & $u < v < w < x$ \\

$r=s=u$ & $t < v < w$ \\

$r=s=t=u$ & $v < w < x$ \\

$r=s$; $t=u$ & $v < w$ \\

$r=s$; $u=v$ & $w < x$ \\

$r=s$; $t=u=v$ & $w < x$ \\

$r=s$; $t=u$; $v = w$ & none \\

$r = 0$ & $s < \text{min}\{t,u\}$; $u < \text{min}\{v, w\}$; $w < x$ \\

$r = 0$; $s=u$ & $t<v$; $w<x$ \\

$r = 0$; $s=u=w$ & $t < v < x$ \\

$r = 0$; $s=u$; $t=w$ & none \\

$r = 0$; $s = u$; $t = w$; $v = x$ & none \\

$r = s = 0$ & $u < v < w< x$ \\

$r=s=0$; $t=u$ & $v < w < x$ \\

$r = s = u = 0$ & $t < v < w < x$ \\
\hline
\end{longtable}

\setlength{\LTleft}{0.17\textwidth}
\begin{longtable}{>{\raggedright\arraybackslash}p{0.06\textwidth - 2\tabcolsep}>{\raggedright\arraybackslash}p{0.6\textwidth-\tabcolsep}@{}}

\caption{Irreducible diagonal subgroups of $\bar{A}_1^4 A_1 = E_8(\#\eeight{174})$ $(p \geq 7)$. \label{barA14A1p7diags}} \\
\hline
ID & Embedding \\
\hline
\eeight{179} & $(a,a^{[r]},b,c,d)$ $(r \neq 0)$  \\

\eeight{180} & $(a^{[r]},b,c,d,a^{[s]})$ $(rs = 0)$ \\

\endfirsthead
\eeight{181} & $(a,a^{[r]},a^{[s]},b,c)$ $(0 < r < s)$  \\

\eeight{182} & $(a^{[r]},a^{[s]},b,c,a^{[t]})$ $(rt =0; r < s)$  \\

\eeight{183} & $(a,a^{[r]},b,b^{[s]},c)$ $(0 < r \leq s)$ \\

\eeight{184} & $(a,a^{[r]},b^{[s]},c,b^{[t]})$ $(r \neq 0 = st)$  \\

\eeight{185} & $(a,a^{[r]},a^{[s]},a^{[t]},b)$ $(0 < r < s < t)$   \\

\eeight{186} & $(a^{[r]},a^{[s]},a^{[t]},b,a^{[u]})$ $(ru=0; r < s < t)$   \\

\eeight{187} & $(a,a^{[r]},b^{[s]},b^{[t]},b^{[u]})$ $(r \neq 0 = su; s < t)$   \\

\eeight{188} & $(a,a^{[r]},a^{[s]},b^{[t]},b^{[u]})$ $(0 < r < s; tu=0)$   \\

\eeight{18} & $(a^{[r]},a^{[s]},a^{[t]},a^{[u]},a^{[v]})$ $(rv=0; r < s < t < u)$ \\
\hline
\end{longtable}

\setlength{\LTleft}{0pt}
\begin{longtable}{>{\raggedright\arraybackslash}p{0.06\textwidth - 2\tabcolsep}>{\raggedright\arraybackslash}p{0.94\textwidth-\tabcolsep}@{}}

\caption{Irreducible diagonal subgroups of $\bar{A}_1^4 A_1 B_2 = E_8(\#\eeight{126})$.  \label{barA14A1B2diags}} \\
\hline
ID & Embedding \\
\hline

\eeight{191} & $(a,a^{[r]},b,c,d,10)$ (if $p=2$ then $r \neq 0)$ \\
\eeight{192} & $(a,b,a^{[r]},c,d,10)$ \\
\eeight{193} & $(a^{[r]},b,c,d,a^{[s]},10)$ $(rs=0)$  \\ \endfirsthead

\eeight{194} & $(a^{[r]},a^{[s]},a^{[t]},b,c,10)$ $(rt=0; r \leq s$; if $p=2$ then $r < s)$ \\
\eeight{195} & $(a^{[r]},a^{[s]},b,c,a^{[t]},10)$ $(rt=0; r \leq s$; if $r=s$ then $r < t$; if $p=2$ then $r < s)$ \\
\eeight{196} & $(a^{[r]},b,a^{[s]},c,a^{[t]},10)$ $(rt=0; r \leq s)$  \\

\eeight{197} & $(a,a^{[r]},b,b^{[s]},c,10)$ ($r \leq s$; if $p=2$ then $rs \neq 0)$ \\
\eeight{198} & $(a,a^{[r]},b^{[s]},c,b^{[t]},10)$ $(st=0$; if $p=2$ then $r \neq 0)$ \\
\eeight{199} & $(a,b^{[s]},a^{[r]},b^{[t]},c,10)$ $(st=0$; $r \leq s+t$; if $r=0$ then $t \neq 0$) \\
\eeight{200} & $(a^{[r]},b^{[t]},a^{[s]},c,b^{[u]},10)$ $(rs=tu=0)$  \\

\eeight{201} & $(a,a^{[r]},a^{[s]},a^{[t]},b,10)$ $(s \leq t$; if $r=0$ then $s < t$; if $s=0$ then $r < t$; if $p=2$ then $r \neq 0$ and $s < t)$ \\
\hline
\end{longtable}

\begin{longtable}{>{\raggedright\arraybackslash}p{0.06\textwidth - 2\tabcolsep}>{\raggedright\arraybackslash}p{0.94\textwidth-\tabcolsep}@{}}

\caption{Irreducible diagonal subgroups of $\bar{A}_1^4 A_1 A_1 = E_8(\#\eeight{189})$ $(p \geq 5)$.  \label{barA14A1A1p5diags}} \\
\hline
ID & Embedding \\
\hline

\eeight{207} & $(a,a^{[r]},b,c,d,e)$ \\
\eeight{208} & $(a^{[r]},b,a^{[s]},c,d,e)$ $(rs=0)$ \\
\eeight{209} & $(a^{[r]},b,c,d,a^{[s]},e)$ $(rs=0)$ \\ \endfirsthead
\eeight{210} & $(a^{[r]},b,c,d,e,a^{[s]})$ $(rs=0)$ \\
\eeight{211} & $(a,b,c,d,e^{[r]},e^{[s]})$ $(rs=0)$ \\

\eeight{212} & $(a^{[r]},a^{[s]},a^{[t]},b,c,d)$ $(rt=0; r \leq s)$  \\
\eeight{213} & $(a^{[r]},a^{[s]},b,c,a^{[t]},d)$ $(rt=0; r \leq s$; if $r=s$ then $r < t$)  \\
\eeight{214} & $(a^{[r]},a^{[s]},b,c,d,a^{[t]})$ $(rt=0; r \leq s)$  \\
\eeight{215} & $(a^{[r]},b,a^{[s]},c,a^{[t]},d)$ $(rt=0; r \leq s)$  \\
\eeight{216} & $(a^{[r]},b,a^{[s]},c,d,a^{[t]})$ $(rt=0; r \leq s)$  \\
\eeight{217} & $(a^{[r]},b,c,d,a^{[s]},a^{[t]})$ $(rst=0)$  \\
\hline
\end{longtable}

\setlength{\LTleft}{0.05\textwidth}
\begin{longtable}{>{\raggedright\arraybackslash}p{0.45\textwidth-2\tabcolsep}>{\raggedright\arraybackslash}p{0.45\textwidth-\tabcolsep}@{}}
\caption{Conditions on field twists for $A_1 = E_8(\#\eeight{23})$.} \label{condsE8A123} \\

\hline \noalign{\smallskip}

Equalities among $r,s,t,u$ and equalities among $v, w$ & Further requirements on $r, \ldots, w$ \\

\hline \noalign{\smallskip}

none & $r < \text{min}\{s,t\}$; $t < u$ \\

$r=s$ & $r < v$; $t < u$ \\
\endfirsthead
$r=t$ & $s < u$ \\

$v=w$ & $r < s < t < u$ \\

$r=s$; $v=w$ & $r < v$; $t < u$ \\

$r=s=t$ & none \\

$r=s=t$; $v=w$ & $r \neq v$ \\
\hline
\end{longtable}

\begin{longtable}{>{\raggedright\arraybackslash}p{0.06\textwidth - 2\tabcolsep}>{\raggedright\arraybackslash}p{0.92\textwidth-\tabcolsep}@{}}

\caption{Irreducible diagonal subgroups of $\bar{A}_1^4 A_1^3 = E_8(\#\eeight{190})$ $(p=2)$.  \label{barA14A13diags}} \\
\hline
ID & Embedding \\

\hline

\eeight{261} & $(a,a^{[r]},b,c,d,e,f)$ $(r \neq 0)$ \\

\eeight{262} & $(a,b,a^{[r]},c,d,e,f)$ \\
\endfirsthead
\eeight{263} & $(a^{[r]},b,c,d,a^{[s]},e,f)$ $(rs = 0)$ \\

\eeight{264}& $(a,b,c,d,e,e^{[r]},f)$ $(r \neq 0)$ \\

\eeight{265} & $(a^{[r]},a^{[s]},a^{[t]},b,c,d,e)$ $(rt = 0; r < s)$ \\

\eeight{266} & $(a^{[r]},a^{[s]},b,c,a^{[t]},d,e)$ $(rt = 0; r < s)$ \\

\eeight{267} & $(a^{[r]},b,a^{[s]},c,a^{[t]},d,e)$ $(rt = 0; r \leq s)$ \\

\eeight{268} & $(a^{[r]},b,c,d,a^{[s]},a^{[t]},e)$ $(rs = 0; s < t)$ \\

\eeight{269} & $(a,b,c,d,e,e^{[r]},e^{[s]})$ $(0 < r < s)$ \\

\eeight{270} & $(a,a^{[r]},b,b^{[s]},c,d,e)$ $(0 < r \leq s)$ \\
\hline
\end{longtable}

\begin{longtable}{>{\raggedright\arraybackslash}p{0.06\textwidth - 2\tabcolsep}>{\raggedright\arraybackslash}p{0.94\textwidth-\tabcolsep}@{}}

\caption{Irreducible diagonal subgroups of $\bar{A}_1^2 A_1^4 = E_8(\#\eeight{365})$ $(p \neq 2)$.  \label{barA12A14diags}} \\
\hline
ID & Embedding \\
\hline

\eeight{374} & $(a,a^{[r]},b,c,d,e)$ \\

\eeight{375} & $(a^{[r]},b,a^{[s]},c,d,e)$ $(rs=0)$ \\
\endfirsthead
\eeight{376} & $(a,b,c,c^{[r]},d,e)$ $(r \neq 0)$ \\

\eeight{377} & $(a^{[r]},a^{[s]},a^{[t]},b,c,d)$ $(rt=0$; $r \leq s$; if $r=s$ then $r < t$)  \\

\eeight{378} & $(a^{[r]},b,a^{[s]},a^{[t]},c,d)$ $(rs=0$; $s < t$)  \\

\eeight{379} & $(a,b,c,c^{[r]},c^{[s]},d)$ $(0 < r < s$)  \\

\eeight{380} & $(a,a^{[r]},b,b^{[s]},c,d)$ $(s \neq 0)$ \\

\eeight{381} & $(a^{[r]},b^{[t]},a^{[s]},b^{[u]},c,d)$ $(rs=tu=0$; $r \leq t$; if $r=t$ then $s \leq u$) \\

\eeight{382} & $(a^{[r]},b,a^{[s]},c,c^{[t]},d)$ $(rs=0 \neq t)$ \\

\eeight{383} & $(a,b,c,c^{[r]},d,d^{[s]})$ $(rs \neq 0; r \leq s)$ \\
\hline
\end{longtable}

\setlength{\LTleft}{0.05\textwidth}
\begin{longtable}{>{\raggedright\arraybackslash}p{0.09\textwidth - 2\tabcolsep}>{\raggedright\arraybackslash}p{0.81\textwidth-\tabcolsep}@{}}

\caption{Irreducible diagonal subgroups of $\bar{A}_1^2 A_1 A_1^2 = E_8(\#\eeight{401})$ $(p \neq 2)$.  \label{barA12A1A12diags}} \\
\hline
ID & Embedding \\

\hline

\eeight{408} & $(a,a^{[r]},b,c,d)$ \\ 

\eeight{409} & $(a^{[r]},b,a^{[s]},c,d)$ $(rs=0)$ \\
\endfirsthead
\eeight{410} & $(a^{[r]},b,c,a^{[s]},d)$ $(rs=0)$ \\

\eeight{411} & $(a,b,c^{[r]},c^{[s]},d)$ $(rs=0)$ \\

\eeight{412} & $(a,b,c,d,d^{[r]})$ $(r \neq 0)$ \\

$\mathit{\eeight{214}^{\{\underline{0}\}}}$ & $(a,b,c,d,d)$ \\

\eeight{413} & $(a^{[r]},a^{[s]},a^{[t]},b,c)$ $(rt=0; r \leq s$; if $r=s$ then $s < t$) \\ 

\eeight{414} & $(a^{[r]},a^{[s]},b,a^{[t]},c)$ $(rt=0; r \leq s)$ \\ 

\eeight{415} & $(a^{[r]},b,a^{[s]},a^{[t]},c)$ $(rst=0)$ \\

\eeight{416} & $(a^{[r]},b,c,a^{[s]},a^{[t]})$ $(rs=0; s < t)$ \\
\hline
\end{longtable}

\setlength{\LTleft}{0pt}
\begin{longtable}{>{\raggedright\arraybackslash}p{0.06\textwidth - 2\tabcolsep}>{\raggedright\arraybackslash}p{0.42\textwidth - 2\tabcolsep}>{\raggedright\arraybackslash}p{0.06\textwidth - 2\tabcolsep}>{\raggedright\arraybackslash}p{0.46\textwidth-\tabcolsep}@{}}

\caption{Irreducible diagonal subgroups of $\bar{A}_1^2 A_1 A_1 = E_8(\#\eeight{402})$ $(p \geq 11)$.  \label{barA12A1A1diags}} \\
\hline
ID & Embedding & ID & Embedding \\

\hline

\eeight{432} & $(a,a^{[r]},b,c)$ & \eeight{437} & $(a^{[r]},a^{[s]},b,a^{[t]})$ $(r \leq s; rt=0)$ \\ 

\eeight{433} & $(a^{[r]},b,a^{[s]},c)$ $(rs=0)$ & \eeight{438} & $(a^{[r]},b,a^{[s]},a^{[t]})$ $(rst=0)$ \\
\endfirsthead
\eeight{434} & $(a^{[r]},b,c,a^{[s]})$ $(rs=0)$ & \eeight{439} & $(a,a^{[r]},b^{[s]},b^{[t]})$ $(st=0)$ \\

\eeight{435} & $(a,b,c^{[r]},c^{[s]})$ $(rs=0)$ & \eeight{440} & $(a^{[r]},b^{[t]},a^{[s]},b^{[u]})$ $(rs=tu=0)$ \\

\eeight{436} & $(a^{[r]},a^{[s]},a^{[t]},b)$ $(r \leq s$; $rt=0$; if $r=s$ then $r < t$) & \eeight{15} & $(a^{[r]},a^{[s]},a^{[t]},a^{[u]})$ ($rst=0$; $r \leq s$; \newline if $r=s$ then $r < t$) \\
\hline
\end{longtable}

\setlength{\LTleft}{0.06\textwidth}
\begin{longtable}{>{\raggedright\arraybackslash}p{0.06\textwidth - 2\tabcolsep}>{\raggedright\arraybackslash}p{0.82\textwidth-\tabcolsep}@{}}

\caption{Irreducible diagonal subgroups of $\bar{A}_1^2 A_1^3 B_2 = E_8(\#\eeight{441})$ $(p=2)$.  \label{barA12A13diags}} \\
\hline
ID & Embedding \\
\hline
\eeight{452} & $(a,a^{[r]},b,c,d,10)$ $(r \neq 0)$  \\

\eeight{453} & $(a^{[r]},b,a^{[s]},c,d,10)$ $(rs=0)$  \\

\endfirsthead
\eeight{454} & $(a,b,c,c^{[r]},d,10)$ $(r \neq 0)$  \\

\eeight{455} & $(a^{[r]},a^{[s]},a^{[t]},b,c,10)$ $(r < s; rt=0)$  \\

\eeight{456} & $(a^{[r]},b,a^{[s]},a^{[t]},c,10)$ $(rs=0$; $s < t$) \\

\eeight{457} & $(a,b,c,c^{[r]},c^{[s]},10)$ $(0 < r < s)$  \\

\eeight{458} & $(a,a^{[r]},b,b^{[s]},c,10)$ $(rs \neq 0)$  \\

\eeight{459} & $(a^{[r]},b^{[t]},a^{[s]},b^{[u]},c,10)$ $(rs=tu=0$; $r \leq t$; if $r=t$ then $s \leq u$)  \\

\eeight{460} & $(a^{[r]},b,a^{[s]},c,c^{[t]},10)$ $(rs=0 \neq t)$  \\ 

\eeight{461} & $(a^{[r]},a^{[s]},a^{[t]},a^{[u]},b,10)$ $(rt=0; r<s; t<u)$  \\

\eeight{462} & $(a^{[r]},b,a^{[s]},a^{[t]},a^{[u]},10)$ $(rs=0; s < t < u)$ \\

\eeight{463} & $(a^{[r]},a^{[s]},a^{[t]},b,b^{[u]},10)$ $(rt=0 \neq u; r < s)$  \\

\eeight{464} & $(a^{[r]},b^{[u]},a^{[s]},a^{[t]},b^{[v]},10)$ $(rs=uv=0; s < t)$  \\

\eeight{465} & $(a,a^{[r]},b,b^{[s]},b^{[t]},10)$ $(r \neq 0; 0 < s < t)$ \\

\eeight{466} & $(a^{[r]},a^{[s]},a^{[t]},a^{[u]},a^{[v]},10)$ $(rt=0; r<s; t<u<v)$ \\

\hline
\end{longtable}

\setlength{\LTleft}{0pt}
\begin{longtable}{>{\raggedright\arraybackslash}p{0.06\textwidth - 2\tabcolsep}>{\raggedright\arraybackslash}p{0.94\textwidth-\tabcolsep}@{}}

\caption{Irreducible diagonal subgroups of $\bar{A}_1^2 A_1^5 = E_8(\#\eeight{451})$ $(p=2)$.  \label{barA12A15diags}} \\
\hline
ID & Embedding \\
\hline

\eeight{467} & $(a,a^{[r]},b,c,d,e,f)$ $(r \neq 0)$ \\
\eeight{468} & $(a^{[r]},b,a^{[s]},c,d,e,f)$ $(rs = 0)$ \\
\eeight{469} & $(a,b,c,c^{[r]},d,e,f)$ $(r \neq 0)$ \\ \endfirsthead

\eeight{470} & $(a^{[r]},a^{[s]},a^{[t]},b,c,d,e)$ $(rt = 0; r <s)$ \\
\eeight{471} & $(a^{[r]},b,a^{[s]},a^{[t]},c,d,e)$ $(rs = 0; s <t)$ \\
\eeight{472} & $(a,b,c,c^{[r]},c^{[s]},d,e)$ $(0 < r < s)$ \\

\eeight{473} & $(a,a^{[r]},b,b^{[s]},c,d,e)$ $(rs \neq 0)$ \\
\eeight{474} & $(a^{[r]},b^{[t]},a^{[s]},b^{[u]},c,d,e)$ $(rs = tu = 0; r \leq t$; if $r=t$ then $s < u$) \\
\eeight{475} & $(a^{[r]},b,a^{[s]},c,c^{[t]},d,e)$ $(rs = 0 \neq t)$ \\
\eeight{476} & $(a,b,c,c^{[r]},d,d^{[s]},e)$ $(rs \neq 0; r \leq s)$ \\

\eeight{477} & $(a^{[r]},a^{[s]},a^{[t]},a^{[u]},b,c,d)$ $(rt = 0; r <s;  t < u)$ \\
\hline
\end{longtable}

\begin{longtable}{>{\raggedright\arraybackslash}p{0.06\textwidth - 2\tabcolsep}>{\raggedright\arraybackslash}p{0.39\textwidth - 2\tabcolsep}>{\raggedright\arraybackslash}p{0.06\textwidth - 2\tabcolsep}>{\raggedright\arraybackslash}p{0.49\textwidth-\tabcolsep}@{}}

\caption{Irreducible diagonal subgroups of $\bar{A}_1^2 A_1 B_2^2 = E_8(\#\eeight{404})$ $(p=2)$.  \label{barA12A1B23diags}} \\
\hline
ID & Embedding & ID & Embedding \\
\hline
\eeight{520} & $(1_a,1_a^{[r]},1_b,10_c,10_d)$ $(r \neq 0)$ & \eeight{526} & $(1_a,1_a^{[r]},1_b,10_c,02_c^{[s]})$ $(r \neq 0)$ \\

\eeight{521} & $(1_a^{[r]},1_b,1_a^{[s]},10_c,10_d)$ $(rs=0)$ & \eeight{527} & $(1_a^{[r]},1_b,1_a^{[s]},10_c,10_c^{[t]})$ $(rs=0 \neq t)$ \\
\endfirsthead
\eeight{522} & $(1_a,1_b,1_c,10_d,10_d^{[r]})$ $(r \neq 0)$ & \eeight{528} & $(1_a^{[r]},1_b,1_a^{[s]},10_c,02_c^{[t]})$ $(rs=0)$ \\

\eeight{523} \newline \eeight{524}  & $(1_a,1_b,1_c,10_d,02_d^{[r]})$ \newline $(1_a^{[r]},1_a^{[s]},1_a^{[t]},10_b,10_c)$  & \eeight{529} & $(1_a^{[r]},1_a^{[s]},1_a^{[t]},10_b,10_b^{[u]})$ $(r < s; rt=0\neq u)$ \\

& $(rt=0; r<s)$ & \eeight{530} & $(1_a^{[r]},1_a^{[s]},1_a^{[t]},10_b,02_b^{[u]})$ $(r<s; rt=0)$ \\

\eeight{525} & $(1_a,1_a^{[r]},1_b,10_c,10_c^{[s]})$ $(rs \neq 0)$ & & \\
\hline
\end{longtable}

\begin{longtable}{>{\raggedright\arraybackslash}p{0.07\textwidth - 2\tabcolsep}>{\raggedright\arraybackslash}p{0.4\textwidth - 2\tabcolsep}>{\raggedright\arraybackslash}p{0.07\textwidth - 2\tabcolsep}>{\raggedright\arraybackslash}p{0.46\textwidth-\tabcolsep}@{}}

\caption{Irreducible diagonal subgroups of $\bar{A}_1^2 A_1 A_1 = E_8(\#\eeight{538})$ $(p \geq 7)$.  \label{barA12A1A1p7diags}} \\
\hline
ID & Embedding & ID & Embedding \\
\hline
\eeight{542} & $(a,a^{[r]},b,c)$ & \eeight{547} & $(a^{[r]},a^{[s]},b,a^{[t]})$ $(r \leq s; rt=0)$ \\ 

\eeight{543} & $(a^{[r]},b,a^{[s]},c)$ $(rs=0)$ & \eeight{548} & $(a^{[r]},b,a^{[s]},a^{[t]})$ $(rst=0)$ \\
\endfirsthead
\eeight{544} & $(a^{[r]},b,c,a^{[s]})$ $(rs=0)$ & \eeight{549} & $(a,a^{[r]},b^{[s]},b^{[t]})$ $(st=0)$ \\

\eeight{545} & $(a,b,c^{[r]},c^{[s]})$ $(rs=0)$ & \eeight{550} & $(a^{[r]},b^{[t]},a^{[s]},b^{[u]})$ $(rs=tu=0)$ \\

\eeight{546} & $(a^{[r]},a^{[s]},a^{[t]},b)$ $(r \leq s$; $rt=0$) & \eeight{17} & $(a^{[r]},a^{[s]},a^{[t]},a^{[u]})$ ($rst=0$; $r \leq s$) \\
\hline
\end{longtable}

\begin{longtable}{>{\raggedright\arraybackslash}p{0.09\textwidth - 2\tabcolsep}>{\raggedright\arraybackslash}p{0.91\textwidth-\tabcolsep}@{}}

\caption{Irreducible diagonal subgroups of $\bar{A}_1^2 A_1^2 A_1 = E_8(\#\eeight{558})$ $(p \neq 2)$.  \label{barA12A12A1diags}} \\
\hline
ID & Embedding \\
\hline

\eeight{566} & $(a,a^{[r]},b,c,d)$ \\
\eeight{567} & $(a^{[r]},b,a^{[s]},c,d)$ $(rs=0)$ \\
\endfirsthead
\eeight{568} & $(a^{[r]},b,c,a^{[s]},d)$ $(rs=0)$ \\
\eeight{569} & $(a^{[r]},b,c,d,a^{[s]})$ $(rs=0)$ \\
\eeight{570} & $(a,b,c,c^{[r]},d)$ $(r \neq 0)$ \\
$\mathit{\eeight{411}^{\{\underline{0}\}}}$ & $(a,b,c,c,d)$ \\
\eeight{571} & $(a,b,c^{[r]},d,c^{[s]})$ $(rs = 0; r \neq s)$ \\
$\mathit{\eeight{227}^{\{\underline{0}\}}}$ & $(a,b,c,d,c)$ \\

\eeight{572} & $(a^{[r]},a^{[s]},a^{[t]},b,c)$ $(rst=0)$ \\
\eeight{573} & $(a^{[r]},a^{[s]},b,c,a^{[t]})$ $(rt=0; r \leq s)$ \\
\hline
\end{longtable}

\begin{longtable}{>{\raggedright\arraybackslash}p{0.09\textwidth - 2\tabcolsep}>{\raggedright\arraybackslash}p{0.34\textwidth - 2\tabcolsep}>{\raggedright\arraybackslash}p{0.06\textwidth - 2\tabcolsep}>{\raggedright\arraybackslash}p{0.51\textwidth-\tabcolsep}@{}}

\caption{Irreducible diagonal subgroups of $\bar{A}_1^2 A_1 A_1 = E_8(\#\eeight{559})$ $(p \geq 7)$.  \label{barA12A1A1p7noS2diags}} \\
\hline
ID & Embedding & ID & Embedding \\
\hline
\eeight{598} & $(a^{[r]},a^{[s]},b,c)$ $(rs=0)$ & \eeight{605} & $(a^{[r]},a^{[s]},b,a^{[t]})$ $(rst=0)$ \\ 

\eeight{599} & $(a^{[r]},b,a^{[s]},c)$ $(rs=0)$ & \eeight{606} & $(a^{[r]},b,a^{[s]},a^{[t]})$ $(rst=0; s \neq t)$ \\
\endfirsthead
\eeight{600} & $(a^{[r]},b,c,a^{[s]})$ $(rs=0)$ & \eeight{607} & $(a,b^{[r]},b^{[s]},b^{[t]})$ $(rst=0; s \neq t)$ \\

\eeight{601} & $(a,b^{[r]},b^{[s]},c)$ $(rs=0)$ & \eeight{608} & $(a^{[r]},a^{[s]},b^{[t]},b^{[u]})$ $(rs=tu=0; t \neq u)$ \\

\eeight{602} & $(a,b^{[r]},c,b^{[s]})$ $(rs=0)$ & \eeight{609} & $(a^{[r]},b^{[t]},a^{[s]},b^{[u]})$ $(rs=tu=0)$ \\

\eeight{603} & $(a,b,c^{[r]},c^{[s]})$ $(r \neq s)$ & \eeight{610} & $(a^{[r]},b^{[t]},b^{[u]},a^{[s]})$ $(rs=tu=0)$ \\

$\mathit{\eeight{545}^{\{\underline{0}\}}}$ & $(a,b,c,c)$ & \eeight{5} & $(a^{[r]},a^{[s]},a^{[t]},a^{[u]})$ ($rstu=0$; $t \neq u$)  \\

\eeight{604} & $(a^{[r]},a^{[s]},a^{[t]},b)$ $(rst=0$) & \\
\hline
\end{longtable}

\begin{longtable}{>{\raggedright\arraybackslash}p{0.06\textwidth - 2\tabcolsep}>{\raggedright\arraybackslash}p{0.94\textwidth-\tabcolsep}@{}}

\caption{Irreducible diagonal subgroups of $A_1^6 = E_8(\#\eeight{641})$ $(p=2)$.  \label{A16diags}} \\
\hline
ID & Embedding \\
\hline

\eeight{648} & $(a,a^{[r]},b,c,d,e)$ $(r \neq 0)$ \\
\eeight{649} & $(a^{[r]},b,c,a^{[s]},d,e)$ $(rs = 0)$ \\\endfirsthead

\eeight{650} & $(a,a^{[r]},a^{[s]},b,c,d)$ $(0 < r < s)$ \\
\eeight{651} & $(a^{[r]},a^{[s]},b,a^{[t]},c,d)$ $(rt = 0; r < s)$ \\

\eeight{652} & $(a,a^{[r]},b^{[s]},b^{[t]},c,d)$ $(r \neq 0 = st)$ \\
\eeight{653} & $(a,a^{[r]},b,c,c^{[s]},d)$ $(rs \neq 0; r \leq s)$ \\
\eeight{654} & $(a,b^{[s]},c,a^{[r]},b^{[t]},d)$ $(st =0; r \leq s+t;$ if $r=0$ then $s \leq t$) \\

\eeight{655} & $(a^{[r]},a^{[s]},a^{[t]},a^{[u]},b,c)$ $(ru = 0; r < s < t)$ \\
\eeight{656} & $(a^{[r]},a^{[s]},b,a^{[t]},a^{[u]},c)$ $(rt = 0; r < s; t < u; r \leq t$; if $r=t$ then $s \leq u$) \\

\eeight{657} & $(a,a^{[r]},a^{[s]},b,b^{[t]},c)$ $(t \neq 0 < r < s)$ \\
\hline
\end{longtable}

\begin{longtable}{>{\raggedright\arraybackslash}p{0.06\textwidth - 2\tabcolsep}>{\raggedright\arraybackslash}p{0.39\textwidth - 2\tabcolsep}>{\raggedright\arraybackslash}p{0.06\textwidth - 2\tabcolsep}>{\raggedright\arraybackslash}p{0.49\textwidth-\tabcolsep}@{}}

\caption{Irreducible diagonal subgroups of $A_1^2 A_1^2 = E_8(\#\eeight{683})$ $(p \geq 5)$.  \label{A1A1A1A1diags1st}} \\
\hline
ID & Embedding & ID & Embedding \\
\hline

\eeight{687} & $(a,a^{[r]},b,c)$ $(r \neq 0)$ & \eeight{691} & $(a^{[r]},b,a^{[s]},a^{[t]})$ $(rs=0; s < t)$ \\
\eeight{688} & $(a^{[r]},b,a^{[s]},c)$ $(rs=0)$ & \eeight{692} & $(a,a^{[r]},b,b^{[s]})$ $(rs \neq 0)$ \\ \endfirsthead
\eeight{689} & $(a,b,c,c^{[r]})$ $(r \neq 0)$ & \eeight{693} & $(a^{[r]},b^{[t]},a^{[s]},b^{[u]})$ $(rs=tu=0; r \leq t$; \newline if $r=t$ then $s \leq u$) \\
\eeight{690} & $(a^{[r]},a^{[s]},a^{[t]},b)$ $(rt=0; r < s)$ & \eeight{21} & $(a^{[r]},a^{[s]},a^{[t]},a^{[u]})$ $(rt=0; r < s; t < u)$ \\
\hline
\end{longtable}

\begin{longtable}{>{\raggedright\arraybackslash}p{0.1\textwidth - 2\tabcolsep}>{\raggedright\arraybackslash}p{0.9\textwidth-\tabcolsep}@{}}

\caption{Irreducible diagonal subgroups of $A_1^2 A_1^2 = E_8(\#\eeight{694})$ $(p \geq 5)$.  \label{A1A1A1A1diags2nd}} \\
\hline
ID & Embedding \\
\hline
\eeight{699} & $(a,a^{[r]},b,c)$ \\
\eeight{700} & $(a^{[r]},b,a^{[s]},c)$ $(rs=0; r \neq s)$ \\
$\mathit{\eeight{688}^{\{\underline{0}\}}}$ &  $(a,b,a,c)$ \\ 
\eeight{701} & $(a^{[r]},b,c,a^{[s]})$ $(rs=0)$ \\ 
\endfirsthead
\eeight{702} & $(a,b,c,c^{[r]})$ \\
\eeight{703} & $(a^{[r]},a^{[s]},a^{[t]},b)$ $(rst=0; r \neq t)$  \\ 
\eeight{704} & $(a^{[r]},b,a^{[s]},a^{[t]})$ $(rst=0; r \neq s)$ \\
\eeight{705} & $(a,a^{[r]},b^{[s]},b^{[t]})$ $(st =0)$ \\
\eeight{706} & $(a^{[r]},b^{[t]},a^{[s]},b^{[u]})$ $(rs=tu=0; r \neq s; t \neq u; r \leq t$; if $r=t$ then $s \leq u$) \\ 
\eeight{707} & $(a^{[r]},b^{[t]},b^{[u]},a^{[s]})$ $(rs=tu=0$; $r \leq t$; if $r=t$ then $s \leq u$) \\
\eeight{22} & $(a^{[r]},a^{[s]},a^{[t]},a^{[u]})$ $(rtu=0; r \neq t; s \neq u; r \leq s$; if $r = s$ then $t \leq u$) \\
\hline
\end{longtable}%

%% file: CompFactorsIrred.tex
\chapter{Composition factors for $G$-irreducible subgroups} \label{tab:compositionfactors}

Let $G$ be a simple exceptional algebraic group. In this section we give the composition factors of the action of each $G$-irreducible connected subgroup $X$ on the minimal and adjoint modules for $G$. However, we do not explicitly list these for diagonal irreducible subgroups $X$ of some semisimple irreducible subgroup $Y$ as it would drastically lengthen the already large tables and, more importantly, because it is easy to recover them from the $Y$-composition factors of the $G$-modules in question. The notation used for describing the composition factors is given in Section \ref{nota}. 

The composition factors of each $G$-irreducible subgroup are found by repeated restrictions from a reductive, maximal connected overgroup $M$. When $X$ is simple this has already been done in \cite[Tables 9--13]{tho2}, \cite[Tables 3--7]{tho1}, where the $X$-composition factors of the minimal and adjoint modules are given.  

The first column of each table gives the identification number $n$ of the irreducible subgroup $X = G(\#n)$. The second column gives the isomorphism type of $X$ and the third gives any characteristic restrictions. In Tables \ref{G2tabcomps}--\ref{E7tabcomps}, the fourth (resp. fifth) column gives the composition factors of the minimal module (resp. adjoint module) for $G$ restricted to $X$. In Table \ref{E8tabcomps} the fourth column gives the composition factors of $L(E_8) \downarrow X$.

\begin{landscape}
\setlength{\LTcapwidth}{\linewidth}
\setlength\LTleft{0.23\linewidth}

\end{landscape}

%% file: CompFactorsLevi.tex
\chapter{Composition factors for the action of Levi subgroups} \label{levicomps}

Let $G$ be a simple exceptional algebraic group. In this section we give the composition factors of the action of proper Levi subgroups $L$ of $G$ on the minimal and adjoint modules for $G$. If $L'$ is simple and has rank at least 2 then these are found in \cite[Tables 8.1--8.7]{LS3}, except for the case where $L < G = F_4$ acting on the minimal module. In all other cases the composition factors are deduced from those of a maximal subsystem subgroup containing $L'$. 

We note that for $G = E_6, E_7$ or $E_8$ any simple factor of a Levi subgroup of type $A_n$, $D_n$ or $E_n$ is generated by long root subgroups of $G$ and so we omit all bars in Tables \ref{levie6}--\ref{levie8}.

\begin{longtable}{p{0.1\textwidth - 2\tabcolsep}>{\raggedright\arraybackslash}p{0.3\textwidth-2\tabcolsep}>{\raggedright\arraybackslash}p{0.6\textwidth-\tabcolsep}@{}}

\caption{The composition factors for the action of Levi subgroups of $G_2$ on $V_{7}$ and $L(G_2)$. \label{levig2}} \\

\hline \noalign{\smallskip}

$L'$ & $V_{7} \downarrow L'$ & $L(G_2) \downarrow L'$ \\

\hline \noalign{\smallskip}

$\bar{A}_1$ & $1^2 /$ $\!\! 0^{3}$ & $W(2) /$ $\!\! 1^{4} /$ $\!\! 0^{3}$ \\

$\tilde{A}_1$ & $W(2) /$ $\!\! 1^2$ & $W(3)^2 /$ $\!\! W(2) /$ $\!\! 0^{3}$ \\

\hline
\endfirsthead
\end{longtable}

\begin{longtable}{p{0.09\textwidth - 2\tabcolsep}>{\raggedright\arraybackslash}p{0.33\textwidth-2\tabcolsep}>{\raggedright\arraybackslash}p{0.58\textwidth-\tabcolsep}@{}}

\caption{The composition factors for the action of Levi subgroups of $F_4$ on $V_{26}$ and $L(F_4)$. \label{levif4}} \\

\hline \noalign{\smallskip}

$L'$ & $V_{26} \downarrow L'$ & $L(F_4) \downarrow L'$ \\

\hline \noalign{\smallskip}

$B_3$ & $W(100) /$ $\!\! 001^2 /$ $\!\! 000^3$ & $W(100)^2 /$ $\!\! W(010) /$ $\!\! 001^2 /$ $\!\! 000$ \\

$B_2$ & $W(10) /$ $\!\! 01^4 /$ $\!\! 00^5$ & $W(10)^4 /$ $\!\! W(02) /$ $\!\! 01^4 /$ $\!\! 00^6$ \\
\endfirsthead
$C_3$ & $100^2 /$ $\!\! W(010)$ & $W(200) /$ $\!\! W(001)^2 /$ $\!\! 000^3$ \\

$\bar{A}_2 \tilde{A}_1$ & $(10,1) /$ $\!\! (10,0) /$ $\!\! (01,1) /$ $\!\! (01,0) /$ $\!\! (00,W(2)) /$ $\!\! (00,1)^2 /$ $\!\! (00,0)$ & $(W(11),0) /$ $\!\! (10,W(2)) /$ $\!\! (10,1) /$ $\!\! (10,0) /$ $\!\! (01,W(2)) /$ $\!\! (01,1) /$ $\!\! (01,0) /$ $\!\! (00,W(2)) /$ $\!\! (00,1)^2 /$ $\!\! (00,0)$  \\

$\tilde{A}_2 \bar{A}_1$ & $(10,1) /$ $\!\! (10,0) /$ $\!\! (01,1) /$ $\!\! (01,0) /$ $\!\! (W(11),0)$ & $(W(20),1) /$ $\!\! (W(20),0) /$ $\!\! (W(11),0) /$ $\!\! (W(02),1) /$ $\!\! (W(02),0) /$ $\!\! (00,W(2)) /$ $\!\! (00,1)^2 /$ $\!\! (00,0)$ \\

$\bar{A}_2$ & $10^3 /$ $\!\! 01^3 /$ $\!\! 00^8$ & $W(11) /$ $\!\! 10^6 /$ $\!\! 01^6 /$ $\!\! 00^8$ \\

$\tilde{A}_2$ & $10^3 /$ $\!\! 01^3 /$ $\!\! W(11)$ & $W(20)^3 /$ $\!\! W(11) /$ $\!\! W(02)^3 /$ $\!\! 00^8$ \\

$\bar{A}_1 \tilde{A}_1$ & $(1,1)^2 /$ $\!\! (1,0)^2 /$ $\!\! (0,W(2)) /$ $\!\! (0,1)^4 /$  $\!\! (0,0)^3$ & $(W(2),0) /$ $\!\! (1,W(2))^2 /$ $\!\! (1,1)^2 /$ $\!\! (1,0)^4 /$ $\!\! (0,W(2))^3 /$  $\!\! (0,1)^4 /$ $\!\! (0,0)^4$ \\

$\bar{A}_1$ & $1^6 /$ $\!\! 0^{14}$ & $W(2) /$ $\!\! 1^{14} /$ $\!\! 0^{21}$ \\

$\tilde{A}_1$ & $W(2) /$ $\!\! 1^8 /$ $\!\! 0^7$ & $W(2)^7 /$ $\!\! 1^8 /$ $\!\! 0^{15}$ \\

\hline

\end{longtable}

\begin{longtable}{p{0.09\textwidth - 2\tabcolsep}>{\raggedright\arraybackslash}p{0.35\textwidth-2\tabcolsep}>{\raggedright\arraybackslash}p{0.56\textwidth-\tabcolsep}@{}}

\caption{The composition factors for the action of Levi subgroups of $E_6$ on $V_{27}$ and $L(E_6)$. \label{levie6}} \\

\hline \noalign{\smallskip}

$L'$ & $V_{27} \downarrow L'$ & $L(E_6) \downarrow L'$ \\

\hline \noalign{\smallskip}

$D_5$ & $\lambda_1 /$ $\!\! \lambda_4 /$ $\!\! 0 $ & $W(\lambda_2) /$ $\!\! \lambda_4 /$ $\!\! \lambda_5 /$ $\!\! 0$ \\

$D_4$ & $\lambda_1 /$ $\!\! \lambda_3 /$ $\!\! \lambda_4 /$ $\!\! 0^3$  & $\lambda_1^2 /$ $\!\! W(\lambda_2) /$  $\!\! \lambda_3^2 /$ $\!\! \lambda_4^2 /$ $\!\! 0^3$ \\
\endfirsthead
$A_5$ & $\lambda_1^2 /$ $\!\! \lambda_4$ & $W(\lambda_1 + \lambda_5) /$ $\!\! \lambda_3^2 /$ $\!\! 0^3$ \\

$A_1 A_4$ & $(1,\lambda_1) /$ $\!\! (1,0) /$ $\!\! (0,\lambda_3) /$ $\!\! (0,\lambda_4) $ & $(W(2),0) /$ $\!\! (1,\lambda_2) /$ $\!\! (1,\lambda_3) /$ $\!\! (0,W(\lambda_1 + \lambda_4)) /$ $\!\! (0,\lambda_1) /$  $\!\! (0,\lambda_4) /$ $\!\! (0,0)$ \\

$A_1 A_2^2$ & $(1,01,00) /$  $\!\! (1,00,10) /$ $\!\! (0,10,01) /$ $\!\! (0,01,00) /$ $\!\! (0,00,10)$  & $(W(2),00,00) /$ $\!\! (1,10,10) /$ $\!\! (1,01,01) /$ $\!\! (1,00,00)^2 /$ $\!\! (0,W(11),00) /$ $\!\! (0,10,10) /$  $\!\! (0,01,01) /$  $\!\! (0,00,W(11)) /$    $\!\! (0,00,00)$ \\

$A_4$ & $\lambda_1^2 /$ $\!\! \lambda_3 /$ $\!\! \lambda_4 /$ $\!\! 0^2$ & $W(\lambda_1 + \lambda_4) /$ $\!\! \lambda_1 /$ $\!\! \lambda_2^2 /$ $\!\! \lambda_3^2 /$ $\!\! \lambda_4 /$ $\!\! 0^4$ \\

$A_1 A_3$ & $(1,100) /$ $\!\! (1,000)^2 /$ $\!\! (0,010) /$ $\!\! (0,001)^2 /$ $\!\! (0,000)$ & $(W(2),000) /$ $\!\! (1,100) /$ $\!\! (1,010)^2 /$ $\!\! (1,001) /$ $\!\! (0,W(101)) /$  $\!\! (0,100)^2 /$ $\!\! (0,001)^2 /$   $\!\! (0,000)^4$ \\

$A_2^2$ & $(10,01) /$ $\!\! (01,00)^3 /$ $\!\! (00,10)^3$  & $(W(11),00) /$ $\!\! (10,10)^3 /$ $\!\! (01,01)^3 /$ $\!\! (00,W(11)) /$ $\!\! (00,00)^8$ \\

$A_1^2 A_2$ & $(1,1,00) /$ $\!\! (1,0,10) /$ $\!\! (1,0,00) /$ $\!\! (0,1,01) /$ $\!\! (0,1,00) /$ $\!\! (0,0,10) /$ $\!\! (0,0,01) /$ $\!\! (0,0,00)$ & $(W(2),0,00) /$ $\!\! (1,1,10) /$ $\!\! (1,1,01) /$  $\!\! (1,0,10) /$ $\!\! (1,0,01) /$ $\!\! (1,0,00)^2 /$ $\!\! (0,W(2),00) /$ $\!\! (0,1,10) /$ $\!\! (0,1,01) /$    $\!\! (0,1,00)^2 /$ $\!\! (0,0,W(11)) /$   $\!\! (0,0,10) /$  $\!\! (0,0,01) /$  $\!\! (0,0,00)^2$ \\ 

$A_3$ & $100^2 /$ $\!\! 010 /$ $\!\! 001^2 /$ $\!\! 000^5$ & $W(101) /$ $\!\! 100^4 /$ $\!\! 010^4 /$ $\!\! 001^4 /$ $\!\! 000^7$ \\

$A_1 A_2$ & $(1,00)^3 /$ $\!\! (1,01) /$ $\!\! (0,10)^3 /$ $\!\! (0,01) /$ $\!\! (0,00)^3$ & $ (W(2),00) /$ $\!\! (1,10)^3 /$ $\!\! (1,01)^3 /$  $\!\! (1,00)^2 /$ $\!\! (0,10)^3 /$ $\!\! (0,W(11)) /$ $\!\! (0,01)^3 /$ $\!\! (0,00)^9$ \\

$A_1^3$ & $(1,1,0) /$ $\!\! (1,0,1) /$ $\!\! (1,0,0)^2 /$ $\!\! (0,1,1) /$  $\!\! (0,1,0)^2 /$ $\!\! (0,0,1)^2 /$ $\!\! (0,0,0)^3$ & $(W(2),0,0) /$ $\!\! (1,1,1)^2 /$ $\!\! (1,1,0)^2 /$ $\!\! (1,0,1)^2 /$ $\!\! (1,0,0)^4 /$ $\!\! (0,W(2),0) /$ $\!\! (0,1,1)^2 /$  $\!\! (0,1,0)^4 /$ $\!\! (0,0,W(2)) /$   $\!\! (0,0,1)^4 /$ $\!\! (0,0,0)^5$ \\ 

$A_2$ & $10^3 /$ $\!\! 01^3 /$ $\!\! 00^9$ & $W(11) /$ $\!\! 10^9 /$ $\!\! 01^9 /$ $\!\! 00^{16}$ \\

$A_1^2$ & $(1,1) /$ $\!\! (1,0)^4 /$ $\!\! (0,1)^4 /$ $\!\! (0,0)^7$ & $(W(2),0) /$ $\!\! (1,1)^6 /$ $\!\! (1,0)^8 /$ $\!\! (0,W(2)) /$  $\!\! (0,1)^8 /$ $\!\! (0,0)^{16}$ \\ 

$A_1$ & $1^6 /$ $\!\! 0^{15}$ & $W(2) /$ $\!\! 1^{20} /$ $\!\! 0^{35}$ \\
\hline
\end{longtable}

\begin{longtable}{p{0.11\textwidth - 2\tabcolsep}>{\raggedright\arraybackslash}p{0.37\textwidth-2\tabcolsep}>{\raggedright\arraybackslash}p{0.52\textwidth-\tabcolsep}@{}}

\caption{The composition factors for the action of Levi subgroups of $E_7$ on $V_{56}$ and $L(E_7)$. \label{levie7}} \\

\hline \noalign{\smallskip}

$L'$ & $V_{56} \downarrow L'$ & $L(E_7) \downarrow L'$ \\

\hline \noalign{\smallskip}

$E_6$ & $\lambda_1 /$ $\!\! \lambda_6 /$ $\!\! 0^2$ & $\lambda_1 /$ $\!\! W(\lambda_2) /$  $\!\! \lambda_6 /$ $\!\! 0$  \\

$D_6$ & $\lambda_1^2 /$ $\!\! \lambda_5$ & $W(\lambda_2) /$ $\!\! \lambda_6^2 /$ $\!\! 0^3$ \\
\endfirsthead
$A_1 D_5$ & $(1, \lambda_1) /$ $\!\! (1,0)^2 /$ $\!\! (0,\lambda_4) /$ $\!\! (0,\lambda_5)$ & $(W(2),0) /$ $\!\! (1,\lambda_4) /$ $\!\! (1,\lambda_5) /$ $\!\! (0,\lambda_1)^2 /$ $\!\! (0,W(\lambda_2)) /$ $\!\! (0,0)$  \\

$D_5$ & $\lambda_1^2 /$ $\!\! \lambda_4 /$ $\!\! \lambda_5 /$ $\!\! 0^4$ & $\lambda_1^2 /$ $\!\! W(\lambda_2) /$  $\!\! \lambda_4^2 /$ $\!\! \lambda_5^2 /$ $\!\! 0^4$  \\

$A_1 D_4$ & $(1,\lambda_1) /$ $\!\! (1,0)^4 /$ $\!\! (0,\lambda_3)^2 /$ $\!\! (0,\lambda_4)^2$ & $(W(2),0) /$ $\!\! (1,\lambda_3)^2 /$ $\!\! (1,\lambda_4)^2 /$ $\!\! (0,\lambda_1)^4 /$ $\!\! (0,W(\lambda_2)) /$ $\!\! (0,0)^6$ \\

$D_4$ & $\lambda_1^2 /$ $\!\! \lambda_3^2 /$ $\!\! \lambda_4^2 /$ $\!\! 0^8$ & $\lambda_1^4 /$ $\!\! W(\lambda_2) /$  $\!\! \lambda_3^4 /$ $\!\! \lambda_4^4 /$ $\!\! 0^9$ \\

$A_6$ & $\lambda_1 /$ $\!\! \lambda_2 /$ $\!\! \lambda_5 /$ $\!\! \lambda_6$ & $W(\lambda_1 + \lambda_6) /$ $\!\! \lambda_1 /$ $\!\! \lambda_3 /$ $\!\! \lambda_4 /$ $\!\! \lambda_6 /$ $\!\! 0$  \\

$A_1 A_5$ & $(1,\lambda_1) /$ $\!\! (1,\lambda_5) /$ $\!\! (0,\lambda_1) /$ $\!\! (0,\lambda_3) /$ $\!\! (0,\lambda_5)$ & $(W(2),0) /$ $\!\! (1,\lambda_2) /$ $\!\! (1,\lambda_4) /$ $\!\! (1,0)^2 /$ $\!\! (0,W(\lambda_1 + \lambda_5)) /$  $\!\! (0,\lambda_2) /$  $\!\! (0,\lambda_4) /$ $\!\! (0,0)$ \\

$A_2 A_4$ & $(10,\lambda_1) /$ $\!\! (10,0) /$ $\!\! (01,\lambda_4) /$ $\!\! (01,0) /$ $\!\! (00,\lambda_2) /$ $\!\! (00,\lambda_3)$ & $(W(11),0) /$ $\!\! (10,\lambda_3) /$ $\!\! (10,\lambda_4) /$ $\!\! (01,\lambda_1) /$ $\!\! (01,\lambda_2) /$  $\!\! (00,W(\lambda_1+\lambda_4)) /$ $\!\! (00,\lambda_1) /$ $\!\! (00,\lambda_4) /$  $\!\! (00,0)$  \\

$A_1 A_2 A_3$ & $(1,10,000) /$ $\!\! (1,01,000) /$ $\!\! (1,00,010) /$ $\!\! (0,10,100) /$ $\!\! (0,01,001) /$ $\!\! (0,00,100) /$  $\!\! (0,00,001)$ & $(W(2),00,000) /$ $\!\! (1,10,001) /$ $\!\! (1,01,100) /$ $\!\! (1,00,100) /$ $\!\! (1,00,001) /$ $\!\! (0,W(11),000) /$  $\!\! (0,10,010) /$ $\!\! (0,10,000) /$ $\!\! (0,01,010) /$  $\!\! (0,01,000) /$ $\!\! (0,00,W(101)) /$  $\!\! (0,00,000)$ \\

$A_5$ & $\lambda_1^3 /$ $\!\! \lambda_3 /$ $\!\! \lambda_5^3$  & $W(\lambda_1 + \lambda_5) /$ $\!\! \lambda_2^3 /$ $\!\! \lambda_4^3 /$ $\!\! 0^8$   \\

$A_5'$ & $\lambda_1^2 /$ $\!\! \lambda_2 /$ $\!\! \lambda_4 /$ $\!\! \lambda_5^2 /$ $\!\! 0^2$ & $W(\lambda_1 + \lambda_5) /$ $\!\! \lambda_1^2 /$ $\!\! \lambda_2 /$ $\!\! \lambda_3^2 /$ $\!\! \lambda_4 /$ $\!\! \lambda_5^2 /$ $\!\! 0^4$ \\

$A_1 A_4$ & $(1,\lambda_1) /$  $\!\! (1,\lambda_4) /$ $\!\! (1,0)^2 /$ $\!\! (0,\lambda_1) /$  $\!\! (0,\lambda_2) /$ $\!\! (0,\lambda_3) /$ $\!\! (0,\lambda_4) /$  $\!\! (0,0)^2$ & $(W(2),0) /$  $\!\! (1,\lambda_1) /$  $\!\! (1,\lambda_2) /$  $\!\! (1,\lambda_3) /$ $\!\! (1,\lambda_4) /$ $\!\! (1,0)^2 /$ $\!\! (0,W(\lambda_1+\lambda_4)) /$ $\!\! (0,\lambda_1)^2 /$ $\!\! (0,\lambda_2) /$ $\!\! (0,\lambda_3) /$ $\!\! (0,\lambda_4)^2 /$  $\!\! (0,0)^2$  \\

$A_2 A_3$ & $(10,100) /$ $\!\! (10,000)^2 /$ $\!\! (01,001) /$ $\!\! (01,000)^2 /$ $\!\! (00,100) /$ $\!\! (00,010)^2 /$  $\!\! (00,001)$ & $(W(11),000) /$  $\!\! (10,010) /$ $\!\! (10,001)^2 /$ $\!\! (10,000) /$ $\!\! (01,100)^2 /$ $\!\! (01,010) /$  $\!\! (01,000) /$ $\!\! (00,W(101)) /$ $\!\! (00,100)^2 /$ $\!\! (00,001)^2 /$  $\!\! (00,000)^4$ \\

$A_1^2 A_3$ & $(1,1,000)^2 /$ $\!\! (1,0,010) /$ $\!\! (1,0,000)^2 /$  $\!\! (0,1,100) /$  $\!\! (0,1,001) /$ $\!\! (0,0,100)^2 /$ $\!\! (0,0,001)^2$ & $(W(2),0,000) /$ $\!\! (1,1,100) /$ $\!\! (1,1,001) /$ $\!\! (1,0,100)^2 /$ $\!\! (1,0,001)^2 /$ $\!\! (0,W(2),000) /$   $\!\! (0,1,010)^2 /$ $\!\! (0,1,000)^4 /$ $\!\! (0,0,W(101)) /$ $\!\! (0,0,010)^2 /$  $\!\! (0,0,000)^4$ \\

$A_4$ & $\lambda_1^3 /$ $\!\! \lambda_2 /$ $\!\! \lambda_3 /$ $\!\! \lambda_4^3 /$ $\!\! 0^6 $ & $W(\lambda_1+\lambda_4) /$ $\!\! \lambda_1^4 /$ $\!\! \lambda_2^3 /$ $\!\! \lambda_3^3 /$ $\!\! \lambda_4^4 /$ $\!\! 0^9$  \\

$A_1 A_3$ & $(1,010) /$  $\!\! (1,000)^6 /$  $\!\! (0,100)^4 /$ $\!\! (0,001)^4$ & $(W(2),000) /$ $\!\! (1,100)^4 /$ $\!\! (1,001)^4 /$ $\!\! (0,W(101)) /$  $\!\! (0,010)^6 /$ $\!\! (0,000)^{15} $ \\

$(A_1 A_3)'$ & $(1,100) /$ $\!\! (1,001) /$  $\!\! (1,000)^4 /$ $\!\! (0,100)^2 /$ $\!\! (0,010)^2 /$   $\!\! (0,001)^2 /$ $\!\! (0,000)^4$ & $(W(2),000) /$ $\!\! (1,100)^2 /$ $\!\! (1,010)^2 /$ $\!\! (1,001)^2 /$ $\!\! (1,000)^4 /$ $\!\! (0,W(101)) /$  $\!\! (0,100)^4 /$ $\!\! (0,010)^2 /$ $\!\! (0,001)^4 /$  $\!\! (0,000)^7$ \\

$A_2^2$ & $(10,10) /$ $\!\! (10,00)^3 /$  $\!\! (01,01) /$ $\!\! (01,00)^3 /$ $\!\! (00,10)^3 /$  $\!\! (00,01)^3 /$  $\!\! (00,00)^2$ & $(W(11),00) /$ $\!\! (10,10)^3 /$ $\!\! (10,01) /$ $\!\! (10,00)^3 /$ $\!\! (01,10) /$ $\!\! (01,01)^3 /$ $\!\! (01,00)^3 /$ $\!\! (00,W(11)) /$   $\!\! (00,10)^3 /$ $\!\! (00,01)^3 /$ $\!\! (00,00)^{9} $ \\

$A_1^2 A_2$ &  $(1,1,00)^2 /$ $\!\! (1,0,10) /$ $\!\! (1,0,01) /$ $\!\! (1,0,00)^2 /$  $\!\! (0,1,10) /$ $\!\! (0,1,01) /$ $\!\! (0,1,00)^2 /$ $\!\! (0,0,10)^2 /$   $\!\! (0,0,01)^2 /$ $\!\! (0,0,00)^4$ & $(W(2),0,00) /$ $\!\! (1,1,10) /$ $\!\! (1,1,01) /$ $\!\! (1,1,00)^2 /$ $\!\! (1,0,10)^2 /$ $\!\! (1,0,01)^2 /$ $\!\! (1,0,00)^4 /$ $\!\! (0,W(2),00) /$   $\!\! (0,1,10)^2 /$ $\!\! (0,1,01)^2 /$  $\!\! (0,1,00)^4 /$  $\!\! (0,0,W(11)) /$   $\!\! (0,0,10)^3 /$ $\!\! (0,0,01)^3 /$ $\!\! (0,0,00)^5$ \\

$A_1^4$ & $(1,1,1,0) /$ $\!\! (1,0,0,1)^2 /$ $\!\! (1,0,0,0)^4 /$ $\!\! (0,1,0,1)^2 /$ $\!\! (0,1,0,0)^4 /$ $\!\! (0,0,1,1)^2 /$   $\!\! (0,0,1,0)^4$ & $(W(2),0,0,0) /$ $\!\! (1,1,0,1)^2 /$  $\!\! (1,1,0,0)^4 /$ $\!\! (1,0,1,1)^2 /$ $\!\! (1,0,1,0)^4 /$ $\!\! (0,W(2),0,0) /$ $\!\! (0,1,1,1)^2 /$  $\!\! (0,1,1,0)^4 /$ $\!\! (0,0,W(2),0) /$ $\!\! (0,0,0,W(2)) /$ $\!\! (0,0,0,1)^8 / (0,0,0,0)^9$ \\

$A_3$ & $100^4 /$ $\!\! 010^2 /$ $\!\! 001^4 /$  $\!\! 000^{12}$ & $W(101) /$ $\!\! 100^8 /$ $\!\! 010^{6} /$ $\!\! 001^8 /$  $\!\! 000^{18}$  \\

$A_1 A_2$ &  $ (1,10) /$ $\!\! (1,01) /$ $\!\! (1,00)^6 /$ $\!\! (0,10)^4 /$ $\!\! (0,01)^4 /$ $\!\! (0,00)^8$ & $(W(2),00) /$    $\!\! (1,10)^4 /$ $\!\! (1,01)^4 /$ $\!\! (1,00)^8 /$ $\!\! (0,W(11)) /$ $\!\! (0,10)^7 /$ $\!\! (0,01)^7 /$ $\!\! (0,00)^{16}$ \\

$A_1^3$ & $(1,1,0)^2 /$ $\!\! (1,0,1)^2 /$ $\!\! (1,0,0)^4 /$ $\!\! (0,1,1)^2 /$  $\!\! (0,1,0)^4 /$ $\!\! (0,0,1)^4 /$ $\!\! (0,0,0)^8$ & $(W(2),0,0) /$ $\!\! (1,1,1)^2 /$ $\!\! (1,1,0)^4 /$ $\!\! (1,0,1)^4 /$ $\!\! (1,0,0)^8 /$ $\!\! (0,W(2),0) /$  $\!\! (0,1,1)^4 /$   $\!\! (0,1,0)^8 /$ $\!\! (0,0,W(2)) /$ $\!\! (0,0,1)^8 /$ $ \!\! (0,0,0)^{12}$ \\

$(A_1^3)'$ & $(1,1,1) /$ $\!\! (1,0,0)^8 /$ $\!\! (0,1,0)^8 /$ $\!\! (0,0,1)^8$ & $(W(2),0,0) /$ $\!\! (1,1,0)^8 /$ $\!\! (1,0,1)^8 /$ $\!\! (0,W(2),0) /$  $\!\! (0,1,1)^8 /$ $\!\! (0,0,W(2)) /$  $\!\! (0,0,0)^{28}$  \\

$A_2$ & $10^6 /$ $\!\! 01^6 /$ $\!\! 00^{20}$ & $W(11) /$ $\!\! 10^{15} /$ $\!\! 01^{15} /$ $\!\! 00^{35}$ \\

$A_1^2$ & $(1,1)^2 /$ $\!\! (1,0)^8 /$ $\!\! (0,1)^8 /$ $\!\! (0,0)^{16}$ & $(W(2),0) /$ $\!\! (1,1)^7 /$ $\!\! (1,0)^{16} /$ $\!\! (0,W(2)) /$  $\!\! (0,1)^{16} /$ $\!\! (0,0)^{31}$ \\

$A_1$ & $1^{12} /$ $\!\! 0^{32}$ & $W(2) /$ $\!\! 1^{32} /$ $\!\! 0^{66}$ \\

\hline

\end{longtable}

\begin{longtable}{p{0.11\textwidth - 2\tabcolsep}>{\raggedright\arraybackslash}p{0.89\textwidth-\tabcolsep}@{}}

\caption{The composition factors for the action of Levi subgroups of $E_8$ on $L(E_8)$. \label{levie8}} \\

\hline \noalign{\smallskip}

$L'$ & $L(E_8) \downarrow L'$ \\

\hline \noalign{\smallskip}

$E_7$ & $W(\lambda_1) /$ $\!\! \lambda_7^2 /$ $\!\! 0^3$ \\

$A_1 E_6$ & $(W(2),0) /$ $\!\! (1,\lambda_1) /$ $\!\! (1,\lambda_6) /$ $\!\! (1,0)^2 /$ $\!\! (0,\lambda_1) /$ $\!\! (0,W(\lambda_2)) /$  $\!\! (0,\lambda_6) /$ $\!\! (0,0)$ \\
\endfirsthead
$E_6$ & $\lambda_1^3 /$ $\!\! W(\lambda_2) /$  $\!\! \lambda_6^3 /$ $\!\! 0^8$ \\

$D_7$ & $\lambda_1^2 /$ $\!\! W(\lambda_2) /$  $\!\! \lambda_6 /$ $\!\! \lambda_7 /$ $\!\! 0$  \\

$A_2 D_5$ & $(W(11),0) /$ $\!\! (10,\lambda_1) /$  $\!\! (10,\lambda_4) /$ $\!\! (10,0) /$ $\!\! (01,\lambda_1) /$ $\!\! (01,\lambda_5) /$ $\!\! (01,0) /$ $\!\! (00,W(\lambda_2)) /$  $\!\! (00,\lambda_4) /$  $\!\! (00,\lambda_5) /$ $\!\! (00,0)$ \\

$D_6$ & $ \lambda_1^4 /$ $\!\! W(\lambda_2) /$  $\!\! \lambda_5^2 /$ $\!\! \lambda_6^2 /$ $\!\! 0^6$ \\

$A_1 D_5$ & $(W(2),0) /$ $\!\! (1,\lambda_1)^2 /$ $\!\! (1,\lambda_4) /$ $\!\! (1,\lambda_5) /$ $\!\! (1,0)^4 /$  $\!\! (0,\lambda_1)^2 /$ $\!\! (0,W(\lambda_2)) /$  $\!\! (0,\lambda_4)^2 /$  $\!\! (0,\lambda_5)^2 /$ $\!\! (0,0)^4$ \\

$A_2 D_4$ & $(W(11),0) /$ $\!\! (10,\lambda_1) /$ $\!\! (10,\lambda_3) /$ $\!\! (10,\lambda_4) /$ $\!\! (10,0)^3 /$ $\!\! (01,\lambda_1) /$ $\!\! (01,\lambda_3) /$ $\!\! (01,\lambda_4) /$ $\!\! (01,0)^3 /$  $\!\! (00,\lambda_1)^2 /$ $\!\! (00,W(\lambda_2)) /$  $\!\! (00,\lambda_3)^2 /$ $\!\! (00,\lambda_4)^2 /$ $\!\! (00,0)^{2}$ \\

$D_5$ & $\lambda_1^6 /$ $\!\! W(\lambda_2) /$  $\!\! \lambda_4^4 /$ $\!\! \lambda_5^4 /$ $\!\! 0^{15}$ \\

$A_1 D_4$ & $(W(2),0) /$  $\!\! (1,\lambda_1)^2 /$ $\!\! (1,\lambda_3)^2 /$ $\!\! (1,\lambda_4)^2 /$ $\!\! (1,0)^8 /$ $\!\! (0,\lambda_1)^4 /$ $\!\! (0,W(\lambda_2)) /$ $\!\! (0,\lambda_3)^4 /$ $\!\! (0,\lambda_4)^4 /$ $\!\! (0,0)^{9}$ \\

$\bar{D}_4$ & $\lambda_1^8 /$ $\!\! W(\lambda_2) /$  $\!\! \lambda_3^8 /$ $\!\! \lambda_4^8 /$ $\!\! 0^{28}$ \\

$A_7$ & $W(\lambda_1 + \lambda_7) /$ $\!\! \lambda_1 /$ $\!\! \lambda_2 /$ $\!\! \lambda_3 /$ $\!\! \lambda_5 /$ $\!\! \lambda_6 /$ $\!\! \lambda_7 /$ $\!\! 0$ \\

$A_3 A_4$ & $(W(101),0) /$ $\!\! (100,\lambda_1) /$ $\!\! (100,\lambda_3) /$ $\!\! (100,0) /$ $\!\! (010,\lambda_1) /$ $\!\! (010,\lambda_4) /$ $\!\! (001,\lambda_2) /$ $\!\! (001,\lambda_4) /$ $\!\! (001,0) /$ $\!\! (000,W(\lambda_1+\lambda_4)) /$  $ \!\! (000,\lambda_2) /$ $\!\! (000,\lambda_3) /$ $\!\! (000,0)$ \\

$A_1 A_6$ & $(W(2),0) /$ $\!\! (1,\lambda_1) /$ $\!\! (1,\lambda_2) /$ $\!\! (1,\lambda_5) /$ $\!\! (1,\lambda_6) /$ $\!\! (0,W(\lambda_1 + \lambda_6)) /$  $\!\! (0,\lambda_1) /$ $\!\! (0,\lambda_3)/$ $\!\! (0,\lambda_4)/$ $\!\! (0,\lambda_6) /$ $\!\! (0,0)$ \\

$A_1 A_2 A_4$ & $(W(2),00,0) /$  $\!\! (1,10,\lambda_4) /$ $\!\! (1,10,0) /$ $\!\! (1,01,\lambda_1) /$ $\!\! (1,01,0) /$ $\!\! (1,00,\lambda_2) /$ $\!\! (1,00,\lambda_3) /$ $\!\! (0,W(11),0) /$  $\!\! (0,10,\lambda_1) /$ $\!\! (0,10,\lambda_2) /$ $\!\! (0,01,\lambda_3) /$ $\!\! (0,01,\lambda_4) /$ $\!\! (0,00,W(\lambda_1+\lambda_4)) /$ $\!\! (0,00,\lambda_1) /$ $\!\! (0,00,\lambda_4) /$ $\!\! (0,00,0)$ \\

$A_6$ & $W(\lambda_1 + \lambda_6) /$ $\!\! \lambda_1^3 /$ $\!\! \lambda_2^2 /$ $\!\! \lambda_3 /$ $\!\! \lambda_4 /$ $\!\! \lambda_5^2 /$ $\!\! \lambda_6^3 /$ $\!\! 0^4$  \\

$A_1 A_5$ & $(W(2),0) /$  $\!\! (1,\lambda_1)^2 /$ $\!\! (1,\lambda_2) /$ $\!\! (1,\lambda_4) /$ $\!\! (1,\lambda_5)^2 /$ $\!\! (1,0)^2 /$ $\!\! (0,W(\lambda_1 + \lambda_5)) /$ $\!\! (0,\lambda_1)^2 /$ $\!\! (0,\lambda_2) /$ $\!\! (0,\lambda_3)^2/$ $\!\! (0,\lambda_4) /$ $\!\! (0,\lambda_5)^2 /$ $\!\! (0,0)^4$ \\

$A_2 A_4$ & $(W(11),0) /$ $\!\! (10,\lambda_1) /$ $\!\! (10,\lambda_2) /$ $\!\! (10,\lambda_4)^2 /$ $\!\! (10,0)^2 /$ $\!\! (01,\lambda_1)^2 /$  $\!\! (01,\lambda_3) /$ $\!\! (01,\lambda_4) /$ $\!\! (01,0)^2 /$ $\!\! (00,W(\lambda_1+\lambda_4)) /$ $\!\! (00,\lambda_1) /$ $\!\! (00,\lambda_2)^2 /$ $\!\! (00,\lambda_3)^2 /$   $\!\! (00,\lambda_4) /$ $\!\! (00,0)^4$ \\

$A_1^2 A_4$ & $(W(2),0,0) /$ $\!\! (1,1,\lambda_1) /$ $\!\! (1,1,\lambda_4) /$ $\!\! (1,1,0)^2 /$ $\!\! (1,0,\lambda_1) /$ $\!\! (1,0,\lambda_2) /$   $\!\! (1,0,\lambda_3) /$ $\!\! (1,0,\lambda_4) /$  $\!\! (1,0,0)^2 /$ $\!\! (0,W(2),0) /$ $\!\! (0,1,\lambda_1) /$ $\!\! (0,1,\lambda_2) /$ $\!\! (0,1,\lambda_3) /$ $\!\! (0,1,\lambda_4) /$ $\!\! (0,1,0)^2 /$ $\!\! (0,0,W(\lambda_1+\lambda_4)) /$ $\!\! (0,0,\lambda_1)^2 /$ $\!\! (0,0,\lambda_2) /$  $\!\! (0,0,\lambda_3) /$ $\!\! (0,0,\lambda_4)^2 /$ $\!\! (0,0,0)^2$ \\

$A_3^2$ & $(W(101),000) /$  $\!\! (100,100) /$ $\!\! (100,010) /$ $\!\! (100,001) /$ $\!\! (100,000)^2 /$ $\!\! (010,100) /$ $\!\! (010,001) /$ $\!\! (010,000)^2  /$ $\!\! (001,100) /$  $\!\! (001,010) /$ $\!\! (001,001) /$ $\!\! (001,000)^2 /$ $\!\! (000,W(101)) /$  $\!\! (000,100)^2 /$ $\!\! (000,010)^2 /$ $\!\! (000,001)^2 /$ $\!\! (000,000)^{2}$ \\

$A_1 A_2 A_3$ & $(W(2),00,000) /$ $\!\! (1,10,001) /$ $\!\! (1,10,000)^2 /$ $\!\! (1,01,100) /$ $\!\! (1,01,000)^2 /$ $\!\! (1,00,100) /$ $\!\! (1,00,010)^2 /$ $\!\! (1,00,001) /$ $\!\! (0,W(11),000) /$   $\!\! (0,10,100)^2 /$ $\!\! (0,10,010) /$ $\!\! (0,10,000) /$ $\!\! (0,01,010) /$ $\!\! (0,01,001)^2 /$ $\!\!(0,01,000) /$ $\!\! (0,00,W(101)) /$ $\!\! (0,00,100)^2 /$ $\!\! (0,00,001)^2 /$ $\!\! (0,00,000)^4$ \\

$A_1^2 A_2^2$ & $(W(2),0,00,00) /$  $\!\!(1,1,10,00) /$ $\!\!(1,1,01,00) /$ $\!\!(1,1,00,10) /$ $\!\!(1,1,00,01) /$ $\!\!(1,0,10,01) /$ $\!\!(1,0,10,00) /$ $\!\!(1,0,01,10) /$ $\!\!(1,0,01,00) /$ $\!\!(1,0,00,10) /$ $\!\!(1,0,00,01) /$ $\!\!(1,0,00,00)^2 /$ $\!\! (0,W(2),00,00) /$ $\!\!(0,1,10,10) /$ $\!\!(0,1,10,00) /$ $\!\!(0,1,01,01) /$ $\!\!(0,1,01,00) /$ $\!\!(0,1,00,10) /$ $\!\!(0,1,00,01) /$ $\!\!(0,1,00,00)^2 /$ $\!\! (0,0,W(11),00) /$ $\!\!(0,0,10,10) /$ $\!\!(0,0,10,01) /$ $\!\!(0,0,10,00) /$ $\!\!(0,0,01,10) /$ $\!\!(0,0,01,01) /$  $\!\!(0,0,01,00) /$ $\!\! (0,0,00,W(11)) /$ $\!\!(0,0,00,10) /$ $\!\!(0,0,00,01) /$ $\!\!(0,0,00,00)^2 $ \\

$A_5$ & $W(\lambda_1 + \lambda_5) /$ $\!\! \lambda_1^6 /$ $\!\! \lambda_2^3 /$ $\!\! \lambda_3^2 /$ $\!\! \lambda_4^3 /$ $\!\! \lambda_5^6 /$ $\!\! 0^{11}$  \\
 
$A_1 A_4$ & $(W(2),0) /$ $\!\! (1,\lambda_1)^3 /$ $\!\! (1,\lambda_2) /$  $\!\! (1,\lambda_3) /$  $\!\! (1,\lambda_4)^3 /$ $\!\! (1,0)^6 /$ $\!\! (0,W(\lambda_1+\lambda_4)) /$ $\!\! (0,\lambda_1)^4 /$ $\!\! (0,\lambda_2)^3 /$ $\!\! (0,\lambda_3)^3 /$ $\!\! (0,\lambda_4)^4 /$  $\!\! (0,0)^9$ \\

$A_2 A_3$ & $(W(11),000) /$ $\!\! (10,100)^2 /$ $\!\! (10,010) /$ $\!\! (10,001)^2 /$ $\!\! (10,000)^5 /$ $\!\! (01,100)^2 /$  $\!\! (01,010) /$ $\!\! (01,001)^2 /$  $\!\! (01,000)^5 /$ $\!\! (00,W(101)) /$ $\!\! (00,100)^4 /$ $\!\! (00,010)^4 /$ $\!\! (00,001)^4 /$ $\!\! (00,000)^7$ \\

$A_1^2 A_3$ & $(W(2),0,000) /$ $\!\! (1,1,100) /$   $\!\! (1,1,001) /$ $\!\! (1,1,000)^4 /$ $\!\! (1,0,100)^2 /$ $\!\! (1,0,010)^2 /$ $\!\! (1,0,001)^2 /$ $\!\! (1,0,000)^4 /$ $\!\! (0,W(2),000) /$ $\!\! (0,1,100)^2 /$ $\!\! (0,1,010)^2 /$ $\!\! (0,1,001)^2 /$ $\!\! (0,1,000)^4 /$  $\!\! (0,0,W(101)) /$ $\!\! (0,0,100)^4 /$ $\!\! (0,0,010)^2 /$ $\!\! (0,0,001)^4 /$ $\!\! (0,0,000)^7$ \\

$A_1 A_2^2$ & $(W(2),00,00) /$ $\!\!(1,10,01) /$ $\!\!(1,10,00)^3 /$ $\!\!(1,01,10) /$ $\!\!(1,01,00)^3 /$ $\!\!(1,00,10)^3 /$ $\!\!(1,00,01)^3 /$   $\!\!(1,00,00)^2 /$ $\!\! (0,W(11),00) /$  $\!\!(0,10,10)^3 /$  $\!\!(0,10,01) /$ $\!\!(0,10,00)^3 /$ $\!\!(0,01,10) /$ $\!\!(0,01,01)^3 /$ $\!\!(0,01,00)^3 /$ $\!\! (0,00,W(11)) /$ $\!\!(0,00,10)^3 /$ $\!\!(0,00,01)^3 /$ $\!\!(0,00,00)^9 $ \\

$A_1^3 A_2$ & $(W(2),0,0,00) /$ $\!\!(1,1,1,00)^2 /$  $\!\!(1,1,0,10) /$ $\!\!(1,1,0,01) /$ $\!\!(1,1,0,00)^2 /$ $\!\!(1,0,1,10) /$ $\!\!(1,0,1,01) /$ $\!\!(1,0,1,00)^2 /$ $\!\!(1,0,0,10)^2 /$ $\!\! (1,0,0,01)^2 /$ $\!\!(1,0,0,00)^4 /$ $\!\! (0,W(2),0,00) /$  $\!\!(0,1,1,10) /$ $\!\!(0,1,1,01) /$ $\!\!(0,1,1,00)^2 /$ $\!\!(0,1,0,10)^2 /$ $\!\!(0,1,0,01)^2 /$  $\!\!(0,1,0,00)^4 /$ $\!\! (0,0,W(2),00) /$  $\!\!(0,0,1,10)^2 /$ $\!\!(0,0,1,01)^2 /$ $\!\!(0,0,1,00)^4 /$ $\!\! (0,0,0,W(11)) /$ $\!\!(0,0,0,10)^3 /$ $\!\!(0,0,0,01)^3 /$ $\!\!(0,0,0,00)^5 $ \\

$A_4$ & $W(\lambda_1+\lambda_4) /$ $\!\! \lambda_1^{10} /$ $\!\! \lambda_2^5 /$ $\!\! \lambda_3^5 /$ $\!\! \lambda_4^{10} /$ $\!\! 0^{24}$  \\

$A_1 A_3$ & $(W(2),000) /$  $\!\! (1,100)^4 /$  $\!\! (1,010)^2 /$ $\!\! (1,001)^4 /$ $\!\! (1,000)^{12} /$ $\!\! (0,W(101)) /$ $\!\! (0,100)^8 /$ $\!\! (0,010)^6 /$ $\!\! (0,001)^8 /$ $\!\! (0,000)^{18}$ \\

$A_2^2$ & $(W(11),00) /$ $\!\!(10,10)^3 /$ $\!\!(10,01)^3 /$ $\!\!(10,00)^9 /$ $\!\!(01,10)^3 /$ $\!\!(01,01)^3 /$ $\!\!(01,00)^9 /$ $\!\! (00,W(11)) /$ $\!\!(00,10)^9 /$ $\!\!(00,01)^9 /$ $\!\!(00,00)^{16} $ \\

$A_1^2 A_2$ & $(W(2),0,00) /$  $\!\!(1,1,10) /$ $\!\!(1,1,01) /$ $\!\!(1,1,00)^6 /$ $\!\!(1,0,10)^4 /$ $\!\!(1,0,01)^4 /$ $\!\!(1,0,00)^8 /$ $\!\! (0,W(2),00) /$  $\!\!(0,1,10)^4 /$ $\!\!(0,1,01)^4 /$ $\!\!(0,1,00)^8 /$ $\!\! (0,0,W(11)) /$ $\!\!(0,0,10)^7 /$ $\!\!(0,0,01)^7 /$ $\!\!(0,0,00)^{16} $ \\

$A_1^4$ & $(W(2),0,0,0) /$    $\!\!(1,1,1,0)^2 /$  $\!\!(1,1,0,1)^2 /$ $\!\!(1,1,0,0)^4 /$ $\!\!(1,0,1,1)^2 /$ $\!\!(1,0,1,0)^4 /$ $\!\!(1,0,0,1)^4 /$ $\!\!(1,0,0,0)^8 /$ $\!\! (0,W(2),0,0) /$ $\!\!(0,1,1,1)^2 /$ $\!\!(0,1,1,0)^4 /$ $\!\!(0,1,0,1)^4 /$ $\!\!(0,1,0,0)^8 /$ $\!\! (0,0,W(2),0) /$ $\!\!(0,0,1,1)^4 /$ $\!\!(0,0,1,0)^8 /$ $\!\! (0,0,0,W(2)) /$  $\!\!(0,0,0,1)^8 /$ $\!\!(0,0,0,0)^{12} $ \\

$A_3$ & $W(101) /$ $\!\! 100^{16} /$ $\!\! 010^{10} /$ $\!\! 001^{16} /$  $\!\! 000^{45}$ \\

$A_1 A_2$ & $(W(2),00) /$  $\!\! (1,10)^{6} /$ $\!\! (1,01)^{6} /$ $\!\! (1,00)^{20} /$  $\!\! (0,W(11)) /$ $\!\! (0,10)^{15} /$ $\!\! (0,01)^{15} /$ $\!\! (0,00)^{35}$  \\

$A_1^3$ & $(W(2),0,0) /$  $\!\! (1,1,1)^{2} /$ $\!\! (1,1,0)^{8} /$ $\!\! (1,0,1)^{8} /$ $\!\! (1,0,0)^{16} /$ $\!\! (0,W(2),0) /$  $\!\! (0,1,1)^{8} /$   $\!\! (0,1,0)^{16} /$ $\!\! (0,0,W(2)) /$  $\!\! (0,0,1)^{16} /$ $\!\! (0,0,0)^{31}$  \\

$A_2$ & $W(11) /$ $\!\! 10^{27} /$ $\!\! 01^{27} /$ $\!\! 00^{78}$ \\

$A_1^2$ & $(W(2),0) /$ $\!\! (1,1)^{12} /$ $\!\! (1,0)^{32} /$ $\!\! (0,W(2)) /$ $\!\! (0,1)^{32} /$  $\!\! (0,0)^{66}$ \\

$A_1$ &  $W(2) / 1^{56} / 0^{133}$ \\

\hline

\end{longtable}

%% file: ThomasIrredSubgp.bbl
\providecommand{\bysame}{\leavevmode\hbox to3em{\hrulefill}\thinspace}
\providecommand{\MR}{\relax\ifhmode\unskip\space\fi MR }
\providecommand{\MRhref}[2]{%
  \href{http://www.ams.org/mathscinet-getitem?mr=#1}{#2}
}
\providecommand{\href}[2]{#2}
\begin{thebibliography}{BMRT13}

\bibitem[ABS90]{ABS}
H.~Azad, M.~Barry, and G.~M. Seitz, \emph{On the structure of parabolic
  subgroups}, Comm. in Alg. \textbf{18} (1990), 551--562.

\bibitem[Ame05]{bon}
B.~Amende, \emph{{$G$}-irreducible subgroups of type {$A_1$}}, Ph.D. thesis,
  University of Oregon, 2005.

\bibitem[BCP97]{magma}
W.~Bosma, J.~Cannon, and C.~Playoust, \emph{The {M}agma algebra system. {I}.
  {T}he user language}, J. Symbolic Comput. \textbf{24} (1997), 235--265.

\bibitem[BGM]{BGM16}
M.~Bate, H.~Geranios, and B.~Martin, \emph{Orbit closures and invariants},
  arXiv:1604.00924 [math.AG].

\bibitem[BMR05]{bmr}
M.~Bate, B.~Martin, and G.~R{\"o}hrle, \emph{A geometric approach to complete
  reducibility}, Invent. Math. \textbf{161} (2005), 177--218.

\bibitem[BMRT11]{BMRT11}
M.~Bate, B.~Martin, G.~R{\"o}hrle, and R.~Tange, \emph{Complete reducibility
  and conjugacy classes of tuples in algebraic groups and {L}ie algebras},
  Math. Z. \textbf{269} (2011), no.~3-4, 809--832.

\bibitem[BMRT13]{BMRT13}
\bysame, \emph{Closed orbits and uniform {$S$}-instability in geometric
  invariant theory}, Trans. Amer. Math. Soc. \textbf{365} (2013), no.~7,
  3643--3673.

\bibitem[Bor91]{borel}
A.~Borel, \emph{Linear algebraic groups}, second ed., Graduate Texts in
  Mathematics, vol. 126, Springer-Verlag, New York, 1991.

\bibitem[Bou68]{bourbaki}
N.~Bourbaki, \emph{Groupes et {A}lgebres de {L}ie \textup{(Chapters 4,5,6)}},
  Hermann, Paris, 1968.

\bibitem[BT71]{BT}
A.~Borel and J.~Tits, \emph{{\'{E}l\'{e}ments unipotents et sous-groupes
  paraboliques de groupes r\'{e}ductifs}}, Invent. Math. \textbf{12} (1971),
  95--104.

\bibitem[Cap09]{Cap09}
P.-E. Caprace, \emph{``{A}bstract'' homomorphisms of split {K}ac-{M}oody
  groups}, Mem. Amer. Math. Soc. \textbf{198} (2009), no.~924.

\bibitem[Car72]{car}
R.~W. Carter, \emph{Conjugacy classes in the {W}eyl group}, Compositio Math
  \textbf{25} (1972), 1--59.

\bibitem[CLSS92]{clss}
A.~M. Cohen, M.~W. Liebeck, J.~Saxl, and G.~M. Seitz, \emph{The local maximal
  subgroups of exceptional groups of {L}ie type, finite and algebraic}, Proc.
  London Math. Soc. \textbf{64} (1992), 21--48.

\bibitem[Daw]{Daw13}
D.~Dawson, \emph{Complete reducibility in {E}uclidean twin buildings},
  arXiv:1110.1048 [math.GR].

\bibitem[GHTar]{GHT16}
R.~M. Guralnick, F.~Herzig, and P.~Tiep, \emph{Adequate subgroups and
  indecomposable modules}, J. Eur. Math. Soc. (to appear).

\bibitem[GLS98]{gls3}
D.~Gorenstein, R.~Lyons, and R.~Solomon, \emph{The classification of the finite
  simple groups. {N}umber 3}, Mathematical Surveys and Monographs, vol. 40.3,
  American Mathematical Society, Providence, RI, 1998.

\bibitem[Gur99]{gur99}
R.~M. Guralnick, \emph{Small representations are completely reducible}, J.
  Algebra \textbf{220} (1999), no.~2, 531--541.

\bibitem[KL90]{KL}
P.~Kleidman and M.~Liebeck, \emph{The subgroup structure of the finite
  classical groups}, Cambridge University Press, 1990.

\bibitem[Kle87]{K}
P.~Kleidman, \emph{The maximal subgroups of the finite 8-dimensional orthogonal
  groups {$P \Omega^+_8(q)$} and of their automorphism groups}, J. Algebra
  \textbf{110} (1987), 173--242.

\bibitem[Litar]{Lit16}
A.~J. Litterick, \emph{On non-generic finite subgroups of exceptional algebraic
  groups}, Mem. Amer. Math. Soc. (to appear).

\bibitem[LS94]{LS6}
M.~W. Liebeck and G.~M. Seitz, \emph{Subgroups generated by root elements in
  groups of {L}ie type}, Annals of Mathematics \textbf{139} (1994), 293--361.

\bibitem[LS96]{LS3}
\bysame, \emph{Reductive subgroups of exceptional algebraic groups}, Mem. Amer.
  Math. Soc. \textbf{121} (1996), no.~580.

\bibitem[LS98]{LS2}
\bysame, \emph{On the subgroup structure of classical groups}, Invent. Math.
  \textbf{134} (1998), 427--453.

\bibitem[LS99]{LS8}
\bysame, \emph{On finite subgroups of exceptional algebraic groups}, J. Reine
  Angew. Math. \textbf{515} (1999), 25--72.

\bibitem[LS03]{LS4}
\bysame, \emph{Variations on a theme of {S}teinberg}, J. Algebra \textbf{260}
  (2003), 261--297.

\bibitem[LS04]{LS1}
\bysame, \emph{The maximal subgroups of positive dimension in exceptional
  algebraic groups}, Mem. Amer. Math. Soc. \textbf{169} (2004), no.~802.

\bibitem[LS12]{LS9}
\bysame, \emph{Unipotent and nilpotent classes in simple algebraic groups and
  {L}ie algebras}, Mathematical surveys and monographs, Amer. Math. Soc., 2012.

\bibitem[LST15]{LST15}
M.~W. Liebeck, G.~M. Seitz, and D.~M. Testerman, \emph{Distinguished unipotent
  elements and multiplicity-free subgroups of simple algebraic groups}, Pacific
  J. Math. \textbf{279} (2015), no.~1-2.

\bibitem[LT99]{Law}
R.~Lawther and D.~M. Testerman, \emph{{$A_1$} subgroups of exceptional
  algebraic groups}, Mem. Amer. Math. Soc. \textbf{141} (1999), no.~674.

\bibitem[LT04]{LT}
M.~W. Liebeck and D.~M. Testerman, \emph{Irreducible subgroups of algebraic
  groups}, Quart. {J}. {M}ath. \textbf{55} (2004), 47--55.

\bibitem[LTar]{littho}
A.~J. Litterick and A.~R. Thomas, \emph{Complete reducibility in good
  characteristic}, Trans. Amer. Math. Soc. (to appear).

\bibitem[L{\"u}b01]{lubeck}
F.~L{\"u}beck, \emph{Small degree representations of finite {C}hevalley groups
  in defining characteristic}, LMS {J}. {C}omput. {M}ath. \textbf{4} (2001),
  135--169.

\bibitem[McN07]{McN07}
G.~McNinch, \emph{Completely reducible {L}ie subalgebras}, Transform. Groups
  \textbf{12} (2007), no.~1, 127--135.

\bibitem[Rot09]{rot}
J.~J. Rotman, \emph{An introduction to homological algebra}, second ed.,
  Universitext, Springer, New York, 2009.

\bibitem[Sei91]{se2}
G.~M. Seitz, \emph{Maximal subgroups of exceptional algebraic groups}, Mem.
  Amer. Math. Soc. \textbf{90} (1991), no.~441.

\bibitem[Ser05]{ser04}
J.-P. Serre, \emph{Compl\`{e}te r\'{e}ductibilit\'{e}}, Ast\'{e}risque
  \textbf{299} (2005), Exp. No. 932, S\'{e}minaire Bourbaki. Vol. 2003-2004.

\bibitem[Ste63]{stein}
R.~Steinberg, \emph{Representations of algebraic groups}, Nagoya Math. J.
  (1963), no.~22, 33--56.

\bibitem[Ste10]{g2dav}
D.~I. Stewart, \emph{The reductive subgroups of {$G_2$}}, J. Group Theory
  \textbf{13} (2010), 117--130.

\bibitem[Ste13]{dav}
\bysame, \emph{The reductive subgroups of {$F_4$}}, Mem. Amer. Math. Soc.
  \textbf{223} (2013), no.~1049.

\bibitem[STar]{StewTho}
D.~I. Stewart and A.~R. Thomas, \emph{The {J}acobson--{M}orozov theorem and
  complete reducibility of {L}ie subalgebras}, Proc. Lond. Math. Soc. (to
  appear).

\bibitem[Tho15]{tho1}
A.~R. Thomas, \emph{Simple irreducible subgroups of exceptional algebraic
  groups}, J. Algebra \textbf{423} (2015), 190--238.

\bibitem[Tho16]{tho2}
\bysame, \emph{Irreducible {$A_1$} subgroups of exceptional algebraic groups},
  J. Algebra \textbf{447} (2016), 240--296.

\end{thebibliography}
